\documentclass{amsart}
\usepackage{amssymb}
\usepackage{tikz-cd}
\usepackage{comment}
\usepackage[mathcal]{euscript}
\usepackage{mathrsfs}
\usepackage[all,cmtip]{xy}
\usepackage[figuresleft]{rotating}

\theoremstyle{plain}

\theoremstyle{definition}

\theoremstyle{remark}

\newcommand{\im}{{\rm Im}\,}
\newcommand{\K}{{\rm Ker}\,}
\newcommand{\Co}{{\rm Coker}\,}
\newcommand{\cl}{{\rm Colim}\,}
\newcommand{\li}{{\rm Lim}\,}

\begin{document}

\title[Remarks on cotorsion theories]
{Remarks on cotorsion theories}

\author{Alexandru Stanculescu}

%\thanks{dffdfd}

%\date{\today}

\email{astanculescu@crcmail.net}

\begin{abstract}
We present some results from classical homological algebra using the language 
of cotorsion theories in abelian categories. The results are a couple of foundational 
facts about homological dimension, the K\"{u}nneth formula and the universal 
coefficients formula. We also present  some ``totally acyclic'' cotorsion theories
and revisit a couple of results from the literature concerning resolutions of 
chain complexes
in abelian categories.
\end{abstract}
\maketitle

We regard cotorsion theories as an organizational principle in homological algebra. 
The main purpose of this paper is to illustrate this point of view by presenting some 
results from classical homological algebra using the language of cotorsion 
theories in abelian categories. The results are a couple of foundational facts 
about homological dimension \cite[\S8 n\textsuperscript{\scriptsize{o}} 1,3]{Bou}, 
the K\"{u}nneth formula [\emph{loc. cit.}, \S4 n\textsuperscript{\scriptsize{o}} 7] and 
the universal coefficients formula [\emph{loc. cit.}, \S5 n\textsuperscript{\scriptsize{o}} 6].
Presenting homological dimension (resp. the universal coefficients formula)
in the language of cotorsion theories means here replacing the pair injective-projective
module by a pair of comparable cotorsion theories (resp. by an arbitrary cotorsion theory)
in an abelian category; the K\"{u}nneth formula in abelian categories roughly means 
the usual K\"{u}nneth formula where the flat cotorsion theory has been extended 
from modules to abelian categories that fit in an abelian Tensor-Hom-Cotensor 
situation \cite{Gra} and the $\mathsf{Tor}$ bifunctors are constructed using 
(complete) flat resolutions. A Tensor-Hom-Cotensor situation is also known under 
the name adjunction of two variables. In linear algebra, two basic examples 
of THC-situations are the one that has as Hom bifunctor the abelian group of 
homomorphisms between two modules and the one that has as 
Tensor bifunctor the tensor product between right and left modules.
The idea of homological algebra in an abelian THC-situation is suggested 
by Bourbaki's functorial treatment of classical homological algebra.
%; nous lui présentons nos hommages.
Overall, extending the triad injective-projective-flat module to more general 
objects is, in the case of the above mentioned results, formal. 

Another purpose of the paper is to present some ``totally acyclic'' cotorsion 
theories in categories of chain complexes and revisit some results from the 
literature concerning resolutions of complexes. The construction of 
resolutions may help establishing that a certain cotorsion theory in 
complexes is complete.

All the arguments use standard properties of abelian categories only and 
we do not use the bifunctors $\mathsf{Ext}$ since the very notion 
of cotorsion theory can be formulated without them. Also, given the level 
of generality that we want to work and with the aim of making the paper 
as simple as possible, we felt that the absence of $\mathsf{Ext}$ would 
not make the results or proofs artificial. The approach to establishing 
that a class of objects is part of a cotorsion theory emphasizes the choice 
of (co)generating classes and adjoint functors. 

Many of our results that involve cotorsion theories assume that these are 
complete. We do not address here the question of completeness of the 
cotorsion theories that we construct as this seems to require different 
techniques. 

In section 1 we introduce and give some basic properties of cotorsion 
theories in abelian categories without using $\mathsf{Ext}$. We spend 
some time studying stability properties of the classes 
$\K(\mathcal{K}\twoheadrightarrow\mathcal{L})$ and 
$\Co(\mathcal{K}^{\perp}\rightarrowtail^{\perp}\mathcal{L})$ 
(whatever they mean).

In section 2 we study the possibility of ``pulling back'' cotorsion 
theories along a functor that has an adjoint (\ref{sec:2.1.2}) and 
of ``pulling back-and-restricting'' cotorsion theories 
(\ref{sec:2.2.2}, \ref{sec:2.2.4}, \ref{sec:2.2.5}) in 
abelian categories that fit in an abelian Tensor-Hom-Cotensor 
situation.

In section 3 we generalize (\ref{sec:3.1.4}) the flat cotorsion theory 
\cite[Lemma 7.1.4]{EJ} from modules to abelian categories that fit in
an abelian THC-situation. The resulting generalized flat cotorsion theory 
has the usual properties (\ref{sec:3.1.5}, \ref{sec:3.1.6}, \ref{sec:3.1.7}).
% We present in \ref{sec:2.2.4} a
%relative version of the generalized flat cotorsion theory and in \ref{sec:2.6.5}
%rudiments of Matlis duality in an abelian THC-situation.
%Completeness of the generalized flat cotorsion theory is not addressed.

In section 4 we recall some facts about complexes in abelian categories.
In section 5 we associate (\ref{sec:5.2.1}, \ref{sec:5.2.4}) to a pair of 
comparable cotorsion theories two new ones; the left class of 
the first consists of objects of ``projective'' dimension $\leqslant n$ 
and the right class of the second consists of objects of 
``injective'' dimension $\leqslant n$. 

In section 6 we revisit the standard and the bar THC-situations for categories of 
complexes in abelian categories \cite[Section 4.3]{Pe}. The former is an enriched 
situation (\ref{sec:6.1.9}). The flat complexes associated to the two 
THC-situations have the usual description (\ref{sec:6.1.3}, \ref{sec:6.2.3}).

In section 7 we very briefly review (\ref{sec:7.1}), following Bourbaki 
\cite[\S3,4]{Bou}, the construction of the torsion product associated 
to the Tensor bifunctor of the standard THC-situation. The only 
difference with \emph{loc. cit.} is that we do not use functorial 
free resolutions but complete flat resolutions; for this reason we 
require that the generalized flat cotorsion theory is complete.
The torsion product is used to prove a K\"{u}nneth type theorem for 
complexes in abelian categories (\ref{sec:7.3.10}, \ref{sec:7.3.11}). The 
proof follows Bourbaki [\emph{loc. cit.}, \S4 n\textsuperscript{\scriptsize{o}} 7],
except that we work with the second argument of the Tensor bifunctor; 
this choice requires some small sign considerations. 
%We don't know whether such a theorem exists in the literature.

In section 8 we revisit some elementary facts about the enriched and 
homotopy orthogonality relations for complexes in abelian categories 
and the formal theory of $K$-flat complexes associated to the standard 
THC-situation (\cite[Sections 1,5]{Sp},\cite[Proposition 3.7]{Gi2}).

In section 9 we prove (\ref{sec:9.1.4}) a universal coefficients theorem 
for complexes in an abelian category. The proof follows Bourbaki 
[\emph{loc. cit.}, \S5 n\textsuperscript{\scriptsize{o}} 6].
%We don't know whether such a theorem exists in the literature.

In section 10 we present (\ref{sec:10.2.1}, \ref{sec:10.2.2}, \ref{sec:10.3.1}) 
a couple of  ``totally acyclic" cotorsion theories and the ``dg-flat'' cotorsion 
theory in categories of complexes in abelian categories.

In section 11 we study the compatibility between some cotorsion
theories and corner morphisms for categories of complexes 
\cite[Theorem 5.1]{Gi4}.

In section 12 we revisit (\ref{sec:12.1.3}, \ref{sec:12.1.9})
a couple of results from the literature \cite{Ha},\cite{Sp} 
concerning resolutions of complexes.
%COVERS AND ADJUNCTIONS, DAY CONVOLUTION AND FLATNESS

\section{Cotorsion theories}

\subsection{} \label{sec:1.1}
Monomorphisms and epimorphisms in abelian categories will be denoted by 
the symbols $\rightarrowtail$ and $\twoheadrightarrow$, respectively.

Let $\mathcal{A}$ be an abelian category. For our purposes a composable 
pair $u:A'\to A$, $v:A\to A''$ of morphisms of $\mathcal{A}$ will be called 
\emph{sequence} and denoted by $(u,v)$. A sequence $(u,v)$ in 
$\mathcal{A}$ is \emph{semi-exact} if $vu=0$, \emph{exact} if $\im u= \K v$.
An exact sequence $A'\overset{u}\rightarrowtail A\to A''$ is \emph{left split} 
if $u$ has a retraction and an exact sequence 
$A'\to A\overset{v}\twoheadrightarrow A''$ is \emph{right split} 
if $v$ has a section.

\subsubsection{} \label{sec:1.1.1}
Suppose that $\mathcal{A}$ has a set $(C_{i})_{i\in I}$ of cogenerators.
Let $A'\overset{u}\rightarrow A\overset{v}\rightarrow A''$ 
be a sequence in $\mathcal{A}$. 

(1) If the sequence $\mathcal{A}(A'',C_{i})\overset{\bar{v}}\to
\mathcal{A}(A,C_{i}) \overset{\bar{u}}\to\mathcal{A}(A',C_{i})$ is 
semi-exact (exact) for all $i\in I$ then $(u,v)$ is semi-exact (exact).

(2) If $\bar{u}$ is a monomorphism for all $i\in I$ then $u$ is an epimorphism and
if $\bar{v}$ is an epimorphism for all $i\in I$ then $v$ is a monomorphism.

%(3) Suppose that $A''/\im v\in \mathcal{K}$ and that $C_{i}\in \mathcal{K}^{\perp}$ for each $i\in  I$. 
%If the sequence $(\bar{v},\bar{u})$ is exact for each $i\in I$ then $(u,v)$ is exact. 
\begin{proof}
(1) Suppose that $(\bar{v},\bar{u})$ is exact for each $i\in  I$. 
Consider the commutative diagram 
\begin{displaymath}
\xymatrix{
&&&{\Co i}\ar[dr]^{k}\\
{A'} \ar[rr]^{u}\ar@{->>}[dr]_{p} && {A} \ar[rr]^{v}\ar@{->>}[dr]_{q} \ar@{->>}[ur]^{r}&& {A''}\\
& {\im u} \ar@{>->}[ur]_{i} && {\im v} \ar@{>->}[ur]_{j}
}
\end{displaymath}
There is a morphism $\delta:\Co i\to \im v$ such that $\delta r=q$ and 
$j\delta=k$. Let $f:B\to \Co i$ be such that $\delta f=0$ but $f\neq 0$. 
There is then $h_{i}:\Co i\to C_{i}$ such that $h_{i}f\neq 0$.
We have $\bar{u}(h_{i}r)=0$, hence there is $g_{i}:A''\to C_{i}$ such 
that $g_{i}v=h_{i}r$. We have $g_{i}kr=h_{i}r$, so $g_{i}k=h_{i}$. 
This implies that $h_{i}f=g_{i}kf=g_{i}j\delta f=0$, contradiction. 
Therefore $\delta$ is an isomorphism, hence $(u,v)$ is exact.
\end{proof}

\subsubsection{} \label{sec:1.1.2}
Suppose that $\mathcal{K}$ is a class of objects of $\mathcal{A}$ 
such that for every $A\in Ob(\mathcal{A})$ there is a 
monomorphism $A\to X$ with $X\in \mathcal{K}$.
Let $A'\overset{u}\rightarrow A\overset{v}\rightarrow A''$ 
be a semi-exact sequence in $\mathcal{A}$. 

(1) If the sequence $\mathcal{A}(A'',X)\overset{\bar{v}}\to
\mathcal{A}(A,X) \overset{\bar{u}}\to\mathcal{A}(A',X)$ is 
exact for all $X\in\mathcal{K}$ then $(u,v)$ is exact.

(2) If $\bar{u}$ is an epimorphism for all $X\in\mathcal{K}$ 
then $u$ is  a monomorphism.
\begin{proof}
(1) Consider the commutative diagram 
\begin{displaymath}
\xymatrix{
{A'} \ar[rrrr]^{u}\ar@{->>}[dr]_{p} &&&& {A} \ar[rr]^{v}\ar@{->>}[dr]_{q} && {A''}\\
& {\im u} \ar@{>->}[r]_{f} & {\K v} \ar@{>->}[urr]_{i} &&& {\im v} \ar@{>->}[ur]_{j}
}
\end{displaymath}
Let $\alpha:\K v\to B$ be a morphism such that $\alpha f=0$. Form the commutative diagram
\[
\xymatrix{
{\K v} \ar@{>->}[r]^{i} \ar[d]_{\alpha}& {A}\ar@{->>}[r]^{q}\ar[d]_{\alpha'} & {\im v}\ar@{=}[d]\\
{B} \ar@{>->}[r]_{i'} & {PO}\ar@{->>}[r] & {\im v}
}
\]
where $PO$ means pushout. We can find a monomorphism $k:PO\to X$, where $X\in \mathcal{K}$.
We have $k\alpha'u=ki'\alpha fp=0$, therefore there is $h:A''\to X$ such that $hv=k\alpha'$.
Then $ki'\alpha=hvi=0$, so $\alpha=0$ and therefore $f$ is an epimorphism.
\end{proof}

\subsection{} \label{sec:1.2}
We denote by $\mathrm{Ab}$ the category of abelian groups. Let $\mathcal{A}$ 
be an abelian category and $\mathcal{K}$ a class of objects of $\mathcal{A}$. 
We define
\[
\mathcal{K}^{\perp}=\{X\in Ob(\mathcal{A}), \ 
\mathcal{A}(-,X):\mathcal{A}^{op}\to \mathrm{Ab} {\rm \ preserves\
epimorphisms\ with\ kernel\ in}\ \mathcal{K}\}
\]
and
\[
^{\perp}\mathcal{K}=\{X\in Ob(\mathcal{A}), \ 
\mathcal{A}(X,-):\mathcal{A}\to \mathrm{Ab} {\rm \ preserves\
epimorphisms\ with\ kernel\ in}\ \mathcal{K}\}
\] 
When $\mathcal{K}=Ob(\mathcal{A})$, $\mathcal{K}^{\perp}$ 
consists of the injective objects of $\mathcal{A}$ 
and $^{\perp}\mathcal{K}$ of the projective objects.

\subsubsection{} \label{sec:1.2.1}
(1) An object $X$ of $\mathcal{A}$ is in $\mathcal{K}^{\perp}$ 
if and only if for every exact sequence
\begin{displaymath}
A'\overset{u}\rightarrow A\overset{v}\rightarrow A''
\end{displaymath}
in $\mathcal{A}$ with $A''/\im v\in \mathcal{K}$ the sequence 
\[
\mathcal{A}(A'',X)\overset{\bar{v}}\to\mathcal{A}(A,X) 
\overset{\bar{u}}\to\mathcal{A}(A',X)
\] 
is exact.

(1bis) An object $X$ of $\mathcal{A}$ is in $^{\perp}\mathcal{K}$ 
if and only if for every exact sequence
\begin{displaymath}
A'\overset{u}\rightarrow A\overset{v}\rightarrow A''
\end{displaymath}
in $\mathcal{A}$ with $\K u\in \mathcal{K}$ the sequence 
\[
\mathcal{A}(X,A')\overset{\bar{u}}\to\mathcal{A}(X,A) 
\overset{\bar{v}}\to\mathcal{A}(X,A'')
\] 
is exact.

(2) If $\mathcal{K},\mathcal{L}$ are two classes of objects of $\mathcal{A}$,

(a) $\mathcal{K}\subset \mathcal{L}$ implies 
$\mathcal{L}^{\perp}\subset\mathcal{K}^{\perp}$;

(b) $\mathcal{K}^{\perp}\cap\mathcal{L}^{\perp}=(\mathcal{K}\cup\mathcal{L})^{\perp}$.

(3) $\mathcal{K}\subset {}^{\perp}(\mathcal{K}^{\perp})$ and 
$({}^{\perp}(\mathcal{K}^{\perp}))^{\perp}=\mathcal{K}^{\perp}$.

(4) \cite[Proposition 4.2]{Ho} If $\mathcal{K},\mathcal{L}$ are two 
classes of objects of $\mathcal{A}$, the following are equivalent:

(a) $\mathcal{L}\subset\mathcal{K}^{\perp}$;

(b) every exact sequence $Y\rightarrowtail A\twoheadrightarrow X$
with $Y\in\mathcal{L},X\in\mathcal{K}$ splits;

(c) $\mathcal{K}\subset ^{\perp}\mathcal{L}$;

(d) every commutative diagram in $\mathcal{A}$
\begin{displaymath}
\xymatrix{
&{Y}\ar@{>->}[d]^{u}\\
{A}\ar[r]^{f}\ar@{>->}[d]_{i}&{C} \ar@{->>}[d]^{v}\\
{B}\ar@{->>}[d]_{p} \ar[r]^{g}&{D}\\
{X}
}
\end{displaymath}
with exact columns and with $Y\in\mathcal{L},X\in\mathcal{K}$ 
has a diagonal filler;

(e) for all exact sequences 
$Y\overset{u}\rightarrowtail C\overset{v}\twoheadrightarrow D$ 
and $A\overset{i}\rightarrowtail B\twoheadrightarrow X$ in 
$\mathcal{A}$ with $Y\in\mathcal{L},X\in\mathcal{K}$, 
the commutative diagram
\begin{displaymath}
\xymatrix{
{\Co\mathcal{A}(B,u)}\ar@{>->}[r] \ar[d]&{\mathcal{A}(B,D)} 
\ar[d]^{\bar{i}}\\
{\Co\mathcal{A}(A,u)} \ar@{>->}[r]&{\mathcal{A}(A,D)}
}
\end{displaymath} 
is a pullback.

(5) The class $\mathcal{K}^{\perp}$ is closed under retracts 
and extensions.

(6) Suppose that $\mathcal{A}$ is has arbitrary products. 
Let $(X_{i})_{i\in I}$ be a family of objects of $\mathcal{A}$. 
Then $\underset{i\in I}\prod X_{i}\in \mathcal{K}^{\perp}$
if and only if $X_{i}\in \mathcal{K}^{\perp}$ for each $i\in I$.

(7) Suppose that $\mathcal{A}$ has arbitrary coproducts. Let 
$X:\mathbb{N}\to\mathcal{A}$ be a functor with transition 
morphisms $f_{m,n}:X_{m}\to X_{n}$ for $m\leqslant n$, such 
that 

(a) $X_{n}\in^{\perp}\mathcal{K}$ for all $n\geqslant 0$;

(b) for all $m\geqslant 0$ there is $n\geqslant m$ such that for all 
$Y\in (^{\perp}\mathcal{K})^{\perp}$ and all $\alpha:X_{n}\to Y$ 
there is $\beta:X_{n+1}\to Y$ such that $\beta f_{m,n+1}=\alpha f_{m,n}$.

Then $\underset{n\in \mathbb{N}}\cl X_{n}\in^{\perp}\mathcal{K}$.

\begin{proof}
(1) Let $X\in \mathcal{K}^{\perp}$. We show that 
$\K \bar{u}\subset\im\bar{v}$. Let $f:A\to X$ be such 
that $\bar{u}(f)=0$. Factor $v$ as $v=jq$, where 
$q:A\twoheadrightarrow \im v$ and $j:\im v \rightarrowtail A''$; 
there is then a morphism $g:\im v\to X$ such that $gq=f$.
Form the commutative solid arrows diagram 
\[
\xymatrix{
{\im v} \ar@{>->}[r]^{j} \ar[d]_{g}& {A''}\ar@{->>}[r]
\ar@{..>}[dl]_{h}\ar[d] & {A''/\im v}\ar@{=}[d]\\
{X} \ar@{>->}[r] & {PO}\ar@{->>}[r] & {A''/\im v}
}
\]
where $PO$ means pushout. By assumption the bottom exact 
sequence splits, therefore there is $h:A''\to X$ such that $hj=g$. Then 
$\bar{v}(h)=hv=hjq=gp=f$. Conversely, let $X\in Ob(\mathcal{A})$
and apply the functor $\mathcal{A}(-,X)$ to the exact sequence 
$0\overset{u}\to A\overset{v}\to A''$ 
with $A''/\im v\in\mathcal{K}$.

(4) It is easy to see that (a)$\Leftrightarrow$(b)$\Leftrightarrow$(c).
That (d)$\Leftrightarrow$(e) follows from the fact that 
$\Co\mathcal{A}(B,u)=\{f:A\to D,{\rm \ there\ is\ g:A\to C\
such\ that }\ vg=f\}$. For (d)$\Rightarrow$(a) we use the 
commutative diagram
\begin{displaymath}
\xymatrix{
&{Y}\ar@{=}[d]\\
{A}\ar[r]\ar@{>->}[d]&{Y} \ar@{->>}[d]\\
{B}\ar@{->>}[d] \ar[r]&{0}\\
{X}
}
\end{displaymath}
where $Y\in\mathcal{L},X\in\mathcal{K}$. We prove (a)$\Rightarrow$(d).
Pulling back the exact sequence $(u,v)$ along $g$ we obtain 
the commutative diagram
\begin{displaymath}
\xymatrix{
&{Y}\ar@{>->}[d]^{u'}\\
{A}\ar@{>->}[r]^{h} \ar@{>->}[d]_{i}&{PB} \ar@{->>}[d]^{v'}\\
{B} \ar@{=}[r]&{B}
}
\end{displaymath}
where $PB$ means pullback and the sequence $(u',v')$ is exact. 
It suffices to show that the commutative square in the previous 
diagram has a diagonal filler $d:B\to PB$. Applying the snake diagram to 
\begin{displaymath}
\xymatrix{
{0}\ar[d]\ar[r]&{Y}\ar@{>->}[d]^{u'}\ar@{>->}[r]^{j}&{\Co h}\ar@{=}[d]\\
{A}\ar@{>->}[r]^{h} \ar@{>->}[d]_{i}&{PB} \ar@{->>}[d]^{v'}
\ar@{->>}[r]^{q}\ar[d]&{\Co h}\ar[d]\\
{B}\ar@{->>}[d]_{p}\ar@{=}[r]&{B}\ar[r]\ar[d]&{0}\\
{X}\ar[r]&{0}
}
\end{displaymath}
we obtain an exact sequence $Y\overset{j}\rightarrowtail 
\Co h\overset{q'}\twoheadrightarrow X$. 
By assumption this sequence splits, so $j$ has a retract $r$. 
Then $rq$ is a retract of $u'$, so $v'$
has a section $s$. We have $v'(h-si)=0$, hence there is 
$l:A\to Y$ such that $u'l=h-si$. 
Form the commutative solid arrows diagram 
\[
\xymatrix{
{A} \ar@{>->}[r]^{i} \ar[d]_{l}& {B}\ar@{->>}[r]^{p}
\ar@{..>}[dl]_{d'}\ar[d] & {X}\ar@{=}[d]\\
{Y} \ar@{>->}[r] & {PO}\ar@{->>}[r] & {X}
}
\]
where $PO$ means pushout. By assumption the bottom exact 
sequence splits, therefore there is $d':B\to Y$ such that $d'i=l$. 
Define $d=s+u'd'$; then $d$ is the required diagonal filler.

(5) We prove the closure under extensions.
Let $Y'\overset{u}\rightarrowtail Y\overset{v}
\twoheadrightarrow Y''$ be an exact sequence with 
$Y',Y''\in \mathcal{K}^{\perp}$ and $A\overset{i}
\rightarrowtail B\overset{p}\twoheadrightarrow X$
an exact sequence with $X\in\mathcal{K}$. Consider 
the commutative diagram
\begin{displaymath}
\xymatrix{
{\mathcal{A}(B,Y')}\ar@{>->}[rr]\ar@{->>}[d]^{\bar{i}}&&
{\mathcal{A}(B,Y)}\ar[rr] \ar[d]
&&{\mathcal{A}(B,Y'')} \ar@{->>}[d]^{\bar{i}}\\
{\mathcal{A}(A,Y')}\ar@{>->}[rr]&&{\mathcal{A}(A,Y)} \ar[rr]
&&{\mathcal{A}(A,Y'')}
}
\end{displaymath}
The two morphisms $\bar{i}$ are epimorphisms by assumption. 
By part 4 the diagram
\begin{displaymath}
\xymatrix{
{\Co\mathcal{A}(B,u)}\ar@{>->}[r]\ar[d]&{\mathcal{A}(B,Y'')} 
\ar@{->>}[d]^{\bar{i}}\\
{\Co\mathcal{A}(A,u)} \ar@{>->}[r]&{\mathcal{A}(A,Y'')}
}
\end{displaymath} 
is a pullback. We obtain the commutative diagram with exact 
rows
\begin{displaymath}
\xymatrix{
{\mathcal{A}(B,Y')}\ar@{>->}[rr]\ar@{->>}[d]^{\bar{i}}&&
{\mathcal{A}(B,Y)}\ar@{->>}[rr] \ar[d]
&&{\Co\mathcal{A}(B,u)} \ar@{->>}[d]^{\bar{i}}\\
{\mathcal{A}(A,Y')}\ar@{>->}[rr]&&{\mathcal{A}(A,Y)} \ar@{->>}[rr]
&&{\Co\mathcal{A}(A,u)}
}
\end{displaymath}
which implies that $Y\in\mathcal{K}^{\perp}$.

(6) Let $Y\rightarrowtail A\twoheadrightarrow B$ be an exact sequence with 
$Y\in ({}^{\perp}\mathcal{K})^{\perp}$; then by (a) we have the exact sequence 
$\mathcal{A}(X_{n},Y)\rightarrowtail \mathcal{A}(X_{n},A)\twoheadrightarrow 
\mathcal{A}(X_{n},B)$ for all $n\geqslant 0$. Consider the inverse system of 
abelian groups $\mathcal{A}(X_{n},Y)_{n\geqslant 0}$
with transition morphisms $\overline{f_{m,n}}:\mathcal{A}(X_{n},Y)
\to\mathcal{A}(X_{m},Y)$; by (b) it satisfies the Mittag-Leffler condition, 
hence the assertion.
\end{proof}

We note that \ref{sec:1.2.1}(7) is a variation on a classic result of Eklof.

\subsubsection{} \label{sec:1.2.2}
A pair $(\mathcal{X},\mathcal{Y})$ of classes of objects of $\mathcal{A}$ 
is a \emph{cotorsion theory} if $\mathcal{X}={}^{\perp}\mathcal{Y}$ and 
$\mathcal{Y}=\mathcal{X}^{\perp}$. We often say that $\mathcal{X}$ is the 
left class and $\mathcal{Y}$ is the right class of the cotorsion theory.

\paragraph{} If $(\mathcal{X},\mathcal{Y})$ is a cotorsion theory in $\mathcal{A}$ 
then $(\mathcal{Y},\mathcal{X})$ is a cotorsion theory in $\mathcal{A}^{op}$.

\paragraph{} Every abelian category has the minimal cotorsion theory 
(projectives, all objects) and the maximal cotorsion theory (all objects, injectives). 
The pair ($Ob(\mathcal{A}),Ob(\mathcal{A}))$ is a cotorsion theory if and only 
if all objects of $\mathcal{A}$ are projective if and only if all objects 
of $\mathcal{A}$ are injective if and only all short exact sequences in 
$\mathcal{A}$ split.

\subsubsection{} \label{sec:1.2.3}
Every class $\mathcal{S}$ of objects of $\mathcal{A}$ \emph{generates} a 
cotorsion theory $(^{\perp}\mathcal{S},(^{\perp}\mathcal{S})^{\perp})$ and 
\emph{cogenerates} a cotorsion theory $(^{\perp}(\mathcal{S}^{\perp}),
\mathcal{S}^{\perp})$ (\ref{sec:1.2.1}(3)).

\paragraph{} Given a class $\mathcal{K}$ of objects of $\mathcal{A}$, the 
pair $(\mathcal{K},\mathcal{K}^{\perp})$ is a cotorsion theory generated by 
$\mathcal{S}$ if and only if $\mathcal{S}\subset \mathcal{K}^{\perp}$ 
and $^{\perp}\mathcal{S}\subset\mathcal{K}$. Dually, the pair 
$({}^{\perp}\mathcal{K},\mathcal{K})$ is a cotorsion theory cogenerated by 
$\mathcal{S}$ if and only if $\mathcal{S}\subset {}^{\perp}\mathcal{K}$ 
and $\mathcal{S}^{\perp}\subset\mathcal{K}$.

\paragraph{} The minimal cotorsion theory is cogenerated by any family 
of projective objects. See \ref{sec:1.2.16} for the finite cogeneration 
of the maximal cotorsion theory.

\paragraph{} If $(^{\perp}\mathcal{K},\mathcal{K})$ and 
$(^{\perp}\mathcal{L},\mathcal{L})$ are cotorsion theories then (\ref{sec:1.2.1}(2))
$(^{\perp}(\mathcal{K}\cap\mathcal{L}),\mathcal{K}\cap\mathcal{L})$
is a cotorsion theory cogenerated by $^{\perp}\mathcal{K}\cup ^{\perp}\mathcal{L}$.

\subsubsection{} \label{sec:1.2.4}
Let $\mathcal{K},\mathcal{L}$ be two classes of objects of $\mathcal{A}$. 
We define $\K(\mathcal{K}\twoheadrightarrow\mathcal{L})$ to consist of 
objects $A$ of $\mathcal{A}$ for which there is an exact sequence 
$A\rightarrowtail K\twoheadrightarrow L$ with $K\in\mathcal{K}, L\in\mathcal{L}$.
We define $\Co(\mathcal{K}\rightarrowtail\mathcal{L})$ to consist of objects 
$A$ of $\mathcal{A}$ for which there is an exact sequence 
$K\rightarrowtail L\twoheadrightarrow A$ with $K\in\mathcal{K}, L\in\mathcal{L}$.

\paragraph{} If $\mathcal{L}$ is closed under extensions then 
\[
\K(\K(\mathcal{K}\twoheadrightarrow\mathcal{L})\twoheadrightarrow\mathcal{L})
\subset\K(\mathcal{K}\twoheadrightarrow\mathcal{L})
\]
and if $\mathcal{K}$ is closed under extensions then
\[
\Co(\mathcal{K}\rightarrowtail\Co(\mathcal{K}\rightarrowtail\mathcal{L}))
\subset\Co(\mathcal{K}\rightarrowtail\mathcal{L})
\]

\subsubsection{} \label{sec:1.2.5}
\cite[Chapter V, Remark 3.2(1)]{BR} 
Let $\mathcal{K}$ be a class of objects of $\mathcal{A}$ that 
is closed under retracts. If $\Co(\mathcal{K}^{\perp}\rightarrowtail\mathcal{K})
=Ob(\mathcal{A})$ then $(\mathcal{K},\mathcal{K}^{\perp})$ is a cotorsion theory.
\begin{proof}
Let $A\in {}^{\perp}(\mathcal{K}^{\perp})$. There is an exact sequence 
$Y\overset{u}\rightarrowtail X\overset{v}\twoheadrightarrow A$ with 
$X\in\mathcal{K}, Y\in\mathcal{K}^{\perp}$. We have then the exact sequence 
$\mathcal{A}(A,Y)\overset{\bar{u}}\rightarrowtail\mathcal{A}(A,X) \overset{\bar{v}}
\twoheadrightarrow\mathcal{A}(A,A)$. It follows that  $A$ is isomorphic to a 
retract of $X$.
\end{proof}

\subsubsection{} \label{sec:1.2.6}
Let $(\mathcal{K},\mathcal{K}^{\perp})$ and 
$(^{\perp}\mathcal{L},\mathcal{L})$ 
be two cotorsion theories in $\mathcal{A}$ with 
$\mathcal{K}^{\perp}\subset \mathcal{L}$.
Then $^{\perp}\mathcal{L}=\mathcal{K}\cap
\Co(\mathcal{K}^{\perp}\rightarrowtail^{\perp}\mathcal{L})$
and $\mathcal{K}^{\perp}=\mathcal{L}\cap
\K(\mathcal{K}^{\perp}\twoheadrightarrow^{\perp}\mathcal{L})$.
\begin{proof}
We show that $^{\perp}\mathcal{L}=\mathcal{K}\cap
\Co(\mathcal{K}^{\perp}\rightarrowtail^{\perp}\mathcal{L})$.
Let $X\in ^{\perp}\mathcal{L}$; then the short exact sequence 
$0\to X=X$ shows that $X\in \Co(\mathcal{K}^{\perp}
\rightarrowtail^{\perp}\mathcal{L})$. Conversely, if
$X\in\mathcal{K}\cap\Co(\mathcal{K}^{\perp}
\rightarrowtail^{\perp}\mathcal{L})$ then
there is a short exact sequence $A\rightarrowtail 
B\twoheadrightarrow X$ with $A\in \mathcal{K}^{\perp}$
and $B\in^{\perp}\mathcal{L}$. It follows that $X$ is a 
retract of $B$, hence $X\in^{\perp}\mathcal{L}$.
\end{proof}

\subsubsection{} \label{sec:1.2.7}
A class of objects of an abelian category is \emph{left exact} 
if it is closed under kernels of epimorphisms between its objects and 
\emph{right exact} if it is left exact in the opposite category. One says that
a cotorsion theory is \emph{left} (\emph{right}) \emph{exact} if its left (right) 
class is left (right) exact and \emph{exact} if it is left and right exact.

In every abelian category the minimal and maximal cotorsion theories are exact.

\subsubsection{} \label{sec:1.2.8}
A commutative square diagram
\begin{displaymath}
\xymatrix{
{\bullet}\ar[r] \ar[d]&{\bullet} \ar[d]\\
{\bullet} \ar[r]&{\bullet}
}
\end{displaymath}
in an abelian category is \emph{exact} if it is a pullback and a pushout.
Given a class $\mathcal{K}$ of objects of $\mathcal{A}$, the following are 
equivalent:

(1) the class $\mathcal{K}^{\perp}$ is right exact ;

(2) for every exact square 
\begin{displaymath}
\xymatrix{
{A_{1}}\ar[r] \ar[d]&{A_{2}} \ar[d]\\
{A_{3}} \ar[r]&{A_{4}}
}
\end{displaymath} 
with $A_{1},A_{2},A_{3}\in \mathcal{K}^{\perp}$ 
we have $A_{4}\in \mathcal{K}^{\perp}$;

(3) for every exact sequences 
$Y'\overset{u}\rightarrowtail Y\twoheadrightarrow Y''$ and
$A\overset{i}\rightarrowtail B\twoheadrightarrow X$
with $Y',Y\in \mathcal{K}^{\perp},X\in\mathcal{K}$, the 
commutative diagram
\begin{displaymath}
\xymatrix{
{\Co\mathcal{A}(B,u)}\ar@{>->}[r] \ar[d]&{\mathcal{A}(B,Y'')} 
\ar[d]^{\bar{i}}\\
{\Co\mathcal{A}(A,u)} \ar@{>->}[r]&{\mathcal{A}(A,Y'')}
}
\end{displaymath} 
is exact.
\begin{proof}
(1)$\Leftrightarrow$(2) This follows from the exact sequence
$A_{1}\rightarrowtail A_{2}\oplus A_{3}\twoheadrightarrow A_{4}$
and \ref{sec:1.2.1}(5). (1)$\Rightarrow$(3) By \ref{sec:1.2.1}(4)
the diagram in (3) is a pullback. Since $Y''\in\mathcal{K}^{\perp}$ 
the morphism $\bar{i}$ is an epimorphism, therefore the diagram 
is a pushout. (3)$\Rightarrow$(1) Since $Y\in\mathcal{K}^{\perp}$ 
the morphism $\Co\mathcal{A}(B,u)\to\Co\mathcal{A}(A,u)$ 
is an epimorphism, therefore so is $\bar{i}$.
\end{proof}

\subsubsection{} \label{sec:1.2.9}
Let $(\mathcal{K},\mathcal{K}^{\perp})$ be a cotorsion theory in $\mathcal{A}$. 
One says that $\mathcal{A}$ has \emph{enough} $\mathcal{K}$ \emph{objects}
if for each $A \in Ob(\mathcal{A})$ there is an epimorphism $X\to A$ with $X\in \mathcal{K}$.
One says that $\mathcal{A}$ has \emph{enough} $\mathcal{K}^{\perp}$ \emph{objects}
if $\mathcal{A}^{op}$ has enough $\mathcal{K}^{\perp}$ objects.

\subsubsection{} \label{sec:1.2.10}
\cite[Lemma 2.6]{YD} Let $(\mathcal{X},\mathcal{Y})$ be a cotorsion theory in $\mathcal{A}$. 
If $\mathcal{X}$ is left exact and $\mathcal{A}$ has enough $\mathcal{X}$ objects 
then $\mathcal{Y}$ is right exact.

\begin{proof}
Let $Y'\rightarrowtail Y\overset{v}\twoheadrightarrow Y''$ be an exact sequence 
with $Y',Y\in \mathcal{Y}$ and $A\overset{u}\rightarrowtail B\twoheadrightarrow X$
an exact sequence with $X\in \mathcal{X}$. Given a morphism $f:A\to Y''$, we show that 
there is $l:B\to Y''$ such that $lu=f$. Form the commutative diagram
\[
\xymatrix{
{A} \ar@{>->}[r]^{u} \ar[d]_{f}& {B}\ar@{->>}[r]\ar[d]& {X}\ar@{=}[d]\\
{Y''} \ar@{>->}[r]^{u'} & {PO}\ar@{->>}[r]^{v'} & {X}
}
\]
where $PO$ means pushout and choose an epimorphism $g:X'\to PO$ with $X'\in\mathcal{X}$.
Form the commutative diagram
\[
\xymatrix{
{Y'}\ar@{>->}[r] \ar@{=}[d]& {PB'}\ar@{->>}[r]\ar[d]& 
{PB} \ar@{>->}[r]^{u''} \ar[d]_{g'}& {X'}\ar@{->>}[r]^{v''}\ar@{->>}[d]_{g}& {X}\ar@{=}[d]\\
{Y'}\ar@{>->}[r] & {Y}\ar@{->>}[r]^{v} & {Y''} \ar@{>->}[r]^{u'} & {PO}\ar@{->>}[r]^{v'} & {X}
}
\]
where $PB,PB'$ are pullbacks. By assumption we have $PB\in \mathcal{X}$, so there is $h:PB\to Y$
such that $vh=g'$. Form the commutative solid arrows diagram
\[
\xymatrix{
{PB} \ar@{>->}[r]^{u''} \ar[d]_{h}& {X'}\ar@{->>}[r]^{v''}\ar[d]_{h'}& {X}\ar@{=}[d]\\
{Y} \ar@{>->}[r]^{u'''} \ar[d]_{v}& {PO'}\ar@{->>}[r]^{v'''} \ar@{..>}[d]^{w}& {X}\ar@{=}[d]\\
{Y''}\ar@{>->}[r]^{u'} & {PO}\ar@{->>}[r]^{v'}& {X}
}
\]
Since $gu''=u'vh$ there is a unique $w:PO'\to PO$ such that $wh'=g$ and $u'v=wu'''$. 
From the universal property of $PO'$ it follows that $v'w=v'''$. Since $Y\in\mathcal{Y}$
the exact sequence $(u''',v''')$ splits; there is then $s:X\to PO$ such that $v's=1_{X}$.
It follows that $u'$ has a retract, therefore there is $l:B\to Y''$ such that $lu=f$.
\end{proof}

\subsubsection{} \label{sec:1.2.11}
\cite[Theorem 1.1]{Gi3} Let $\mathcal{K},\mathcal{L}$ be two classes of objects of $\mathcal{A}$. 
Suppose that

(1) $\mathcal{K}^{\perp}$ is right exact;

(2) $\mathcal{K}^{\perp}\subset (^{\perp}\mathcal{L})^{\perp}$;

(3) $^{\perp}\mathcal{L}\subset \K(\mathcal{K}^{\perp}\twoheadrightarrow^{\perp}\mathcal{L})$.

Then $\K(\mathcal{K}^{\perp}\twoheadrightarrow^{\perp}\mathcal{L})$ is right exact 
and $\Co(\mathcal{K}^{\perp}\rightarrowtail^{\perp}\mathcal{L})\subset
\K(\mathcal{K}^{\perp}\twoheadrightarrow^{\perp}\mathcal{L})$.
If $\Co(\mathcal{K}^{\perp}\rightarrowtail^{\perp}\mathcal{L})=
\K(\mathcal{K}^{\perp}\twoheadrightarrow^{\perp}\mathcal{L})$ then
$\K(\mathcal{K}^{\perp}\twoheadrightarrow^{\perp}\mathcal{L})$ is thick.
\begin{proof}
We first consider three preparatory steps.
\paragraph{Step 1} Consider a short exact sequence $Y\rightarrowtail \bullet\twoheadrightarrow X$
with $Y\in\mathcal{K}^{\perp}$ and $X\in^{\perp}\mathcal{L}$. By (2) we have a short exact
sequence $X\rightarrowtail \bullet\twoheadrightarrow Y$. By (3) we have a short exact sequence 
$X\rightarrowtail Y'\twoheadrightarrow X'$ with $Y'\in\mathcal{K}^{\perp}$ and $X'\in^{\perp}\mathcal{L}$.
The commutative diagram
\[
\xymatrix{
{X} \ar@{>->}[r] \ar[d]& {\bullet}\ar@{->>}[r]\ar[d]& {Y}\ar@{=}[d]\\
{Y'} \ar@{>->}[r] \ar@{->>}[d] & {PO}\ar@{->>}[r] \ar@{->>}[d]& {Y}\\
{X'}\ar@{=}[r] & {X'}
}
\]
where $PO$ means pushout, shows that $PO\in\mathcal{K}^{\perp}$, hence 
$\bullet\in \K(\mathcal{K}^{\perp}\twoheadrightarrow^{\perp}\mathcal{L})$.
\paragraph{Step 2} Consider a short exact sequence $\bullet\rightarrowtail A\twoheadrightarrow X$
with $A\in \K(\mathcal{K}^{\perp}\twoheadrightarrow^{\perp}\mathcal{L})$ and $X\in^{\perp}\mathcal{L}$.
By (3) we have a short exact sequence 
$A\rightarrowtail Y'\twoheadrightarrow X'$ with $Y'\in\mathcal{K}^{\perp}$ and $X'\in^{\perp}\mathcal{L}$.
The commutative diagram
\[
\xymatrix{
{\bullet} \ar@{>->}[r] \ar@{=}[d]& {A}\ar@{->>}[r]\ar@{>->}[d]& {X}\ar@{>->}[d]\\
{\bullet} \ar@{>->}[r] & {Y'}\ar@{->>}[r] \ar@{->>}[d]& {PO}\ar@{->>}[d]\\
& {X'}\ar@{=}[r] & {X'}
}
\]
where $PO$ means pushout, shows that $PO\in{\perp}\mathcal{L}$, hence 
$\bullet\in \K(\mathcal{K}^{\perp}\twoheadrightarrow^{\perp}\mathcal{L})$.
\paragraph{Step 3} Consider a short exact sequence $Y\rightarrowtail A\twoheadrightarrow \bullet$
with $Y\in\mathcal{K}^{\perp}$ and $A\in\K(\mathcal{K}^{\perp}\twoheadrightarrow^{\perp}\mathcal{L})$.
By (3) we have a short exact sequence 
$A\rightarrowtail Y'\twoheadrightarrow X'$ with $Y'\in\mathcal{K}^{\perp}$ and $X'\in^{\perp}\mathcal{L}$.
Form commutative diagram
\[
\xymatrix{
{Y} \ar@{>->}[r] \ar@{=}[d]& {A}\ar@{->>}[r]\ar@{>->}[d]& {\bullet}\ar@{>->}[d]\\
{Y} \ar@{>->}[r] & {Y'}\ar@{->>}[r] \ar@{->>}[d]& {PO}\ar@{->>}[d]\\
& {X'}\ar@{=}[r] & {X'}
}
\]
where $PO$ means pushout. By (1) $PO\in\mathcal{K}^{\perp}$, hence 
$\bullet\in \K(\mathcal{K}^{\perp}\twoheadrightarrow^{\perp}\mathcal{L})$.

\paragraph{} Consider now a short exact sequence $A'\rightarrowtail A\twoheadrightarrow \bullet$ with 
$A',A\in\K(\mathcal{K}^{\perp}\twoheadrightarrow^{\perp}\mathcal{L})$.
By (3) we have short exact sequences
$A'\rightarrowtail Y'\twoheadrightarrow X'$ and 
$A\rightarrowtail Y\twoheadrightarrow X$ with $Y',Y\in\mathcal{K}^{\perp}$ and $X',X\in^{\perp}\mathcal{L}$.
Form the commutative diagram
\[
\xymatrix{
{A'} \ar@{>->}[rrrr] \ar[dr]& & &&{A}\ar@{->>}[rrrr]\ar[dr]\ar@{>->}'[d]'[dd][ddd] & &&&{\bullet}\ar@{=}[dr]\\
&{Y'} \ar@{>->}[rrrr]\ar@{->>}[dr] & &&&{PO}\ar@{->>}[rrrr]\ar@{->>}[dr]\ar@{>->}[ddd]&& &&{\bullet}\\
&&{X'}\ar@{=}[rrrr] &&&&{X'}\ar@{=}[ddd]\\
&&&& {Y}\ar@{->>}[ddd]\ar@{>->}[dr]\\
&&&&&{PO'}\ar@{->>}[dr]\ar@{->>}[ddd]\\
&&&&&&{X'}\\
&&&&{X}\ar@{=}[dr]\\
&&&&&{X}
}
\]
where $PO,PO'$ are pushouts. Then $PO'\in\K(\mathcal{K}^{\perp}\twoheadrightarrow^{\perp}\mathcal{L})$
by step 1, $PO\in\K(\mathcal{K}^{\perp}\twoheadrightarrow^{\perp}\mathcal{L})$ by step 2 and 
$\bullet\in\K(\mathcal{K}^{\perp}\twoheadrightarrow^{\perp}\mathcal{L})$ by step 3.
\paragraph{} We now show that $\Co(\mathcal{K}^{\perp}\twoheadrightarrow^{\perp}\mathcal{L})\subset
\K(\mathcal{K}^{\perp}\twoheadrightarrow^{\perp}\mathcal{L})$. Let 
$Y\rightarrowtail X\twoheadrightarrow \bullet$ be an exact sequence with 
$Y\in\mathcal{K}^{\perp}$ and $X\in^{\perp}\mathcal{L}$. By (3) we have a short exact sequence 
$X\rightarrowtail Y'\twoheadrightarrow X'$ with $Y'\in\mathcal{K}^{\perp}$ and $X'\in^{\perp}\mathcal{L}$.
Form the commutative diagram
\[
\xymatrix{
{Y} \ar@{>->}[r] \ar@{=}[d]& {X}\ar@{->>}[r]\ar@{>->}[d]& {\bullet}\ar@{>->}[d]\\
{Y} \ar@{>->}[r] & {Y'}\ar@{->>}[r] \ar@{->>}[d]& {PO}\ar@{->>}[d]\\
& {X'}\ar@{=}[r] & {X'}
}
\]
where $PO$ means pushout. By (1) $PO\in\mathcal{K}^{\perp}$, hence 
$\bullet\in \K(\mathcal{K}^{\perp}\twoheadrightarrow^{\perp}\mathcal{L})$.
\paragraph{} We next show that $\K(\mathcal{K}^{\perp}\twoheadrightarrow^{\perp}\mathcal{L})$ 
is left exact. Let $\bullet\rightarrowtail A\twoheadrightarrow A''$ be an exact sequence with 
$A,A''\in\K(\mathcal{K}^{\perp}\twoheadrightarrow^{\perp}\mathcal{L})$.
By (3) we have short exact sequences $A\rightarrowtail Y\twoheadrightarrow X$ and 
$A''\rightarrowtail Y''\twoheadrightarrow X''$ with $Y,Y''\in\mathcal{K}^{\perp}$ and $X,X''\in^{\perp}\mathcal{L}$.
Form the commutative diagram
\[
\xymatrix{
{\bullet}\ar@{>->}[rrrr]\ar@{=}[dr] &&&& {A}\ar@{->>}[rrrr] \ar@{>->}[dr]
&&&&{A''}\ar@{>->}[dr]\ar@{>->}'[d]'[dd][ddd]\\
&{\bullet} \ar@{>->}[rrrr]&&&&{Y} \ar@{->>}[rrrr]\ar@{->>}[dr] &&&&{PO}\ar@{->>}[dr]\ar@{>->}[ddd]\\
&&&&&&{X}\ar@{=}[rrrr] &&&&{X}\ar@{=}[ddd]\\
&&&&&&&& {Y''}\ar@{->>}[ddd]\ar@{>->}[dr]\\
&&&&&&&&&{PO'}\ar@{->>}[dr]\ar@{->>}[ddd]\\
&&&&&&&&&&{X}\\
&&&&&&&& {X''}\ar@{=}[dr]\\
&&&&&&&&&{X''}
}
\]
where $PO,PO'$ are pushouts. Then $PO'\in\K(\mathcal{K}^{\perp}\twoheadrightarrow^{\perp}\mathcal{L})$
by step 1 and $PO\in\K(\mathcal{K}^{\perp}\twoheadrightarrow^{\perp}\mathcal{L})$ by step 2.
We can find an exact sequence $Y'\rightarrowtail X'\twoheadrightarrow PO$ with 
$Y'\in\mathcal{K}^{\perp}$ and $X'\in^{\perp}\mathcal{L}$. Form the diagram
\[
\xymatrix{
&{Y'}\ar@{>->}[d]\ar@{=}[r]&{Y'}\ar@{>->}[d]\\
{\bullet} \ar@{>->}[r] \ar@{=}[d]& {PB}\ar@{->>}[r]\ar@{->>}[d]& {X'}\ar@{->>}[d]\\
{\bullet} \ar@{>->}[r] & {Y}\ar@{->>}[r] & {PO}
}
\]
where $PB$ means pullback. It follows that $PB\in\mathcal{K}^{\perp}$, hence
$\bullet\in\K(\mathcal{K}^{\perp}\twoheadrightarrow^{\perp}\mathcal{L})$.
\paragraph{} We show that $\K(\mathcal{K}^{\perp}\twoheadrightarrow^{\perp}\mathcal{L})$ 
is closed under extensions. Let $A'\rightarrowtail \bullet\twoheadrightarrow A''$ be an exact 
sequence with $A',A''\in\K(\mathcal{K}^{\perp}\twoheadrightarrow^{\perp}\mathcal{L})$.
We can find short exact sequences $A'\rightarrowtail Y'\twoheadrightarrow X'$ and 
$A''\rightarrowtail Y''\twoheadrightarrow X''$ with $Y',Y''\in\mathcal{K}^{\perp}$ 
and $X',X''\in^{\perp}\mathcal{L}$. Form the commutative diagram
\[
\xymatrix{
&&&&&{Y''}\ar@{=}[rrrr]\ar@{>->}[dd] &&&&{Y''}\ar@{>->}[dd]\\
\\
&{Y'}\ar@{=}[dd]\ar@{>->}[rrrr] & &&&{PB}\ar@{->>}[rrrr]\ar@{->>}[dd]&& &&{X''}\ar@{->>}[dd]\\
{A'} \ar@{>->}[rrrr] \ar[dr]& & &&{\bullet}\ar@{->>}[rrrr]\ar[dr]& &&&{A''}\ar@{=}[dr]\\
&{Y'} \ar@{>->}[rrrr]\ar@{->>}[dr] & &&&{PO}\ar@{->>}[rrrr]\ar@{->>}[dr]&& &&{A''}\\
&&{X'}\ar@{=}[rrrr] &&&&{X'}
}
\]
Then $PB\in\K(\mathcal{K}^{\perp}\twoheadrightarrow^{\perp}\mathcal{L})$
by step 1, $PO\in\K(\mathcal{K}^{\perp}\twoheadrightarrow^{\perp}\mathcal{L})$ by step 3
and $\bullet\in\K(\mathcal{K}^{\perp}\twoheadrightarrow^{\perp}\mathcal{L})$ by step 2.
\end{proof}

\subsubsection{} \label{sec:1.2.12}
Let $(\mathcal{K},\mathcal{K}^{\perp})$ be a right exact cotorsion theory in
$\mathcal{A}$. If $\mathcal{K}\subset\K(\mathcal{K}^{\perp}\twoheadrightarrow\mathcal{K})$
then $\K(\mathcal{K}^{\perp}\twoheadrightarrow\mathcal{K})$ is right exact and 
$\Co(\mathcal{K}^{\perp}\rightarrowtail\mathcal{K})\subset
\K(\mathcal{K}^{\perp}\twoheadrightarrow\mathcal{K})$.

\begin{proof}
This is a consequence of \ref{sec:1.2.11} with $\mathcal{L}=\mathcal{K}^{\perp}$.
\end{proof}

\subsubsection{} \label{sec:1.2.13}
\cite[Lemma 12.3.3]{Pe} 
 Let $\mathcal{K},\mathcal{L}$ be two classes of objects of $\mathcal{A}$ with 
$\mathcal{K}\subset\mathcal{L}$.

(1) $^{\perp}\mathcal{K}\cap\mathcal{K}\subset ^{\perp}\mathcal{L}$ and
$^{\perp}\mathcal{L}\cap\mathcal{L}\subset \mathcal{K}$ if and only if 
$^{\perp}\mathcal{K}\cap\mathcal{K}=^{\perp}\mathcal{L}\cap\mathcal{L}$.

(2) Suppose $^{\perp}\mathcal{L}\subset \mathcal{K}$, $\mathcal{K}$ is left exact
and $\mathcal{A}$ has enough $^{\perp}\mathcal{L}$ objects, meaning that for every
$A\in Ob(\mathcal{A})$ there is an epimorphism $X\to A$ with $X\in^{\perp}\mathcal{L}$.
Then $^{\perp}\mathcal{K}\cap\mathcal{K}=^{\perp}\mathcal{L}\cap\mathcal{L}$.
 
(2bis) Suppose $\mathcal{L}^{\perp}\subset \mathcal{K}$, $\mathcal{K}$ is right exact
and $\mathcal{A}$ has enough $\mathcal{L}^{\perp}$ objects, meaning that for every
$A\in Ob(\mathcal{A})$ there is a monomorphism $A\to X$ with $X\in\mathcal{L}^{\perp}$.
Then $\mathcal{K}\cap\mathcal{K}^{\perp}=\mathcal{L}\cap\mathcal{L}^{\perp}$.
\begin{proof}
(2) Using \ref{sec:1.2.1}(2) we have $^{\perp}\mathcal{L}\cap\mathcal{L}\subset
^{\perp}\mathcal{K}\cap\mathcal{K}$. Let $X''\in ^{\perp}\mathcal{K}\cap\mathcal{K}$.
We can find an epimorphism $f:X\to X''$ with $X\in ^{\perp}\mathcal{L}$; then $\K f\in
\mathcal{K}$ and so by \ref{sec:1.2.1}(1bis) $X''$ is a retract of $X$.
\end{proof}

\subsubsection{} \label{sec:1.2.14}
Let $(\mathcal{X},\mathcal{Y})$ be a cotorsion theory in $\mathcal{A}$.
We denote by $\underline{\mathcal{Y}}$ the class of objects $Y\in \mathcal{Y}$ 
such that every exact sequence $Y\rightarrowtail A\rightarrow X$ with $X\in 
\mathcal{X}$ is left split (\ref{sec:1.1}). We denote by $\underline{\mathcal{X}}$ 
the class of objects $X\in \mathcal{X}$ such that every exact sequence 
$Y\rightarrow A\twoheadrightarrow X$ with $Y\in \mathcal{Y}$ is right split.
The class $\underline{\mathcal{Y}}$ contains the injective objects of 
$\mathcal{A}$ and $\underline{\mathcal{X}}$ contains the projective 
objects.

\subsubsection{} \label{sec:1.2.15}
Let $(\mathcal{X},\mathcal{Y})$ be a cotorsion theory in $\mathcal{A}$.
The class $\mathcal{Y}$ is closed under quotients if and only if
$\underline{\mathcal{X}}=\mathcal{X}$. Dually, the class $\mathcal{X}$ 
is closed under subobjects if and only if $\underline{\mathcal{Y}}=\mathcal{Y}$.
\begin{proof}
Suppose $\mathcal{Y}$ is closed under quotients. Let $X\in\mathcal{X}$ and 
$Y\overset{u}\rightarrow A\overset{v}\twoheadrightarrow X$ be an exact sequence 
with $Y\in \mathcal{Y}$. Factor $u$ as $u=ip$, where $p:Y\twoheadrightarrow \im u$ 
and $i:\im u \rightarrowtail A$. Then $\im u\in \mathcal{Y}$, hence the exact sequence
$\im u\overset{i}\rightarrowtail A\overset{v}\twoheadrightarrow X$ splits, so $v$ has 
a section. 
Conversely, let $f:Y\twoheadrightarrow Y'$
with $Y\in\mathcal{Y}$ and let $A\overset{u}\rightarrowtail B\overset{v}
\twoheadrightarrow X$ be an exact sequence with $X\in\mathcal{X}$. 
Let $g:A\to Y'$. Form the commutative diagram
\[
\xymatrix{
& {A} \ar@{>->}[r]^{u} \ar[d]_{g}& {B}\ar@{->>}[r]^{v}\ar[d]& {X}\ar@{=}[d]\\
{Y}\ar@{->>}[r]^{f}& {Y'} \ar@{>->}[r] & {PO}\ar@{->>}[r] & {X}
}
\]
where $PO$ means pushout. Then we have the exact sequence 
$Y\to PO\twoheadrightarrow X$ which is right split by assumption, 
hence $Y'\rightarrowtail PO$ has a retraction.
It follows that the morphism $\bar{u}:\mathcal{A}(B,Y')\to \mathcal{A}(A,Y')$ 
is an epimorphism, hence $Y'\in\mathcal{Y}$.
\end{proof}

\subsubsection{} \label{sec:1.2.16}
Suppose that $\mathcal{A}$ is a Grothendieck category with a generator 
$U$. Let $Sub(U)$ be the set of subobjects of $U$ and let $\mathcal{G}$ 
be a subset of $Sub(U)$ satisfying the following conditions:

(a) for every commutative diagram
\begin{displaymath}
\xymatrix{
{R}\ar[rr] \ar[dr]_{r}&&{S} \ar[dl]^{s}\\
&{U} 
}
\end{displaymath}
in $\mathcal{A}$, if $r:R\to U$ belongs to $\mathcal{G}$ then so does
$s:S\to U$.

(b) if $s:S\to U$ belongs to $\mathcal{G}$ then for all morphisms $f:U\to U/S$,
$k:\K f\to U$ belongs to $\mathcal{G}$.

Recall \cite{Ga} that every class of objects of $\mathcal{A}$ that is closed under 
subquotients gives rise to a subset $\mathcal{G}$ of $Sub(U)$ satisfying (a) and (b).
Also, when $\mathcal{A}$ is the category of left modules over a ring $R$, condition
(b) means: for every left ideal $I$ of $R$ that belongs to $\mathcal{G}$ and every 
$r\in R$, the left ideal $(I:r)=\{x\in R, xr\in I\}$ belongs to $\mathcal{G}$. 
Define now
\[
\mathcal{G}_{t}=\{X\in Ob(\mathcal{A}), \ k:\K f\to U {\rm \ belongs\ to} \ U 
\ {\rm for\ all}\ f:U\to X\}
\]
The cotorsion theory $(^{\perp}(\mathcal{G}_{t}^{\perp}),
\mathcal{G}_{t}^{\perp})$ is cogenerated by a set.

\begin{proof}
We put $\mathcal{S}=\{U/S, \ s:S\to U {\rm \ belongs\ to} \ \mathcal{G}\}$;
by (b) we have $\mathcal{S}\subset \mathcal{G}_{t}$ 
and we show that $\mathcal{S}^{\perp}\subset \mathcal{G}_{t}^{\perp}$.
Let $Y\overset{u}\rightarrowtail A\overset{v}\twoheadrightarrow X$ be an 
exact sequence with $Y\in \mathcal{S}^{\perp},X\in\mathcal{G}_{t}$. 
\paragraph{Step 1} Let $s:S\rightarrowtail X,l:S\to A$ be morphisms such
that $vl=s$. Suppose that $v$ is not an isomorphism. Then there is 
$\alpha:X\to B$ such that $\alpha s=0$ and $\alpha\neq 0$. In turn,
there is $f:U\to X$ such that $\im f\neq 0$. Let $U\twoheadrightarrow
\im f\overset{r}\rightarrowtail X$ be the image factorization of $f$.
Form the commutative diagram
\begin{displaymath}
\xymatrix{
{\K f}\ar@{>->}[d]\ar@{=}[r]&{\K f}\ar@{>->}[d]\\
{PB}\ar@{>->}[r] \ar@{->>}[d]&{U} \ar@{->>}[d]
\ar@{->>}[r]&{U/PB}\ar@{=}[d]\\
{S\cap \im f}\ar@{>->}[d]\ar@{>->}[r]&{\im f}
\ar@{->>}[r]\ar@{>->}[d]^{r}&{\im f/S\cap \im f}\\
{S}\ar@{>->}[r]^{s}&{X}
}
\end{displaymath}
where $PB$ means pullback. By (a) we have that 
$PB\to U$ belongs to $\mathcal{G}$. Since 
$\im f/S\cap \im f\cong S\cup \im f/S$, by (b) we
have that  $S\cup \im f/S\in \mathcal{S}$. Consider 
the commutative diagram
\begin{displaymath}
\xymatrix{
{S\cap \im f}\ar[r] \ar[d]&{\im f} \ar[d]\ar[ddr]^{r}\\
{S} \ar[r]^{i}\ar[drr]_{s}&{S\cup\im f}\ar[dr]^{s'}\\
&&{X}
}
\end{displaymath} 
where $s'$ is defined by the universal property of 
pushout. We obtain then the commutative solid arrows
diagram 
\begin{displaymath}
\xymatrix{
&{Y}\ar@{>->}[d]^{u}\\
{S}\ar[r]^{l}\ar@{>->}[d]_{i}&{A} \ar@{->>}[d]^{v}\\
{S\cup\im f}\ar@{->>}[d]\ar[r]_{s'}\ar@{..>}[ur]^{l'}&{X}\\
{S\cup \im f/S}
}
\end{displaymath}
By \ref{sec:1.2.1}(4) there is $l':S\cup\im f\to A$ such that
$l'i=l$ and $vl'=s'$.
\paragraph{Step 2} Let $P$ be the set of all ordered pairs
$((S,s),l)$, where $s:S\rightarrowtail X,l:S\to A$ and $vl=s$.
The set $P$ is nonempty since $((0,0),0)\in P$. The relation
$((S,s),l)\leqslant ((S',s'),l')$ if there is $i:S\to S'$ such that
$l'i=l$ defines an ordering on $P$. We claim that $P$ is inductive
with respect to this ordering. Let $((S_{i},s_{i}),l_{i})_{i\in  I}$
be a totally ordered set of elements of $P$. The element
$((\underset{i\in I}\sum S_{i},s),l)$, where 
$\underset{i\in I}\bigoplus S_{i}\twoheadrightarrow
\underset{i\in I}\sum S_{i}\overset{s}\rightarrowtail X$ 
is the image factorization of the natural morphism 
$\underset{i\in I}\bigoplus S_{i}\to X$ and 
$l:\underset{i\in I}\sum S_{i}\to A$ is defined by the 
universal property of colimits, is an upper bound. The
the claim is proved. Therefore by Zorn's lemma there
exists a maximal element $((S,s),l)$ in $P$. If $s$ is
an isomorphism then $v$ has a section and we are done.
If not, from Step 1 there is an element $((S\cup \im f,s'),l')$
of $P$ such that $((S,s),l)< ((S\cup \im f,s'),l')$,
contradicting the maximality of $((S,s),l)$.
\end{proof}

\section{Cotorsion theories and adjoint functors}

\subsection{} \label{sec:2.1} 
Let $\mathcal{A},\mathcal{A'}$ be abelian categories and 
\begin{displaymath}
\xymatrix{
F:{\mathcal{A}}\ar@<0.8ex>[r]&{\mathcal{A'}}:G\ar@<0.4ex>[l]\\
}
\end{displaymath}
an adjoint pair. Let $\mathcal{K}$ be a class of objects of $\mathcal{A}$ 
and $\mathcal{K}'$ a class of objects of $\mathcal{A}'$.

\subsubsection{} \label{sec:2.1.1}

(1) Suppose that $\mathcal{A}$ has enough $\mathcal{K}$ objects, meaning that
for every $A\in Ob(\mathcal{A})$ there is an epimorphism $X\to A$ with $X\in\mathcal{K}$.
If $F(\mathcal{K})\subset \mathcal{K}'$ then $G$ preserves epimorphisms with kernel in 
$\mathcal{K}'^{\perp}$. 

(1bis) Suppose that $\mathcal{A}'$ has enough 
$\mathcal{K'}^{\perp}$ objects, meaning that for every $A'\in Ob(\mathcal{A}')$ 
there is a monomorphism $A'\to Y'$ with $Y'\in\mathcal{K}'^{\perp}$. If
$G(\mathcal{K'}^{\perp})\subset \mathcal{K}^{\perp}$ then 
$F$ preserves monomorphisms with cokernel in $\mathcal{K}$.

(2) If $G$ preserves epimorphisms with kernel in $\mathcal{K'}$ and 
$G(\mathcal{K'})\subset \mathcal{K}$ then $F({}^{\perp}\mathcal{K})\subset {}^{\perp}(\mathcal{K'})$.

(2bis) If $F$ preserves monomorphisms with cokernel in $\mathcal{K}$ 
and $F(\mathcal{K})\subset \mathcal{K}'$ then $G(\mathcal{K'}^{\perp})\subset \mathcal{K}^{\perp}$.

(3) Suppose that $(\mathcal{K},\mathcal{K}^{\perp})$ and $(\mathcal{K}',\mathcal{K}'^{\perp})$
are cotorsion theories, $F$ preserves monomorphisms with cokernel in $\mathcal{K}$
and $G$ preserves epimorphisms with kernel in $\mathcal{K'}^{\perp}$.
Then $F(\mathcal{K})\subset\mathcal{K}'\Leftrightarrow 
G(\mathcal{K'}^{\perp})\subset \mathcal{K}^{\perp}$.
\begin{proof}
(1) Let $Y'\rightarrowtail A'\twoheadrightarrow B'$ be an exact sequence 
with $Y'\in \mathcal{K'}^{\perp}$. We can find an epimorphism $f:X\to G(B')$ with $X\in \mathcal{K}$.
Form the commutative solid arrows diagram
\begin{displaymath}
\xymatrix{
{Y'}\ar@{>->}[r] \ar@{=}[d] & {PB}\ar@{->>}[r] \ar[d] & {F(X)} \ar[d]^{f^{\#}} \ar@{..>}[dl]\\
{Y'} \ar@{>->}[r] & {A'}\ar@{->>}[r] & {B'}
}
\end{displaymath}
where $f^{\#}$ is the adjoint transpose of $f$ and $PB$ means pullback. 
Since $F(X)\in \mathcal{K}'$ the top exact sequence splits, therefore 
there is a dotted arrow that makes the lower triangular diagram involving 
$f^{\#}$ commute. By adjunction, there is a dotted arrow in the 
diagram
\begin{displaymath}
\xymatrix{
& & {X} \ar@{->>}[d]^{f} \ar@{..>}[dl]\\
{G(Y')} \ar@{>->}[r] & {G(A')}\ar[r] & {G(B')}
}
\end{displaymath}
that makes the triangular diagram commute, hence the 
morphism $G(A')\to G(B')$ is an epimorphism.

(2) Let $X\in{}^{\perp}\mathcal{K}$ and $X'\rightarrowtail A'\twoheadrightarrow B'$ be an exact 
sequence with $X'\in \mathcal{K'}$. Applying the functor $\mathcal{A}'(F(X),-)$ to the previous 
exact sequence and using adjunction we obtain the exact sequence
\begin{displaymath}
\mathcal{A}(X,G(X'))\rightarrowtail \mathcal{A}(X,G(A'))\rightarrow \mathcal{A}(X,G(B'))
\end{displaymath}
By hypothesis the right morphism of the previous exact sequence is an epimorphism, 
therefore $F(X)\in{}^{\perp}\mathcal{K}'$.

(3) This follows from parts 2 and 2bis.
\end{proof}
The previous result contains the familiar statements

(1bis) if the target category of a left adjoint has enough injective objects and the right adjoint
preserves injective objects then the left adjoint is exact, and

(2bis) if a left adjoint is exact then the right adjoint preserves injective objects.

\subsubsection{} \label{sec:2.1.2}

Suppose that $(\mathcal{K},\mathcal{K}^{\perp})$ and $(\mathcal{K}',\mathcal{K}'^{\perp})$
are cotorsion theories.

(1) If $F$ preserves monomorphisms with cokernel in $F^{-1}(\mathcal{K'})$ and 
$G$ preserves epimorphisms with kernel in $\mathcal{K}'^{\perp}$ then 
$(F^{-1}(\mathcal{K}'),F^{-1}(\mathcal{K}')^{\perp})$ is a cotorsion theory in $\mathcal{A}$.
If $(F^{-1}(\mathcal{K}'),F^{-1}(\mathcal{K}')^{\perp})$ is cogenerated (\ref{sec:1.2.3})
by $\mathcal{S}$ then $(\mathcal{K}',\mathcal{K}'^{\perp})$ is cogenerated by $F(\mathcal{S})$.

(1bis) If $G$ preserves epimorphisms with kernel in $G^{-1}(\mathcal{K}^{\perp})$ 
and $F$ preserves monomorphisms with cokernel in $\mathcal{K}$ then 
$({}^{\perp}G^{-1}(\mathcal{K}^{\perp}),G^{-1}(\mathcal{K}^{\perp}))$ 
is a cotorsion theory in $\mathcal{A'}$.

\begin{proof}
(1) We shall use \ref{sec:1.2.3}. We put $\mathcal{S}=\{G(Y')\}$, where $Y'\in\mathcal{K}'^{\perp}$. 
We show that $\mathcal{S}\subset F^{-1}(\mathcal{K}')^{\perp}$. Let $Y'\in\mathcal{K}'^{\perp}$ 
and $A\rightarrowtail B\twoheadrightarrow X$ be an exact sequence with $F(X)\in \mathcal{K}'$.
By hypothesis and assumption we have the exact sequence 
\begin{displaymath}
\mathcal{A}'(F(X),Y')\rightarrowtail \mathcal{A}'(F(A),Y')\twoheadrightarrow \mathcal{A}'(F(B),Y')
\end{displaymath}
By adjunction we have the exact sequence 
\begin{displaymath}
\mathcal{A}(X,G(Y'))\rightarrowtail \mathcal{A}(B,G(Y'))\twoheadrightarrow \mathcal{A}(A,G(Y'))
\end{displaymath}
which implies that $G(Y')\in F^{-1}(\mathcal{K}')^{\perp}$. We now show that 
${}^{\perp}\mathcal{S}\subset F^{-1}(\mathcal{K}')$. Let $X\in {}^{\perp}\mathcal{S}$
and let $Y'\rightarrowtail A'\twoheadrightarrow B'$ be an exact sequence with 
$Y'\in \mathcal{K}'^{\perp}$. Applying the functor $\mathcal{A}'(F(X),-)$ to the 
previous exact sequence and using adjunction we obtain the exact sequence
\begin{displaymath}
\mathcal{A}(X,G(Y'))\rightarrowtail \mathcal{A}(X,G(A'))\rightarrow \mathcal{A}(X,G(B'))
\end{displaymath}
By hypothesis and assumption on $X$ the right morphism of the previous exact sequence 
is an epimorphism, therefore $F(X)\in\mathcal{K}'$. 

Suppose now that $(F^{-1}(\mathcal{K}'),F^{-1}(\mathcal{K}')^{\perp})$ is 
cogenerated by $\mathcal{S}$. Let $X\in \mathcal{S}$ and 
$Y'\rightarrowtail A'\overset{v}\twoheadrightarrow B'$ be 
an exact sequence in $\mathcal{A}'$ with $Y'\in \mathcal{K}'^{\perp}$. 
Form the diagram with exact rows 
\[
\xymatrix{
{G(Y')} \ar@{>->}[r] \ar@{=}[d]& {PB}\ar@{->>}[r]\ar[d]& {X}\ar[d]_{\eta_{X}}\\
{G(Y')} \ar@{>->}[r] & {G(A')}\ar@{->>}[r]^{G(v)} & {GF(X)}
}
\]
where $\eta_{X}$ is the unit of the adjunction and $PB$ means pullback.
By \ref{sec:2.1.1}(2bis) we have $G(Y')\in F^{-1}(\mathcal{K}')^{\perp}$
hence the top exact sequence splits. There is then $f:X\to G(A')$ such that
$G(v)f=\eta_{X}$. The adjoint transpose $f^{\#}$ of $f$ is a section of $v$.
\end{proof}

The following observation \cite{Gi2} is ``often rediscovered folklore''.

\subsubsection{} \label{sec:2.1.3}

(1) Suppose that $F$ preserves monomorphisms with cokernel in $\mathcal{K}$. For each 
$X\in \mathcal{K}$ and $A'\in Ob(\mathcal{A'})$ we have a natural morphism 
\begin{displaymath}
\Theta_{X,A'}:\mathsf{Ext}_{\mathcal{A}}^{1}(X,G(A')) \to \mathsf{Ext}_{\mathcal{A'}}^{1}(F(X),A')
\end{displaymath}
If, moreover, $F$ is full and faithful then $\Theta_{X,A'}$ is a monomorphism.

(1bis) Suppose that $G$ preserves epimorphisms with kernel in $\mathcal{K'}$. Then for each 
$A\in Ob(\mathcal{A})$ and $X'\in \mathcal{K'}$ we have a natural morphism 
\begin{displaymath}
\Psi_{A,X'}:\mathsf{Ext}_{\mathcal{A'}}^{1}(F(A),X')\to\mathsf{Ext}_{\mathcal{A}}^{1}(A,G(X'))
\end{displaymath}
If, moreover, $G$ is full and faithful then $\Psi_{A,X'}$ is a monomorphism.

(2) Suppose that $F$ preserves monomorphisms with cokernel in $\mathcal{K}$ and 
$G$ preserves epimorphisms with kernel in $\mathcal{K'}$. For each $X\in \mathcal{K}$ and 
$X'\in \mathcal{K'}$ the morphism $\Theta_{X,X'}$ is an isomorphism with inverse $\Psi_{X,X'}$.

\subsection{} \label{sec:2.2}
Let $\mathcal{A}_{i}$ $(1\leqslant i\leqslant 3)$ be a category. 
A \emph{THC--situation} \cite{Gra}, or adjunction of two variables, 
consists of three functors
\begin{displaymath}
T:\mathcal{A}_{1}\times \mathcal{A}_{2}\to \mathcal{A}_{3}, \
H:\mathcal{A}_{2}^{op}\times \mathcal{A}_{3}\to \mathcal{A}_{1}, \
C:\mathcal{A}_{1}^{op}\times \mathcal{A}_{3}\to \mathcal{A}_{2}
\end{displaymath}
and natural isomorphisms
\begin{displaymath}
\mathcal{A}_{3}(T(X_{1},X_{2}),X_{3})\cong \mathcal{A}_{1}(X_{1},H(X_{2},X_{3}))
\cong \mathcal{A}_{2}(X_{2},C(X_{1},X_{3}))
\end{displaymath}
For each object $X_{i}$ of $\mathcal{A}_{i}$ we have then adjoint pairs
\begin{displaymath}
T(-,X_{2}):\mathcal{A}_{1}\rightleftarrows \mathcal{A}_{3}:H(X_{2},-)
\end{displaymath}
\begin{displaymath}
T(X_{1},-):\mathcal{A}_{2}\rightleftarrows \mathcal{A}_{3}:C(X_{1},-)
\end{displaymath}
\begin{displaymath}
C(-,X_{3}):\mathcal{A}_{1}\rightleftarrows \mathcal{A}_{2}^{op}:H(-,X_{3})
\end{displaymath}
An \emph{abelian THC--situation} is a THC--situation as above with 
$\mathcal{A}_{i}$ abelian. In this case $T,H$ and $C$ are biadditive, 
meaning additive in each variable separately.

\subsubsection{} \label{sec:2.2.1}
We recall \cite{Fr},\cite{Ga} that every complete and cocomplete abelian 
category $\mathcal{A}$ fits into the natural abelian THC-situation 
\begin{displaymath}
T:\mathrm{Ab}\times \mathcal{A}\to \mathcal{A}, \
\mathcal{A}(-,-):\mathcal{A}^{op}\times \mathcal{A}\to \mathrm{Ab}, \
C:\mathrm{Ab}^{op}\times \mathcal{A}\to \mathcal{A}
\end{displaymath}
where $\mathrm{Ab}$ is the category of abelian groups.
The Tensor functor can be constructed as follows. Let $M$ be an abelian 
group and $A\in Ob(\mathcal{A})$. Let $\sigma_{m}:A\to\underset{M}\bigoplus A$
and $\sigma_{m,m'}:A\to\underset{M\times M}\bigoplus A$ be coprojection
morphisms. We define 
\[
a(M,A):\underset{M\times M}\bigoplus A\to
\underset{M}\bigoplus A
\] to be the unique morphism such that 
$a(M,A)\sigma_{m,m'}=\sigma_{m+m'}-(\sigma_{m}+\sigma_{m'})$
for all $m,m'\in M$; then $T(M,A)=\Co a(M,A)$. The construction of
Cotensor is dual.

\subsubsection{} \label{sec:2.2.2}
Let 
\[
T:\mathcal{A}_{1}\times \mathcal{A}_{2}\to \mathcal{A}_{3}, \
H:\mathcal{A}_{2}^{op}\times \mathcal{A}_{3}\to \mathcal{A}_{1}, \
C:\mathcal{A}_{1}^{op}\times \mathcal{A}_{3}\to \mathcal{A}_{2}
\]
be an abelian THC--situation. Let $\mathcal{K}_{1}\subset Ob(\mathcal{A}_{1})$
be a class of objects, $(\mathcal{X}_{2},\mathcal{Y}_{2})$ 
a cotorsion theory in $\mathcal{A}_{2}$ and $(\mathcal{X}_{3},\mathcal{Y}_{3})$ 
a cotorsion theory in $\mathcal{A}_{3}$. We define
\[
T(\mathcal{K}_{1},-)^{-1}(\mathcal{X}_{3})=
\{A_{2}\in Ob(\mathcal{A}_{2}),\ T(X_{1},A_{2})\in \mathcal{X}_{3} {\rm \ for\ all}\  X_{1}
\in\mathcal{K}_{1}\}
\]

Suppose that

(1) $(\mathcal{X}_{2},\mathcal{Y}_{2})$ is generated 
by $\mathcal{S}_{2}$ and $(\mathcal{X}_{3},\mathcal{Y}_{3})$ 
is generated by $\mathcal{S}_{3}$;

(2) for each object $X_{1}$ of $\mathcal{K}_{1}$ the functor $T(X_{1},-)$ 
preserves monomorphisms with cokernel in $\mathcal{X}_{2}$ and the 
functor $C(X_{1},-)$ preserves epimorphisms with kernel in $\mathcal{S}_{3}$.

Then $(T(\mathcal{K}_{1},-)^{-1}(\mathcal{X}_{3})\cap \mathcal{X}_{2},
(T(\mathcal{K}_{1},-)^{-1}(\mathcal{X}_{3})\cap \mathcal{X}_{2})^{\perp})$ 
is a cotorsion theory. A sufficient condition for the functor $C(X_{1},-)$ 
to preserve the epimorphisms with kernel in $\mathcal{S}_{3}$ is that the category
$\mathcal{A}_{2}$ has enough $T(\mathcal{K}_{1},-)^{-1}(\mathcal{X}_{3})\cap \mathcal{X}_{2}$
objects, meaning that for every $A_{2}\in Ob(\mathcal{A}_{2})$ there is an epimorphism
$X_{2}\to A_{2}$ with $X_{2}\in T(\mathcal{K}_{1},-)^{-1}(\mathcal{X}_{3})\cap \mathcal{X}_{2}$.

\begin{proof}
We shall use \ref{sec:1.2.3}. We put $\mathcal{S}=\{Y_{2},C(X_{1},Y_{3})\}$, 
where $Y_{2}\in\mathcal{S}_{2},X_{1}\in\mathcal{K}_{1},Y_{3}\in\mathcal{S}_{3}$.
\paragraph{Step 1} We show that $\mathcal{S}\subset
(T(\mathcal{K}_{1},-)^{-1}(\mathcal{X}_{3})\cap \mathcal{X}_{2})^{\perp}$.
Let $A_{2}\rightarrowtail B_{2}\twoheadrightarrow X_{2}$ be an exact sequence with
$X_{2}\in T(\mathcal{K}_{1},-)^{-1}(\mathcal{X}_{3})\cap \mathcal{X}_{2}$.
Let $Y_{2}\in\mathcal{S}_{2}$; then we have the exact sequence 
$\mathcal{A}_{2}(X_{2},Y_{2})\rightarrowtail \mathcal{A}_{2}(B_{2},Y_{2})
\twoheadrightarrow \mathcal{A}_{2}(A_{2},Y_{2})$ since $Y_{2}\in\mathcal{Y}_{2}$
and $X_{2}\in\mathcal{X}_{2}$. Let $X_{1}\in\mathcal{K}_{1},Y_{3}\in\mathcal{S}_{3}$.
By assumption we have the exact sequence 
\begin{displaymath}
\xymatrix{
{T(X_{1},A_{2})} \ar@{>->}[r] & {T(X_{1},B_{2})} 
\ar@{->>}[r]  & {T(X_{1},X_{2})}
}
\end{displaymath}
and therefore the exact sequence 
\begin{displaymath}
\mathcal{A}_{3}(T(X_{1},X_{2}),Y_{3})\rightarrowtail \mathcal{A}_{3}(T(X_{1},B_{2}),Y_{3})
\twoheadrightarrow \mathcal{A}_{3}(T(X_{1},A_{2}),Y_{3})
\end{displaymath}
By adjunction we have the exact sequence 
\begin{displaymath}
\mathcal{A}_{2}(X_{2},C(X_{1},Y_{3}))\rightarrowtail \mathcal{A}_{2}(B_{2},C(X_{1},Y_{3}))
\twoheadrightarrow \mathcal{A}_{2}(A_{2},C(X_{1},Y_{3}))
\end{displaymath}
\paragraph{Step 2} We show that ${}^{\perp}\mathcal{S}\subset
T(\mathcal{K}_{1},-)^{-1}(\mathcal{X}_{3})\cap \mathcal{X}_{2}$. 
Let $X_{2}\in {}^{\perp}\mathcal{S}$. Since 
$\mathcal{S}_{2}\subset \mathcal{S}$ we have (\ref{sec:1.2.1}(2)) 
$^{\perp}\mathcal{S}\subset ^{\perp}\mathcal{S}_{2}=\mathcal{X}_{2}$, hence 
it suffices to show that $T(X_{1},X_{2})\in ^{\perp}\mathcal{S}_{3}$ whenever
$X_{1}\in\mathcal{K}_{1}$. Let $Y_{3}\rightarrowtail A_{3}\twoheadrightarrow B_{3}$ 
be an exact sequence with $Y_{3}\in \mathcal{S}_{3}$. By assumption we have 
the exact sequence
\begin{displaymath}
\xymatrix{
{C(X_{1},Y_{3})} \ar@{>->}[r] & {C(X_{1},A_{3})} 
\ar@{->>}[r]  & {C(X_{1},B_{3})}
}
\end{displaymath}
and therefore by the definition of $\mathcal{S}$ we have the exact sequence 
\begin{displaymath}
\mathcal{A}_{2}(X_{2},C(X_{1},Y_{3}))\rightarrowtail \mathcal{A}_{2}(X_{2},C(X_{1},A_{3})))
\twoheadrightarrow \mathcal{A}_{2}(X_{2},C(X_{1},B_{3}))
\end{displaymath}
By adjunction we have the exact sequence 
\begin{displaymath}
\mathcal{A}_{3}(T(X_{1},X_{2}),Y_{3})\rightarrowtail \mathcal{A}_{3}(T(X_{1},X_{2}),A_{3})
\twoheadrightarrow \mathcal{A}_{3}(T(X_{1},X_{2}),B_{3})
\end{displaymath}
hence $T(X_{1},X_{2})\in ^{\perp}\mathcal{S}_{3}$ and 
$(T(\mathcal{K}_{1},-)^{-1}(\mathcal{X}_{3})\cap \mathcal{X}_{2},
(T(\mathcal{K}_{1},-)^{-1}(\mathcal{X}_{3})\cap \mathcal{X}_{2})^{\perp})$ 
is a cotorsion theory.
Let now $X_{1}\in\mathcal{K}_{1}$ and $Y_{3}\rightarrowtail A_{3}\twoheadrightarrow B_{3}$ 
be an exact sequence with $Y_{3}\in \mathcal{S}_{3}$. Applying $C(X_{1},-)$ to the sequence
we obtain the exact sequence $C(X_{1},Y_{3})\rightarrowtail C(X_{1},A_{3})
\rightarrow C(X_{1},B_{3})$. Let $f:X_{2}'\to C(X_{1},B_{3})$ be an epimorphism
with $X_{2}'\in T(\mathcal{K}_{1},-)^{-1}(\mathcal{X}_{3})\cap\mathcal{X}_{2}$.
Applying $\mathcal{A}_{3}(T(X_{1},X_{2}'),-)$ to $Y_{3}\rightarrowtail A_{3}\twoheadrightarrow B_{3}$ 
we obtain the exact sequence 
\begin{displaymath}
\mathcal{A}_{3}(T(X_{1},X_{2}'),Y_{3})\rightarrowtail \mathcal{A}_{3}(T(X_{1},X_{2}'),A_{3})
\twoheadrightarrow \mathcal{A}_{3}(T(X_{1},X_{2}'),B_{3})
\end{displaymath}
By adjunction we have the exact sequence 
\begin{displaymath}
\mathcal{A}_{2}(X_{2}',C(X_{1},Y_{3}))\rightarrowtail \mathcal{A}_{2}(X_{2}',C(X_{1},A_{3}))
\twoheadrightarrow \mathcal{A}_{2}(X_{2}',C(X_{1},B_{3}))
\end{displaymath}
Therefore there is $g:X_{2}'\to C(X_{1},A_{3})$ such that the diagram
\begin{displaymath}
\xymatrix{
& {X_{2}'} \ar@{->>}[d]^{f} \ar[dl]_{g}\\
{C(X_{1},A_{3})}\ar[r] & {C(X_{1},B_{3})}
}
\end{displaymath}
commutes, hence its horizontal arrow is an epimorphism.
\end{proof}

\subsubsection{} \label{sec:2.2.3}

Under the assumptions \ref{sec:2.2.2}, if $\mathcal{X}_{2},\mathcal{X}_{3}$ 
are left exact (\ref{sec:1.2.7}) then so is
$T(\mathcal{K}_{1},-)^{-1}(\mathcal{X}_{3})\cap \mathcal{X}_{2}$.
\begin{proof}
Let $X_{2}'\rightarrowtail X_{2}\twoheadrightarrow X_{2}''$ be an exact sequence with
$X_{2},X_{2}''\in T(\mathcal{K}_{1},-)^{-1}(\mathcal{X}_{3})\cap \mathcal{X}_{2}$.
Then $X_{2}'\in \mathcal{X}_{2}$ since $\mathcal{X}_{2}$ is left exact.
Let $X_{1}\in \mathcal{K}_{1}$; then we have by \ref{sec:2.2.2}(2) the exact sequence 
\begin{displaymath}
\xymatrix{
{T(X_{1},X_{2}')} \ar@{>->}[r] & {T(X_{1},X_{2})} 
\ar@{->>}[r]  & {T(X_{1},X_{2}'')}
}
\end{displaymath}
It follows that $T(X_{1},X_{2}')\in \mathcal{X}_{3}$ since $\mathcal{X}_{3}$ is left exact.
\end{proof}

\subsubsection{} \label{sec:2.2.4}
Let 
\[
T:\mathcal{A}_{1}\times \mathcal{A}_{2}\to \mathcal{A}_{3}, \
H:\mathcal{A}_{2}^{op}\times \mathcal{A}_{3}\to \mathcal{A}_{1}, \
C:\mathcal{A}_{1}^{op}\times \mathcal{A}_{3}\to \mathcal{A}_{2}
\]
be an abelian THC--situation. Let $\mathcal{K}_{2}\subset Ob(\mathcal{A}_{2})$
be a class of objects, $(\mathcal{X}_{3},\mathcal{Y}_{3})$ 
a cotorsion theory in $\mathcal{A}_{3}$ and $(\mathcal{X}_{1},\mathcal{Y}_{1})$ 
a cotorsion theory in $\mathcal{A}_{1}$. We define
\[
H(\mathcal{K}_{2},-)^{-1}(\mathcal{Y}_{1})=
\{A_{3}\in Ob(\mathcal{A}_{3}),\ H(K_{2},A_{3})\in \mathcal{Y}_{1} {\rm \ for\ all}\  K_{2}
\in\mathcal{K}_{2}\}
\]

Suppose that

(1) $(\mathcal{X}_{3},\mathcal{Y}_{3})$ is cogenerated by $\mathcal{S}_{3}$ and
$(\mathcal{X}_{1},\mathcal{Y}_{1})$ is cogenerated by $\mathcal{S}_{1}$;

(2) for each object $K_{2}$ of $\mathcal{K}_{2}$ the functor $H(K_{2},-)$ preserves 
epimorphisms with kernel in $\mathcal{Y}_{3}$;

(3) the category $\mathcal{A}_{3}$ has enough 
$H(\mathcal{K}_{2},-)^{-1}(\mathcal{Y}_{1})\cap \mathcal{Y}_{3}$ objects,
meaning that for every $A_{3}\in Ob(\mathcal{A}_{3})$ there is a monomorphism
$A_{3}\to Y_{3}$ with $Y_{3}\in H(\mathcal{K}_{2},-)^{-1}(\mathcal{Y}_{1})\cap \mathcal{Y}_{3}$.

Then $(^{\perp}(H(\mathcal{K}_{2},-)^{-1}(\mathcal{Y}_{1})\cap \mathcal{Y}_{3}),
H(\mathcal{K}_{2},-)^{-1}(\mathcal{Y}_{1})\cap \mathcal{Y}_{3})$ is a cotorsion theory.

\begin{proof}
The proof is similar to the proof of \ref{sec:2.2.2}, with $\mathcal{S}=\{X_{3},T(X_{1},K_{2})\}$,
where $X_{3}\in\mathcal{S}_{3},X_{1}\in\mathcal{S}_{1},K_{2}\in\mathcal{K}_{2}$.
\end{proof}

\subsubsection{} \label{sec:2.2.5}
Let 
\[
T:\mathcal{A}_{1}\times \mathcal{A}_{2}\to \mathcal{A}_{3}, \
H:\mathcal{A}_{2}^{op}\times \mathcal{A}_{3}\to \mathcal{A}_{1}, \
C:\mathcal{A}_{1}^{op}\times \mathcal{A}_{3}\to \mathcal{A}_{2}
\]
be an abelian THC--situation. Let $\mathcal{K}_{3}\subset Ob(\mathcal{A}_{3})$
be a class of objects, $(\mathcal{X}_{2},\mathcal{Y}_{2})$ 
a cotorsion theory in $\mathcal{A}_{2}$ and $(\mathcal{X}_{1},\mathcal{Y}_{1})$ 
a cotorsion theory in $\mathcal{A}_{1}$. We define
\[
H(-,\mathcal{K}_{3})^{-1}(\mathcal{Y}_{1})=
\{A_{2}\in Ob(\mathcal{A}_{2}),\ H(A_{2},K_{3})\in 
\mathcal{Y}_{1} {\rm \ for\ all}\  K_{3}
\in\mathcal{K}_{3}\}
\]

Suppose that

(1) $(\mathcal{X}_{2},\mathcal{Y}_{2})$ is generated by $\mathcal{S}_{2}$ and
$(\mathcal{X}_{1},\mathcal{Y}_{1})$ is cogenerated by $\mathcal{S}_{1}$;

(2) for each object $K_{3}$ of $\mathcal{K}_{3}$ the functor $H(-,K_{3})$ preserves 
epimorphisms with kernel in $\mathcal{X}_{2}$;

(3) the category $\mathcal{A}_{2}$ has enough 
$H(-,\mathcal{K}_{3})^{-1}(\mathcal{Y}_{1})\cap \mathcal{X}_{2}$ objects,
meaning that for every $A_{2}\in Ob(\mathcal{A}_{2})$ there is an epimorphism
$X_{2}\to A_{2}$ with $X_{2}\in H(-,\mathcal{K}_{3})^{-1}(\mathcal{Y}_{1})\cap \mathcal{X}_{2}$.

Then $(H(-,\mathcal{K}_{3})^{-1}(\mathcal{Y}_{1})\cap \mathcal{X}_{2},
(H(-,\mathcal{K}_{3})^{-1}(\mathcal{Y}_{1})\cap \mathcal{X}_{2})^{\perp})$ is a cotorsion theory.

\begin{proof}
The proof is similar to the proof of \ref{sec:2.2.2}, with $\mathcal{S}=\{Y_{2},C(X_{1},K_{3})\}$,
where $Y_{2}\in\mathcal{S}_{2},X_{1}\in\mathcal{S}_{1},K_{3}\in\mathcal{K}_{3}$.
\end{proof}

\section{The flat cotorsion theory associated to an abelian tensor-hom-cotensor situation}

\subsection{} \label{sec:3.1}
Let 
\[
T:\mathcal{A}_{1}\times \mathcal{A}_{2}\to \mathcal{A}_{3}, \
H:\mathcal{A}_{2}^{op}\times \mathcal{A}_{3}\to \mathcal{A}_{1}, \
C:\mathcal{A}_{1}^{op}\times \mathcal{A}_{3}\to \mathcal{A}_{2}
\]
be an abelian THC--situation. An object $P_{2}$ (resp. $P_{1}$) of 
$\mathcal{A}_{2}$ (resp. $\mathcal{A}_{1}$) is called \emph{flat} if 
$T(-,P_{2})$ (resp. $T(P_{1},-)$) preserves monomorphisms. 
We denote by $\mathcal{F}$ the class of flat objects of $\mathcal{A}_{2}$ 
and by $\mathcal{F}'$ the class of flat objects of $\mathcal{A}_{1}$.

\subsubsection{} \label{sec:3.1.1}
An object $P_{2}$ of $\mathcal{A}_{2}$ is flat if and only if 
for every exact sequence
\[
A_{1}'\overset{u}\rightarrow A_{1}\overset{v}\rightarrow A_{1}''
\]
in $\mathcal{A}_{1}$ the sequence 
\begin{displaymath}
\xymatrix{
{T(A'_{1},P_{2})} \ar[rr]^{T(u,P_{2})} && {T(A_{1},P_{2})} 
\ar[rr]^{T(v,P_{2})}  && {T(A''_{1},P_{2})}
}
\end{displaymath}
is exact.
\begin{proof}
Suppose $P_{2}$ is flat. Consider the commutative diagram
\begin{displaymath}
\xymatrix{
{A'_{1}} \ar[rr]^{u}\ar@{->>}[dr]_{p} && {A_{1}} \ar[rr]^{v}\ar@{->>}[dr]_{q} 
&& {A''_{1}}\ar@{->>}[r]^{p''} & {A''_{1}/\im v}\\
& {\im u} \ar@{>->}[ur]_{i} && {\im v} \ar@{>->}[ur]_{j}
}
\end{displaymath}
Since $T(-,P_{2})$ is a left adjoint we have the exact sequence
\begin{displaymath}
\xymatrix{
{T(A'_{1},P_{2})} \ar@{>->}[rr]^{T(u,P_{2})} && {T(A_{1},P_{2})} 
\ar@{->>}[rr]^{T(q,P_{2})}  && {T(\im v,P_{2})}
}
\end{displaymath}
By assumption we have the exact sequence
\begin{displaymath}
\xymatrix{
{T(\im v,P_{2})} \ar@{>->}[rr]^{T(j,P_{2})} && {T(A''_{1},P_{2})} 
\ar@{->>}[rr]^{T(p'',P_{2})}  && {T(A''/\im v,P_{2})}
}
\end{displaymath}
By splicing we have the exact sequence 
\begin{displaymath}
\xymatrix{
{T(A'_{1},P_{2})}\ar[rr]^{T(u,P_{2})}&& {T(A_{1},P_{2})}
\ar[rr]^{T(jq,P_{2})}&& {T(\im v,P_{2})}\ar[rr]^{T(p'',P_{2})}&&{T(A''/\im v,P_{2})}
}
\end{displaymath}
hence the claim. Conversely, one uses the exact sequence 
$0\overset{u}\to A\overset{v}\to B$.
\end{proof}
\subsubsection{} \label{sec:3.1.2}
(Lambek) (1) Let $P_{2}\in Ob(\mathcal{A}_{2})$.

(i) If $P_{2}$ is flat and $E_{3}$ is an injective object of $\mathcal{A}_{3}$
then  $H(P_{2},E_{3})$ is injective.

(ii) Suppose that $\mathcal{A}_{3}$ has a set $(E_{i})_{i\in I}$ 
of injective cogenerators. If $H(P_{2},E_{i})$ is injective for each $i\in I$
then $P_{2}$ is flat.

(2) Let $P_{1}\in Ob(\mathcal{A}_{1})$.

(i) If $P_{1}$ is flat and $E_{3}$ is an injective object of $\mathcal{A}_{3}$
then  $C(P_{1},E_{3})$ is injective.

(ii) Suppose that $\mathcal{A}_{3}$ has a set $(E_{i})_{i\in I}$ 
of injective cogenerators. If $C(P_{1},E_{i})$ is injective for each $i\in I$
then $P_{1}$ is flat.
\begin{proof}
(1) Let $A_{1}\rightarrowtail B_{1}\twoheadrightarrow C_{1}$ be 
an exact sequence in $\mathcal{A}_{1}$. We prove (i). 
Since $E_{3}$ is injective we have the exact sequence
\begin{displaymath}
\mathcal{A}_{3}(T(C_{1},P_{2}),E_{3})\rightarrowtail 
\mathcal{A}_{3}(T(B_{1},P_{2}),E_{3})
\twoheadrightarrow \mathcal{A}_{3}(T(A_{1},P_{2}),E_{3})
\end{displaymath}
By adjunction we have the exact sequence 
\begin{displaymath}
\mathcal{A}_{1}(C_{1},H(P_{2},E_{3}))\rightarrowtail 
\mathcal{A}_{1}(B_{1},H(P_{2},E_{3}))
\twoheadrightarrow \mathcal{A}_{1}(A_{1},H(A_{2},E_{3}))
\end{displaymath}
showing that $H(P_{2},E_{3})$ is injective. We prove (ii). 
The sequence 
\[
T(A_{1},P_{2})\to T(B_{1},P_{2})\twoheadrightarrow T(C_{1},P_{2})
\] is short exact if
\begin{equation} \label{eq:3.1.2}
\mathcal{A}_{3}(T(C_{1},P_{2}),E_{i})\rightarrow 
\mathcal{A}_{3}(T(B_{1},P_{2}),E_{i})
\rightarrow \mathcal{A}_{3}(T(A_{1},P_{2}),E_{i})
\end{equation}
is short exact for each $i\in I$ (\ref{sec:1.1.1}). 
By adjunction and assumption the sequence 
(\ref{eq:3.1.2}) is short exact. The proof of part 2 is similar.
\end{proof}

\subsubsection{} \label{sec:3.1.3}
Let $E_{3}$ be an injective object of $\mathcal{A}_{3}$.

(1) If $\mathcal{A}_{1}$ has enough flat objects, meaning that for each 
$A_{1}\in Ob(\mathcal{A}_{1})$ there is an epimorphism $P_{1}\to A_{1}$ 
with $P_{1}$ flat, then $H(-,E_{3})$ is exact.

(2) If $\mathcal{A}_{2}$ has enough flat objects, meaning that for each 
$A_{2}\in Ob(\mathcal{A}_{2})$ there is an epimorphism $P_{2}\to A_{2}$ 
with $P_{2}$ flat, then $C(-,E_{3})$ is exact.
\begin{proof}
(1) We use \ref{sec:2.1.1}(1) with $F=C(-,E_{3}), \mathcal{K}=\mathcal{F}$
and $\mathcal{K}'$ equal to the class of injective objects of $\mathcal{A}_{2}$.
It suffices to show that if $P_{1}$ is a flat object of $\mathcal{A}_{1}$ then 
$C(P_{1},E_{3})$ is injective. But this is \ref{sec:3.1.2}(2).
The proof of part 2 is similar.
%Let $A_{2}\rightarrowtail B_{2}\twoheadrightarrow C_{2}$ be an 
%exact sequence in $\mathcal{A}_{2}$.
%Consider the solid arrows diagram
%\begin{displaymath}
%\xymatrix{
%& & {P_{1}} \ar@{->>}[d]^{f} \ar@{..>}[dl]\\
%{H(C_{2},E_{3})} \ar@{>->}[r] & {H(B_{2},E_{3})}\ar[r] & {H(A_{2},E_{3})}
%}
%\end{displaymath}
%with $P_{1}$ flat. By adjunction we have a commutative solid arrows diagram
%\begin{displaymath}
%\xymatrix{
%{T(P_{1},A_{2})} \ar@{>->}[r] \ar[d]_{f^{\#}}& 
%{T(P_{1},B_{2})}\ar@{->>}[r]\ar@{..>}[dl]\ar[d] & {T(P_{1},C_{2})}\ar@{=}[d]\\
%{E_{3}} \ar@{>->}[r] & {PO}\ar@{->>}[r] & {T(P_{1},C_{2})}
%}
%\end{displaymath}
%where $f^{\#}$ is the adjoint of $f$ and $PO$ means pushout. Since the 
%bottom exact sequence splits, there is a dotted arrow in the previous 
%d%iagram that makes commutative the triangular diagram containing 
%$f%^{\#}$. By adjunction there is a dotted arrow in the first diagram 
%that makes commutative the triangular diagram containing $f$, hence 
%the morphism $H(B_{2},E_{3})\to H(A_{2},E_{3})$ is an epimorphism. 
%The argument for part 2 is dual.
\end{proof}

\subsubsection{} \label{sec:3.1.4}
Suppose that

(1) the category $\mathcal{A}_{i}$ $(1\leqslant i\leqslant 2)$ has enough 
flat objects, meaning that for each $A_{i}\in Ob(\mathcal{A}_{i})$ there 
is an epimorphism $P_{i}\to A_{i}$ with $P_{i}$ flat;

(2) the category $\mathcal{A}_{3}$ has enough injective objects.

Then $(\mathcal{F},\mathcal{F}^{\perp})$ is a cotorsion theory. The 
objects of $\mathcal{F}^{\perp}$ are usually called cotorsion objects.
\begin{proof}
We shall use \ref{sec:1.2.3}. We put $\mathcal{S}=\{C(A_{1},E_{3})\}$, where 
$A_{1}\in Ob(\mathcal{A}_{1})$ and $E_{3}$ is an injective object of $\mathcal{A}_{3}$.
\paragraph{Step 1} We show that $\mathcal{S}\subset \mathcal{F}^{\perp}$.
Let $A_{2}\rightarrowtail B_{2}\twoheadrightarrow P_{2}$ be
an exact sequence with $P_{2}$ flat, let $A_{1}\in Ob(\mathcal{A}_{1})$ and $E_{3}$ an 
injective object of $\mathcal{A}_{3}$. Applying the functor $\mathcal{A}_{2}(-,C(A_{1},E_{3}))$
to the previous exact sequence and using adjunction we have the exact sequence
\begin{displaymath}
\mathcal{A}_{1}(A_{1},H(P_{2},E_{3}))\rightarrowtail \mathcal{A}_{1}(A_{1},H(B_{2},E_{3}))
\rightarrow \mathcal{A}_{1}(A_{1},H(A_{2},E_{3}))
\end{displaymath}
Since the functor $T(-,P_{2})$ is exact, its right adjoint $H(P_{2},-)$ preserves injective objects
(\ref{sec:2.1.1}(2bis)), hence $H(P_{2},E_{3})$ is injective. Since $H(-,E_{3})$ is exact  
(\ref{sec:3.1.3}) the right arrow in the previous exact sequence is an epimorphism, 
showing that $C(A_{1},E_{3})\in \mathcal{F}^{\perp}$.

\paragraph{Step 2} We show that ${}^{\perp}\mathcal{S}\subset \mathcal{F}$.
Let $P_{2}\in {}^{\perp}\mathcal{S}$ and $A_{1}\rightarrowtail B_{1}\twoheadrightarrow C_{1}$
be an exact sequence in $\mathcal{A}_{1}$. There is a monomorphism $T(A_{1},P_{2})\to E_{3}$
with $E_{3}$ injective. Since $C(-,E_{3})$ is exact (\ref{sec:3.1.3}) we have the exact sequence
$C(C_{1},E_{3})\rightarrowtail C(B_{1},E_{3})\twoheadrightarrow C(A_{1},E_{3})$
hence by the assumption on $P_{2}$ the exact sequence
\begin{displaymath}
\mathcal{A}_{2}(P_{2},C(C_{1},E_{3}))\rightarrowtail 
\mathcal{A}_{2}(P_{2},C(B_{1},E_{3}))\twoheadrightarrow 
\mathcal{A}_{2}(P_{2},C(A_{1},E_{3}))
\end{displaymath}
By adjunction we have the exact sequence
\begin{displaymath}
\mathcal{A}_{3}(T(A_{1},P_{2}),E_{3})\rightarrowtail 
\mathcal{A}_{3}(T(B_{1},P_{2}),E_{3})\twoheadrightarrow 
\mathcal{A}_{3}(T(C_{1},P_{2}),E_{3})
\end{displaymath}
This implies that $T(A_{1},P_{2})\to T(B_{1},P_{2})$ is a monomorphism.
\end{proof}

\subsubsection{} \label{sec:3.1.5}
The class $\mathcal{F}$ is left exact (\ref{sec:1.2.7}).
\begin{proof}
Let $P_{2}'\rightarrowtail P_{2}\twoheadrightarrow P_{2}''$
be an exact sequence with $P_{2},P_{2}''\in \mathcal{F}$ and 
$A_{1}\rightarrowtail B_{1}\twoheadrightarrow C_{1}$
an exact sequence in $\mathcal{A}_{1}$. There is a monomorphism 
$T(A_{1},P_{2}')\to E_{3}$ with $E_{3}$ injective.
By \ref{sec:3.1.3} and \ref{sec:2.1.1}(2bis) we have the exact sequence
\[H(P_{2}'',E_{3})\rightarrowtail H(P_{2},E_{3})\twoheadrightarrow H(P_{2}',E_{3})
\]
with $H(P_{2}'',E_{3}),H(P_{2},E_{3})$ injective. Then 
$H(P_{2}',E_{3})$ is injective and therefore we have the exact sequence
\begin{displaymath}
\mathcal{A}_{1}(C_{1},H(P_{2}',E_{3}))\rightarrowtail 
\mathcal{A}_{1}(B_{1},H(P_{2}',E_{3}))\twoheadrightarrow 
\mathcal{A}_{1}(A_{1},H(P_{2}',E_{3}))
\end{displaymath}
By adjunction we have the exact sequence
\begin{displaymath}
\mathcal{A}_{3}(T(C_{1},P_{2}'),E_{3})\rightarrowtail 
\mathcal{A}_{3}(T(B_{1},P_{2}'),E_{3})\twoheadrightarrow 
\mathcal{A}_{3}(T(A_{1},P_{2}'),E_{3})
\end{displaymath}
This implies that $T(A_{1},P_{2}')\to T(B_{1},P'_{2})$ is a monomorphism.
\end{proof}

\subsubsection{} \label{sec:3.1.6}
(1) For every $A_{1}\in Ob(\mathcal{A}_{1})$ the functor $T(A_{1},-)$ 
preserves the monomorphisms with flat cokernel.

(2) Let $P_{2}\in Ob(\mathcal{A}_{2})$. Suppose that for every 
$A_{1}\in Ob(\mathcal{A}_{1})$ the functor $T(A_{1},-)$ preserves 
the monomorphisms with cokernel $P_{2}$. Then $P_{2}$ is flat.

(3) Suppose that $\mathcal{A}_{3}$ has a set $(E_{i})_{i\in I}$ 
of injective cogenerators and that the class of injective objects of 
$\mathcal{A}_{1}$ is cogenerated by a set $\mathcal{S}_{1}$. 
Then $\mathcal{F}$ is generated by the set 
$\mathcal{S}'=\{C(X_{1},E_{i})\}$, where $X_{1}\in\mathcal{S}_{1}$.

\begin{proof}
(1) This is a consequence of the definition of $\mathcal{S}$ and 
\ref{sec:2.1.1}(1bis) applied to the adjunction $(T(A_{1},-),C(A_{1},-))$ 
and the classes $\mathcal{K}=\mathcal{F}$,
$\mathcal{K}'=Ob(\mathcal{A}_{3})$.

(2) Let $P_{2}\in Ob(\mathcal{A}_{2})$ and 
$A_{1}\rightarrowtail B_{1}\twoheadrightarrow C_{1}$ 
an exact sequence in $\mathcal{A}_{1}$.
There is a monomorphism 
$f:T(A_{1},P_{2})\to E_{3}$ with $E_{3}$ injective. Since $C(-,E_{3})$ 
is exact we have the commutative solid arrows diagram
\begin{displaymath}
\xymatrix{
{C(C_{1},E_{3})} \ar@{>->}[r] \ar@{=}[d]& {PB}\ar@{->>}[r]\ar[d] & 
{P_{2}}\ar@{..>}[dl]\ar[d]^{f^{\flat}}\\
{C(C_{1},E_{3})} \ar@{>->}[r] & {C(B_{1},E_{3})}
\ar@{->>}[r] & {C(A_{1},E_{3})}
}
\end{displaymath}
where $f^{\flat}$ is the adjoint of $f$ and $PB$ means pullback. 
By the first assumption we have the exact sequence 
$T(C_{1},C(C_{1},Y_{3}))\rightarrowtail T(C_{1},PB)\twoheadrightarrow T(C_{1},P_{2})$.
Since $E_{3}$ is injective we have the exact sequence
\begin{displaymath}
\mathcal{A}_{3}(T(C_{1},P_{2}),E_{3})\rightarrowtail \mathcal{A}_{3}(T(C_{1},PB),E_{3})
\twoheadrightarrow \mathcal{A}_{3}(T(C_{1},C(C_{1},E_{3})),E_{3})
\end{displaymath}
By adjunction we have the exact sequence
\begin{displaymath}
\mathcal{A}_{2}(P_{2},C(C_{1},E_{3}))\rightarrowtail \mathcal{A}_{2}(PB,C(C_{1},E_{3}))
\twoheadrightarrow \mathcal{A}_{2}(C(C_{1},E_{3}),C(C_{1},E_{3}))
\end{displaymath}
This implies that the map $PB\to P_{2}$ in the previous commutative diagram has a section, 
therefore there is a dotted arrow that makes commutative the lower triangular diagram  
containing $f^{\flat}$. An adjunction argument implies then 
that $T(A_{1},P_{2})\to T(B_{1},P_{2})$ is a monomorphism.

(3) From the proof of \ref{sec:3.1.4} we have $\mathcal{S}'\subset \mathcal{S}$,
so $\mathcal{F}=^{\perp}\mathcal{S}\subset ^{\perp}\mathcal{S}'$. 
Let $P_{2}\in^{\perp}\mathcal{S}'$. By \ref{sec:3.1.2} it suffices to show that $H(P_{2},E_{i})$
is injective. Let $A_{1}\rightarrowtail B_{1}\twoheadrightarrow X_{1}$ be
an exact sequence in $\mathcal{A}_{1}$ with $X_{1}\in\mathcal{S}_{1}$. By \ref{sec:3.1.3}
we have the exact sequence $C(X_{1},E_{i})\rightarrowtail C(B_{1},E_{i})\twoheadrightarrow 
C(A_{1},E_{i})$. Applying $\mathcal{A}_{2}(P_{2},-)$ to it and using an adjunction
argument we obtain that $H(P_{2},E_{i})$ is injective.
\end{proof}

\subsubsection{} \label{sec:3.1.7}
For a morphism $u_{1}:A_{1}\to B_{1}$ of $\mathcal{A}_{1}$ and a morphism 
$u_{2}:A_{2}\to B_{2}$ of $\mathcal{A}_{2}$, we denote by $u_{1}Tu_{2}$ 
the natural morphism
\[
T(A_{1},B_{2})\bigcup_{T(A_{1},A_{2})} T(B_{1},A_{2})\to T(B_{1},B_{2})
\]
If $u_{1},u_{2}$ are monomorphisms with flat cokernel then 
$u_{1}Tu_{2}$ is a monomorphism.
\begin{proof}
We shall need the fact that \ref{sec:3.1.1}--\ref{sec:3.1.6} also hold, with suitable 
modifications in proofs, for the class of flat objects of $\mathcal{A}_{1}$. 
Using this and the exact sequence
$A_{i}\overset{u_{i}}\rightarrowtail B_{i}\twoheadrightarrow 
P_{i}$ ($1\leqslant i\leqslant 2$) we have the commutative diagram
\begin{displaymath}
\xymatrix{
{T(A_{1},A_{2})}\ar@{>->}[rr]\ar@{>->}[d]&&{T(A_{1},B_{2})}\ar@{>->}[d]\\
{T(B_{1},A_{2})}\ar@{>->}[rr]&&
{T(B_{1},B_{2})}\ar@{->>}[rr] \ar@{->>}[d]
&&{T(B_{1},P_{2})} \ar@{->>}[d]\\
& &{T(P_{1},B_{2})} \ar@{->>}[rr]
&&{T(P_{1},P_{2})}
}
\end{displaymath}
The arrows that form of the upper left square diagram are
monomorphisms by \ref{sec:3.1.6}(1). It follows then from 
\cite[Lemma 1.10.5 and Proposition 1.7.4]{Bo2} that
\[
T(A_{1},B_{2})\bigcup_{T(A_{1},A_{2})} T(B_{1},A_{2})
\cong T(B_{1},A_{2})\cup T(A_{1},B_{2})
\]
hence $u_{1}Tu_{2}$ is a monomorphism.
\end{proof}

Here is another proof that $(\mathcal{F},\mathcal{F}^{\perp})$ 
is a cotorsion theory.

\subsubsection{} \label{sec:3.1.8}
Suppose that

(1) the category $\mathcal{A}_{i}$ $(1\leqslant i\leqslant 2)$ has 
enough flat objects, meaning that for each $A_{i}\in Ob(\mathcal{A}_{i})$ 
there is an epimorphism $P_{i}\to A_{i}$ with $P_{i}$ flat;

(2) the category $\mathcal{A}_{3}$ has an injective cogenerator $E_{3}$.

Then $(\mathcal{F},\mathcal{F}^{\perp})$ is a cotorsion theory.

\begin{proof}
Consider the adjoint pair 
\begin{displaymath}
C(-,E_{3}):\mathcal{A}_{1}\rightleftarrows \mathcal{A}_{2}^{op}:H(-,E_{3})
\end{displaymath}
By \ref{sec:3.1.3} both functors are exact. By \ref{sec:2.1.2}(1bis) applied to 
$\mathcal{K}=Ob(\mathcal{A}_{1})$ we then obtain a cotorsion theory in 
$\mathcal{A}_{2}^{op}$ whose right class is the class of flat objects by \ref{sec:3.1.2}.
\end{proof}

\subsection{} \label{sec:3.2}
We shall repeat most of \ref{sec:3.1} replacing ``preserves monomorphisms" by 
``preserves monomorphisms with cokernel in a given class". Let 
\[
T:\mathcal{A}_{1}\times \mathcal{A}_{2}\to \mathcal{A}_{3}, \
H:\mathcal{A}_{2}^{op}\times \mathcal{A}_{3}\to \mathcal{A}_{1}, \
C:\mathcal{A}_{1}^{op}\times \mathcal{A}_{3}\to \mathcal{A}_{2}
\]
be an abelian THC--situation. Let $(\mathcal{X}_{3},\mathcal{Y}_{3})$ 
be a cotorsion theory in $\mathcal{A}_{3}$.

Let $\mathcal{K}_{1}\subset Ob(\mathcal{A}_{1})$ 
be a class of objects. An object $P_{2}$ of $\mathcal{A}_{2}$ is 
called $\mathcal{K}_{1}-$\emph{flat} if $T(-,P_{2})$ 
preserves monomorphisms with cokernel in $\mathcal{K}_{1}$. We denote by 
$\mathcal{K}_{1}-\mathcal{F}$ the class of $\mathcal{K}_{1}-$flat objects.

Let $\mathcal{K}_{2}$ be a class of objects of $\mathcal{A}_{2}$ that contains 
$\mathcal{K}_{1}-\mathcal{F}$. An object $P_{1}$ of $\mathcal{A}_{1}$ is called 
$\mathcal{K}_{2}-$\emph{flat} if $T(P_{1},-)$ preserves monomorphisms with 
cokernel in $\mathcal{K}_{2}$. We denote by 
$\mathcal{K}_{2}-\mathcal{F}$ the class of $\mathcal{K}_{2}-$flat objects.

\subsubsection{} \label{sec:3.2.1}
An object $P_{2}$ of $\mathcal{A}_{2}$ is $\mathcal{K}_{1}-$flat if and only if 
for every exact sequence
\[
A_{1}'\overset{u}\rightarrow A_{1}\overset{v}\rightarrow A_{1}''
\]
in $\mathcal{A}_{1}$ with $A''/\im v\in \mathcal{K}_{1}$ the sequence 
\begin{displaymath}
\xymatrix{
{T(A'_{1},P_{2})} \ar[rr]^{T(u,P_{2})} && {T(A_{1},P_{2})} 
\ar[rr]^{T(v,P_{2})}  && {T(A''_{1},P_{2})}
}
\end{displaymath}
is exact.
\begin{proof}
If $P_{2}$ is $\mathcal{K}_{1}-$flat the proof is the same as the first part of the proof of 
\ref{sec:3.1.1}. Conversely, use the exact sequence 
$0\overset{u}\to A\overset{v}\to B\twoheadrightarrow X_{1}$ with $X_{1}\in\mathcal{K}_{1}$.
\end{proof}

\subsubsection{} \label{sec:3.2.2}
Suppose that the category $\mathcal{A}_{3}$ has a set $(Y_{i})_{i\in I}$ of 
cogenerators that belong to $\mathcal{Y}_{3}$. Moreover, suppose that 
$T(\mathcal{K}_{1},\mathcal{K}_{1}-\mathcal{F})\subset \mathcal{X}_{3}$.
Then an object $P_{2}$ of $\mathcal{A}_{2}$ is $\mathcal{K}_{1}-$flat if 
and only if $H(P_{2},Y_{i})\in\mathcal{K}_{1}^{\perp}$ for each $i\in I$.
\begin{proof}
The proof is the same as for \ref{sec:3.1.2}. Note that it also shows 
that $H(P_{2},Y_{3})\in\mathcal{K}_{1}^{\perp}$ for all $P_{2}$ $\mathcal{K}_{1}-$flat
and all $Y_{3}\in\mathcal{Y}_{3}$.
\end{proof}

\subsubsection{} \label{sec:3.2.3}
Suppose that

(1) $T(\mathcal{K}_{1},\mathcal{K}_{1}-\mathcal{F})\subset \mathcal{X}_{3}$;

(2) the category $\mathcal{A}_{2}$ has enough $\mathcal{K}_{1}-$flat objects, 
meaning that for each $A_{2}\in Ob(\mathcal{A}_{2})$ there is an 
epimorphism $P_{2}\to A_{2}$ with $P_{2}$ $\mathcal{K}_{1}-$flat;

(3) $T(\mathcal{K}_{2}-\mathcal{F},\mathcal{K}_{1}-\mathcal{F})\subset \mathcal{X}_{3}$;

(4) the category $\mathcal{A}_{1}$ has enough $\mathcal{K}_{2}-$flat objects, 
meaning that for each $A_{1}\in Ob(\mathcal{A}_{1})$ there is an 
epimorphism $P_{1}\to A_{1}$ with $P_{1}$ $\mathcal{K}_{2}-$flat.

Then for every $Y_{3}\in\mathcal{Y}_{3}$ the functor $C(-,Y_{3})$ 
preserves monomorphisms  with cokernel in $\mathcal{K}_{1}$ and 
the functor $H(-,Y_{3})$ preserves epimorphisms 
with kernel in $\mathcal{K}_{1}-\mathcal{F}$.
\begin{proof}
Let $A_{1}\rightarrowtail B_{1}\twoheadrightarrow X_{1}$ be an 
exact sequence with $X_{1}\in \mathcal{K}_{1}$. Consider the solid arrows 
diagram
\begin{displaymath}
\xymatrix{
& & {P_{2}} \ar@{->>}[d]^{f} \ar@{..>}[dl]\\
{C(X_{1},Y_{3})} \ar@{>->}[r] & {C(B_{1},Y_{3})}\ar[r] & {C(A_{1},Y_{3})}
}
\end{displaymath}
with $P_{2}$ flat. By adjunction we have a commutative solid arrow diagram
\begin{displaymath}
\xymatrix{
{T(A_{1},P_{2})} \ar@{>->}[r] \ar[d]_{f^{\#}}& 
{T(B_{1},P_{2})}\ar@{->>}[r]\ar@{..>}[dl]\ar[d] & {T(X_{1},P_{2})}\ar@{=}[d]\\
{Y_{3}} \ar@{>->}[r] & {PO}\ar@{->>}[r] & {T(X_{1},P_{2})}
}
\end{displaymath}
where $f^{\#}$ is the adjoint of $f$ and $PO$ means pushout. The 
bottom exact sequence splits by (1), so there is a dotted arrow 
that makes commutative the triangular diagram containing $f^{\#}$. By adjunction 
there is a dotted arrow in the first diagram that makes commutative the triangular
diagram containing $f$. The argument for $H(-,Y_{3})$ is dual, using (3) and (4).
\end{proof}

\subsubsection{} \label{sec:3.2.4}
Suppose that

(1) $T(\mathcal{K}_{1},\mathcal{K}_{1}-\mathcal{F})\subset \mathcal{X}_{3}$;

(2) the category $\mathcal{A}_{2}$ has enough $\mathcal{K}_{1}-$flat objects, 
meaning that for each $A_{2}\in Ob(\mathcal{A}_{2})$ there is an epimorphism 
$P_{2}\to A_{2}$ with $P_{2}$ $\mathcal{K}_{1}-$flat;

(3) $T(\mathcal{K}_{2}-\mathcal{F},\mathcal{K}_{1}-\mathcal{F})\subset \mathcal{X}_{3}$;

(4) the category $\mathcal{A}_{1}$ has enough $\mathcal{K}_{2}-$flat objects, 
meaning that for each $A_{1}\in Ob(\mathcal{A}_{1})$ there is an epimorphism 
$P_{1}\to A_{1}$ with $P_{1}$ $\mathcal{K}_{2}-$flat.

(5) the category $\mathcal{A}_{3}$ has enough $\mathcal{Y}_{3}$ objects.

Then $(\mathcal{K}_{1}-\mathcal{F},\mathcal{K}_{1}-\mathcal{F}^{\perp})$ 
is a cotorsion theory.
\begin{proof}
We shall use \ref{sec:1.2.3}. We put $\mathcal{S}=\{C(X_{1},Y_{3})\}$, where 
$X_{1}\in \mathcal{K}_{1},Y_{3}\in\mathcal{Y}_{3}$.
\paragraph{Step 1} We show that  $\mathcal{S}\subset \mathcal{K}_{1}-\mathcal{F}^{\perp}$.
Let $A_{2}\rightarrowtail B_{2}\twoheadrightarrow P_{2}$ be
an exact sequence with $P_{2}$ $\mathcal{K}_{1}-$flat and let 
$X_{1}\in \mathcal{K}_{1}$, $Y_{3}\in\mathcal{Y}_{3}$. 
By \ref{sec:3.2.3} we have the exact sequence
$H(P_{2},Y_{3})\rightarrowtail H(B_{2},Y_{3})\twoheadrightarrow H(A_{2},Y_{3})$. 
By the proof of \ref{sec:3.2.2} 
we have the exact sequence
\begin{displaymath}
\mathcal{A}_{1}(X_{1},H(P_{2},Y_{3}))\rightarrowtail 
\mathcal{A}_{1}(X_{1},H(B_{2},Y_{3}))\twoheadrightarrow 
\mathcal{A}_{1}(X_{1},H(A_{2},Y_{3}))
\end{displaymath}
By adjunction we have the exact sequence
\begin{displaymath}
\mathcal{A}_{2}(P_{2},C(X_{1},Y_{3}))\rightarrowtail 
\mathcal{A}_{2}(B_{2},C(X_{1},Y_{3}))\twoheadrightarrow 
\mathcal{A}_{2}(A_{2},C(X_{1},Y_{3}))
\end{displaymath}
hence $C(X_{1},Y_{3})\in \mathcal{K}_{1}-\mathcal{F}^{\perp}$.
\paragraph{Step 2} We show that ${}^{\perp}\mathcal{S}\subset \mathcal{K}_{1}-\mathcal{F}$.
Let $P_{2}\in {}^{\perp}\mathcal{S}$ and $A_{1}\rightarrowtail B_{1}\twoheadrightarrow X_{1}$
be an exact sequence with $X_{1}\in \mathcal{K}_{1}$. There is a monomorphism 
$T(A_{1},P_{2})\to Y_{3}$ with $Y_{3}\in\mathcal{Y}_{3}$. By \ref{sec:3.2.3} 
we have the exact sequence $C(X_{1},Y_{3})\rightarrowtail C(B_{1},Y_{3})
\twoheadrightarrow C(A_{1},Y_{3})$
hence by the assumption on $P_{2}$ the exact sequence
\begin{displaymath}
\mathcal{A}_{2}(P_{2},C(X_{1},Y_{3}))\rightarrowtail 
\mathcal{A}_{2}(P_{2},C(B_{1},Y_{3}))\twoheadrightarrow 
\mathcal{A}_{2}(P_{2},C(A_{1},Y_{3}))
\end{displaymath}
By adjunction we have the exact sequence
\begin{displaymath}
\mathcal{A}_{3}(T(X_{1},P_{2}),Y_{3})\rightarrowtail 
\mathcal{A}_{3}(T(B_{1},P_{2}),Y_{3})\twoheadrightarrow 
\mathcal{A}_{3}(T(A_{1},P_{2}),Y_{3})
\end{displaymath}
This implies that $T(A_{1},P_{2})\to T(B_{1},P_{2})$ is a monomorphism.
\end{proof}

\subsubsection{} \label{sec:3.2.5}
If $\mathcal{K}_{1}^{\perp}$ is left exact 
then $\mathcal{K}_{1}-\mathcal{F}$ is left exact.
\begin{proof}
Let $P_{2}'\rightarrowtail P_{2}\twoheadrightarrow P_{2}''$
be an exact sequence with $P_{2},P_{2}''\in \mathcal{K}_{1}-\mathcal{F}$ and 
$A_{1}\rightarrowtail B_{1}\twoheadrightarrow X_{1}$
an exact sequence with $X_{1}\in \mathcal{K}_{1}$. There is a monomorphism 
$T(A_{1},P_{2}')\to Y_{3}$ with $Y_{3}\in\mathcal{Y}_{3}$.
By \ref{sec:3.2.3} and the proof of \ref{sec:3.2.2} we have the exact sequence
\[H(P_{2}'',Y_{3})\rightarrowtail H(P_{2},Y_{3})\twoheadrightarrow H(P_{2}',Y_{3})
\]
with $H(P_{2}'',Y_{3}),H(P_{2},Y_{3})\in\mathcal{K}_{1}^{\perp}$. Then 
$H(P_{2}',Y_{3})\in\mathcal{K}_{1}^{\perp}$ and therefore we have the exact sequence
\begin{displaymath}
\mathcal{A}_{1}(X_{1},H(P_{2}',Y_{3}))\rightarrowtail 
\mathcal{A}_{1}(B_{1},H(P_{2}',Y_{3}))\twoheadrightarrow 
\mathcal{A}_{1}(A_{1},H(P_{2}',Y_{3}))
\end{displaymath}
By adjunction we have the exact sequence
\begin{displaymath}
\mathcal{A}_{3}(T(X_{1},P_{2}'),Y_{3})\rightarrowtail 
\mathcal{A}_{3}(T(B_{1},P_{2}'),Y_{3})\twoheadrightarrow 
\mathcal{A}_{3}(T(A_{1},P_{2}'),Y_{3})
\end{displaymath}
This implies that $T(A_{1},P_{2}')\to T(B_{1},P_{2}')$ is a monomorphism.
\end{proof}

\subsubsection{} \label{sec:3.2.6}
(1) For every $X_{1}\in \mathcal{K}_{1}$ the functor $T(X_{1},-)$ preserves 
the monomorphisms with cokernel in $\mathcal{K}_{1}-\mathcal{F}$.

(2) Let $P_{2}\in \mathcal{K}_{2}$. Suppose that for every 
$X_{1}\in \mathcal{K}_{1}$ the functor $T(X_{1},-)$ preserves the 
monomorphisms with cokernel $P_{2}$. If 
$T(\mathcal{K}_{1},\mathcal{K}_{2})\subset \mathcal{X}_{3}$ then 
$P_{2}$ is $\mathcal{K}_{1}-$flat.
\begin{proof}
(1) This is a consequence of the definition of $\mathcal{S}$ and 
\ref{sec:2.1.1}(1bis) applied to the adjunction $(T(X_{1},-),C(X_{1},-))$ 
and the classes $\mathcal{K}=\mathcal{K}_{1}-\mathcal{F}$, 
$\mathcal{K}'=\mathcal{X}_{3}$.

(2) Let $P_{2}\in\mathcal{K}_{2}$ and 
$A_{1}\rightarrowtail B_{1}\twoheadrightarrow X_{1}$ 
be an exact sequence with $X_{1}\in \mathcal{K}_{1}$.
There is a monomorphism $f:T(A_{1},P_{2})\to Y_{3}$ with 
$Y_{3}\in\mathcal{Y}_{3}$. By \ref{sec:3.2.3} we have the exact sequence
$C(X_{1},Y_{3})\rightarrowtail C(B_{1},Y_{3})\twoheadrightarrow C(A_{1},Y_{3})$.
Form the commutative solid arrows diagram
\begin{displaymath}
\xymatrix{
{C(X_{1},Y_{3})} \ar@{>->}[r] \ar@{=}[d]& {PB}\ar@{->>}[r]\ar[d] & {P_{2}}\ar@{..>}[dl]
\ar[d]^{f^{\flat}}\\
{C(X_{1},Y_{3})} \ar@{>->}[r] & {C(B_{1},Y_{3})}\ar@{->>}[r] & {C(A_{1},Y_{3})}
}
\end{displaymath}
where $f^{\flat}$ is the adjoint of $f$ and $PB$ means pullback. 
By the first assumption we have the exact sequence 
$T(X_{1},C(X_{1},Y_{3}))\rightarrowtail T(X_{1},PB)\twoheadrightarrow T(X_{1},P_{2})$.
By the second assumption we have the exact sequence
\begin{displaymath}
\mathcal{A}_{3}(T(X_{1},P_{2}),Y_{3})\rightarrowtail \mathcal{A}_{3}(T(X_{1},PB),Y_{3})
\twoheadrightarrow \mathcal{A}_{3}(T(X_{1},C(X_{1},Y_{3})),Y_{3})
\end{displaymath}
By adjunction we have the exact sequence
\begin{displaymath}
\mathcal{A}_{2}(P_{2},C(X_{1},Y_{3}))\rightarrowtail \mathcal{A}_{2}(PB,C(X_{1},Y_{3}))
\twoheadrightarrow \mathcal{A}_{2}(C(X_{1},Y_{3}),C(X_{1},Y_{3}))
\end{displaymath}
This implies that the map $PB\to P_{2}$ in the previous commutative diagram has a section, 
therefore there is a dotted arrow that makes commutative the lower triangular 
containing $f^{\flat}$. An adjunction argument implies then that 
$T(A_{1},P_{2})\to T(B_{1},P_{2})$ is a monomorphism.
\end{proof}

\subsubsection{} \label{sec:3.2.7} 
In \ref{sec:3.2} we take $\mathcal{K}_{1}=\mathcal{I}_{1}$, the class of 
injective objects of $\mathcal{A}_{1}$, $\mathcal{K}_{2}=Ob(\mathcal{A}_{2})$ 
and $\mathcal{X}_{3}=Ob(\mathcal{A}_{3})$. We assume that

(1) the category $\mathcal{A}_{i}$ $(1\leqslant i\leqslant 2)$ has enough flat objects, 
meaning that for each $A_{i}\in Ob(\mathcal{A}_{i})$ there is an epimorphism 
$P_{i}\to A_{i}$ with $P_{i}$ flat;

(2) the category $\mathcal{A}_{3}$ has enough injective objects.

Then by \ref{sec:3.2.4} $(\mathcal{I}_{1}-\mathcal{F},
\mathcal{I}_{1}-\mathcal{F}^{\perp})$ is a cotorsion theory.

\subsection{} \label{sec:3.3}
Let $\mathcal{A}$ be an abelian category and
\[
T:\mathcal{A}\times \mathcal{A}\to \mathcal{A}, \
H:\mathcal{A}^{op}\times \mathcal{A}\to \mathcal{A}, \
C:\mathcal{A}^{op}\times \mathcal{A}\to \mathcal{A}
\]
a THC--situation. We say that a class $\mathcal{K}$ of objects of 
$\mathcal{A}$ is $T$-\emph{closed} if if $T(X,X')\in\mathcal{K}$ 
whenever $X,X'\in\mathcal{K}$. We say that a cotorsion theory 
$(\mathcal{X},\mathcal{Y})$ in $\mathcal{A}$ is $T$-closed if
$\mathcal{X}$ is $T$-closed.

\subsubsection{} \label{sec:3.3.1}
A sufficient condition for $(\mathcal{F},\mathcal{F}^{\perp})$ 
to be $T$-closed is to have, for all $X,Y,Z\in Ob(\mathcal{A})$, a 
natural morphism
\[
T(X,T(Y,Z))\to T(T(X,Y),Z)
\]
that is an isomorphism.

\subsubsection{} \label{sec:3.3.2}
Suppose that $\mathcal{X}\subset \mathcal{F}$ and 
that $\mathcal{A}$ has enough $\mathcal{X}$ objects. 
The following are equivalent:

(1) $(\mathcal{X},\mathcal{Y})$ is $T$-closed;

(2) $C(X,Y)\in \mathcal{Y}$ for all $X\in\mathcal{X},Y\in \mathcal{Y}$.
\begin{proof}
$(1)\Rightarrow (2)$ Let $X\in \mathcal{X}$. The claim follows from 
\ref{sec:2.1.1}(2bis) applied to the adjunction $(T(X,-),C(X,-))$ and 
$\mathcal{K}=\mathcal{K}'=\mathcal{X}$, together with the assumption that
$\mathcal{X}\subset \mathcal{F}$. $(2)\Rightarrow (1)$ Let $X'\in \mathcal{X}$.
The claim follows from \ref{sec:2.1.1}(2) applied to the adjunction 
$(T(X',-),C(X',-))$ and $\mathcal{K}=\mathcal{K}'=\mathcal{Y}$, 
together with the assumption that $\mathcal{A}$ has enough $\mathcal{X}$ 
objects.
\end{proof}

\subsection{} \label{sec:3.4}
Let 
\[
T:\mathcal{A}_{1}\times \mathcal{A}_{2}\to \mathcal{A}_{3}, \
H:\mathcal{A}_{2}^{op}\times \mathcal{A}_{3}\to \mathcal{A}_{1}, \
C:\mathcal{A}_{1}^{op}\times \mathcal{A}_{3}\to \mathcal{A}_{2}
\]
be an abelian THC--situation. An object $P_{2}$ (resp. $P_{1}$) of 
$\mathcal{A}_{2}$ ( resp. $\mathcal{A}_{1}$)
is called \emph{faithfully flat} if $T(-,P_{2})$ (resp. $T(P_{1},-)$)
preserves monomorphisms and is faithful. 
A monomorphism $s:B_{2}\to A_{2}$ of $\mathcal{A}_{2}$ is \emph{pure}
if $T(A_{1},s)$ is a monomorphism for all $A_{1}\in Ob(\mathcal{A}_{1})$.

\subsubsection{} \label{sec:3.4.1}
Let $P_{2}\in Ob(\mathcal{A}_{2})$. The following are equivalent:

(1) $P_{2}$ is faithfully flat;

(2) a sequence $(u_{1},v_{1})$ in $\mathcal{A}_{1}$ is exact if and 
only if the sequence $(T(u_{1},P_{2}),T(v_{1},P_{2}))$ is exact.
\begin{proof}
The proof is standard, using \ref{sec:3.1.1}.
\end{proof}

\subsubsection{} \label{sec:3.4.2}
Suppose that $\mathcal{A}_{3}$ has an injective cogenerator $E_{3}$ and that 
$\mathcal{A}_{1}$ has a faithfully flat object $P_{1}$. Suppose moreover that 
$\mathcal{A}_{2}$ is complete. Then $\mathcal{A}_{2}$ has enough injectives.
\begin{proof}
The object $C(P_{1},E_{3})$ is injective since $C(P_{1},-)$ has a left adjoint 
that is exact (\ref{sec:2.1.1}(2bis)). The object $C(P_{1},E_{3})$ is a 
cogenerator since $C(P_{1},-)$ is the right adjoint of a faithful functor,
so it preserves cogenerators. Let $A_{2}\in Ob(\mathcal{A}_{2})$. 
Consider the natural morphism
\[
c_{A_{2}}:A_{2}\to \underset{\mathcal{A}_{2}(A_{2},C(P_{1},E_{3}))}
\prod C(P_{1},E_{3}) 
\]
We have $\mathcal{A}_{2}(A_{2},C(P_{1},E_{3}))\cong 
\mathcal{A}_{3}(T(P_{1},A_{2}),E_{3})\neq 0$
if $A_{2}\neq 0$ since $T(P_{1},-)$ and $\mathcal{A}_{3}(-,E_{3})$ 
are faithful. The target of $c_{A_{2}}$ is injective and $c_{A_{2}}$ is 
a monomorphism since $C(P_{1},E_{3})$ is a cogenerator.
\end{proof}
If $\mathcal{A}_{2}$ is not complete, the following is a non-functorial 
variant of \ref{sec:3.4.2}.

\subsubsection{} \label{sec:3.4.3}
We first show that $H(-,E_{3})$ is faithful. Let $A_{2}\in Ob(\mathcal{A}_{2})$, 
$A_{2}\neq 0$. Suppose that $H(A_{2},E_{3})=0$. Then 
$0=\mathcal{A}_{1}(P_{1},H(A_{2},E_{3}))\cong
\mathcal{A}_{3}(T(P_{1},A_{2}),E_{3})$, which implies that 
$T(P_{1},A_{2})=0$ since $\mathcal{A}_{3}(-,E_{3})$ is faithful. 
But $T(P_{1},A_{2})\neq 0$ by the assumption on $P_{1}$,
contradiction. Therefore $H(-,E_{3})$ is faithful. It follows that for 
each $A_{2}\in Ob(\mathcal{A}_{2})$ the natural morphism 
$\epsilon_{A_{2}}:A_{2}\to C(H(A_{2},E_{3}),E_{3})$ is a monomorphism. 
Let $P_{1}'\to H(A_{2},E_{3})$ be an epimorphism with $P_{1}'$ flat. 
Then $C(H(A_{2},E_{3}),E_{3})\to C(P_{1}',E_{3})$ is a monomorphism 
since $C(-,E_{3})$ is a left adjoint. The object $C(P_{1}',E_{3})$ is 
injective since $C(P_{1}',-)$ has a left adjoint that is exact (\ref{sec:2.1.1}(2bis)). 
The composite monomorphism
\[
A_{2}\to C(H(A_{2},E_{3}),E_{3})\to C(P_{1}',E_{3})
\]
shows that $\mathcal{A}_{2}$ has enough injectives.

\subsubsection{} \label{sec:3.4.4}
Suppose that $\mathcal{A}_{1}$ has a faithfully flat 
object $P_{1}$ and that $\mathcal{A}_{3}$ has a set $(E^{i}_{3})_{i\in  I}$
of injective cogenerators. Let $(u_{2},v_{2})$ be a sequence in 
$\mathcal{A}_{2}$. If the sequence $(H(v_{2},E^{i}_{3}),H(u_{2},E^{i}_{3}))$
is exact and $\K H(v_{2},E^{i}_{3})\in\mathcal{F}'^{\perp}$ for all $i\in I$
then $(u_{2},v_{2})$ is exact.
\begin{proof}
Using \ref{sec:3.4.1}, \ref{sec:1.1.1} and adjunction we have
the equivalences 
\[
(u_{2},v_{2})\ {\rm exact} \Leftrightarrow (\mathcal{A}_{3}(T(P_{1},v_{2}),E^{i}_{3}),
\mathcal{A}_{3}(T(P_{1},u_{2}),E^{i}_{3})))\ {\rm exact\ for\ all}\ i\in  I
\]
\[
\Leftrightarrow (\mathcal{A}_{1}(P_{1},H(v_{2},E^{i}_{3})),
\mathcal{A}_{1}(P_{1},H(u_{2},E^{i}_{3}))\ {\rm exact\ for\ all}\ i\in  I
\]
The last sequence is exact by \ref{sec:1.2.1}(1bis).
\end{proof}

\subsubsection{} \label{sec:3.4.5}
Suppose that $\mathcal{A}_{1}$ has a faithfully flat 
object $P_{1}$ and that $\mathcal{A}_{3}$ has enough injectives. 
Let $(u_{2},v_{2})$ be a semi-exact sequence in $\mathcal{A}_{2}$. 
If the sequence $(H(v_{2},E_{3}),H(u_{2},E_{3}))$
is exact for all injective objects $E_{3}$ of $\mathcal{A}_{3}$
and $\K H(v_{2},E^{i}_{3})\in\mathcal{F}'^{\perp}$ 
for all $i\in I$ then $(u_{2},v_{2})$ is exact.
\begin{proof}
The proof is the same as the proof of \ref{sec:3.4.4}, 
using \ref{sec:1.1.2} instead of \ref{sec:1.1.1}.
\end{proof}

\subsubsection{} \label{sec:3.4.6}
Consider in $\mathcal{A}_{2}$ the commutative diagram 
with top row exact 
\[
\xymatrix{
{A_{2}} \ar[r]^{u_{2}} \ar[d]_{f}& {B_{2}}\ar@{->>}[r]^{v_{2}}
\ar[d]_{g} & {C_{2}}\ar[d]^{h}\\
{A'_{2}} \ar[r]_{u'_{2}} & {B'_{2}}\ar[r]_{v'_{2}} & {C'_{2}}
}
\]

(1) If $u'_{2},f$ and $h$ are pure monomorphisms then so is $g$.

(2) If $(u'_{2},v'_{2})$ is exact, $f,v'_{2}$ are epimorphisms 
and $g$ is a pure monomorphism then $h$ is a pure monomorphism.
\begin{proof}
Let $A_{1}\in Ob(\mathcal{A}_{1})$. For part 1 we apply $T(A_{1},-)$
to the diagram in the statement; we obtain the commutative diagram
with top row exact 
\[
\xymatrix{
{T(A_{1},A_{2})} \ar[r]\ar@{>->}[d]& {T(A_{1},B_{2})}\ar@{->>}[r]
\ar[d] & {T(A_{1},C_{2})}\ar@{>->}[d]\\
{T(A_{1},A'_{2})} \ar@{>->}[r] & {T(A_{1},B'_{2})}\ar[r]& {T(A_{1},C'_{2})}
}
\]
Now the assertion follows from the kernel exact sequence. Part 2 is similar. 
\end{proof}

\subsubsection{} \label{sec:3.4.7}
Let $S_{2}\to P_{2}$ be a monomorphism in $\mathcal{A}_{2}$
with $P_{2}$ flat.

(1) The monomorphism $S_{2}\to P_{2}$ is pure if and only if $P_{2}/S_{2}$ is flat.

(2) Suppose that $\mathcal{A}_{1},\mathcal{A}_{2},\mathcal{A}_{3}$ are
Grothendieck categories and $\mathcal{A}_{1}$ has a flat generator $U_{1}$.
The monomorphism $S_{2}\to P_{2}$ is pure if and only if for all monomorphisms
$S_{1}\to U_{1}$ the natural morphism $f:T(S_{1},S_{2})\to T(U_{1},S_{2})
\cap T(S_{1},P_{2})$ is an isomorphism.
\begin{proof}
(1) If $P_{2}/S_{2}$ is flat then we use \ref{sec:3.1.6}(1). Conversely, let
$A_{1}\rightarrowtail B_{1}\twoheadrightarrow C_{1}$ be an exact 
sequence in $\mathcal{A}_{1}$. The assertion follows upon applying 
the snake diagram to 
\[
\xymatrix{
{T(A_{1},S_{2})}\ar@{>->}[r]\ar@{>->}[d]& {T(B_{1},S_{2})}\ar@{->>}[r]
\ar@{>->}[d] & {T(C_{1},S_{2})}\ar@{>->}[d]\\
{T(A_{1},P_{2})} \ar@{>->}[r] & {T(B_{1},P_{2})}\ar@{->>}[r]& {T(C_{1},P_{2})}
}
\]
(2) The morphism $f:T(S_{1},S_{2})\to T(U_{1},S_{2})
\cap T(S_{1},P_{2})$ is obtained from the universal property 
of pullback applied to the diagram
\[
\xymatrix{
{T(S_{1},S_{2})}\ar[r]\ar[d]& {T(S_{1},P_{2})}\ar@{>->}[d]\\
{T(U_{1},S_{2})} \ar@{>->}[r] & {T(U_{1},P_{2})}
}
\]
Suppose $S_{2}\to P_{2}$ is pure. Then $P_{2}/S_{2}$ is flat by
part 1 and it follows that the previous diagram is a pullback, hence 
$f$ is an isomorphism. Conversely, consider the commutative diagram
with exact rows
\[
\xymatrix{
{T(S_{1},S_{2})}\ar@{>->}[r]\ar@{>->}[d]& {T(S_{1},P_{2})}\ar@{->>}[r]
\ar@{>->}[d] & {T(S_{1},P_{2}/S_{2})}\ar@{=}[d]\\
{T(U_{1},S_{2})} \ar@{>->}[r] \ar@{=}[d]& {T(U_{1},S_{2})\cup T(S_{1},P_{2})}
\ar@{>->}[d]\ar@{->>}[r]& {T(S_{1},P_{2}/S_{2})}\ar[d]\\
{T(U_{1},S_{2})}\ar@{>->}[r]&{T(U_{1},P_{2})}\ar@{->>}[r]&{T(U_{1},P_{2}/S_{2})}
}
\]
The lower right square diagram is a pushout, hence 
$T(S_{1},P_{2}/S_{2})\to T(U_{1},P_{2}/S_{2})$ is a monomorphism.
Using \ref{sec:3.1.2} and \ref{sec:1.2.16} with $\mathcal{G}=Sub(U_{1})$
we then obtain that $P_{2}/S_{2}$ is flat, hence $S_{2}\to P_{2}$ is pure
by part 1.
\end{proof}

\subsubsection{} \label{sec:3.4.8}
Suppose that $\mathcal{A}_{3}$ has an injective cogenerator $E_{3}$ 
and that $\mathcal{A}_{1}$ has a faithfully flat object $P_{1}$. For all 
$A_{2}\in Ob(\mathcal{A}_{2})$ the monomorphism 
$\epsilon_{A_{2}}:A_{2}\to C(H(A_{2},E_{3}),E_{3})$ (\ref{sec:3.4.3}) 
is pure. 
\begin{proof}
Since $(C(-,E_{3}),H(-,E_{3}))$ is an adjoint pair the composite
\[
\xymatrix{
{H(A_{2},E_{3}}\ar[r]& {H(C(H(A_{2},E_{3}),E_{3}),E_{3})}
\ar[rrr]^{H(\epsilon_{A_{2}},E_{3})}
&&&{H(A_{2},E_{3})}
}
\]
is the identity. Therefore for all $A_{1}\in Ob(\mathcal{A}_{1})$
the morphism $\mathcal{A}_{1}(A_{1},H(\epsilon_{A_{2}},E_{3}))$
is an epimorphism. Using adjunctions we have the commutative diagram
\[
\xymatrix{
{\mathcal{A}_{3}(T(A_{1},C(H(A_{2},E_{3}),E_{3})),E_{3})}
\ar[rrr]^{\mathcal{A}_{3}(T(A_{1},\epsilon_{A_{2}}),E_{3})}
\ar[d]&&&{\mathcal{A}_{3}(T(A_{1},A_{2}),E_{3})}\ar[d]\\
{\mathcal{A}_{1}(A_{1},H(C(H(A_{2},E_{3}),E_{3}),E_{3}))} 
\ar[rrr]_{\mathcal{A}_{1}(A_{1},H(\epsilon_{A_{2}},E_{3}))}
 &&&{\mathcal{A}_{1}(A_{1},H(A_{2},E_{3}))}
}
\]
in which the vertical arrows are isomorphisms. It follows that
$\mathcal{A}_{3}(T(A_{1},\epsilon_{A_{2}}),E_{3})$ is an 
epimorphism, hence (\ref{sec:1.1.1}(2)) $T(A_{1},\epsilon_{A_{2}})$
is a monomorphism.
\end{proof}

\section{Recollections on complexes}

In this section we recall some facts about (chain) complexes in abelian categories; 
for more details we refer to \cite{Bou},\cite{EM},\cite{Gi1},\cite{KP},\cite{Pe}. 
We adopt Bourbaki's notations, terminology and sign conventions, except that 
\emph{homotopisme} is translated as homotopy equivalence and 
\emph{homologisme} as quasi-isomorphism.

\subsection{} \label{sec:4.1}
Let $\mathcal{A}$ be an abelian category. We denote by $Gr(\mathcal{A})$
the category of $\mathbb{Z}$-graded objects in $\mathcal{A}$  and by 
$Ch(\mathcal{A})$ the category of complexes in $\mathcal{A}$.

Let $p$ be an integer. For a $\mathbb{Z}$-graded object $\mathbf{Z}$ in 
$\mathcal{A}$ we define the $\mathbb{Z}$-graded object $\mathbf{Z}(p)$
as $\mathbf{Z}(p)_{n}=\mathbf{Z}_{n+p}$, and if $\mathbf{Z},\mathbf{Z}'$ 
are $\mathbb{Z}$-graded  objects in $\mathcal{A}$ and 
$u:\mathbf{Z}\to\mathbf{Z}'$ is a morphism of $\mathbb{Z}$-graded 
objects we define $u(p)_{n}=u_{n+p}$. For a complex $(\mathbf{C},d)$ in 
$\mathcal{A}$ we denote by $(\mathbf{C}(p),d(p))$ the $p$-th 
translate of $(\mathbf{C},d)$, where $d(p)_{n}=(-1)^{p}d_{n+p}$. 
The differential $d$ is a morphism $d:\mathbf{C}\to\mathbf{C}(-1)$ and
a morphism $u:(\mathbf{C},d)\to (\mathbf{C}',d')$ is a morphism of 
$\mathbb{Z}$-graded objects $u:\mathbf{C}\to\mathbf{C}'$ such that 
$u(-1)d=d'u$.

Let $(\mathbf{C},d)$ be a complex in $\mathcal{A}$. We denote by 
$d^{\#}:\mathbf{C}(1)\to \mathbf{C}$ the adjoint transpose of $d$,
by $B(\mathbf{C})$ the complex of boundaries of $\mathbf{C}$ and by 
$Z(\mathbf{C})$ the complex of cycles. There are natural morphisms
$\delta(1):\mathbf{C}(1)\to B(\mathbf{C}), i:B(\mathbf{C})\to Z(\mathbf{C}), 
j:Z(\mathbf{C})\to \mathbf{C}, p:Z(\mathbf{C})\to H(\mathbf{C})$ and 
exact sequences
\[
Z(\mathbf{C})\overset{j}\rightarrowtail \mathbf{C}\overset{\delta}
\twoheadrightarrow B(\mathbf{C})(-1)\ \ \ \ H(\mathbf{C})\rightarrowtail 
\mathbf{C}/B(\mathbf{C}) \overset{\bar{\delta}}
\twoheadrightarrow B(\mathbf{C})(-1)
\]
For each integer $n$ we denote by $(-)_{n}:Ch(\mathcal{A})\to \mathcal{A}$
the evaluation at $n$ functor and by $D^{n}:\mathcal{A}\to Ch(\mathcal{A})$ 
the $n$th disk functor; the functor $(-)_{n}$ is left adjoint to $D^{n+1}$ 
and right adjoint to $D^{n}$.

The complex of cycles is a functor $Z:Ch(\mathcal{A})\to Ch(\mathcal{A})$ 
that has a left adjoint $-/B$ which sends a complex $\mathbf{C}$ to the
quotient object $\mathbf{C}/B(\mathbf{C})$. The functor $Z$ preserves 
epimorphisms with exact kernel and the functor $-/B$ preserves 
monomorphisms with exact cokernel. We have $-/BD^{n}=ZD^{n+1}$ 
and the two composites are the $n$th sphere functor $S^{n}$, therefore
$S^{n}$ has $(-)_{n}/B_{n}$ as left adjoint and $Z_{n}$ as right adjoint.

The underlying $\mathbb{Z}$-graded object faithful functor $U:Ch(\mathcal{A})
\to Gr(\mathcal{A})$ has a left adjoint $L$ defined as $L(\mathbf{Z})=
\mathbf{Z}(1)\oplus\mathbf{Z}$, with differential 
$\left[
\begin{array}{rr}
0 &  1_{\mathbf{Z}} \\
0   &  0 \\
\end{array}
\right]$.
The functor $U$ has also a right adjoint $R$ defined as $R(\mathbf{Z})=
\mathbf{Z}(-1)\oplus\mathbf{Z}$, with differential 
$\left[
\begin{array}{rr}
0 &  0 \\
1_{\mathbf{Z}(-1)}  &  0 \\
\end{array}
\right]$. By construction, $L$ and $R$ are exact and the complexes
$L(\mathbf{Z}),R(\mathbf{Z})$ are exact with $B(L(\mathbf{Z}))=\mathbf{Z}(1),
B(R(\mathbf{Z}))=\mathbf{Z}$.

\subsubsection{} \label{sec:4.1.1}
Let $u:\mathbf{C}'\to \mathbf{C}$ be a morphism. Then $H(u)$ is a  
monomorphism if and only if the diagram
\begin{displaymath}
\xymatrix{
{B(\mathbf{C}')}\ar[r]^{i'} \ar[d]&{Z(\mathbf{C}')} \ar[d]\\
{B(\mathbf{C})} \ar[r]^{i}&{Z(\mathbf{C})}
}
\end{displaymath} 
is a pullback and $H(u)$ is an isomorphism if and only if the 
previous diagram is exact (\ref{sec:1.2.8}). 
If $u$ is a monomorphism then $H(u)$ is an isomorphism
if and only if $B(\mathbf{C}')=B(\mathbf{C})\cap Z(\mathbf{C}')$
and $Z(\mathbf{C})=B(\mathbf{C})\cup Z(\mathbf{C}')$.
Dually, $H(u)$ is an epimorphism if and only if the diagram
\begin{displaymath}
\xymatrix{
{\mathbf{C}'/B(\mathbf{C}')}\ar[r]^{\bar{\delta'}} 
\ar[d]&{B(\mathbf{C}')(-1)} \ar[d]\\
{\mathbf{C}/B(\mathbf{C})} \ar[r]^{\bar{\delta}}&{B(\mathbf{C})(-1)}
}
\end{displaymath} 
is a pushout and $H(u)$ is an isomorphism if and only if the 
previous diagram is exact.

\subsubsection{} \label{sec:4.1.2}
\cite[Proposition 3.2]{Gi2} For a class $\mathcal{K}$ of objects of 
$\mathcal{A}$ we denote by $Gr(\mathcal{K})$ (resp. $Ch(\mathcal{K})$) 
the class of $\mathbb{Z}$-graded objects (resp. complexes) that are 
degreewise in $\mathcal{K}$. We denote by $ex[\mathcal{K}]$ the 
class of exact complexes in $\mathcal{A}$ that have cycles in $\mathcal{K}$.

Let $(\mathcal{K},\mathcal{K}^{\perp})$ be a cotorsion theory in $\mathcal{A}$.
It is clear that $(Gr(\mathcal{K}),Gr(\mathcal{K}^{\perp}))$ is a cotorsion theory 
in $Gr(\mathcal{A})$. Then, using the adjoint pairs $(L,U)$ and $(U,R)$ (\ref{sec:4.1})
we obtain from \ref{sec:2.1.2}(1) and \ref{sec:2.1.2}(1bis) the cotorsion theories
$(Ch(\mathcal{K}),Ch(\mathcal{K})^{\perp})$ and 
$(^{\perp}Ch(\mathcal{K}^{\perp}),Ch(\mathcal{K}^{\perp}))$ in $Ch(\mathcal{A})$.
See \ref{sec:8.2.5}(1) for a description of $Ch(\mathcal{K})^{\perp}$. The 
cotorsion theory $(Ch(\mathcal{K}),Ch(\mathcal{K})^{\perp})$ is generated 
(\ref{sec:1.2.3}) by $\{D^{n}(Y)\}$, where $Y\in \mathcal{K}^{\perp}$ and 
$n$ is an integer; the cotorsion theory 
$(^{\perp}Ch(\mathcal{K}^{\perp}),Ch(\mathcal{K}^{\perp}))$ is cogenerated 
by $\{D^{n}(X)\}$, where $X\in \mathcal{K}$ and $n$ is an integer.
Suppose, moreover, that $\mathcal{A}$ has enough $\mathcal{K}$ objects.
Then $Gr(\mathcal{A})$ has enough $Gr(\mathcal{K})$ objects 
and $Ch(\mathcal{A})$ has enough $ex[\mathcal{K}]$ objects, meaning that
for each complex $\mathbf{C}$ in $\mathcal{A}$ there is an epimorphism 
$\mathbf{X}\to\mathbf{C}$ with $\mathbf{X}\in ex[\mathcal{K}]$. In particular,
$Ch(\mathcal{A})$ has enough $Ch(\mathcal{K})$ objects.

\subsubsection{} \label{sec:4.1.3}
Let $(\mathbf{C},d),(\mathbf{C}',d')$ be complexes in $\mathcal{A}$.
Two morphisms $f,g$ from $\mathbf{C}$ to $\mathbf{C}'$ are homotopic,
written $f\sim g$, if there is a morphism of $\mathbb{Z}$-graded objects
$s:\mathbf{C}\to\mathbf{C}'(1)$ such that $g-f=d's+sd$. A complex is 
homotopic to zero if its identity morphism is homotopic to the zero morphism.

Let $u:\mathbf{C}'\to\mathbf{C}$ be a morphism. 
We denote by $Con(u)$ (resp. $Cyl(u)$) the cone 
(resp. cylinder) of $u$. There are exact sequences
\[
\mathbf{C}'\overset{\tilde{u}}\rightarrowtail 
Cyl(u)\overset{\tilde{\pi}}\twoheadrightarrow Con(u)\ {\rm and}\ 
\mathbf{C}\overset{\pi}\rightarrowtail 
Con(u)\overset{\delta}\twoheadrightarrow\mathbf{C}'(-1)
\]
The morphism $u$ factors into $\mathbf{C}'\overset{\tilde{u}}\to 
Cyl(u)\overset{\beta}\to \mathbf{C}$, where the morphism $\tilde{u}$ 
is a monomorphism and the morphism $\beta$  is a homotopy 
equivalence that has a right inverse $\alpha$ such that 
$\pi=\tilde{\pi}\alpha$.

If the complex $Con(u)$ homotopic to zero then $u$ is a 
homotopy equivalence.

For every integer $p$ we have a natural isomorphism 
$Con(u(p))\cong Con(u)(p)$ defined as
\[
((-1)^{p}1_{\mathbf{C}'_{p+n-1}},1_{\mathbf{C}_{p+n}}):
\mathbf{C}'_{p+n-1}\oplus \mathbf{C}_{p+n}\to
\mathbf{C}'_{p+n-1}\oplus\mathbf{C}_{p+n}
\]

\subsection{} \label{sec:4.2}
Let $\mathcal{A}$ be an abelian category and let 
$(\mathbf{C},d),(\mathbf{C}',d')$ be complexes in $\mathcal{A}$.
We define the complex $\mathrm{Homgr}_{\mathcal{A}}(\mathbf{C}',\mathbf{C})$
in $\mathrm{Ab}$ as $\mathrm{Homgr}_{\mathcal{A}}(\mathbf{C}',\mathbf{C})_{n}=
\underset{p}\prod \mathcal{A}(\mathbf{C}'_{p},\mathbf{C}_{p+n})$ with differential 
$D(\mathbf{C}',\mathbf{C})_{n}=D_{n}$ defined as $D_{n}(f)=df-(-1)^{n}fd'$.
There are natural isomorphisms
\[
Ch(\mathrm{Ab})(S^{n}(\mathbb{Z}),\mathrm{Homgr}_{\mathcal{A}}(\mathbf{C}',\mathbf{C}))
\cong Z_{n}(\mathrm{Homgr}_{\mathcal{A}}(\mathbf{C}',\mathbf{C}))
\cong Ch(\mathcal{A})(\mathbf{C}',\mathbf{C}(n))
\]
For every integer $r$ we have
\[
(\mathrm{Homgr}_{\mathcal{A}}(\mathbf{C}',\mathbf{C})(r),D(\mathbf{C}',\mathbf{C})(r))=
(\mathrm{Homgr}_{\mathcal{A}}(\mathbf{C}',\mathbf{C}(r)),D(\mathbf{C}',\mathbf{C}(r)))
\]
and a natural isomorphism
\[
\Sigma(\mathbf{C}',\mathbf{C}):\mathrm{Homgr}_{\mathcal{A}}(\mathbf{C}'(r),\mathbf{C})
\cong \mathrm{Homgr}_{\mathcal{A}}(\mathbf{C}',\mathbf{C})(-r)
\]
constructed as follows. Let $n,p$ be integers and $f_{p}:\mathbf{C}'_{p+r}\to\mathbf{C}_{p+n}$ 
a morphism. We define $\Sigma(\mathbf{C}',\mathbf{C})_{n}((f_{p}))=((-1)^{nr}f_{p-r}))$. 
One can check that $\Sigma(\mathbf{C}',\mathbf{C})$ is an isomorphism.

For every integer $n$ we have a natural isomorphism 
\[
Ch(\mathcal{A})(\mathbf{C}',\mathbf{C}(n))/\sim \cong 
H_{n}(\mathrm{Homgr}_{\mathcal{A}}(\mathbf{C}',\mathbf{C}))
\]
Let $u:\mathbf{C}'\to\mathbf{C}$ be a morphism and $\mathbf{C}''$ 
another complex. We have 
\[
\mathrm{Homgr}_{\mathcal{A}}(\mathbf{C''},Con(u))=
Con(\mathrm{Homgr}_{\mathcal{A}}(\mathbf{C}'',u))
\]
and a natural isomorphism 
\[
\Sigma(u,\mathbf{C}''):\mathrm{Homgr}_{\mathcal{A}}(Con(u)(1),\mathbf{C}'')
\cong Con(\mathrm{Homgr}_{\mathcal{A}}(u,\mathbf{C}''))
\]
constructed as follows. We have 
\[
\mathrm{Homgr}_{\mathcal{A}}(Con(u)(1),\mathbf{C}'')_{n}=
\underset{p}\prod \mathcal{A}(\mathbf{C}'_{p}\oplus\mathbf{C}_{p+1},\mathbf{C}''_{p+n})
\]
and
\[
Con(\mathrm{Homgr}_{\mathcal{A}}(u,\mathbf{C}''))_{n}=
\underset{i}\prod \mathcal{A}(\mathbf{C}_{i},\mathbf{C}''_{i+n-1})\oplus
\underset{j}\prod \mathcal{A}(\mathbf{C}'_{j},\mathbf{C}''_{j+n})
\]
We define 
\[
\Sigma(u,\mathbf{C}'')_{n}((g_{p})_{p})=(((-1)^{n}g_{i-1}\sigma_{2})_{i},(g_{j}\sigma_{1})_{j})
\]
where $\sigma_{1}:\mathbf{C}'_{p}\to\mathbf{C}'_{p}\oplus\mathbf{C}_{p+1}$
and $\sigma_{2}:\mathbf{C}_{p+1}\to\mathbf{C}'_{p}\oplus\mathbf{C}_{p+1}$
are coprojection morphisms. One can check that $\Sigma(u,\mathbf{C}'')$ is a morphism.

Let $r$ be an integer and $A\in Ob(\mathcal{A})$. Then
\[
\mathrm{Homgr}_{\mathcal{A}}(S^{r}(A),\mathbf{C})
=\mathcal{A}(A,\mathbf{C}_{r+})
\]
where $\mathcal{A}(A,\mathbf{C}_{r+})_{n}=
\mathcal{A}(A,\mathbf{C}_{r+n})$ has differential $\mathcal{A}(A,d_{r+n})$,
and 
\[
\mathrm{Homgr}_{\mathcal{A}}(\mathbf{C}',S^{r}(A))=
\mathcal{A}(\mathbf{C}'_{r-},A)
\]
where $\mathcal{A}(\mathbf{C}'_{r-},A)_{n}=
\mathcal{A}(\mathbf{C}'_{r-n},A)$ 
has differential $\mathcal{A}((-1)^{n+1}d'_{r-n+1},A)$.

\subsubsection{} \label{sec:4.2.1}
The functors $\mathrm{Homgr}_{\mathcal{A}}(\mathbf{C}',-)$ 
and $\mathrm{Homgr}_{\mathcal{A}}(-,\mathbf{C})$ 
send degreewise split short exact sequences to short exact sequences.

\subsubsection{} \label{sec:4.2.2}
Let $u:\mathbf{C}'\to\mathbf{C}$ be a morphism and $\mathbf{C}''$
a complex. Consider the degreewise split exact sequence (\ref{sec:4.1.3}) 
$\mathbf{C}\overset{\pi}\rightarrowtail Con(u)\overset{\delta}
\twoheadrightarrow\mathbf{C}'(-1)$ and translate it (\ref{sec:4.1})
we obtain the exact sequence $\mathbf{C}(1)\overset{\pi(1)}
\rightarrowtail Con(u)(1)\overset{\delta(1)}\twoheadrightarrow\mathbf{C}'$. 
Applying $\mathrm{Homgr}_{\mathcal{A}}(-,\mathbf{C}'')$ to the 
previous sequence we obtain  (\ref{sec:4.2.1}, \ref{sec:4.2})
the commutative diagram
\begin{displaymath}
\xymatrix{
{\mathrm{Homgr}_{\mathcal{A}}(\mathbf{C}',\mathbf{C}'')}
\ar@{>->}[rr]^{\overline{\delta(1)}}
\ar@{=}[d]&&{\mathrm{Homgr}_{\mathcal{A}}(Con(u)(1),\mathbf{C}'')}
\ar[d]_{\Sigma(u,\mathbf{C}'')}\ar@{->>}[rr]^{\overline{\pi(1)}}
&&{\mathrm{Homgr}_{\mathcal{A}}(\mathbf{C}(1),\mathbf{C}'')}
\ar[d]_{\Sigma(\mathbf{C},\mathbf{C}'')}\\
{\mathrm{Homgr}_{\mathcal{A}}(\mathbf{C}',\mathbf{C}'')}\ar@{>->}[rr]^{\pi'}&&
{Con(\mathrm{Homgr}_{\mathcal{A}}(u,\mathbf{C}''))}\ar@{->>}[rr]^{\delta'}&&
{\mathrm{Homgr}_{\mathcal{A}}(\mathbf{C},\mathbf{C}'')(-1)}
}
\end{displaymath}
where $\overline{\delta(1)}=\mathrm{Homgr}_{\mathcal{A}}(\delta(1),\mathbf{C}'')$,
and $\overline{\pi(1)}=\mathrm{Homgr}_{\mathcal{A}}(\pi(1),\mathbf{C}'')$.
By the naturality of the connecting morphism 
\cite[\S2 n\textsuperscript{\scriptsize{o}} 3 prop. 2]{Bou}
we have the commutative diagram
\begin{displaymath}
\xymatrix{
{H(\mathrm{Homgr}_{\mathcal{A}}(\mathbf{C}(1),\mathbf{C}''))}
\ar[rr]^{\partial(\overline{\delta(1)},\overline{\pi(1)})}
\ar[d]_{H(\Sigma(\mathbf{C},\mathbf{C}''))} && 
{H(\mathrm{Homgr}_{\mathcal{A}}(\mathbf{C}',\mathbf{C}''))(-1)}\ar@{=}[d]\\
{H(\mathrm{Homgr}_{\mathcal{A}}(\mathbf{C},\mathbf{C}'')(-1))}\ar[rr]^{\partial(\pi',\delta')}
&& {H(\mathrm{Homgr}_{\mathcal{A}}(\mathbf{C}',\mathbf{C}''))(-1)}
}
\end{displaymath}
We have $\partial(\pi',\delta')=-H(\mathrm{Homgr}_{\mathcal{A}}(u,\mathbf{C}'')(-1))$ 
\cite[\S2 n\textsuperscript{\scriptsize{o}} 6 lemme 3a)]{Bou} and 
$\partial(\pi',\delta')$ is an isomorphism if and only if 
$H(\mathrm{Homgr}_{\mathcal{A}}(Con(u)(1),\mathbf{C}''))=0$
\cite[\S2 n\textsuperscript{\scriptsize{o}} 3 cor. 1]{Bou}.

\subsubsection{} \label{sec:4.2.3}
\cite[Lemma 3.4]{Gi1} Suppose that $\mathbf{C}'$ is right bounded, 
that is, there is an integer $i$ such that $\mathbf{C}'_{n}=0$ for $n<i$. 
For each integer $r\geqslant i$ we define the subcomplex 
$\mathrm{Homgr}_{\mathcal{A}}^{(r)}(\mathbf{C}',\mathbf{C})$ 
of $\mathrm{Homgr}_{\mathcal{A}}(\mathbf{C}',\mathbf{C})$ as
\[
\mathrm{Homgr}_{\mathcal{A}}^{(r)}(\mathbf{C}',\mathbf{C})_{n}=
\underset{i\leqslant p\leqslant r}\prod \mathcal{A}(\mathbf{C}'_{p},\mathbf{C}_{p+n})
\]
with differential induced by $D(\mathbf{C}',\mathbf{C})$. The complex
$\mathrm{Homgr}_{\mathcal{A}}^{(i)}(\mathbf{C}',\mathbf{C})=
\mathcal{A}(\mathbf{C}'_{i},\mathbf{C})$ has differential 
$\mathcal{A}(\mathbf{C}'_{i},d)$. Let 
$t_{r}:\mathrm{Homgr}_{\mathcal{A}}^{(r)}(\mathbf{C}',\mathbf{C})\to
\mathrm{Homgr}_{\mathcal{A}}^{(r-1)}(\mathbf{C}',\mathbf{C})$ be the 
natural projection epimorphism $t_{r}((f_{p})_{i\leqslant p\leqslant r})
=(f_{p})_{i\leqslant p\leqslant r-1}$; we have $\K t_{r}=
\mathcal{A}(\mathbf{C}'_{r},\mathbf{C}_{r+})$, where 
$\mathcal{A}(\mathbf{C}'_{r},\mathbf{C}_{r+})_{n}=
\mathcal{A}(\mathbf{C}'_{r},\mathbf{C}_{r+n})$ has 
differential $\mathcal{A}(\mathbf{C}'_{r},d_{r+n})$.
We have a natural isomorphism 
\[
\mathrm{Homgr}_{\mathcal{A}}(\mathbf{C}',\mathbf{C})\cong
\underset{r\in(\mathbb{Z}_{\geqslant i})^{op}}\li 
\mathrm{Homgr}_{\mathcal{A}}^{(r)}(\mathbf{C}',\mathbf{C})
\]
hence an exact sequence 
\[
\mathrm{Homgr}_{\mathcal{A}}(\mathbf{C}',\mathbf{C})\rightarrowtail 
\underset{r\geqslant i}\prod\mathrm{Homgr}_{\mathcal{A}}^{(r)}
(\mathbf{C}',\mathbf{C})\twoheadrightarrow \underset{r\geqslant i}
\prod\mathrm{Homgr}_{\mathcal{A}}^{(r)}(\mathbf{C}',\mathbf{C})
\]
Suppose, in addition, that the complexes $\mathcal{A}(\mathbf{C}'_{r},
\mathbf{C}_{r+})$ and $\mathcal{A}(\mathbf{C}'_{i},\mathbf{C})$ 
are exact. Induction over $r\geqslant i$ implies then that 
$\mathrm{Homgr}_{\mathcal{A}}^{(r)}(\mathbf{C}',\mathbf{C})$ 
is exact for all $r\geqslant i$, therefore 
$\mathrm{Homgr}_{\mathcal{A}}(\mathbf{C}',\mathbf{C}) $ is exact.
In particular, let $\mathcal{K}$ be a class of objects of $\mathcal{A}$ 
and suppose $\mathbf{C}'_{n}\in \mathcal{K}$ for $n\geqslant i$ and 
$\mathbf{C}\in ex[\mathcal{K}^{\perp}]$. By \ref{sec:1.2.1}(1bis) 
the complexes $\mathcal{A}(\mathbf{C}'_{r},\mathbf{C}_{r+})$ 
and $\mathcal{A}(\mathbf{C}'_{i},\mathbf{C})$ are exact, therefore
$\mathrm{Homgr}_{\mathcal{A}}(\mathbf{C}',\mathbf{C}) $ is exact.

\subsubsection{} \label{sec:4.2.4}
Let $k$ be an integer. We denote by $\tau_{k}(\mathbf{C}')$ the complex
\[
...\to\mathbf{C}'_{k+2}\to\mathbf{C}'_{k+1}\to Z_{k}(\mathbf{C}')\to 0\to ...
\]
where $Z_{k}(\mathbf{C}')$ is in degree $k$. There are natural 
morphisms $\tau_{k}(\mathbf{C}')\to\tau_{k-1}(\mathbf{C}')$
and for every integer $n$ the direct system $\tau_{k}(\mathbf{C}')_{n}
\to\tau_{k-1}(\mathbf{C}')_{n}\to \tau_{k-2}(\mathbf{C}')_{n}\to ...$
is eventually equal to $\mathbf{C}'_{n}$; hence 
$\underset{i\in(\mathbb{Z}_{\leqslant k})^{op}}
\cl \tau_{i}(\mathbf{C}')=\mathbf{C}'$ and we have 
\[
\mathrm{Homgr}_{\mathcal{A}}(\mathbf{C}',\mathbf{C})=
\underset{i\in\mathbb{Z}_{\leqslant k}}
\li \mathrm{Homgr}_{\mathcal{A}}(\tau_{i}(\mathbf{C}'),
\mathbf{C})
\]
For every integer $n$ the inverse system 
$(\mathrm{Homgr}_{\mathcal{A}}(\tau_{i}(\mathbf{C}'),
\mathbf{C})_{n})_{i\in\mathbb{Z}_{\leqslant k}}$
satisfies the Mittag-Leffler condition.
\begin{proof}
Let $l\leqslant i\leqslant k$ and let 
\[
t_{i,l}:\mathrm{Homgr}_{\mathcal{A}}(\tau_{l}(\mathbf{C}'),
\mathbf{C})_{n}\to\mathrm{Homgr}_{\mathcal{A}}(\tau_{i}(\mathbf{C}'),
\mathbf{C})_{n}
\]
be the natural morphism. We have $\mathrm{Homgr}_{\mathcal{A}}(\tau_{i}(\mathbf{C}'),
\mathbf{C})_{n}=\underset{p\geqslant i}\prod 
\mathcal{A}(\tau_{i}(\mathbf{C}')_{p},\mathbf{C}_{p+n})$ and 
\[
t_{i,l}((f^{l}_{p})_{p})=(f^{l}_{i}j'_{i},f^{l}_{i+1},f^{l}_{i+2},...)
\]
It follows that $\im t_{i,i-1}=\im t_{i,i-2}=...$, hence the assertion.
\end{proof}

\subsubsection{} \label{sec:4.2.5}
\cite[Prop. 13.2.3]{GD} Let $F:\mathcal{A}'\to\mathcal{A}$ be a left exact 
functor between abelian categories. For every complex $\mathbf{D}$
in $\mathcal{A}'$ there is a natural morphism
\begin{equation}\label{eq:4.2.5}
\xymatrix{
{H(F(\mathbf{D}))}\ar[r]^{\lambda(\mathbf{D})}& {F(H(\mathbf{D}))}
}
\end{equation} 
Let $\mathcal{K}'$ be a class of objects of $\mathcal{A}'$ such that 
$F$ preserves epimorphisms with kernel in $\mathcal{K}'$. If
$B(\mathbf{D}),Z(\mathbf{D})\in Ch(\mathcal{K}')$ then
(\ref{eq:4.2.5}) is an isomorphism. We note that 

(i) $B(\mathbf{D})\in Ch(\mathcal{K}')$ if $\mathbf{D}\in Ch(\mathcal{K}')$ 
and $\mathcal{K}'$ is closed under quotients, and 

(ii) $Z(\mathbf{D})\in Ch(\mathcal{K}')$ if $H(\mathbf{D})\in Ch(\mathcal{K}')$ 
and $\mathcal{K}'$ is closed under extensions.

\subsubsection{} \label{sec:4.2.6}
\cite[Lemma 3.3]{EGR},\cite[Lemma 3.9]{Gi1}
Let $k$ be an integer. Supose that $\mathcal{K}$ is a class of 
objects of $\mathcal{A}$, $\mathbf{C}',Z(\mathbf{C}')\in Ch(\mathcal{K})$ 
and $\mathbf{C}\in ex[\mathcal{K}^{\perp}]$. By \ref{sec:4.2.3}
applied to $\tau_{i}(\mathbf{C}')$ instead of $\mathbf{C}'$
we obtain that $\mathrm{Homgr}_{\mathcal{A}}(\tau_{i}(\mathbf{C}'),\mathbf{C}) 
$ is exact for all $i\leqslant k$. We apply \ref{sec:4.2.5} to 
$\mathcal{A}'=\mathrm{Ab}^{\mathbb{N}^{op}}, \mathcal{A}=\mathrm{Ab},
F=\underset{i\in\mathbb{Z}_{\leqslant k}}\li, \mathbf{D}_{i}=
\mathrm{Homgr}_{\mathcal{A}}(\tau_{i}(\mathbf{C}'),\mathbf{C})$
and $\mathcal{K}'$ equal to the class of inverse systems that 
satisfy the Mittag-Leffler condition; we obtain from 
\cite[prop. 13.2.1 et 13.2.2]{GD} and \ref{sec:4.2.4} 
that $\mathrm{Homgr}_{\mathcal{A}}
(\mathbf{C}',\mathbf{C})$ is exact. 

See \ref{sec:9.1.5}(2) for a different approach to the exactness of 
$\mathrm{Homgr}_{\mathcal{A}}(\mathbf{C}',\mathbf{C})$.
%$\underset{i\in\mathbb{Z}_{\leqslant k}}\li
%\mathrm{Homgr}_{\mathcal{A}}(\mathbf{X}^{i},\mathbf{D})$
\subsubsection{} \label{sec:4.2.7}\cite{EM}
Suppose that $\mathcal{A}$ has countable coproducts.
Let $k$ be an integer and 
\[
\mathbf{C}^{k}\overset{u^{k}}\rightarrow \mathbf{C}^{k-1}\overset{u^{k-1}}
\rightarrow \mathbf{C}^{k-2}\rightarrow ...
\]
a direct system of complexes in $\mathcal{A}$. Let $\sigma^{i}:
\mathbf{C}^{i}\to\underset{i\leqslant k}\bigoplus\mathbf{C}^{i}$ 
be coprojection morphisms and let $\Phi(\mathbf{C}):
\underset{i\leqslant k}\bigoplus\mathbf{C}^{i}\to
\underset{i\leqslant k}\bigoplus\mathbf{C}^{i}$ be the 
unique morphism such that $\Phi(\mathbf{C})\sigma^{i}=
\sigma^{i-1}u^{i}-\sigma^{i}$ for $i\leqslant k$. 
We have the exact sequence
\[
\xymatrix{
{\underset{i\leqslant k}\bigoplus\mathbf{C}^{i}}\ar[rr]^{\Phi(\mathbf{C})}
&&{\underset{i\leqslant k}\bigoplus\mathbf{C}^{i}}\ar@{->>}[rr]^{\pi(\mathbf{C})}
&&{\underset{i\in(\mathbb{Z}_{\leqslant k})^{op}}\cl \mathbf{C}^{i}}
}
\]
Let $\sigma^{H,i}:H(\mathbf{C}^{i})\to\underset{i\leqslant k}\bigoplus 
H(\mathbf{C}^{i})$ be coprojection morphisms and let $\Phi(H(\mathbf{C})):
\underset{i\leqslant k}\bigoplus H(\mathbf{C}^{i})\to
\underset{i\leqslant k}\bigoplus H(\mathbf{C}^{i})$ be the 
unique morphism such that $\Phi(H(\mathbf{C}))\sigma^{H,i}=
\sigma^{H,i-1}H(u^{i})-\sigma^{H,i}$ for $i\leqslant k$. 
We have the exact sequence
\[
\xymatrix{
{\underset{i\leqslant k}\bigoplus H(\mathbf{C}^{i})}\ar[rr]^{\Phi(H(\mathbf{C}))}
&&{\underset{i\leqslant k}\bigoplus H(\mathbf{C}^{i})}\ar@{->>}[rr]^{\pi(H(\mathbf{C}))}
&&{\underset{i\in(\mathbb{Z}_{\leqslant k})^{op}}\cl H(\mathbf{C}^{i})}
}
\]
Let $h:\underset{i\leqslant k}\bigoplus H(\mathbf{C}^{i})\to
H(\underset{i\leqslant k}\bigoplus \mathbf{C}^{i})$ be the unique
morphism such that $h\sigma^{H,i}=H(\sigma^{i})$ for $i\leqslant k$.
We have the commutative solid arrows diagram
\[
\xymatrix{
{\underset{i\leqslant k}\bigoplus H(\mathbf{C}^{i})}\ar[rr]^{\Phi(H(\mathbf{C}))}
\ar[d]^{h}&&{\underset{i\leqslant k}\bigoplus H(\mathbf{C}^{i})}
\ar@{->>}[rr]^{\pi(H(\mathbf{C}))}\ar[d]^{h} && 
{\underset{i\in(\mathbb{Z}_{\leqslant k})^{op}}\cl H(\mathbf{C}^{i})}
\ar@{..>}[d]^{h'}\\
{H(\underset{i\leqslant k}\bigoplus \mathbf{C}^{i})} \ar[rr]^{H(\Phi(\mathbf{C}))}
&& {H(\underset{i\leqslant k}\bigoplus \mathbf{C}^{i})}\ar[rr]^{H(\pi(\mathbf{C}))}
&& {H(\underset{i\in(\mathbb{Z}_{\leqslant k})^{op}}\cl \mathbf{C}^{i})}
}
\]
from which we obtain the unique dotted arrow $h'$ such that
$h'\pi(H(\mathbf{C}))=H(\pi(\mathbf{C}))h$. Suppose now that
$\mathcal{A}$ satisfies AB4 and that $u^{i}$ is a degreewise split 
monomorphism for $i\leqslant k$. Then $\Phi(\mathbf{C})$ is a 
monomorphism. If, moreover, $H(u^{i})$ is a degreewise split 
monomorphism for $i\leqslant k$ then $h'$ is an isomorphism.
\begin{proof}
Let $n$ be an integer and $r^{i}_{n}:\mathbf{C}^{i-1}_{n}\to
\mathbf{C}^{i}_{n}$ a retract of $u^{i}_{n}$, $i\leqslant k$. 
We inductively define 
$t^{i}_{n}:\mathbf{C}^{i}_{n}\to\underset{i\leqslant k}
\bigoplus \mathbf{C}^{i}_{n}$ as $t^{k}_{n}=0, t^{i-1}_{n}=
(\sigma^{i}_{n}+t^{i}_{n})r^{i}_{n}$ for $i<k$. Let 
$r(\mathbf{C})_{n}:\underset{i\leqslant k}\bigoplus \mathbf{C}^{i}
\to\underset{i\leqslant k}\bigoplus \mathbf{C}^{i}$ be the 
unique morphism such that $r(\mathbf{C})_{n}\sigma^{i}_{n}=t^{i}_{n}$
for $i\leqslant k$. Then $r(\mathbf{C})_{n}
\Phi(\mathbf{C})_{n}\sigma^{i}_{n}=
r(\mathbf{C})_{n}\sigma^{i-1}_{n}u^{i}_{n}-r(\mathbf{C})_{n}
\sigma^{i}_{n}=t^{i-1}_{n}u^{i}_{n}-t^{i}_{n}=\sigma^{i}_{n}$, 
therefore $r(\mathbf{C})_{n}$ is a retract of $\Phi(\mathbf{C})_{n}$.
Suppose that $H(u^{i})$ is a degreewise split monomorphism; 
then the same reasoning as above shows that $\Phi(H(\mathbf{C}))$ 
is a monomorphism. Let 
\[
\partial(\Phi(\mathbf{C}),\pi(\mathbf{C})):
H(\underset{i\in(\mathbb{Z}_{\leqslant k})^{op}}\cl \mathbf{C}^{i})
\to \underset{i\leqslant k}H(\bigoplus\mathbf{C}^{i})(-1)
\]
be the connecting morphism associated to 
$(\Phi(\mathbf{C}),\pi(\mathbf{C}))$. Since 
$h$ is an isomorphism and $\Phi(H(\mathbf{C}))$ is a 
monomorphism we have that $H(\Phi(\mathbf{C}))$ is a 
monomorphism, hence $\partial(\Phi(\mathbf{C}),\pi(\mathbf{C}))=0$
and $H(\pi(\mathbf{C}))$ is an epimorphism. It follows that $h'$
is an isomorphism.
\end{proof}

\section{Cotorsion theories and dimension}

\subsection{} \label{sec:5.1}
Let $\mathcal{A}$ be an abelian category and 
$\mathcal{K}$ a class of objects of $\mathcal{A}$. 

\subsubsection{} \label{sec:5.1.1}
\cite[\S3 n\textsuperscript{\scriptsize{o}} 1 lemme 1]{Bou} (1) 
Consider the solid arrows diagram 
\[
\xymatrix{
&{X} \ar[r]^{u'} \ar[d]_{f'} \ar@{..>}[dl]_{k'}& 
{B}\ar[r]^{v'} \ar[d]_{f} \ar[dl]_{k}& {B''}\ar[dl]_{k''}\\
{A'} \ar[r]^{u} & {A}\ar[r]^{v} & {A''}
}
\]
with $X\in ^{\perp}\mathcal{K}$, $(u',v')$ semi-exact, $\K u\in\mathcal{K}$, $(u,v)$ exact, 
$fu'=vf'$ and $f=vk+k''v'$. Then there is a morphism $k':X\to A'$ such that $f'=uk'+ku'$.

\cite[\S3 n\textsuperscript{\scriptsize{o}} 1 lemme 2]{Bou} (2) 
Consider the commutative solid arrows diagram 
\[
\xymatrix{
{X} \ar[r]^{u'} \ar@{..>}[d]_{f'}& {B}\ar[r]^{v'} \ar[d]_{f} & {B''}\ar[d]_{f''}\\
{A'} \ar[r]^{u} & {A}\ar[r]^{v} & {A''}
}
\]
with $X\in ^{\perp}\mathcal{K}$, $(u',v')$ semi-exact, $\K u\in\mathcal{K}$ and $(u,v)$ exact.
Then there is a morphism $f':X\to A'$ such that $uf'=fu'$.
\begin{proof}
We prove 1, part 2 is similar. We apply the functor $\mathcal{A}(X,-)$ to the given diagram;
we obtain the commutative diagram
\[
\xymatrix{
&{\mathcal{A}(X,X)} \ar[r]^{\bar{u'}} \ar[d]_{\bar{f'}} & 
{\mathcal{A}(X,B)}\ar[r]^{\bar{v'}} \ar[d]_{\bar{f}} & {\mathcal{A}(X,B'')}\\
{\mathcal{A}(X,A')} \ar[r]^{\bar{u}} & {\mathcal{A}(X,A)}\ar[r]^{\bar{v}} & {\mathcal{A}(X,A'')}
}
\]
with top row semi-exact and bottom row exact (\ref{sec:1.2.1}(1bis)).
We have  $\bar{v}\overline{f'-ku'}=\overline{k''v'u'}=0$, hence there is 
$l:\mathcal{A}(X,X)\to\im \bar{u}$ such that $il=\overline{f'-ku'}$, where
$i:\im \bar{u}\to\mathcal{A}(X,A)$ is the inclusion. Taking  $1_{X}$ we obtain
that there is $k':X\to A'$ such that $\bar{u}(k')=f'-ku'$, that is, $f'=uk'+ku'$.
\end{proof}

\subsubsection{} \label{sec:5.1.2}
\cite[\S3 n\textsuperscript{\scriptsize{o}} 1 prop. 1a)]{Bou} Let $(\mathbf{X},d'),(\mathbf{C},d)$ be 
complexes in $\mathcal{A}$ and $r$ an integer. Let $(u_{i}:\mathbf{X}_{i}\to\mathbf{C}_{i})_{i\leqslant r}$ 
be a family of morphisms such that $d_{i}u_{i}=u_{i-1}d'_{i}$ for $i\leqslant r$. Suppose that 
$\mathbf{X}_{i}\in ^{\perp}\mathcal{K}$ for $i>r$, $H_{i}(\mathbf{C})=0$ for $i\geqslant r$ and 
$B_{i}(\mathbf{C})\in \mathcal{K}$ for $i>r$. Then the family $(u_{i})$ extends to a morphism
$v:\mathbf{X}\to\mathbf{C}$; any two such extensions are homotopic.

\subsubsection{} \label{sec:5.1.3}
\cite[\S5 n\textsuperscript{\scriptsize{o}} 2 lemme 1a)]{Bou}
Let $\mathbf{X},\mathbf{C}$ be complexes in $\mathcal{A}$ with 
$\mathbf{X}\in Ch(^{\perp}\mathcal{K})$, $\mathbf{C}$ exact with cycles 
in $\mathcal{K}$ and $\mathbf{C}_{n}=0$ for $n<0$. Then 
$\mathrm{Homgr}_{\mathcal{A}}(\mathbf{X},\mathbf{C})$ is exact.

\subsubsection{} \label{sec:5.1.4}
Suppose that $(\mathcal{K},\mathcal{K}^{\perp})$ is a cotorsion theory.
A $\mathcal{K}$-resolution of an object $A$ of $\mathcal{A}$ is a pair 
$(\mathbf{C},p)$, where $\mathbf{C}$ is a 
complex that is degreewise in $\mathcal{K}$, $\mathbf{C}_{n}=0$ 
for $n<0$ and $p:\mathbf{C}\to S^{0}(A)$ is a quasi-isomorphism. 
A $\mathcal{K}^{\perp}$-resolution of an object $A$ of $\mathcal{A}$ is a 
pair $(\mathbf{C},e)$, where $\mathbf{C}$ is a 
complex that is degreewise in $\mathcal{K}^{\perp}$, $\mathbf{C}_{n}=0$ 
for $n>0$ and $e:S^{0}(A)\to\mathbf{C}$ is a quasi-isomorphism. 

\subsubsection{} \label{sec:5.1.5}
\cite[\S3 n\textsuperscript{\scriptsize{o}} 2 prop. 3]{Bou}
Let $f:A\to B$ be a morphism in $\mathcal{A}$. Let $p':\mathbf{X}\to S^{0}(A)$
be a morphism with $\mathbf{X}_{n}=0$ for $n<0$ and $\mathbf{X}$ 
degreewise in $\mathcal{K}$, and $p:\mathbf{C}\to S^{0}(B)$ a
$\mathcal{K}$-resolution of $B$ with $B(\mathbf{C})$ degreewise in 
$\mathcal{K}^{\perp}$. Then there is a morphism $\tilde{f}:\mathbf{X}\to
\mathbf{C}$, unique up to homotopy, such that $p\tilde{f}=S^{0}(f)p'$.

\subsubsection{} \label{sec:5.1.6}
Suppose that $(\mathcal{K},\mathcal{K}^{\perp})$ is a left complete cotorsion 
theory, meaning that $\Co(\mathcal{K}^{\perp}\rightarrowtail\mathcal{K})
=Ob(\mathcal{A})$ (\ref{sec:1.2.4}). A complete $\mathcal{K}$-resolution 
of an object $A$ of $\mathcal{A}$ is a pair $(\mathbf{C},p)$, where 
$(\mathbf{C},p)$ is a $\mathcal{K}$-resolution of $A$ and $B(\mathbf{C})$ is 
degreewise in $\mathcal{K}^{\perp}$. Every object of $\mathcal{A}$ has a 
complete $\mathcal{K}$-resolution.

\paragraph{} Suppose that $(\mathcal{K},\mathcal{K}^{\perp})$ is a right complete 
cotorsion theory, meaning that $\K(\mathcal{K}^{\perp}\twoheadrightarrow\mathcal{K})
=Ob(\mathcal{A})$. A complete $\mathcal{K}^{\perp}$-resolution 
of an object $A$ of $\mathcal{A}$ is a pair $(\mathbf{C},e)$, where 
$(\mathbf{C},e)$ is a $\mathcal{K}^{\perp}$-resolution of $A$ and $B(\mathbf{C})$ is 
degreewise in $\mathcal{K}$. Every object of $\mathcal{A}$ has a 
complete $\mathcal{K}^{\perp}$-resolution.

\subsubsection{} \label{sec:5.1.7}
\cite[\S3 n\textsuperscript{\scriptsize{o}} 2 cor. \`{a} la prop. 3]{Bou}
Let $(\mathbf{C},p)$ and $(\mathbf{C}',p')$ be two complete
$\mathcal{K}$-resolutions of $A$. There is a homotopy equivalence 
$\alpha:\mathbf{C}'\to\mathbf{C}$, unique up to homotopy,
such that $p\alpha=p'$.

\subsubsection{} \label{sec:5.1.8}
Let $(\mathbf{C},p)$ and $(\mathbf{X},p')$ be two 
$\mathcal{K}$-resolutions of $A$ with $(\mathbf{C},p)$
complete. Then there is a morphism $\tilde{f}:\mathbf{X}\to\mathbf{C}$,
unique up to homotopy, such that $p\tilde{f}=p'$.

\subsubsection{} \label{sec:5.1.9}
Let $\mathbf{C},\mathbf{C}'$ be exact complexes in $\mathcal{A}$
and $n$ an integer. Every morphism $u:\mathbf{C}\to\mathbf{C}'$ with $B_{n-1}(u)$ 
an isomorphism factors as $\mathbf{C}\overset{u'}\to\mathbf{C}''\overset{v}\to\mathbf{C}$, 
where: $\mathbf{C}''$ is exact, $B_{n}(u')$ is an isomorphism, $B_{n-1}(v)=B_{n-1}(u)$ 
and the diagram 
\begin{displaymath}
\xymatrix{
{\mathbf{C}''_{n+1}}\ar[r] \ar[d]&{\mathbf{C}''_{n}} \ar[d]\\
{\mathbf{C}'_{n+1}} \ar[r]&{\mathbf{C}'_{n}}
}
\end{displaymath}
is exact (\ref{sec:1.2.8}).
\begin{proof}
The commutative diagram
\begin{displaymath}
\xymatrix{
{\mathbf{C}_{n+2}}\ar[rr] \ar[d]_{u_{n+2}}&&{\mathbf{C}_{n+1}}\ar[d]_{u'_{n+1}}
\ar[r] & {\mathbf{C}_{n}} \ar@{=}[d]\ar@{->>}[r] & {B_{n-1}(\mathbf{C})}\ar@{=}[d]\\
{\mathbf{C}'_{n+2}}\ar@{->>}[r]\ar@{=}[d]&{B_{n+1}(\mathbf{C}')}\ar@{>->}[r]\ar@{=}[d] & 
{\mathbf{C}'_{n+1}\times_{\mathbf{C}'_{n}}\mathbf{C}_{n}}\ar[d]_{v_{n+1}}
\ar[r]&{\mathbf{C}_{n}}\ar@{->>}[r]\ar[d]_{u_{n}}&{B_{n-1}(\mathbf{C})}\ar[d]^{B_{n-1}(u)}\\
{\mathbf{C}'_{n+2}} \ar@{->>}[r]&{B_{n+1}(\mathbf{C}')}\ar@{>->}[r] & {\mathbf{C}'_{n+1}}
\ar[r]&{\mathbf{C}'_{n}}\ar@{->>}[r]&{B_{n-1}(\mathbf{C}')}
}
\end{displaymath}
shows that we can put
\begin{displaymath}
\mathbf{C}''_{k} =
\begin{cases}
\mathbf{C}'_{k}, & k\geqslant n+2,\\
\mathbf{C}'_{n+1}\times_{\mathbf{C}'_{n}}\mathbf{C}_{n}, & k=n+1,\\
\mathbf{C}_{n}, & k\leqslant n
\end{cases} 
\end{displaymath}
\begin{displaymath}
u'_{k} =
\begin{cases}
u_{k}, & k\geqslant n+2,\\
1_{\mathbf{C}_{k}}, & k\leqslant n
\end{cases} 
\end{displaymath}
and 
\begin{displaymath}
v_{k} =
\begin{cases}
1_{\mathbf{C}_{k}}, & k\geqslant n+2,\\
u_{k}, & k\leqslant n
\end{cases} 
\end{displaymath}
\end{proof}

\subsubsection{} \label{sec:5.1.10}
Suppose that $(\mathcal{K},\mathcal{K}^{\perp})$ is a left exact 
cotorsion theory (\ref{sec:1.2.7}) and let $u:\mathbf{C}\to\mathbf{C}'$
be a morphism as in \ref{sec:5.1.9}. If $\mathbf{C}_{k},\mathbf{C}'_{k}\in\mathcal{K}$
for all $k\geqslant n$ then $\mathbf{C}''_{k}\in\mathcal{K}$ for all $k\geqslant n$.
\begin{proof}
This follows from the construction of $\mathbf{C}''_{k}$ and \ref{sec:1.2.8}.
\end{proof}

\subsubsection{} \label{sec:5.1.11}
Suppose that $(\mathcal{K},\mathcal{K}^{\perp})$ is a left 
exact and left complete cotorsion theory. Let $A\in Ob(\mathcal{A})$ 
and $n$ an integer $\geqslant 0$. The following are equivalent:

(1) there is a $\mathcal{K}$-resolution $(\mathbf{C},p)$ of $A$ with 
$B_{n-1}(\mathbf{C})\in \mathcal{K}$;

(2) any $\mathcal{K}$-resolution $(\mathbf{D},q)$ of $A$ has 
$B_{n-1}(\mathbf{D})\in \mathcal{K}$.

\begin{proof}
(1)$\Rightarrow$(2) Let $(\mathbf{D},q)$ be an arbitrary $\mathcal{K}$-resolution of $A$
and let $(\mathbf{C}',p')$ be a complete $\mathcal{K}$-resolution of $A$. By \ref{sec:5.1.8}
there are morphisms $u:\mathbf{C}\to\mathbf{C}'$ and $w:\mathbf{D}\to\mathbf{C}'$ 
such that $p'u=p$ and $p'w=q$. The morphism $u$ is a morphism between exact complexes
with $B_{-1}(u)=1_{A}$, hence (\ref{sec:5.1.9}) we can factor it as $u=vu'$. We continue
this factorisation procedure with $u'$ in place of $u$; after $n-1$ such factorizations we obtain 
an exact diagram 
\begin{displaymath}
\xymatrix{
{B_{n-1}(\mathbf{C})}\ar@{>->}[r] \ar[d]&{\mathbf{C}_{n-1}} \ar[d]\\
{B_{n-1}(\mathbf{C}')} \ar@{>->}[r]&{\mathbf{X}_{n-1}}
}
\end{displaymath}
Applying \ref{sec:5.1.10} each time we factor we have that $\mathbf{X}_{n-1}\in\mathcal{K}$.
Since $\mathcal{K}$ is closed under extensions and retracts we obtain that
$B_{n-1}(\mathbf{C}')\in\mathcal{K}$. Now, the morphism $w$ is a morphism 
between exact complexes with $B_{-1}(w)=1_{A}$, so after applying to it
the same procedure as to $u$ we obtain an exact diagram 
\begin{displaymath}
\xymatrix{
{B_{n-1}(\mathbf{D})}\ar@{>->}[r] \ar[d]&{\mathbf{D}_{n-1}} \ar[d]\\
{B_{n-1}(\mathbf{C}')} \ar@{>->}[r]&{\mathbf{X}_{n-1}}
}
\end{displaymath}
with $\mathbf{X}_{n-1}\in\mathcal{K}$. By \ref{sec:1.2.8} we then obtain that 
$B_{n-1}(\mathbf{D})\in\mathcal{K}$.
\end{proof}

\subsection{} \label{sec:5.2}

Let $\mathcal{A}$ be an abelian category and $(\mathcal{K},\mathcal{K}^{\perp}), 
(^{\perp}\mathcal{L},\mathcal{L})$ be cotorsion theories in 
$\mathcal{A}$ with $\mathcal{K}^{\perp}\subset \mathcal{L}$.

\subsubsection{} \label{sec:5.2.1}
Suppose that $\mathcal{A}$ has enough $\mathcal{K}^{\perp}$ objects 
and $(^{\perp}\mathcal{L},\mathcal{L})$ is left exact and left complete. 
We inductively define an ascending sequence 
$(n-^{\perp}\mathcal{L})_{n\geqslant 0}$ and a descending sequence 
$(\mathcal{S}_{n})_{n\geqslant 0}$ of 
classes of objects of $\mathcal{A}$ such that 
$n-^{\perp}\mathcal{L}\subset \mathcal{K}$, 
$\mathcal{K}^{\perp}\subset \mathcal{S}_{n}$
and $n-^{\perp}\mathcal{L}=^{\perp}\mathcal{S}_{n}$. 
The argument will show that 
$(n-^{\perp}\mathcal{L}, (n-^{\perp}\mathcal{L})^{\perp})$ 
is a cotorsion theory.

We put $0-^{\perp}\mathcal{L}=^{\perp}\mathcal{L}$ and $\mathcal{S}_{0}=\mathcal{L}$.
Having defined $n-^{\perp}\mathcal{L}$ and $\mathcal{S}_{n}$, we set
$(n+1)-^{\perp}\mathcal{L}=\mathcal{K}\cap \Co(n-^{\perp}\mathcal{L}
\rightarrowtail^{\perp}\mathcal{L})$ and $\mathcal{S}_{n+1}=
\Co(\mathcal{S}_{n}\rightarrowtail\mathcal{K}^{\perp})$.

We show that $\mathcal{S}_{n+1}\subset ((n+1)-^{\perp}\mathcal{L})^{\perp}$.
Let $C\rightarrowtail B\twoheadrightarrow A$ be an exact sequence with
$C\in\mathcal{S}_{n+1}$ and $A\in (n+1)-^{\perp}\mathcal{L}$.
We can find an exact sequence $S_{n}\rightarrowtail Y\twoheadrightarrow C$
with $S_{n}\in \mathcal{S}_{n}$ and $Y\in\mathcal{K}^{\perp}$.
Form the commutative diagram
\[
\xymatrix{
{S_{n}}\ar@{>->}[r] \ar@{=}[d]& {PB'}\ar@{->>}[r]\ar[d]& 
{PB} \ar@{>->}[r] \ar[d]_{f} & {X_{0}}\ar@{->>}[r]\ar@{->>}[d]_{g} & {A}\ar@{=}[d]\\
{S_{n}}\ar@{>->}[r] & {Y}\ar@{->>}[r] & {C} \ar@{>->}[r] & {B}\ar@{->>}[r] & {A}
}
\]
where $X_{0}\in ^{\perp}\mathcal{L}$ and $PB,PB'$ are pullbacks. 
Then (\ref{sec:5.1.11}) $PB\in n-^{\perp}\mathcal{L}$, hence the 
sequence $S_{n}\rightarrowtail PB'\twoheadrightarrow PB$
splits by the inductive hypothesis. Form now the commutative 
solid arrows diagram
\[
\xymatrix{
{PB} \ar@{>->}[r] \ar[d]& {X_{0}}\ar@{->>}[r]\ar[d]& {A}\ar@{=}[d]\\
{Y} \ar@{>->}[r]^{u} \ar[d]& {PO}\ar@{->>}[r]^{v} \ar@{..>}[d]& {A}\\
{C}\ar@{>->}[r] & {B}
}
\]
where the composite $PB\to Y\to C$ is $f$ and $PO$ means pushout. 
There is then the dotted arrow that makes commutative the resulting 
square diagram and such that the composite $X_{0}\to PO\to B$ is $g$. 
It follows that the diagram
\begin{displaymath}
\xymatrix{
{Y}\ar@{>->}[r]^{u} \ar[d]&{PO} \ar[d]\\
{C} \ar@{>->}[r]&{B}
}
\end{displaymath}
is a pushout. Since $A\in \mathcal{K}$ the sequence $(u,v)$ splits, hence
$u$ has a retract and therefore so does $C\to B$.

We show that $^{\perp}\mathcal{S}_{n+1}\subset (n+1)-^{\perp}\mathcal{L}$.
Let $A\in ^{\perp}\mathcal{S}_{n+1}$. We first show that $A\in\mathcal{K}$.
Let $Y\rightarrowtail B\twoheadrightarrow A$ be an exact sequence with 
$Y\in\mathcal{K}^{\perp}$; then $Y\in\mathcal{S}_{n+1}$ since we have the 
exact sequence $0\to Y=Y$ and $0\in\mathcal{S}_{n}$. Consequently the first
short exact sequence splits, so $A\in\mathcal{K}$. Now, we can find an exact sequence
$K\rightarrowtail X_{0}\twoheadrightarrow A$ where $X_{0}\in ^{\perp}\mathcal{L}$.
We show that $K\in n-^{\perp}\mathcal{L}$. Given an exact sequence 
$S_{n}\rightarrowtail B\twoheadrightarrow K$ with $S_{n}\in\mathcal{S}_{n}$,
form the commutative diagram
\[
\xymatrix{
{S_{n}}\ar@{>->}[r] \ar@{=}[d]& {B}\ar@{->>}[r]\ar@{>->}[d]_{g}& 
{K} \ar@{>->}[r] \ar@{>->}[d]_{f} & {X_{0}}\ar@{->>}[r]\ar@{>->}[d] & {A}\ar@{=}[d]\\
{S_{n}}\ar@{>->}[r] & {Y}\ar@{->>}[r] & {PO} \ar@{>->}[r] & {PO'}\ar@{->>}[r] & {A}
}
\]
where $PO,PO'$ are pushouts and $Y\in \mathcal{K}^{\perp}$. Then $PO\in\mathcal{S}_{n+1}$,
hence the sequence $PO\rightarrowtail PO'\twoheadrightarrow A$ splits.
 We can form then the commutative solid arrows diagram
\[
\xymatrix{
& {B}\ar@{->>}[r]\ar@{..>}[d] & {K}\ar[d]\\
{S_{n}} \ar@{>->}[r]^{u} \ar@{=}[d]& {PB}\ar@{->>}[r]^{v}\ar[d]& {X_{0}}\ar[d]\\
{S_{n}} \ar@{>->}[r]& {Y}\ar@{->>}[r] & {PO}
}
\]
where $PB$ means pullback and the composite $K\to X_{0}\to PO$ is $f$.
There is then the dotted arrow that makes commutative the resulting square diagram 
and such that the composite $B\to PB\to Y$ is $g$. It follows that the diagram
\begin{displaymath}
\xymatrix{
{B}\ar@{->>}[r] \ar[d]&{K} \ar[d]\\
{PB} \ar@{->>}[r]^{v}&{X_{0}}
}
\end{displaymath}
is a pullback. Since $X_{0}\in ^{\perp}\mathcal{L}\subset n-^{\perp}\mathcal{L}$ 
the sequence $(u,v)$ splits, hence $v$ has a section and therefore so does $B\to K$.

%Summing up, we obtain that $(n+1)-^{\perp}\mathcal{L}=^{\perp}\mathcal{S}_{n+1}$.

\subsubsection{} \label{sec:5.2.2}
Under the assumptions \ref{sec:5.2.1}, if $\mathcal{K}$ is left exact then 
so is $n-^{\perp}\mathcal{L}$.
\begin{proof}
We use induction over $n\geqslant 0$. The case $n=0$ holds by assumption.
Suppose that $n-^{\perp}\mathcal{L}$ is left exact and let 
$A'\rightarrowtail A\twoheadrightarrow A''$ be an exact sequence with 
$A,A''\in (n+1)-^{\perp}\mathcal{L}$. Then $A'\in\mathcal{K}$ by assumption.
We can find an exact sequence $A''_{0}\rightarrowtail X''_{0}\twoheadrightarrow A''$
with $X''_{0}\in^{\perp}\mathcal{L},A''_{0}\in n-^{\perp}\mathcal{L}$ (\ref{sec:5.1.11}).
Form the commutative diagram
\begin{displaymath}
\xymatrix{
{K'}\ar@{>->}[d]\ar@{=}[r]&{K'}\ar@{>->}[d]\\
{PB'}\ar@{>->}[r] \ar@{->>}[d]& {X_{0}}\ar@{->>}[r]\ar@{->>}[d] & {X''_{0}}\ar@{=}[d]\\
{A'} \ar@{>->}[r] \ar@{=}[d]& {PB}\ar@{->>}[r]\ar@{->>}[d] & {X''_{0}}\ar@{->>}[d]\\
{A'} \ar@{>->}[r] & {A}\ar@{->>}[r] & {A''}
}
\end{displaymath}
where $PB$ and $PB'$ mean pullback, $X_{0}\in^{\perp}\mathcal{L}$ and 
$K'$ is the kernel of $X_{0}\to PB$. Then $PB'\in ^{\perp}\mathcal{L}$
and we consider the commutative solid arrows diagram
\begin{displaymath}
\xymatrix{
{K'}\ar@{>..>}[r] \ar@{>->}[d]& {K_{0}}\ar@{..>>}[r]\ar@{>->}[d] & {A''_{0}}\ar@{>->}[d]\\
{PB'} \ar@{>->}[r] \ar@{->>}[d]& {X_{0}}\ar@{->>}[r]\ar@{->>}[d] & {X''_{0}}\ar@{->>}[d]\\
{A'} \ar@{>->}[r] & {A}\ar@{->>}[r] & {A''}
}
\end{displaymath}
where $K_{0}$ is the kernel of the composite morphism $X_{0}\to PB\to A$.
By the universal property of kernel we have the induced dotted arrows that make 
the resulting square diagrams commute; moreover, by the snake diagram 
the dotted sequence of the previous diagram is a short exact sequence. 
By the assumption on $A$ and \ref{sec:5.1.11} we have that 
$K_{0}\in  n-^{\perp}\mathcal{L}$, therefore $K'\in n-^{\perp}\mathcal{L}$ 
by the inductive hypothesis, hence the assertion.
\end{proof}

%\subsubsection{} \label{sec:5.2.3}
% ascending sequence $(n-^{\perp}\mathcal{L})_{n\geqslant 0}$ is stationary.

\subsubsection{} \label{sec:5.2.3}
Suppose that $\mathcal{A}$ has enough injectives and projectives. 
Putting $\mathcal{K}=\mathcal{L}=Ob(\mathcal{A})$ 
in \ref{sec:5.2.1} and \ref{sec:5.2.2} we obtain that the class of objects 
of $\mathcal{A}$ of projective dimension $\leqslant n$ is the left class 
of a left exact cotorsion theory.

\subsubsection{} \label{sec:5.2.4}
Suppose that $\mathcal{A}$ has enough $^{\perp}\mathcal{L}$
objects and $(\mathcal{K},\mathcal{K}^{\perp})$ is right exact 
and right complete. We inductively define an ascending sequence 
$(n-\mathcal{K}^{\perp})_{n\geqslant 0}$ and a descending sequence 
$(\mathcal{T}_{n})_{n\geqslant 0}$ of classes of objects of $\mathcal{A}$ 
such that $n-\mathcal{K}^{\perp}\subset \mathcal{L}$, 
$^{\perp}\mathcal{L}\subset \mathcal{T}_{n}$
and $n-\mathcal{K}^{\perp}=\mathcal{T}_{n}^{\perp}$. 

We put $0-\mathcal{K}^{\perp}=\mathcal{K}^{\perp}$ and $\mathcal{T}_{0}
=\mathcal{K}$. Having defined $n-\mathcal{K}^{\perp}$ and $\mathcal{T}_{n}$, 
we set $(n+1)-\mathcal{K}^{\perp}=\mathcal{L}\cap \K (\mathcal{K}^{\perp}
\twoheadrightarrow n-\mathcal{K}^{\perp})$ and $\mathcal{T}_{n+1}=
\K (^{\perp}\mathcal{L}\twoheadrightarrow\mathcal{T}_{n})$.

An argument similar to \ref{sec:5.2.1} and \ref{sec:5.2.2} shows that 
$(^{\perp}(n-\mathcal{K}^{\perp}),n-\mathcal{K}^{\perp})$ is a cotorsion 
theory that is right exact if $\mathcal{L}$ is moreover right exact. 
If $(\mathcal{K},\mathcal{K}^{\perp})$ is cogenerated (\ref{sec:1.2.3})
by $\mathcal{S}$ and  $(^{\perp}\mathcal{L},\mathcal{L})$ is
cogenerated by $\mathcal{T}$ then $(^{\perp}(n-\mathcal{K}^{\perp}),
n-\mathcal{K}^{\perp})$ is cogenerated by $\mathcal{C}_{n}$, where
$\mathcal{C}_{0}=\mathcal{S}, \mathcal{C}_{1}=
\K (\mathcal{T}\twoheadrightarrow\mathcal{C}_{0}),...,
\mathcal{C}_{n+1}=\K (\mathcal{T}\twoheadrightarrow\mathcal{C}_{n})$ 
for $n\geqslant 0$.

Suppose that $\mathcal{A}$ has enough injectives and projectives. 
Putting $\mathcal{K}=\mathcal{L}=Ob(\mathcal{A})$ we obtain that the 
class of objects of $\mathcal{A}$ of injective dimension $\leqslant n$ is 
the right class of a right exact cotorsion theory.

\subsubsection{} \label{sec:5.2.5}
Recall that $\mathcal{S}_{0}=\mathcal{L}$, $\mathcal{S}_{n+1}=
\Co(\mathcal{S}_{n}\rightarrowtail\mathcal{K}^{\perp})$, $\mathcal{T}_{0}=\mathcal{K}$
and $\mathcal{T}_{n+1}=\K (^{\perp}\mathcal{L}\twoheadrightarrow\mathcal{T}_{n})$.

(1) Suppose that $\mathcal{A}$ has enough $\mathcal{K}^{\perp}$ objects 
and $(^{\perp}\mathcal{L},\mathcal{L})$ is left exact and left complete. Then
\[
n-^{\perp}\mathcal{L}=\mathcal{K}\Leftrightarrow\mathcal{T}_{n}\subset ^{\perp}\mathcal{L}
\Rightarrow n-\mathcal{K}^{\perp}=\mathcal{L}
\]
(1bis) Suppose that $\mathcal{A}$ has enough $^{\perp}\mathcal{L}$ objects 
and $(\mathcal{K},\mathcal{K}^{\perp})$ is right exact and right complete. Then
\[
n-\mathcal{K}^{\perp}=\mathcal{L}\Leftrightarrow\mathcal{S}_{n}\subset \mathcal{K}^{\perp}
\Rightarrow n-^{\perp}\mathcal{L}=\mathcal{K}
\]
\begin{proof}
(1) The assertion is clear when $n=0$. Suppose that $n>0$ and 
$n-^{\perp}\mathcal{L}=\mathcal{K}$
and let $K_{n}\in\mathcal{T}_{n}$. We can find an exact sequence
\[
\xymatrix{
{K_{n}}\ar@{>->}[r] & {X_{n-1}}\ar[r] & {...} \ar[r] & {X_{0}}\ar@{->>}[r] & {K}
}
\]
with $K\in\mathcal{K}$ and $X_{0},...,X_{n-1}\in ^{\perp}\mathcal{L}$. 
Since $K\in n-^{\perp}\mathcal{L}$ it follows from \ref{sec:5.1.11} 
that $K_{n}\in ^{\perp}\mathcal{L}$. Then 
$\mathcal{L}\subset \mathcal{T}_{n}^{\perp}=n-\mathcal{K}^{\perp}$, 
hence $\mathcal{L}=n-\mathcal{K}^{\perp}$. Conversely, 
let $X\in \mathcal{K}$ and let 
\[
\xymatrix{
{...}\ar[r]&{X_{n}}\ar[r] & {X_{n-1}}\ar[r]^{d_{n-1}}
& {...} \ar[r] & {X_{0}}\ar@{->>}[r] & {X}
}
\]
be a $^{\perp}\mathcal{L}$-resolution of $X$. Then 
$\K d_{n-1}\in\mathcal{T}_{n}\subset^{\perp}\mathcal{L}$,
hence $X\in n-^{\perp}\mathcal{L}$.
\end{proof}

\subsubsection{} \label{sec:5.2.6}
Putting $\mathcal{K}=Ob(\mathcal{A})$ in \ref{sec:5.2.5} and 
assumming that $\mathcal{A}$ has enough injectives we obtain 
that $n-^{\perp}\mathcal{L}=Ob(\mathcal{A})$ if and only if 
all objects of $\mathcal{L}$ have injective dimension $\leqslant n$.
Putting $\mathcal{L}=Ob(\mathcal{A})$ in \ref{sec:5.2.5} and 
assumming that $\mathcal{A}$ has enough projectives we obtain 
that $n-\mathcal{K}^{\perp}=Ob(\mathcal{A})$ if and only if 
all objects of $\mathcal{K}$ have projective dimension $\leqslant n$.
Putting $\mathcal{K}=\mathcal{L}=Ob(\mathcal{A})$ in \ref{sec:5.2.5} 
and assumming that $\mathcal{A}$ has enough projectives and injectives 
we obtain that all objects of $\mathcal{A}$ have projective dimension 
$\leqslant n$ if and only if all objects of $\mathcal{A}$ have injective 
dimension $\leqslant n$.

%\subsubsection{} \label{sec:5.2.8} adjunction
\subsubsection{} \label{sec:5.2.7}
Suppose that there is $n\geqslant 0$ such that 
$n-^{\perp}\mathcal{L}=n-\mathcal{K}^{\perp}$.
This happens if $\underset{n\geqslant 0}\bigcup n-^{\perp}\mathcal{L}
=\underset{n\geqslant 0}\bigcup n-\mathcal{K}^{\perp}$ and the 
sequences $(n-^{\perp}\mathcal{L})_{n\geqslant 0}$ and 
$(n-\mathcal{K}^{\perp})_{n\geqslant 0}$ are eventually constant.
In this case (\ref{sec:1.2.13}) 
$n-^{\perp}\mathcal{L}\cap (n-^{\perp}\mathcal{L})^{\perp}=
\mathcal{K}\cap\mathcal{K}^{\perp}$ and 
$^{\perp}(n-\mathcal{K}^{\perp})\cap n-\mathcal{K}^{\perp}=
^{\perp}\mathcal{L}\cap\mathcal{L}$.

\subsection{} \label{sec:5.3}
Let 
\[
T:\mathcal{A}'\times \mathcal{A}\to \mathcal{A}'', \
H:\mathcal{A}^{op}\times \mathcal{A}''\to \mathcal{A}', \
C:\mathcal{A}'^{op}\times \mathcal{A}''\to \mathcal{A}
\]
be an abelian THC--situation.
We denote by $\mathcal{F}$ the class of flat objects (\ref{sec:3.1}) of $\mathcal{A}$.
We assume that the cotorsion theory $(\mathcal{F},\mathcal{F}^{\perp})$ 
(\ref{sec:3.1.4}) is left complete and that $\mathcal{A}$ has enough injective objects. 
Then, using \ref{sec:5.2.1} and \ref{sec:5.2.2} with $\mathcal{K}=Ob(\mathcal{A})$
and $^{\perp}\mathcal{L}=\mathcal{F}$, we obtain that $n-\mathcal{F}$ is the left class
of a left exact cotorsion theory, as already noticed in \cite[Proposition 10.2.9]{Pe}. 
The objects of $n-\mathcal{F}$ are said to have flat dimension $\leqslant n$, or 
to be $n$-flat [\emph{loc. cit.}, Definition 3.3.3]. From \ref{sec:5.2.4} we obtain 
that the class of cotorsion objects of $\mathcal{A}$ of injective dimension 
$\leqslant n$ is the right class of a right exact cotorsion theory. From \ref{sec:5.2.5} 
we obtain that all objects of $\mathcal{A}$ have flat dimension $\leqslant n$ if 
and only if all cotorsion objects of $\mathcal{A}$ have injective dimension $\leqslant n$.
If $\mathcal{A}$ has enough projectives we obtain from \emph{loc. cit.} that all
flat objects of $\mathcal{A}$ have projective dimension $\leqslant n$ if and only
if $Ob(\mathcal{A})=n-\mathcal{F}^{\perp}$.

We denote by $\mathcal{F}'$ the class of flat objects of $\mathcal{A}'$
and we further assume that the cotorsion theory $(\mathcal{F}',\mathcal{F}'^{\perp})$ 
is left complete and that $\mathcal{A}'$ has enough injective objects. 

\subsubsection{} \label{sec:5.3.1}
\cite[Proposition 14.2.3(2)]{Pe} Suppose, moreover, that $\mathcal{A}''$ has 
a set $(E''_{i})_{i\in I}$ of injective cogenerators. Then an object $A$ of 
$\mathcal{A}$ is $n$-flat if and only if $H(A,E''_{i})$ is (a cotorsion object)
of injective dimension $\leqslant n$ for all $i\in I$.
\begin{proof}
Suppose that $A$ is $n$-flat. We can find an exact sequence
\[
\xymatrix{
{K}\ar@{>->}[r] & {P_{n-1}}\ar[r] & {...} \ar[r] & {P_{0}}\ar@{->>}[r] & {A}
}
\]
with $P_{0},...,P_{n-1},K$ flat. Since $H(-,E''_{i})$ is exact (\ref{sec:3.1.3})
we have the exact sequence 
\[
\xymatrix{
{H(A,E''_{i})}\ar@{>->}[r] & {H(P_{0},E''_{i})}\ar[r] & {...} \ar[r] & 
{H(P_{n-1},E''_{i})}\ar@{->>}[r] & {H(K,E''_{i})}
}
\]
where $H(P_{0},E''_{i}),...,H(P_{n-1},E''_{i}),H(K,E''_{i})$ are injective by \ref{sec:3.1.2}.
This implies that $H(A,E''_{i})$ is $n$-injective. Conversely, let 
\[
\xymatrix{
{...}\ar[r]&{P_{n}}\ar[r]^{d_{n}} & {P_{n-1}}\ar[r]^{d_{n-1}} & {...} 
\ar[r] & {P_{0}}\ar@{->>}[r] & {A}
}
\]
be an $\mathcal{F}$-resolution of $A$. Then we have the exact sequence 
\[
\xymatrix{
{H(A,E''_{i})}\ar@{>->}[r] & {H(P_{0},E''_{i})}\ar[r] & {...} \ar[r] & 
{H(P_{n-1},E''_{i})}\ar[r]^{H(d_{n},E''_{i})} & {H(P_{n},E''_{i})}\ar[r]&{...}
}
\]
where $H(P_{n},E''_{i})$ is injective for all $n\geqslant 0$. Since 
$H(A,E''_{i})$ is $n$-injective it follows from \ref{sec:5.1.11} that
$\im H(d_{n},E''_{i})=H(\K d_{n-1},E''_{i})$ is injective, hence
$\K d_{n-1}$ is flat by \ref{sec:3.1.2}.
\end{proof}

\section{Tensor-hom-cotensor situations for categories of complexes}

\subsection{} \label{sec:6.1}
Let 
\[
T:\mathcal{A}'\times \mathcal{A}\to \mathcal{A}'', \
H:\mathcal{A}^{op}\times \mathcal{A}''\to \mathcal{A}', \
C:\mathcal{A}'^{op}\times \mathcal{A}''\to \mathcal{A}
\]
be an abelian THC--situation. We assume that 
$\mathcal{A}',\mathcal{A}$ have all limits and 
$\mathcal{A}''$ has all colimits. We revisit the 
standard THC--situation 
\[
T:Ch(\mathcal{A}')\times Ch(\mathcal{A})\to Ch(\mathcal{A}'')
\]
\[
H:Ch(\mathcal{A})^{op}\times Ch(\mathcal{A}'')\to Ch(\mathcal{A}')
\]
\[
C:Ch(\mathcal{A}')^{op}\times Ch(\mathcal{A}'')\to Ch(\mathcal{A})
\]
associated to the given one. Let $(\mathbf{C}',d')$ (resp. $(\mathbf{C},d),(\mathbf{C}'',d'')$)
be complexes in $\mathcal{A}'$ (resp. $\mathcal{A},\mathcal{A}''$).
\paragraph{(Tensor)} We define $T(\mathbf{C}',\mathbf{C})_{n}=
\underset{p+q=n}\bigoplus T(\mathbf{C}'_{p},\mathbf{C}_{q})$. 
Let $\sigma_{p,q}:T(\mathbf{C}'_{p},\mathbf{C}_{q})\to T(\mathbf{C}',\mathbf{C})_{n}$
be coprojection morphisms. The differential 
\[
D(\mathbf{C}',\mathbf{C})_{n}=D_{n}:T(\mathbf{C}',\mathbf{C})_{n}\to T(\mathbf{C}',\mathbf{C})_{n-1}
\]
is the unique morphism that satisfies
\[
D_{n}\sigma_{p,q}=\sigma_{p-1,q}T(d'_{p},\mathbf{C}_{q})+\sigma_{p,q-1}(-1)^{p}T(\mathbf{C}'_{p},d_{q})
\]
whenever $p+q=n$. The universal property that defines $D$ implies that $D^{2}=0$.
We omit the details that $T$ assembles into a bifunctor.
\paragraph{(Hom)} We define $H(\mathbf{C},\mathbf{C}'')_{n}=
\underset{p}\prod H(\mathbf{C}_{p},\mathbf{C}''_{p+n})$. Let $\pi_{p}:H(\mathbf{C},\mathbf{C}'')_{n}\to
H(\mathbf{C}_{p},\mathbf{C}''_{p+n})$ be projection morphisms. The differential 
\[
D(\mathbf{C},\mathbf{C}'')_{n}=D_{n}:H(\mathbf{C},\mathbf{C}'')_{n}\to H(\mathbf{C},\mathbf{C}'')_{n-1}
\]
is the unique morphism that satisfies
\[
\pi_{p}D_{n}=H(\mathbf{C}_{p},d''_{p+n})\pi_{p}-(-1)^{n}H(d_{p},\mathbf{C}''_{p+n-1})\pi_{p-1}
\]
for all integers $p$. The universal property that defines $D$ implies that $D^{2}=0$.
\paragraph{(Cotensor)} We define $C(\mathbf{C}',\mathbf{C}'')_{n}=
\underset{p}\prod C(\mathbf{C}'_{p},\mathbf{C}''_{p+n})$. Let $\pi_{p}:C(\mathbf{C}',\mathbf{C}'')_{n}\to
C(\mathbf{C}'_{p},\mathbf{C}''_{p+n})$ be projection morphisms. The differential 
\[
D(\mathbf{C}',\mathbf{C}'')_{n}=D_{n}:C(\mathbf{C}',\mathbf{C}'')_{n}\to C(\mathbf{C}',\mathbf{C}'')_{n-1}
\]
is the unique morphism that satisfies
\[
\pi_{p}D_{n}=C(\mathbf{C}'_{p},d''_{p+n})\pi_{p}-(-1)^{n}C(d'_{p},\mathbf{C}''_{p+n-1})\pi_{p-1}
\]
for all integers $p$. The universal property that defines $D$ implies that $D^{2}=0$.
\paragraph{(Tensor $\Leftrightarrow$ Hom)} To give a morphism
$u:T(\mathbf{C}',\mathbf{C})\to \mathbf{C}''$ is to give
morphisms $u_{p,q}: T(\mathbf{C}'_{p},\mathbf{C}_{q})\to\mathbf{C}''_{p+q}$ such that
\[
d''_{p+q}u_{p,q}=u_{p-1,q}T(d'_{p},\mathbf{C}_{q})+(-1)^{p}u_{p,q-1}T(\mathbf{C}'_{p},d_{q})
\]
To give a morphism $v:\mathbf{C}'\to H(\mathbf{C},\mathbf{C}'')$ is to give
morphisms $v_{n,p}:\mathbf{C}'_{n}\to H(\mathbf{C}_{p},\mathbf{C}''_{p+n})$ such that
\[
v_{n-1,p}d'_{n}=H(\mathbf{C}_{p},d''_{p+n})v_{n,p}-(-1)^{n}H(d_{p},\mathbf{C}''_{p+n-1})v_{n,p-1}
\]
Given a morphism $u:T(\mathbf{C}',\mathbf{C})\to \mathbf{C}''$ 
we define $v_{n,p}$ to be the adjoint transpose of $u_{n,p}$; conversely, 
given a morphism $v:\mathbf{C}'\to H(\mathbf{C},\mathbf{C}'')$ 
we define $u_{n,p}$ to be the adjoint transpose of $v_{n,p}$. 
\paragraph{(Tensor $\Leftrightarrow$ Cotensor)} 
To give a morphism $w:\mathbf{C}\to C(\mathbf{C}',\mathbf{C}'')$ is to give
morphisms $w_{n,p}:\mathbf{C}_{n}\to C(\mathbf{C}'_{p},\mathbf{C}''_{p+n})$ such that
\[
w_{n-1,p}d_{n}=C(\mathbf{C}'_{p},d''_{p+n})w_{n,p}-(-1)^{n}C(d'_{p},\mathbf{C}''_{p+n-1})w_{n,p-1}
\]
Given a morphism 
$u:T(\mathbf{C}',\mathbf{C})\to \mathbf{C}''$ we define $w_{n,p}=(-1)^{np}u_{p,n}^{\flat}$,
where $u_{p,n}^{\flat}:\mathbf{C}_{n}\to C(\mathbf{C}'_{p},\mathbf{C}''_{p+n})$ 
is the adjoint transpose of $u_{p,n}$; conversely, given a morphism 
$w:\mathbf{C}\to C(\mathbf{C}',\mathbf{C}'')$ we define $u_{p,n}=(-1)^{pn}w_{n,p}^{\#}$, 
where $w_{n,p}^{\#}:T(\mathbf{C}'_{p},\mathbf{C}_{n})\to\mathbf{C}''_{p+n}$ is the 
adjoint transpose of $w_{n,p}$.

\subsubsection{} \label{sec:6.1.1}
Consider the bifunctors $T,H$ from \ref{sec:6.1}. 
For every integer $r$ we have
\[
(T(\mathbf{C}',\mathbf{C})(r),D(\mathbf{C}',\mathbf{C})(r))=
(T(\mathbf{C}'(r),\mathbf{C}),D(\mathbf{C}'(r),\mathbf{C}))
\]
and
\[
(H(\mathbf{C},\mathbf{C}'')(r),D(\mathbf{C},\mathbf{C}'')(r))=
(H(\mathbf{C},\mathbf{C}''(r)),D(\mathbf{C},\mathbf{C}''(r)))
\]
For every integer $r$ we have a natural isomorphism
\[
\Sigma(\mathbf{C}',\mathbf{C}):T(\mathbf{C}',\mathbf{C}(r))\cong T(\mathbf{C}',\mathbf{C})(r)
\]
constructed as follows. For each integer $n$ let $\sigma_{p,q}:T(\mathbf{C}'_{p},\mathbf{C}_{q+r})
\to T(\mathbf{C}',\mathbf{C}(r))_{n}$, where $p+q=n$, and 
$\sigma'_{a,b}:T(\mathbf{C}'_{a},\mathbf{C}_{b})\to
T(\mathbf{C}',\mathbf{C})(r)_{n}$, where $a+b=n+r$, be coproduct coprojections.
Consider the morphism $(-1)^{pr}1_{T(\mathbf{C}'_{p},\mathbf{C}_{q+r})}:
T(\mathbf{C}'_{p},\mathbf{C}_{q+r})\to T(\mathbf{C}'_{p},\mathbf{C}_{q+r})$;
there is then a unique morphism 
\begin{displaymath}
\Sigma(\mathbf{C}',\mathbf{C})_{n}:
T(\mathbf{C}',\mathbf{C}(r))_{n}\to T(\mathbf{C}',\mathbf{C})(r)_{n}
\end{displaymath}
such that $\Sigma(\mathbf{C}',\mathbf{C})_{n}\sigma_{p,q}=(-1)^{pr}\sigma'_{p,q+r}$.
One can check that $\Sigma(\mathbf{C}',\mathbf{C})$ is a morphism. The inverse of 
$\Sigma(\mathbf{C}',\mathbf{C})_{n}$ is the unique morphism
\begin{displaymath}
\Sigma^{-1}(\mathbf{C}',\mathbf{C})_{n}:
T(\mathbf{C}',\mathbf{C})(r)_{n}\to T(\mathbf{C}',\mathbf{C}(r))_{n}
\end{displaymath}
such that $\Sigma^{-1}(\mathbf{C}',\mathbf{C})_{n}\sigma'_{a,b}=(-1)^{ar}\sigma_{a,b-r}$.
Naturality of $\Sigma$ in both arguments is easy to see.

We note that if $\mathbf{C}'$ has zero differential then the following diagram commutes
\begin{displaymath}
\xymatrix{
{T(\mathbf{C}',\mathbf{C})}\ar[rr]^{T(\mathbf{C}',d)}\ar[drr]_{D(\mathbf{C}',\mathbf{C})}
&&{T(\mathbf{C}',\mathbf{C}(-1))} \ar[d]^{\Sigma(\mathbf{C}',\mathbf{C})}\\
&&{T(\mathbf{C}',\mathbf{C})(-1)} 
}
\end{displaymath}

\subsubsection{} \label{sec:6.1.2} 
Let $\mathbf{C}'$ be a complex in $\mathcal{A}'$ and
$u:\mathbf{D}\to\mathbf{C}$ a morphism of complexes in $\mathcal{A}$.
We define
\[
\Sigma(\mathbf{C}',u):T(\mathbf{C}',Con(u))\to Con(T(\mathbf{C}',u))
\]
as $\Sigma(\mathbf{C}',u)=\Sigma(\mathbf{C}',\mathbf{D})\oplus 
1_{T(\mathbf{C}',\mathbf{C})}$, where $\Sigma(\mathbf{C}',\mathbf{D}):
T(\mathbf{C}',\mathbf{D}(-1))\cong T(\mathbf{C}',\mathbf{D})(-1)$ 
has been defined in \ref{sec:6.1.1}. One can check that $\Sigma(\mathbf{C}',u)$ 
is a morphism, hence an isomorphism. Moreover, we have the following 
commutative diagram 
\begin{displaymath}
\xymatrix{
{T(\mathbf{C}',\mathbf{C})}\ar@{>->}[rr]^{T(\mathbf{C}',\pi)}\ar@{=}[d]&&
{T(\mathbf{C}',Con(u))}\ar[d]_{\Sigma(\mathbf{C}',u)}\ar@{->>}[rr]^{T(\mathbf{C}',\delta)}
&&{T(\mathbf{C}',\mathbf{D}(-1))}\ar[d]_{\Sigma(\mathbf{C}',\mathbf{D})}\\
{T(\mathbf{C}',\mathbf{C})}\ar@{>->}[rr]^{\pi}&&
{Con(T(\mathbf{C}',u))}\ar@{->>}[rr]^{\delta}&&{T(\mathbf{C}',\mathbf{D})(-1)}
}
\end{displaymath}

If $u':\mathbf{D}'\to\mathbf{C}'$ is a morphism of complexes 
in $\mathcal{A}'$ and $\mathbf{C}$ is a complex in $\mathcal{A}$ 
we similarly have an isomorphism
$T(Con(u'),\mathbf{C})\cong Con(T(u',\mathbf{C}))$.

\subsubsection{} \label{sec:6.1.3}
\cite[Proposition 4.5.2]{Pe} Let $\mathbf{C}$ be a complex in 
$\mathcal{A}$. If $T(-,\mathbf{C})$ preserves monomorphisms 
then $\mathbf{C}$ is degreewise flat (\ref{sec:3.1}). Conversely, 
if $\mathcal{A}''$ satisfies AB4 and $\mathbf{C}$ is degreewise flat
then $T(-,\mathbf{C})$ preserves monomorphisms.

\subsubsection{} \label{sec:6.1.4}
Suppose $\mathcal{A}''$ satisfies AB4. Then \ref{sec:6.1.3} says that 
the class of flat objects of $Ch(\mathcal{A})$ is $Ch(\mathcal{F})$, hence 
(\ref{sec:4.1.2}) in this case they form the left class of a cotorsion theory.

\subsubsection{} \label{sec:6.1.5}
\cite[\S4 n\textsuperscript{\scriptsize{o}} 2 cor. 1]{Bou}
Suppose that $\mathcal{A}''$ satisfies AB4 and that $\mathbf{C}$ is a complex in 
$\mathcal{A}$ that is degreewise flat with zero differential. Then 
$T(-,\mathbf{C})$ preserves exact complexes.
\begin{proof}
Let $(\mathbf{C}',d')$ be an exact complex in $\mathcal{A}'$. The differential 
of $T(\mathbf{C}',\mathbf{C})$ is $T(d',\mathbf{C})$; we have then 
a natural isomorphism $ T(B(\mathbf{C}'),\mathbf{C})\cong B(T(\mathbf{C}',\mathbf{C}))$
and by \ref{sec:6.1.3} a natural isomorphism 
$T(Z(\mathbf{C}'),\mathbf{C})\cong Z(T(\mathbf{C}',\mathbf{C}))$.
The assertion follows from \ref{sec:6.1.3}.
\end{proof}

\subsubsection{} \label{sec:6.1.6}
\cite[\S4 n\textsuperscript{\scriptsize{o}} 3  lemme 1 et prop. 4]{Bou} 
Suppose that $\mathcal{A}''$ satisfies AB4.
If $\mathbf{C}$ is a degreewise flat complex in $\mathcal{A}$ with 
$\mathbf{C}_{n}=0$ for $n<0$, then $T(-,\mathbf{C})$ preserves 
exact complexes and quasi-isomorphisms.
Similarly, if $\mathbf{C}'$ is a complex in $\mathcal{A}'$ that is degreewise 
flat with $\mathbf{C}'_{n}=0$ for $n<0$, then $T(\mathbf{C}',-)$  
preserves exact complexes and quasi-isomorphisms.
\begin{proof}
Let $(\mathbf{C}',d')$ be an exact complex in $\mathcal{A}'$.
We have $T(\mathbf{C}',\mathbf{C})_{n}=
\underset{q\geqslant 0}\bigoplus T(\mathbf{C}'_{n-q},\mathbf{C}_{q})$.
For each integer $r\geqslant 0$ we define the subcomplex 
$T^{(r)}(\mathbf{C}',\mathbf{C})$ of $T(\mathbf{C}',\mathbf{C})$ by
\[
T^{(r)}(\mathbf{C}',\mathbf{C})_{n}=
\underset{0\leqslant q\leqslant r}\bigoplus T(\mathbf{C}'_{n-q},\mathbf{C}_{q})
\]
Let $u^{r}:T^{(r)}(\mathbf{C}',\mathbf{C})\to 
T^{(r+1)}(\mathbf{C}',\mathbf{C})$
be the natural monomorphism; we have $\Co u^{r}=
T(\mathbf{C}',S^{r+1}(\mathbf{C}_{r+1}))$.
By \ref{sec:6.1.5} the complex $\Co u^{r}$ is exact, hence $u^{r}$ is 
a quasi-isomorphism. Since $T^{(0)}(\mathbf{C}',\mathbf{C})$ is exact 
(\ref{sec:3.1.1}), induction on $r\geqslant 0$ shows that 
$T^{(r)}(\mathbf{C}',\mathbf{C})$ is exact for each 
$r\geqslant 0$. We have a natural isomorphism 
\[
\underset{r\in\mathbb{N}}\cl T^{(r)}(\mathbf{C}',\mathbf{C})\cong 
T(\mathbf{C}',\mathbf{C})
\]
Let $\sigma_{r}:T^{(r)}(\mathbf{C}',\mathbf{C})\to
\underset{r\geqslant 0}\bigoplus T^{(r)}(\mathbf{C}',\mathbf{C})$
be coprojection morphisms. Let 
\[
c:\underset{r\geqslant 0}\bigoplus T^{(r)}(\mathbf{C}',\mathbf{C})\to
\underset{r\geqslant 0}\bigoplus T^{(r)}(\mathbf{C}',\mathbf{C})
\]
be the unique morphism such that $c\sigma_{r}=\sigma_{r+1}u^{r}-\sigma_{r}$
for $r\geqslant 0$; then $\Co c=\underset{r\in\mathbb{N}}
\cl T^{(r)}(\mathbf{C}',\mathbf{C})$. Since $u^{r}_{n}$ is a split 
monomorphism for $r\geqslant 0$ and all integers $n$, the morphism
$c$ is a monomorphism (\ref{sec:4.2.7}). Since $\mathcal{A}''$ satisfies AB4 
the homology exact sequence associated to 
\[
\underset{r\geqslant 0}\bigoplus T^{(r)}(\mathbf{C}',\mathbf{C})\overset{c}
\rightarrowtail\underset{r\geqslant 0}\bigoplus T^{(r)}(\mathbf{C}',\mathbf{C})
\twoheadrightarrow T(\mathbf{C}',\mathbf{C})
\]
 implies that 
$T(\mathbf{C}',\mathbf{C})$ is exact. Let now $u':\mathbf{D}'\to\mathbf{C}'$ 
be a quasi-isomorphism in $Ch(\mathcal{A}')$.
We factor (\ref{sec:4.1.3}) $u'$ into $\mathbf{D}'\overset{\tilde{u'}}
\to Cyl(u')\overset{\beta'}\to \mathbf{C}'$, where $\tilde{u'}$ is a monomorphism 
and $\beta'$ is a homotopy equivalence. Since $T(-,\mathbf{C})$ preserves 
homotopy equivalences, we can assume that $u'$ is a monomorphism. 
Then $\Co u'$ is exact and we have an exact sequence 
\begin{displaymath}
\xymatrix{
{T(\mathbf{D}',\mathbf{C})}\ar@{>->}[rr]^{T(u',\mathbf{C})}&&
{T(\mathbf{C}',\mathbf{C})}
\ar@{->>}[r]&{T(\Co u',\mathbf{C})}
}	
\end{displaymath}
By the first part $ T(\Co u',\mathbf{C})$ is exact, therefore $T(u',\mathbf{C})$
is a quasi-isomorphism. The case of $T(\mathbf{C}',-)$ is dealt with similarly.
\end{proof}

\subsubsection{} \label{sec:6.1.7} 
We denote by $\mathfrak{TC}$ the class of complexes $\mathbf{C}$ in $\mathcal{A}$
for which $T(-,\mathbf{C})$ preserves monomorphisms with exact cokernel. 
Suppose that $\mathcal{A}',\mathcal{A}$ are cocomplete and $\mathcal{A}''$
satisfies AB4. Then $(\mathfrak{TC},\mathfrak{TC}^{\perp})$ is a cotorsion theory.

\begin{proof}
We shall use \ref{sec:3.2.4}. We take $\mathcal{K}_{1}=ex(\mathcal{A}')$, the class
of exact complexes in $\mathcal{A}'$, $\mathcal{K}_{2}=Ob(Ch(\mathcal{A}))$
and $\mathcal{K}_{3}=Ob(Ch(\mathcal{A}''))$. It suffices to show that 
$Ch(\mathcal{A})$ and $Ch(\mathcal{A}')$ have enough flat objects and that
$Ch(\mathcal{A}'')$ has enough injective objects. This follows from \ref{sec:4.1.2}.
\end{proof}

\subsubsection{} \label{sec:6.1.8}
The category $Ch(\mathrm{Ab})$ is a closed symmetric monoidal category 
with monoidal product $\otimes=T$ (\ref{sec:6.1}) and unit $S^{0}(\mathbb{Z})$. 
A category enriched over $(Ch(\mathrm{Ab}),\otimes,S^{0}(\mathbb{Z}))$ 
will be called dg-category and a 
$(Ch(\mathrm{Ab}),\otimes,S^{0}(\mathbb{Z}))$-functor will be called 
dg-functor. If $\mathcal{A}$ is an abelian category then $Ch(\mathcal{A})$ 
is a dg-category whose dg-hom is $\mathrm{Homgr}_{\mathcal{A}}$ (\ref{sec:4.2}).

\subsubsection{} \label{sec:6.1.9}
The standard THC-situation is an enriched THC-situation. 
By this we mean that $T,H,C$ are dg-functors and there are isomorphisms
\[
\mathrm{Homgr}_{\mathcal{A}''}(T(\mathbf{C}',\mathbf{C}),\mathbf{C}'')\cong 
\mathrm{Homgr}_{\mathcal{A}'}(\mathbf{C}',H(\mathbf{C},\mathbf{C}''))
\]
\[
\mathrm{Homgr}_{\mathcal{A}''}(T(\mathbf{C}',\mathbf{C}),\mathbf{C}'')\cong 
\mathrm{Homgr}_{\mathcal{A}}(\mathbf{C},C(\mathbf{C}',\mathbf{C}''))
\]
which are dg-natural in $\mathbf{C}',\mathbf{C}$ and $\mathbf{C}''$ 
\cite[Chapter 6, definition 6.7.1]{Bo2}.
\begin{proof}
We will show in detail only a part of the assertion.

\paragraph{Step 1} We sketch that $T$ is a dg-functor. For 
complexes $(\mathbf{C}',d^{\mathbf{C}'}),(\mathbf{D}',d^{\mathbf{D}'})$
in $\mathcal{A}'$ and $(\mathbf{C},d^{\mathbf{C}}),(\mathbf{D},d^{\mathbf{D}})$ 
in $\mathcal{A}$ we need to define a morphism
\[
T_{\mathbf{C}'\mathbf{C}\mathbf{D}'\mathbf{D}}:
\mathrm{Homgr}_{\mathcal{A}'}(\mathbf{C}',\mathbf{D}')\otimes
\mathrm{Homgr}_{\mathcal{A}}(\mathbf{C},\mathbf{D})
\to \mathrm{Homgr}_{\mathcal{A}''}
(T(\mathbf{C}',\mathbf{C}),T(\mathbf{D}',\mathbf{D}))
\]
This is equivalent to defining a morphism
\[
u_{p,q}:\mathrm{Homgr}_{\mathcal{A}'}(\mathbf{C}',\mathbf{D}')_{p}
\otimes_{\mathbb{Z}}\mathrm{Homgr}_{\mathcal{A}}(\mathbf{C},\mathbf{D})_{q}
\to \mathrm{Homgr}_{\mathcal{A}''}
(T(\mathbf{C}',\mathbf{C}),T(\mathbf{D}',\mathbf{D}))_{p+q}
\]
for all integers $p,q$ such that
\[
D(T(\mathbf{C}',\mathbf{C}),T(\mathbf{D}',\mathbf{D}))_{p+q}u_{p,q}=
u_{p-1,q}(D(\mathbf{D}',\mathbf{C}')\otimes 1)+(-1)^{p}u_{p,q-1}(1\otimes
D(\mathbf{C},\mathbf{D})_{q})
\]
Let $\pi_{i}:\mathrm{Homgr}_{\mathcal{A}''}
(T(\mathbf{C}',\mathbf{C}),T(\mathbf{D}',\mathbf{D}))_{p+q}\to
\mathcal{A}''(T(\mathbf{C}',\mathbf{C})_{i},T(\mathbf{D}',\mathbf{D})_{i+p+q})$
be projection morphisms. We define a morphism 
\[
\mu_{i,p,q}:\mathrm{Homgr}_{\mathcal{A}'}(\mathbf{C}',\mathbf{D}')_{p}
\otimes_{\mathbb{Z}}\mathrm{Homgr}_{\mathcal{A}}(\mathbf{C},\mathbf{D})_{q}
\to\mathcal{A}''(T(\mathbf{C}',\mathbf{C})_{i},T(\mathbf{D}',\mathbf{D})_{i+p+q})
\]
as follows. Let $(f'_{k})_{k}\otimes (f_{l})_{l}\in
\mathrm{Homgr}_{\mathcal{A}'}(\mathbf{C}',\mathbf{D}')_{p}
\otimes_{\mathbb{Z}}\mathrm{Homgr}_{\mathcal{A}}(\mathbf{C},\mathbf{D})_{q}$;
then $\mu_{i,p,q}((f'_{k})_{k}\otimes (f_{l})_{l})$ is the unique morphism 
$T(\mathbf{C}',\mathbf{C})_{i}\to T(\mathbf{D}',\mathbf{D})_{i+p+q}$ such that 
\[
\mu_{i,p,q}((f'_{k})_{k}\otimes (f_{l})_{l})\sigma_{x,y}=(-1)^{qx}\sigma_{x+p,y+q}
T(f'_{x},f_{y})
\]
for all $x+y=i$, where $\sigma_{x,y}:T(\mathbf{C}'_{x},\mathbf{C}_{y})
\to T(\mathbf{C}',\mathbf{C})_{i}$
and $\sigma_{x+p,y+q}:T(\mathbf{D}'_{x+p},\mathbf{D}_{y+q})
\to T(\mathbf{D}',\mathbf{D})_{i+p+q}$
are coprojection morphisms. We now define $u_{p,q}$ to be the 
unique morphism such that $\pi_{i}u_{p,q}=\mu_{i,p,q}$.
\paragraph{} We calculate $u_{p-1,q}(D(\mathbf{D}',\mathbf{C}')\otimes 1)$.
For $x+y=i$ we have
\[
\begin{aligned}
M\overset{def}=\pi_{i}u_{p-1,q}(D(\mathbf{D}',\mathbf{C}')\otimes 1)
((f'_{k})_{k}\otimes (f_{l})_{l})\sigma_{x,y}=\\
\mu_{i,p-1,q}((d^{\mathbf{D}'}_{k+p}f'_{k}-
(-1)^{p}f'_{k-1}d^{\mathbf{C}'}_{k})_{k}\otimes (f_{l})_{l})
\sigma_{x,y}=\\
(-1)^{qx}\sigma_{x+p-1,y+q}T(d^{\mathbf{D}'}_{x+p}f'_{x}-
(-1)^{p}f'_{x-1}d^{\mathbf{C}'}_{x},f_{y})=\\
\end{aligned}
\]
\[
(-1)^{qx}\sigma_{x+p-1,y+q}T(d^{\mathbf{D}'}_{x+p}f'_{x},f_{y})-
(-1)^{qx+p}\sigma_{x+p-1,y+q}T(f'_{x-1}d^{\mathbf{C}'}_{x},f_{y})
\]
\paragraph{} We calculate $(-1)^{p}u_{p,q-1}(1\otimes D(\mathbf{C},\mathbf{D})_{q})$.
For $x+y=i$ we have
\[
\begin{aligned}
N\overset{def}=(-1)^{p}\pi_{i}u_{p,q-1}(1\otimes D(\mathbf{C},\mathbf{D})_{q})
((f'_{k})_{k}\otimes (f_{l})_{l})\sigma_{x,y}=\\
(-1)^{p}\mu_{i,p,q-1}((f'_{k})_{k}\otimes 
(d^{\mathbf{D}}_{l+q}f_{l}-(-1)^{q}f_{l-1}d^{\mathbf{C}}_{l})_{l})
\sigma_{x,y}=\\
\end{aligned}
\]
\[
(-1)^{p+(q-1)x}\sigma_{x+p,y+q-1}T(f'_{x},d^{\mathbf{D}}_{y+q}f_{y})-
(-1)^{p+q+(q-1)x}\sigma_{x+p,y+q-1}T(f'_{x},f_{y-1}d^{\mathbf{C}}_{y})
\]
\paragraph{} We calculate $D(T(\mathbf{C}',\mathbf{C}),
T(\mathbf{D}',\mathbf{D}))_{p+q}u_{p,q}$. For $x+y=i$ we have
\[
P\overset{def}=\pi_{i}D(T(\mathbf{C}',\mathbf{C}),
T(\mathbf{D}',\mathbf{D}))_{p+q}u_{p,q}((f'_{k})_{k}\otimes (f_{l})_{l})\sigma_{x,y}=
\]
\[
\mathcal{A}''(1,D(\mathbf{D}',\mathbf{D})_{i+p+q})
(\mu_{i,p,q}((f'_{k})_{k}\otimes (f_{l})_{l}))\sigma_{x,y}-
\]
\[
(-1)^{p+q}\mathcal{A}''(D(\mathbf{C}',\mathbf{C})_{i},1)
(\mu_{i-1,p,q}((f'_{k})_{k}\otimes (f_{l})_{l}))\sigma_{x,y}
\]
Now 
\[
\begin{aligned}
\mathcal{A}''(1,D(\mathbf{D}',\mathbf{D})_{i+p+q})
(\mu_{i,p,q}((f'_{k})_{k}\otimes (f_{l})_{l}))\sigma_{x,y}=\\
D(\mathbf{D}',\mathbf{D})_{i+p+q}\mu_{i,p,q}((f'_{k})_{k}\otimes (f_{l})_{l})\sigma_{x,y}=\\
\end{aligned}
\]
\[
(-1)^{qx}\sigma_{x+p-1,y+q}T(d^{\mathbf{D}'}_{x+p}f'_{x},f_{y})+(-1)^{qx+x+p}
\sigma_{x+p,y+q-1}T(f'_{x},d^{\mathbf{D}}_{y+q}f_{y})
\]
and
\[
\begin{aligned}
(-1)^{p+q}\mathcal{A}''(D(\mathbf{C}',\mathbf{C})_{i},1)
(\mu_{i-1,p,q}((f'_{k})_{k}\otimes (f_{l})_{l}))\sigma_{x,y}=\\
(-1)^{p+q}\mu_{i-1,p,q}((f'_{k})_{k}\otimes (f_{l})_{l})
D(\mathbf{C}',\mathbf{C})_{i}\sigma_{x,y}=
\end{aligned}
\]
\[
(-1)^{p+q+q(x-1)}\sigma_{x+p-1,y+q}T(f'_{x-1}d^{\mathbf{C}'}_{x},f_{y})
-(-1)^{p+q+x+qx}\sigma_{x+p,y+q-1}T(f'_{x},f_{y-1}d^{\mathbf{C}}_{y})
\]
We obtain $M+N=P$, as needed.
\paragraph{Step 2} We have 
\[
\mathrm{Homgr}_{\mathcal{A}''}(T(\mathbf{C}',\mathbf{C}),\mathbf{C}'')_{n}
=\underset{i}\prod\mathcal{A}''(\underset{p+q=i}\bigoplus
T(\mathbf{C}'_{p},\mathbf{C}_{q}),\mathbf{C}''_{i+n})
\]
and
\[
\mathrm{Homgr}_{\mathcal{A}'}(\mathbf{C}',H(\mathbf{C},\mathbf{C}''))_{n}
=\underset{p}\prod\mathcal{A}'(\mathbf{C}'_{p},\underset{q}\prod
H(\mathbf{C}_{q},\mathbf{C}''_{p+q+n}))
\]
We define $\Theta_{n}:\mathrm{Homgr}_{\mathcal{A}''}
(T(\mathbf{C}',\mathbf{C}),\mathbf{C}'')_{n}\to
\mathrm{Homgr}_{\mathcal{A}'}(\mathbf{C}',H(\mathbf{C},\mathbf{C}''))_{n}$
as 
\[
\Theta_{n}((u_{p,q}:T(\mathbf{C}'_{p},\mathbf{C}_{q})\to\mathbf{C}''_{p+q+n}))
=(v_{p}:\mathbf{C}'_{p}\to H(\mathbf{C},\mathbf{C}'')_{p+n})
\]
where $v_{p}$ is defined as the unique morphism such that $\pi_{q}v_{p}$
is the adjoint transpose $v_{p,q}:\mathbf{C}'_{p}\to 
H(\mathbf{C}_{q},\mathbf{C}''_{p+q+n})$ of $u_{p,q}$ and where $\pi_{q}$
is the projection morphism. It is clear that $\theta_{n}$ is an isomorphism of 
graded objects. We show that $\Theta$ is a morphism. We denote by 
$D(\mathbf{C}',\mathbf{C};\mathbf{C}'')$ the differential of 
$\mathrm{Homgr}_{\mathcal{A}''}(T(\mathbf{C}',\mathbf{C}),\mathbf{C}'')$
and by $D(\mathbf{C}';\mathbf{C},\mathbf{C}'')$ the differential of
$\mathrm{Homgr}_{\mathcal{A}'}(\mathbf{C}',H(\mathbf{C},\mathbf{C}''))$.
We calculate $D(\mathbf{C}',\mathbf{C};\mathbf{C}'')$.
Let $f_{i}: T(\mathbf{C}',\mathbf{C})_{i}\to \mathbf{C}''_{i+n}$ be a 
morphism and $\sigma_{p,q}:T(\mathbf{C}'_{p},\mathbf{C}_{q})
\to T(\mathbf{C}',\mathbf{C})_{i}$, where $p+q=i$, be coprojection
morphisms. If $u_{p,q}=f_{p+q}\sigma_{p,q}$ then 
\[
D(\mathbf{C}',\mathbf{C};\mathbf{C}'')_{n}((u_{p,q}))=
d''_{p+q+n}u_{p,q}-(-1)^{n}u_{p-1,q}T(d'_{p},1)-(-1)^{n+p}u_{p,q-1}T(1,d_{q})
\]
We calculate $D(\mathbf{C}';\mathbf{C},\mathbf{C}'')$. Let $g_{p}:\mathbf{C}'_{p}
\to H(\mathbf{C},\mathbf{C}'')_{p+n}$ be a morphism and 
$\pi_{q}: H(\mathbf{C},\mathbf{C}'')_{p+n}\to H(\mathbf{C}_{q},\mathbf{C}''_{p+q+n})$
be projection morphisms. If $v_{p,q}=\pi_{q}g_{p}$ then
\[
D(\mathbf{C}';\mathbf{C},\mathbf{C}'')_{n}((v_{p,q}))=
H(1,d''_{p+q+n})v_{p,q}-(-1)^{n}v_{p-1,q}d'_{p}-(-1)^{p+n}H(d_{q},1)v_{p,q-1}
\]
We have
\[
D(\mathbf{C}';\mathbf{C},\mathbf{C}'')_{n}\Theta_{n}((u_{p,q}))=
(H(1,d''_{p+q+n})v_{p,q}-(-1)^{n}v_{p-1,q}d'_{p}-(-1)^{p+n}H(d_{q},1)v_{p,q-1})
\]
and
\[
\Theta_{n-1}D(\mathbf{C}',\mathbf{C};\mathbf{C}'')_{n}((u_{p,q}))=
\]
\[
\Theta_{n-1}((d''_{p+q+n}u_{p,q}-(-1)^{n}u_{p-1,q}T(d'_{p},1)-(-1)^{n+p}u_{p,q-1}T(1,d_{q}))=(v'_{p})
\]
where  $\pi_{q}v'_{p}$ is the adjoint transpose of 
\[
d''_{p+q+n}u_{p,q}-(-1)^{n}u_{p-1,q}T(d'_{p},1)-(-1)^{n+p}u_{p,q-1}T(1,d_{q})
\]
The adjoint transpose of $d''_{p+q+n}u_{p,q}$ is $H(1,d''_{p+q+n})v_{p,q}$,
the adjoint transpose of $u_{p-1,q}T(d'_{p},1)$ is $v_{p-1,q}d'_{p}$ and 
the adjoint transpose of $u_{p,q-1}T(1,d_{q})$ is $H(d_{q},1)v_{p,q-1}$.
The desired equality follows.
\end{proof}

\subsection{} \label{sec:6.2}
Let 
\[
T:\mathcal{A}'\times \mathcal{A}\to \mathcal{A}'', \
H:\mathcal{A}^{op}\times \mathcal{A}''\to \mathcal{A}', \
C:\mathcal{A}'^{op}\times \mathcal{A}''\to \mathcal{A}
\]
be an abelian THC--situation. We assume that $\mathcal{A}',\mathcal{A}$
have all limits and $\mathcal{A}''$ has all colimits.
We revisit the bar THC--situation 
\[
\overline{T}:Ch(\mathcal{A}')\times Ch(\mathcal{A})\to Ch(\mathcal{A}'')
\]
\[
\overline{H}:Ch(\mathcal{A})^{op}\times Ch(\mathcal{A}'')\to Ch(\mathcal{A}')
\]
\[
\overline{C}:Ch(\mathcal{A}')^{op}\times Ch(\mathcal{A}'')\to Ch(\mathcal{A})
\]
associated to the given one \cite[Section 4.3]{Pe}. Let $(\mathbf{C}',d')$ (resp. $(\mathbf{C},d),
(\mathbf{C}'',d'')$) be complexes in $\mathcal{A}'$ (resp. $\mathcal{A},\mathcal{A}''$).
\paragraph{(Tensor)}
Consider the bifunctor $T$ from \ref{sec:6.1}. Since $T(d',\mathbf{C})$ is a morphism
we have (\ref{sec:4.1}) the commutative diagram 
\[
\xymatrix{
{T(\mathbf{C}',\mathbf{C})(1)}\ar[rr]^{T(d',\mathbf{C})(1)} \ar[d]_{D(\mathbf{C}',\mathbf{C})^{\#}}
&&{T(\mathbf{C}'(-1),\mathbf{C})(1)} \ar[d]^{D(\mathbf{C}'(-1),\mathbf{C})^{\#}}\\
{T(\mathbf{C}',\mathbf{C})} \ar[rr]_{T(d',\mathbf{C})}
&&{T(\mathbf{C}'(-1),\mathbf{C})}
}
\]
where $D(\mathbf{C}',\mathbf{C})^{\#}$ is the adjoint transpose of $D(\mathbf{C}',\mathbf{C})$.
We have then an induced morphism
\[
B(T(d',\mathbf{C})):B(T(\mathbf{C}',\mathbf{C}))\to B(T(\mathbf{C}'(-1),\mathbf{C}))
\]
and a commutative solid arrows diagram
\[
\xymatrix{
{B(T(\mathbf{C}',\mathbf{C}))} \ar@{>->}[r]
\ar[d]_{B(T(d',\mathbf{C}))}& 
{T(\mathbf{C}',\mathbf{C})}\ar@{->>}[rr]^{\pi}
\ar[d]^{T(d',\mathbf{C})} & &
{T(\mathbf{C}',\mathbf{C})/B(T(\mathbf{C}',\mathbf{C}))}
\ar@{..>}[d]^{\overline{D}(\mathbf{C}',\mathbf{C})}\\
{B(T(\mathbf{C}'(-1),\mathbf{C}))} \ar@{>->}[r] & 
{T(\mathbf{C}'(-1),\mathbf{C})}\ar@{->>}[rr]^{\pi'} & & 
{T(\mathbf{C}'(-1),\mathbf{C})/B(T(\mathbf{C}'(-1),\mathbf{C}))}
}
\]
There is a unique morphism $\overline{D}(\mathbf{C}',\mathbf{C})=\overline{D}$
such that $\overline{D}(\mathbf{C}',\mathbf{C})\pi=\pi'T(d',\mathbf{C})$.
By \ref{sec:6.1.1} we have 
\[
T(\mathbf{C}'(-1),\mathbf{C})/B(T(\mathbf{C}'(-1),\mathbf{C}))=
(T(\mathbf{C}',\mathbf{C})/B(T(\mathbf{C}',\mathbf{C})))(-1)
\]
hence if we define $\overline{T}(\mathbf{C}',\mathbf{C})=
T(\mathbf{C}',\mathbf{C})/B(T(\mathbf{C}',\mathbf{C}))$ we have 
\[
\overline{D}(\mathbf{C}',\mathbf{C}):\overline{T}(\mathbf{C}',\mathbf{C})
\to \overline{T}(\mathbf{C}',\mathbf{C})(-1)
\]
Since $d'^{2}=0$ we have by construction that $\overline{D}^{2}=0$.
We omit the details that $\overline{T}$ assembles into a bifunctor.

\paragraph{(Hom)} 
Consider the bifunctor $H$ from \ref{sec:6.1}. From the commutative diagram
\[
\xymatrix{
{H(\mathbf{C},\mathbf{C}'')}\ar[rr]^{H(\mathbf{C},d'')} \ar[d]_{D(\mathbf{C},\mathbf{C}'')}
&&{H(\mathbf{C},\mathbf{C}''(-1))} \ar[d]^{D(\mathbf{C},\mathbf{C}''(-1))}\\
{H(\mathbf{C},\mathbf{C}'')(-1)} \ar[rr]_{H(\mathbf{C},d'')(-1)}
&&{H(\mathbf{C},\mathbf{C}'')(-1)}
}
\]
we have an induced morphism 
\[
Z(H(\mathbf{C},d'')):Z(H(\mathbf{C},\mathbf{C}''))\to Z(H(\mathbf{C},\mathbf{C}''(-1)))
\]
By \ref{sec:6.1.1} we have $Z(H(\mathbf{C},\mathbf{C}''(-1)))=
Z(H(\mathbf{C},\mathbf{C}''))(-1)$ hence if we define $\overline{H}(\mathbf{C},\mathbf{C}'')=
Z(H(\mathbf{C},\mathbf{C}''))$ and $\overline{D}=\overline{D}(\mathbf{C},\mathbf{C}'')=
Z(H(\mathbf{C},d''))$ we have 
\[
\overline{D}(\mathbf{C},\mathbf{C}''):\overline{H}(\mathbf{C},\mathbf{C}'')
\to \overline{H}(\mathbf{C},\mathbf{C}'')(-1)
\]
Since $d''^{2}=0$ we have by construction that $\overline{D}^{2}=0$.
We omit the details that $\overline{H}$ assembles into a bifunctor.
\paragraph{(Cotensor)}
Consider the following modified version $C$ of the cotensor bifunctor from \ref{sec:6.1}.
We define $C(\mathbf{C}',\mathbf{C}'')_{n}=
\underset{p}\prod C(\mathbf{C}'_{p},\mathbf{C}''_{p+n})$ with differential
$D(\mathbf{C}',\mathbf{C}'')_{n}=D_{n}$ the unique morphism that satisfies
\[
\pi_{p}D_{n}=C(\mathbf{C}'_{p},d''_{p+n})\pi_{p}-(-1)^{n-1}C(d'_{p},\mathbf{C}''_{p+n-1})\pi_{p-1}
\]
for all integers $p$. Then for all integers $r$ we have
\begin{equation}\label{eq:4.3.0}
(C(\mathbf{C}',\mathbf{C}'')(r),D(\mathbf{C}',\mathbf{C}'')(r))=
(C(\mathbf{C}',\mathbf{C}''(r)),D(\mathbf{C}',\mathbf{C}''(r)))
\end{equation}
From the commutative diagram
\[
\xymatrix{
{C(\mathbf{C}',\mathbf{C}'')}\ar[rr]^{C(\mathbf{C}',d'')} \ar[d]_{D(\mathbf{C}',\mathbf{C}'')}
&&{C(\mathbf{C}',\mathbf{C}''(-1))} \ar[d]^{D(\mathbf{C}',\mathbf{C}''(-1))}\\
{C(\mathbf{C}',\mathbf{C}'')(-1)} \ar[rr]_{C(\mathbf{C}',d'')(-1)}
&&{C(\mathbf{C}',\mathbf{C}'')(-1)}
}
\]
we have an induced morphism 
\[
Z(C(\mathbf{C}',d'')):Z(C(\mathbf{C}',\mathbf{C}''))\to Z(C(\mathbf{C},\mathbf{C}''(-1)))
\]
By (\ref{eq:4.3.0}) we have $Z(C(\mathbf{C}',\mathbf{C}''(-1)))=
Z(C(\mathbf{C}',\mathbf{C}''))(-1)$ hence if we define $\overline{C}(\mathbf{C}',\mathbf{C}'')=
Z(C(\mathbf{C}',\mathbf{C}''))$ and $\overline{D}=\overline{D}(\mathbf{C}',\mathbf{C}'')=
Z(C(\mathbf{C},d''))$ we have 
\[
\overline{D}(\mathbf{C}',\mathbf{C}''):\overline{C}(\mathbf{C}',\mathbf{C}'')
\to \overline{C}(\mathbf{C}',\mathbf{C}'')(-1)
\]
Since $d''^{2}=0$ we have by construction that $\overline{D}^{2}=0$.
We omit the details that $\overline{C}$ assembles into a bifunctor.

\paragraph{(Tensor $\Leftrightarrow$ Hom)}
To give a morphism $\overline{T}(\mathbf{C}',\mathbf{C})\to \mathbf{C}''$ is to give
a morphism of graded objects  $u:T(\mathbf{C}',\mathbf{C})\to \mathbf{C}''$ such that
$uD(\mathbf{C}',\mathbf{C})(1)=0$, where 
\begin{displaymath}
\xymatrix{
{T(\mathbf{C}',\mathbf{C})(1)}\ar[rr]^{D(\mathbf{C}',\mathbf{C})(1)} 
&& {T(\mathbf{C}',\mathbf{C})} \ar[r]^{u} & {\mathbf{C}''}
}
\end{displaymath}
with $D(\mathbf{C}',\mathbf{C})(1)_{n}=D(\mathbf{C}',\mathbf{C})_{n+1}$
and such that the diagram 
\begin{displaymath}
\xymatrix{
{T(\mathbf{C}',\mathbf{C})}\ar[rr]^{u} \ar[d]_{T(d',\mathbf{C})}
&&{\mathbf{C}''} \ar[d]^{d''}\\
{T(\mathbf{C}'(-1),\mathbf{C})} \ar[rr]_{u(-1)}
&&{\mathbf{C}''(-1)}
}
\end{displaymath}
commutes. This data is equivalent to giving morphisms
$u_{p,q}:T(\mathbf{C}'_{p},\mathbf{C}_{q})\to\mathbf{C}''_{p+q}$
such that
\begin{equation}\label{eq:4.3}
u_{p-1,q}T(d'_{p},\mathbf{C}_{q})+(-1)^{p}u_{p,q-1}T(\mathbf{C}'_{p},d_{q})=0
\end{equation}
and 
\begin{equation}\label{eq:4.3'}
d''_{p+q}u_{p,q}=u_{p-1,q}T(d'_{p},\mathbf{C}_{q})
\end{equation}
for all integers $p,q$. To give a morphism $\mathbf{C}'\to \overline{H}(\mathbf{C},\mathbf{C}'')$ 
is to give a morphism of graded objects $v:\mathbf{C}'\to H(\mathbf{C},\mathbf{C}'')$ such that
$D(\mathbf{C},\mathbf{C}'')v=0$ and such that the diagram 
\begin{displaymath}
\xymatrix{
{\mathbf{C}'}\ar[rr]^{v} \ar[d]_{d'}
&&{H(\mathbf{C},\mathbf{C}'')} \ar[d]^{H(\mathbf{C},d'')}\\
{\mathbf{C}'(-1)} \ar[rr]_{v(-1)}
&&{H(\mathbf{C},\mathbf{C}''(-1))}
}
\end{displaymath}
commutes. This data is equivalent to giving morphisms
$v_{n,p}:\mathbf{C}'_{n}\to H(\mathbf{C}_{p},\mathbf{C}''_{p+n})$ such that
\begin{equation}\label{eq:4.3''}
H(\mathbf{C}_{n},d''_{n+p})v_{n,p}-(-1)^{n}H(d_{p},\mathbf{C}''_{n+p-1})v_{n,p-1}=0
\end{equation}
and 
\begin{equation}\label{eq:4.3'''}
v_{n-1,p}d'_{n}=H(\mathbf{C}_{n},d''_{n+p})v_{n,p}
\end{equation}
for all integers $n,p$. Given a morphism 
$\overline{T}(\mathbf{C}',\mathbf{C})\to \mathbf{C}''$ we define $v_{n,p}=(-1)^{p}u_{n,p}^{\flat}$,
where $u_{n,p}^{\flat}:\mathbf{C}'_{n}\to H(\mathbf{C}_{p},\mathbf{C}''_{p+n})$ 
is the adjoint transpose of $u_{n,p}$; then (\ref{eq:4.3''}) and (\ref{eq:4.3'''}) follow
from (\ref{eq:4.3}) and (\ref{eq:4.3'}) combined. Conversely, 
given a morphism $\mathbf{C}'\to \overline{H}(\mathbf{C},\mathbf{C}'')$ 
we define $u_{p,q}=(-1)^{q}v_{p,q}^{\#}$, where $v_{p,q}^{\#}$ is the adjoint transpose of $v_{p,q}$. 
\paragraph{(Tensor $\Leftrightarrow$ Cotensor)}
To give a morphism $\mathbf{C}\to \overline{C}(\mathbf{C}',\mathbf{C}'')$ is to give
a morphism $w:\mathbf{C}\to \overline{C}(\mathbf{C}',\mathbf{C}'')$ of graded objects 
such that $D(\mathbf{C}',\mathbf{C}'')w=0$ and such that the diagram
\begin{displaymath}
\xymatrix{
{\mathbf{C}}\ar[rr]^{w} \ar[d]_{d}
&&{C(\mathbf{C}',\mathbf{C}'')} \ar[d]^{C(\mathbf{C}',d'')}\\
{\mathbf{C}(-1)} \ar[rr]_{w(-1)}
&&{C(\mathbf{C}',\mathbf{C}''(-1))}
}
\end{displaymath}
commutes. This data is equivalent to giving morphisms
$w_{n,p}:\mathbf{C}_{n}\to C(\mathbf{C}'_{p},\mathbf{C}''_{p+n})$ such that
\begin{equation}\label{eq:4.3''''}
C(\mathbf{C}_{p},d''_{n+p})w_{n,p}-(-1)^{n-1}C(d'_{p},\mathbf{C}''_{n+p-1})w_{n,p-1}=0
\end{equation}
and 
\begin{equation}\label{eq:4.3''''''}
C(\mathbf{C}_{p},d''_{n+p})w_{n,p}=w_{n-1,p}d_{n}
\end{equation}
for all integers $n,p$. Given a morphism 
$\overline{T}(\mathbf{C}',\mathbf{C})\to \mathbf{C}''$ we define 
$w_{n,p}=(-1)^{\alpha(n,p)}u_{p,n}^{\flat}$, where 
$u_{p,n}^{\flat}:\mathbf{C}_{n}\to C(\mathbf{C}'_{p},\mathbf{C}''_{p+n})$ 
is the adjoint transpose of $u_{p,n}$ and
\[
\alpha(n,p)=\frac{n(n-1)}{2}+\frac{p(p-1)}{2}+\frac{(n+p)(n+p+1)}{2}
\]
see \cite[Proposition 2.1.2]{EGR}. With this definition one can check that 
(\ref{eq:4.3''''}) and (\ref{eq:4.3''''''})  follow from (\ref{eq:4.3}) and (\ref{eq:4.3'}) 
combined. Conversely, given a morphism $\mathbf{C}\to \overline{C}(\mathbf{C}',\mathbf{C}'')$ 
we define $u_{p,q}=(-1)^{\alpha(p,q)}w_{q,p}^{\#}$, where $w_{q,p}^{\#}$ 
is the adjoint transpose of 
$w_{q,p}:\mathbf{C}_{q}\to C(\mathbf{C}'_{p},\mathbf{C}''_{p+q})$.

\subsubsection{} \label{sec:6.2.1}
For all integers $r$ we have natural isomorphisms
\[
\overline{T}(D^{r}(A'),\mathbf{C})\cong T(S^{r}(A'),\mathbf{C})  {\rm \ and} \ 
\overline{H}(\mathbf{C},D^{r}(A''))\cong H(\mathbf{C},S^{r-1}(A''))
\] 

\subsubsection{} \label{sec:6.2.2}
Suppose that $\mathcal{A}''$ has a set $(E''_{i})_{i\in I}$ of injective cogenerators.
Let $\mathbf{C}$ be a complex in $\mathcal{A}$. Then $\overline{T}(-,\mathbf{C})$ 
preserves monomorphisms if and only if $H(\mathbf{C},S^{r}(E''_{i}))$ is an injective
complex for all integers $r$ and all $i\in  I$.
\begin{proof}
The set $(D^{r}(E''_{i}))_{i\in  I}$, where $r$ is an integer, is a set of injective cogenerators 
of $Ch(\mathcal{A}'')$. The assertion follows then from \ref{sec:3.1.2} and \ref{sec:6.2.1}.
\end{proof}

Recall (\ref{sec:3.1}) that $\mathcal{F}'$ denotes the class of flat objects of $\mathcal{A}'$.

\subsubsection{} \label{sec:6.2.3} 
\cite[Proposition 2.4]{EGR}
Suppose that $\mathcal{A}''$ has a set $(E''_{i})_{i\in I}$ of injective cogenerators.
For a complex $(\mathbf{C},d)$ in $\mathcal{A}$ consider the following statements:

(1) $\mathbf{C}\in ex[\mathcal{F}]$ (\ref{sec:4.1.2});

(2) $H(\mathbf{C},S^{r}(E''))$ is injective for all integers $r$ and all injective objects 
$E''$ of $\mathcal{A}''$;

(3) $H(\mathbf{C},S^{r}(E''_{i}))$ is exact and $Z(H(\mathbf{C},S^{r}(E''_{i})))$ 
is degreewise in $\mathcal{F}'^{\perp}$ for all integers $r$ and all $i\in  I$.

We have $(1)\Rightarrow (2)\Rightarrow (3)$. If $\mathcal{A}'$ has a faithfully flat 
object (\ref{sec:3.4}) then part 3 implies that $\mathbf{C}$ is exact. If, moreover, 
$Z(H(\mathbf{C},S^{r}(E''_{i})))$ is degreewise injective for all integers $r$ and all 
$i\in  I$ then $(3)\Rightarrow (1)$.
\begin{proof}
We have  $H(\mathbf{C},S^{r}(E''))_{n}=H(\mathbf{C}_{r-n},E'')$ and the differential of
$H(\mathbf{C},S^{r}(E''))$ is $D(\mathbf{C},S^{r}(E''))_{n}=(-1)^{n-1}H(d_{r-n+1},E'')$.

$(1)\Rightarrow (2)$ Consider the commutative diagram
\[
\xymatrix{
{\mathbf{C}_{r-n+1}} \ar[rr]^{(-1)^{n-1}d_{r-n+1}}\ar@{->>}[dr] && 
{\mathbf{C}_{r-n}} \ar[rr]^{(-1)^{n}d_{r-n}}\ar@{->>}[dr]_{p}
& & {\mathbf{C}_{r-n-1}}\\
& {B_{r-n}(\mathbf{C})} \ar@{>->}[ur]_{i} && {B_{r-n-1}(\mathbf{C})} \ar@{>->}[ur]
}
\]
where the sequence $(i,p)$ is exact. Applying the functor $H(-,E'')$ to the sequence $(i,p)$
we obtain (\ref{sec:3.1.3}) the exact sequence
\[
H(B_{r-n-1}(\mathbf{C}),E'')\rightarrowtail H(\mathbf{C}_{r-n},E'')\twoheadrightarrow
H(B_{r-n}(\mathbf{C}),E'')
\]
hence $H(\mathbf{C},S^{r}(E''))$ is exact. In particular we have the commutative diagram
\[
\xymatrix{
{H(\mathbf{C}_{r-n},E'')}\ar[rr]^{(-1)^{n-1}H(d_{r-n+1},E'')}\ar@{->>}[dr] 
&& {H(\mathbf{C}_{r-n+1},E'')}\\
& {H(B_{r-n}(\mathbf{C}),E'')} \ar@{>->}[ur]
}
\]
Since $B_{r-n}(\mathbf{C})$ is flat, by \ref{sec:3.1.2} we have that 
$H(B_{r-n}(\mathbf{C}),E'')$ is injective, hence $H(\mathbf{C},S^{r}(E''))$ is injective.
$(2)\Rightarrow (3)$ is clear. Let now $P'$ be a faithfully flat object of $\mathcal{A}'$.
By \ref{sec:3.4.4} applied to the sequence $(d_{r-n+1},d_{r-n})$ we obtain that 
$\mathbf{C}$ is exact. Then $Z_{n}(H(\mathbf{C},S^{r}(E'')))\cong 
H(B_{r-n-1}(\mathbf{C}),E''_{i})$ (cf. $(1)\Rightarrow (2)$), hence if
$Z_{n}(H(\mathbf{C},S^{r}(E''_{i})))$ is injective it follows from \ref{sec:3.1.2}
that $B_{r-n-1}(\mathbf{C})$ is flat.
\end{proof}

%\subsubsection{} \label{sec:4.5.4}
%characterization of flatness for the Gabriel THC-situation.

\subsubsection{} \label{sec:6.2.4}
%The category $Ch(\mathrm{Ab})$ is a closed (symmetric) monoidal category with monoidal
%product $\overline{\otimes}$ (\ref{sec:4.5}) and unit $D^{0}(\mathbb{Z})$. A category 
%enriched over $(Ch(\mathrm{Ab}),\overline{\otimes},D^{0}(\mathbb{Z}))$ will be called 
%$\overline{\text{dg}}$-category and a $(Ch(\mathrm{Ab}),\overline{\otimes},
%D^{0}(\mathbb{Z}))$-functor will be called $\overline{\text{dg}}$-functor.

Let $\mathcal{A}$ be an abelian category. For any two complexes $\mathbf{C},\mathbf{C}'$
in $\mathcal{A}$ we define the complex 
$\overline{\mathrm{Homgr}}_{\mathcal{A}}(\mathbf{C}',\mathbf{C})$ of abelian 
groups as 
\[
\overline{\mathrm{Homgr}}_{\mathcal{A}}(\mathbf{C}',\mathbf{C})_{n}=
Ch(\mathrm{Ab})(S^{n}(\mathbb{Z}),\mathrm{Homgr}_{\mathcal{A}}(\mathbf{C}',\mathbf{C}))
\]
with differential defined as $\overline{D}(\mathbf{C}',\mathbf{C})_{n}(u)=
\mathrm{Homgr}_{\mathcal{A}}(\mathbf{C}',d)u$,
where $d$ is the differential of $\mathbf{C}$. The morphism 
$\overline{D}(\mathbf{C}',\mathbf{C})_{n}$ is well defined since 
$S^{n}(\mathbb{Z})(1)=S^{n-1}(\mathbb{Z})$. We have
$\overline{\mathrm{Homgr}}_{\mathcal{A}}(\mathbf{C}',\mathbf{C})_{n}\cong
Ch(\mathcal{A})(\mathbf{C}',\mathbf{C}(n))$, hence 
$\overline{D}(\mathbf{C}',\mathbf{C})_{n}(u)=d(n)u$, where $d(n)$
is the differential of $\mathbf{C}(n)$. Thus,
\[
\overline{\mathrm{Homgr}}_{\mathcal{A}}(\mathbf{C}',\mathbf{C})\cong
\mathrm{Homgr}_{Ch(\mathcal{A})}(S^{0}(\mathbf{C}'),\mathbf{C}(-))
\]
where $(\mathbf{C}(-),d(-))$ is the complex in $Ch(\mathcal{A})$ defined
as $\mathbf{C}(-)_{n}=\mathbf{C}(n), d(-)_{n}=d(n)$.
%Then $Ch(\mathcal{A})$ is a 
%$\overline{\text{dg}}$-category whose $\overline{\text{dg}}$-hom we shall denote by 
%$\overline{\mathrm{Homgr}}_{\mathcal{A}}(\mathbf{C}',\mathbf{C})$; one has 

\subsubsection{} \label{sec:6.2.5}
\cite[Proposition 2.1.1)]{EGR} 
There are natural isomorphisms
\[
\overline{\mathrm{Homgr}}_{\mathcal{A}''}(\overline{T}(\mathbf{C}',\mathbf{C}),\mathbf{C}'')\cong 
\overline{\mathrm{Homgr}}_{\mathcal{A}'}(\mathbf{C}',\overline{H}(\mathbf{C},\mathbf{C}''))
\]
\[
\overline{\mathrm{Homgr}}_{\mathcal{A}''}(\overline{T}(\mathbf{C}',\mathbf{C}),\mathbf{C}'')\cong 
\overline{\mathrm{Homgr}}_{\mathcal{A}}(\mathbf{C},\overline{C}(\mathbf{C}',\mathbf{C}''))
\]
%which are $\overline{\text{dg}}$-natural in $\mathbf{C}',\mathbf{C}$ and $\mathbf{C}''$.
\begin{proof}
We prove the first isomorphism. We note that from \ref{sec:6.1.1} and the 
construction of $\overline{H}$ we have 
\[
(\overline{H}(\mathbf{C},\mathbf{C}'')(r),\overline{D}(\mathbf{C},\mathbf{C}'')(r))=
(\overline{H}(\mathbf{C},\mathbf{C}''(r)),\overline{D}(\mathbf{C},\mathbf{C}''(r)))
\]
Moreover, 
\begin{equation} \label{eq:4.6.1}
\overline{D}(\mathbf{C},\mathbf{C}''(r)))=\overline{H}(\mathbf{C},d''(r))
\end{equation} 
where $d''$ is the differential of $\mathbf{C}''$ and $d''(r)$ is the differential of 
$\mathbf{C}''(r)$. We then have 
\[
\begin{aligned}
\overline{\mathrm{Homgr}}_{\mathcal{A}''}(\overline{T}(\mathbf{C}',\mathbf{C}),\mathbf{C}'')_{n} & =
Ch(\mathcal{A}'')(\overline{T}(\mathbf{C}',\mathbf{C}),\mathbf{C}''(n))\\
& \cong Ch(\mathcal{A}')(\mathbf{C}',\overline{H}(\mathbf{C},\mathbf{C}''(n)) \\
& = Ch(\mathcal{A}')(\mathbf{C}',\overline{H}(\mathbf{C},\mathbf{C}'')(n)) \\
& =\overline{\mathrm{Homgr}}_{\mathcal{A}'}(\mathbf{C}',\overline{H}(\mathbf{C},\mathbf{C}''))_{n}
\end{aligned}
\]
By (\ref{eq:4.6.1}) this isomorphism is compatible with the differentials, 
hence we have the required isomorphism.
The second isomorphism is proved similarly.
\end{proof}

\section{The K\"{u}nneth theorem for complexes in abelian categories}

\subsection{} \label{sec:7.1}
Consider the standard THC--situation (\ref{sec:6.1}).
We assume that $\mathcal{A}''$ satisfies AB4. Recall (\ref{sec:3.1}) 
that $\mathcal{F}$ (resp. $\mathcal{F}'$) denotes the class of flat 
objects of $\mathcal{A}$ (resp. $\mathcal{A}'$). We assume
that $(\mathcal{F},\mathcal{F}^{\perp})$ and $(\mathcal{F}',
\mathcal{F}'^{\perp})$ are left complete (\ref{sec:5.1.6}). 
We simply say flat resolution (resp. complete flat resolution) 
for an $\mathcal{F}$-resolution (resp. complete $\mathcal{F}$-resolution) 
of an object of $\mathcal{A}$; same for the objects of $\mathcal{A}'$ 
with $\mathcal{F}'$ in place of $\mathcal{F}$.

\subsubsection{} \label{sec:7.1.1}
Let $A\in Ob(\mathcal{A}')$ and $A\in Ob(\mathcal{A})$.
For all integers $n$ we define $\mathsf{Tor}_{n}(A',A)=
H_{n}(T(\mathbf{C}'_{A'},\mathbf{C}_{A}))$,
where $\mathbf{C}'_{A'}$ (resp. $\mathbf{C}_{A}$) is a 
complete flat resolution of $A'$ (resp. $A$). 

Let $\mathbf{P}'$ (resp. $\mathbf{P}$) be a flat 
resolution of $A'$ (resp. $A$). By \ref{sec:5.1.5} applied to 
$1_{A'}$ and $1_{A}$ we have commutative diagrams
\begin{displaymath}
\xymatrix{
{\mathbf{P}'}\ar[r]\ar[d]& {S^{0}(A')}&&{{\rm and}}&&
 {\mathbf{P}}\ar[r]\ar[d]& {S^{0}(A)}\\
{\mathbf{C}'_{A'}}\ar[ur]&&&&&{\mathbf{C}_{A}}\ar[ur]
}
\end{displaymath}
in which the vertical arrows are quasi-isomorphisms. Applying
the bifunctor $T$ to the two diagrams we obtain the 
commutative diagram
\begin{displaymath}
\xymatrix{
{T(\mathbf{P}',\mathbf{P})}\ar[d]\ar[r]&
{T(\mathbf{P}',\mathbf{C}_{A})}\ar[d]\ar[r]&{T(\mathbf{P}',S^{0}(A))}\\
{T(\mathbf{C}'_{A'},\mathbf{P})}\ar[r] \ar[d]&
{T(\mathbf{C}'_{A'},\mathbf{C}_{A})} \ar[d]\ar[r]\ar[d]&{T(\mathbf{C}'_{A'},S^{0}(A))}\\
{T(S^{0}(A'),\mathbf{P})}&{T(S^{0}(A'),\mathbf{C}_{A})}\\
}
\end{displaymath}
By \ref{sec:6.1.6} all arrows of the previous diagram are quasi-isomorphisms. 
In particular, $\mathsf{Tor}_{n}(A',A)$ can be computed using flat resolutions.

\subsection{} \label{sec:7.2}
For an exact sequence 
$\mathbf{D}\overset{u}\rightarrowtail \mathbf{C}\overset{v}
\twoheadrightarrow \mathbf{E}$ 
of complexes in an abelian category we denote by 
\[
\partial(u,v):H(\mathbf{E})\to H(\mathbf{D})(-1)
\]
the connecting morphism associated to $(u,v)$. We construct here 
$\partial(u,v)$ when $\mathbf{D}$ and $\mathbf{E}$ have zero differential. 
We follow a method due to Lambek \cite[Section 5.5]{La}, which constructs 
the connecting morphism for arbitrary short exact sequences by giving its 
image factorization. In this way the naturality of the connecting morphism 
is immediate. We refer to \ref{sec:4.1} for notations.

\subsubsection{} \label{sec:7.2.1}
We have the commutative diagram
\begin{displaymath}
\xymatrix{
{B(\mathbf{C})}\ar@{>->}[r]^{i}& {Z(\mathbf{C})} \ar@{>->}[d]^{j}
\ar@{->>}[r]^{p}&{H(\mathbf{C})}\ar[d]^{\bar{j}}\\
{\mathbf{C}(1)} \ar[r]^{d^{\#}} \ar@{->>}[u]_{\delta(1)} & 
{\mathbf{C}}\ar@{->>}[r]^{\pi} & {\mathbf{C}/B(\mathbf{C})}
}
\end{displaymath}
where $\bar{j}$ is induced by $j$. The sequence $(u,v)$ gives 
rise to the commutative diagram
\begin{displaymath}
\xymatrix{
{B(\mathbf{C})}\ar@{>->}[r]^{i}& {Z(\mathbf{C})} \ar@{>->}[d]^{j} 
\ar[dr]^{Z(v)} \ar@{->>}[r]^{p}&
 H(\mathbf{C})\ar[d]^{H(v)}\\
{\mathbf{D}} \ar@{>->}[r]^{u} \ar[ur]^{Z(u)} \ar[dr]_{\bar{u}}&
{\mathbf{C}}\ar@{->>}[r]^{v} \ar@{->>}[d]_{\pi} & {\mathbf{E}}\\
& {\mathbf{C}/B(\mathbf{C})} \ar[ur]_{\bar{v}}
}
\end{displaymath}
where $\bar{u}$ and $\bar{v}$ are induced by $u$ and $v$; we have $H(u)=pZ(u)$.
The (commutative, solid arrows) diagram that gives rise to the connecting morphism 
$\partial(u,v)$ is then
\begin{displaymath}
\xymatrix{
&{\mathbf{D}}\ar@{=}[d]\ar[r]^{H(u)}&{H(\mathbf{C})}
\ar@{>->}[d]^{\bar{j}}\ar[r]^{H(v)}&
{\mathbf{E}}\ar@{=}[d]\ar@{..>>}[r]^{e}&{\Co H(v)}\\
&{\mathbf{D}}\ar[r]^{\pi u} \ar[d]_{0}&
{\mathbf{C}/B(\mathbf{C})} \ar[d]^{i(-1)\bar{\delta}}
\ar@{->>}[r]^{\bar{v}}\ar[d]&{\mathbf{E}}\ar[d]^{0}\\
&{\mathbf{D}(-1)}\ar@{=}[d]\ar[r]^{Z(u)(-1)}&
{Z(\mathbf{C})(-1)}\ar[r]^{Z(v)(-1)}\ar[d]&{\mathbf{E}(-1)}\ar@{=}[d]\\
{\K H(\mathbf{C})(-1)}\ar@{>..>}[r]^{m(-1)}&{\mathbf{D}(-1)}
\ar[r]^{H(u)(-1)}&{H(\mathbf{C})(-1)}\ar[r]^{H(v)(-1)}&{\mathbf{E}(-1)}
}
\end{displaymath}
cf. \cite[\S2 n\textsuperscript{\scriptsize{o}} 3 p. 30]{Bou}.
\paragraph{Claim 1} $\Co H(v)=B(\mathbf{C})(-1)$. We have 
$\Co H(v)=\Co H(v)p=\Co vj$. Consider the morphism 
$\delta:\mathbf{C}\to B(\mathbf{C})(-1)$; we have $\delta u=0$ so there 
is a unique $e:\mathbf{E}\to B(\mathbf{C})(-1)$ such that $ev=\delta$. 
From the commutative diagram
\begin{displaymath}
\xymatrix{
{Z(\mathbf{C})}\ar@{>->}[r]^{j}& {\mathbf{C}} 
\ar@{->>}[d]^{v}\ar@{->>}[r]^{\delta}&
{B(\mathbf{C})(-1)}\ar@{=}[d]\\
{Z(\mathbf{C})} \ar[r]^{vj} \ar@{=}[u]& 
{\mathbf{E}}\ar@{->>}[r]^{e} & {B(\mathbf{C})(-1)} 
}
\end{displaymath}
it follows that $\Co vj=B(\mathbf{C})(-1)$.
\paragraph{Claim 2} $\K H(u)=B(\mathbf{C})$. We have $\K H(u)=\K \bar{j}H(u)=\K \pi u$.
Since $vji=H(v)pi=0$ there is a unique $m:B(\mathbf{C})\to \mathbf{D}$ such that $um=ji$.
From the commutative diagram
\begin{displaymath}
\xymatrix{
{B(\mathbf{C})}\ar@{>->}[r]^{m}\ar@{=}[d]& {\mathbf{D}} 
\ar@{>->}[d]^{u}\ar[r]^{\pi u}&{\mathbf{C}/B(\mathbf{C})}\ar@{=}[d]\\
{B(\mathbf{C})} \ar@{>->}[r]^{ji}& {\mathbf{C}}\ar@{->>}[r]^{\pi} & 
{\mathbf{C}/B(\mathbf{C})}
}
\end{displaymath}
it follows that $\K \pi u=B(\mathbf{C})$. 

From the two claims we obtain that $\partial(u,v)=m(-1)e$.

\subsubsection{} \label{sec:7.2.2}
When the sequence $(u,v)$ is $(j,\delta)$ we have $H(j)=p$ 
and $H(\delta)=0$, so $e$ is the identity and $m=i$, 
therefore $\partial(j,\delta)=i(-1)$.

\subsection{} \label{sec:7.3}
Consider the standard THC--situation (\ref{sec:6.1}).
Let $(\mathbf{C}',d')$ (resp. $(\mathbf{C},d)$)
be a complex in $\mathcal{A}'$ (resp. $\mathcal{A}$).
Recall that we denote by $D(\mathbf{C}',\mathbf{C})$ 
the differential of $T(\mathbf{C}',\mathbf{C})$.
We denote by 
\[
\begin{aligned}
i(\mathbf{C}',\mathbf{C}):B(T(\mathbf{C}',\mathbf{C}))\to Z(T(\mathbf{C}',\mathbf{C}))\\
j(\mathbf{C}',\mathbf{C}):Z(T(\mathbf{C}',\mathbf{C}))\to T(\mathbf{C}',\mathbf{C})\\
p(\mathbf{C}',\mathbf{C}):Z(T(\mathbf{C}',\mathbf{C}))\to H(T(\mathbf{C}',\mathbf{C}))\\
\delta(\mathbf{C}',\mathbf{C}):T(\mathbf{C}',\mathbf{C})\to B(T(\mathbf{C}',\mathbf{C}))(-1)
\end{aligned}
\]
the natural morphisms.

\subsubsection{} \label{sec:7.3.1}
We have $D(\mathbf{C}',\mathbf{C})T(j',j)=0$.
\begin{proof}
Let $n$ be an integer and let $\sigma'_{p,q}:T(Z_{p}(\mathbf{C}'),Z_{q}(\mathbf{C}))
\to T(Z(\mathbf{C}'),Z(\mathbf{C}))_{n}$ be coproduct coprojections, where $p+q=n$.
We have 
\[
D(\mathbf{C}',\mathbf{C})_{n}T(j',j)_{n}\sigma'_{p,q}=
D(\mathbf{C}',\mathbf{C})_{n}\sigma_{p,q}T(j'_{p},j_{q})=
\]
\[
(\sigma_{p-1,q}T(d'_{p},\mathbf{C}_{q})+
\sigma_{p,q-1}(-1)^{p}T(\mathbf{C}'_{p},d_{q}))T(j'_{p},j_{q})=0
\]
\end{proof}

\subsubsection{} \label{sec:7.3.2}
From \ref{sec:7.3.1} it follows that there is a unique morphism 
\[
z(\mathbf{C}',\mathbf{C}):T(Z(\mathbf{C}'),Z(\mathbf{C}))\to Z(T(\mathbf{C}',\mathbf{C}))
\]
such that $j(\mathbf{C}',\mathbf{C})z(\mathbf{C}',\mathbf{C})=T(j',j)$.
One can check that this morphism is natural in both arguments, 
so that we obtain a natural transformation
\begin{displaymath}
z:T(Z(-),Z(-))\Rightarrow Z(T(-,-)):Ch(\mathcal{A}')\times Ch(\mathcal{A})\to Ch(\mathcal{A}'')
\end{displaymath}

\subsubsection{} \label{sec:7.3.3}
If $f:X\to A,g:Y\to A$ are morphisms in an abelian category, we denote by $(f,g)$
the induced morphism $X\oplus Y\to A$. One has $\im (f,g)\cong \im f \cup \im g$.

\subsubsection{} \label{sec:7.3.4}
We claim that the diagram
\begin{displaymath}
\xymatrix{
{T(\mathbf{C}'(1),Z(\mathbf{C}))\oplus T(Z(\mathbf{C}'),\mathbf{C}(1))}
\ar[rrrrr]^{(T(i'\delta'(1),Z(\mathbf{C})),T(Z(\mathbf{C}'),i\delta(1)))} 
\ar[d]^{(T(\mathbf{C}'(1),j),T(j',\mathbf{C})(1)\Sigma(Z(\mathbf{C}'),\mathbf{C}))}&&&&&
{T(Z(\mathbf{C}'),Z(\mathbf{C}))} \ar[d]^{T(j',j)}\\
{T(\mathbf{C}',\mathbf{C})(1)} \ar[rrrrr]_{D(\mathbf{C}',\mathbf{C})^{\#}}&&&&&
{T(\mathbf{C}',\mathbf{C})}
}
\end{displaymath}
where $\Sigma(Z(\mathbf{C}'),\mathbf{C})$ has been defined in \ref{sec:6.1.1}, 
is commutative.
\paragraph{Step 1} We show that $D(\mathbf{C}',\mathbf{C})^{\#}
T(\mathbf{C}'(1),j)=T(j',j)T(i'\delta'(1),Z(\mathbf{C}))$. It suffices to show 
that $D(\mathbf{C}',\mathbf{C})T(\mathbf{C}',j)=T(d',j)$.
Since $T(\mathbf{C}',j)$ is a morphism we have (\ref{sec:4.1})
$D(\mathbf{C}',\mathbf{C})T(\mathbf{C}',j)=
T(\mathbf{C}',j)(-1)D(\mathbf{C}',Z(\mathbf{C}))$.
Since $Z(\mathbf{C})$ has zero differential we have 
$D(\mathbf{C}',Z(\mathbf{C}))=T(d',Z(\mathbf{C}))$
and the equality follows. 
\paragraph{Step 2} We show that $D(\mathbf{C}',\mathbf{C})^{\#}
T(j',\mathbf{C})(1)\Sigma(Z(\mathbf{C}'),\mathbf{C})=
T(j',j)T(Z(\mathbf{C}'),i\delta(1))$. Since $T(j,\mathbf{C})$ is a morphism 
we have $D(\mathbf{C}',\mathbf{C})^{\#}T(j',\mathbf{C})(1)=
T(j',\mathbf{C})D(Z(\mathbf{C}'),\mathbf{C})^{\#}$. Since $Z(\mathbf{C}')$ 
has zero differential we have $D(Z(\mathbf{C}'),\mathbf{C})^{\#}
\Sigma(Z(\mathbf{C}'),\mathbf{C})=T(Z(\mathbf{C}'),d^{\#})$
and the equality follows. Combining now the two steps we obtain the 
commutativity of the diagram.

\subsubsection{} \label{sec:7.3.5}
Let $u':A'\to B'$ (resp. $u:A\to B$) be an epimorphism
in $\mathcal{A}'$ (resp. $\mathcal{A}$). Then the diagram
\begin{displaymath}
\xymatrix{
{T(A',A)}\ar@{->>}[r]^{T(A',u)} \ar@{->>}[d]_{T(u',A)}
&{T(A',B)} \ar@{->>}[d]^{T(u',B)}\\
{T(B',A)} \ar@{->>}[r]_{T(B',u)}
&{T(B',B)}
}
\end{displaymath}
is a pushout. Consequently, if $i':\K u'\to A'$ and $i:\K u\to A$
are the natural morphisms, there is an exact sequence
\begin{equation}\label{eq:7.3.5}
\xymatrix{
{\im T(i',A)\cup \im T(A',i)}\ar@{>->}[rr]&&{T(A',A)}
\ar@{->>}[rr]^{T(u',u)}&&{T(B',B)}
}
\end{equation}
\begin{proof}
Let $f:T(B',A)\to A''$ and $g:T(A',B)\to A''$ be 
such that $fT(u',A)=gT(A',u)$. By adjunction we 
have a commutative solid arrows diagram
\begin{displaymath}
\xymatrix{
{A'}\ar[rr]^{g^{\flat}} \ar@{->>}[d]_{u'}
&&{H(B,A'')} \ar@{>->}[d]^{H(u,A'')}\\
{B'} \ar@{..>}[urr]^{h^{\flat}}\ar[rr]_{f^{\flat}}
&&{H(A,A'')}
}
\end{displaymath}
where $f^{\flat}$ and $g^{\flat}$ are the adjoints of $f$ and $g$, respectively.
The map $H(u,A'')$ is a monomorphism since $H(-,A'')$ is a right adjoint.
Since in an abelian category we have the factorization system 
(epimorphisms, monomorphisms) we have a unique dotted diagonal filler $h^{\flat}$ 
in the previous diagram. The adjoint of this diagonal filler is a morphism 
$h:T(B',B)\to A''$ having the property that $f=hT(B',u)$ and $g=hT(u',B)$. 
Since $h^{\flat}$ is unique, so is $h$ with this property. Therefore 
the required diagram is a pushout. For the second part, consider the diagram
\begin{displaymath}
\xymatrix{
&&{T(\K u',A)}\ar[d]^{T(i',A)}\\
{T(A',\K u)}\ar[rr]^{T(A',i)}&&{T(A',A)}
\ar@{->>}[rr]^{T(A',u)} \ar@{->>}[d]_{T(u',A)}
&&{T(A',B)} \ar@{->>}[d]^{T(u',B)}\\
& &{T(B',A)} \ar@{->>}[rr]_{T(B',u)}
&&{T(B',B)}
}
\end{displaymath}
We have 
\[ T(A',B)=T(A',A)/\im T(A',i) \ {\rm and} \ 
T(B',A)=T(A',A)/\im T(i',A)
\]
The existence of the exact sequence (\ref{eq:7.3.5}) is then a 
consequence of the first part and of \cite[Corollary 1.7.5]{Bo2}.
\end{proof}
The exact sequence \ref{sec:7.3.5}(\ref{eq:7.3.5}) is a familiar one in 
the theory of modules over a ring where it says that, if $A$ is a ring, 
$M'\overset{u}\rightarrow M\overset{v}\twoheadrightarrow M''$
an exact sequence of right $A$-modules and 
$N'\overset{s}\rightarrow N\overset{t}\twoheadrightarrow N''$
an exact sequence of left $A$-modules,
then $v\otimes t$ is surjective and its kernel is equal to 
$\im (u\otimes 1_{N})+\im (1_{M}\otimes s)$.

\subsubsection{} \label{sec:7.3.6}
We apply the naturality of the image factorization to the 
commutative diagram \ref{sec:7.3.4}; from \ref{sec:7.3.3} 
we obtain a unique morphism
\[
b(\mathbf{C}',\mathbf{C}):\im T(i',Z(\mathbf{C}))\cup 
\im T(Z(\mathbf{C}'),i)\to B(T(\mathbf{C}',\mathbf{C}))
\]
such that $j(\mathbf{C}',\mathbf{C})i(\mathbf{C}',\mathbf{C})b(\mathbf{C}',\mathbf{C})
=T(j',j)k(\mathbf{C}',\mathbf{C})$, where 
\[
k(\mathbf{C}',\mathbf{C}):\im T(i',Z(\mathbf{C}))\cup \im T(Z(\mathbf{C}'),i)
\rightarrowtail T(Z(\mathbf{C}'),Z(\mathbf{C}))
\]
is the monomorphism of the image factorization of the top arrow of diagram \ref{sec:7.3.4}.
From \ref{sec:7.3.5} the top row of the commutative solid arrow diagram
\begin{displaymath}
\xymatrix{
{\im T(i',Z(\mathbf{C}))\cup \im T(Z(\mathbf{C}'),i)}\ar@{>->}[rr]^{k(\mathbf{C}',\mathbf{C})}
\ar[d]_{b(\mathbf{C}',\mathbf{C})}&&{T(Z(\mathbf{C}'),Z(\mathbf{C}))}
\ar[d]_{z(\mathbf{C}',\mathbf{C})}\ar@{->>}[rr]^{T(p',p)}&&
{T(H(\mathbf{C}'),H(\mathbf{C}))}\ar@{..>}[d]_{\gamma(\mathbf{C}',\mathbf{C})}\\
{B(T(\mathbf{C}',\mathbf{C}))}\ar@{>->}[rr]^{i(\mathbf{C}',\mathbf{C})}&&
{Z(T(\mathbf{C}',\mathbf{C}))}\ar@{->>}[rr]^{p(\mathbf{C}',\mathbf{C})}&&{H(T(\mathbf{C}',\mathbf{C}))}
}
\end{displaymath}
is exact, therefore there is a unique morphism 
\begin{equation}\label{eq:7.3.6}
\gamma(\mathbf{C}',\mathbf{C}):T(H(\mathbf{C}'),H(\mathbf{C}))\to H(T(\mathbf{C}',\mathbf{C}))
\end{equation}
such that $\gamma(\mathbf{C}',\mathbf{C})T(p',p)=p(\mathbf{C}',\mathbf{C})z(\mathbf{C}',\mathbf{C})$.
This morphism is natural in both arguments. For example, let $u:\mathbf{D}\to\mathbf{C}$ be a morphism.
Form the cube diagram
\[
\xymatrix{
{T(Z(\mathbf{C}'),Z(\mathbf{D}))} \ar[rrr]^{T(p',p_{\mathbf{D}})} \ar[dr]^{z(\mathbf{C}',\mathbf{D})} 
\ar[dd]|-{T(Z(\mathbf{C}'),Z(u)))} &&& {T(H(\mathbf{C}'),H(\mathbf{D}))}\ar[dr]^{\gamma(\mathbf{C}',\mathbf{D})} 
\ar[dd]|-{T(H(\mathbf{C}'),H(u)))}\\
& {Z(T(\mathbf{C}',\mathbf{D}))} \ar[rrr]_{p(\mathbf{C}',\mathbf{D})} \ar[dd] &&& {H(T(\mathbf{C}',\mathbf{D}))} 
\ar[dd]|-{H(T(\mathbf{C}',u))}\\
{T(Z(\mathbf{C}'),Z(\mathbf{C}))} \ar[rrr]^{T(p',p_{\mathbf{C}})} \ar[dr]^{z(\mathbf{C}',\mathbf{C})} &&&
{T(H(\mathbf{C}'),H(\mathbf{C}))} \ar[dr]^{\gamma(\mathbf{C}',\mathbf{C})}\\
& {Z(T(\mathbf{C}',\mathbf{C}))} \ar[rrr]_{p(\mathbf{C}',\mathbf{C})} &&& {H(T(\mathbf{C}',\mathbf{C}))}\\
 }
\]
The left face of the cube commutes by the naturality of $z$ (\ref{sec:7.3.2}), the top and 
bottom faces commute by the construction of $\gamma(\mathbf{C}',\mathbf{D})$ and 
$\gamma(\mathbf{C}',\mathbf{C})$, the back face commutes because $T$ is a bifunctor. 
Since $T(p',p_{\mathbf{D}})$ is an epimorphism, it follows that the right face commutes. 
Thus, we obtain a natural transformation
\begin{displaymath}
\gamma:T(H(-),H(-))\Rightarrow H(T(-,-)):Ch(\mathcal{A}')\times Ch(\mathcal{A})\to Ch(\mathcal{A}'')
\end{displaymath}
For every integer $r$ the diagram 
\begin{displaymath}
\xymatrix{
{T(H(\mathbf{C}'),H(\mathbf{C}(r)))}\ar[rr]^{\gamma(\mathbf{C}',\mathbf{C}(r))}
\ar[d]_{\Sigma(H(\mathbf{C}'),H(\mathbf{C}))}
&&{H(T(\mathbf{C}',\mathbf{C}(r)))} \ar[d]^{H(\Sigma(\mathbf{C}',\mathbf{C}))}\\
{T(H(\mathbf{C}'),H(\mathbf{C}))(r)}\ar[rr]^{\gamma(\mathbf{C}',\mathbf{C})(r)}
&&{H(T(\mathbf{C}',\mathbf{C}))(r)} 
}
\end{displaymath}
commutes, where $\Sigma(\mathbf{C}',\mathbf{C})$ has been defined in \ref{sec:6.1.1}.

\subsubsection{} \label{sec:7.3.7}
Suppose that $\mathcal{A}''$ satisfies AB4. The morphism 
$\gamma(\mathbf{C}',\mathbf{C})$ (\ref{sec:7.3.6}(\ref{eq:7.3.6})) is 
an isomorphism if $\mathbf{C}$ is degreewise flat with zero differential.
\begin{proof}
From the proof of \ref{sec:6.1.5} we have that the morphisms $z(\mathbf{C}',\mathbf{C})$
(\ref{sec:7.3.2}) and $b(\mathbf{C}',\mathbf{C})$ (\ref{sec:7.3.6}) are isomorphisms, 
hence so is $\gamma(\mathbf{C}',\mathbf{C})$ from the way it was constructed.
\end{proof}

\subsubsection{} \label{sec:7.3.8}
Suppose that $\mathcal{A}''$ satisfies AB4, $\mathbf{C}'$ has zero differential 
and there is an exact sequence $\mathbf{D}\overset{u}\rightarrowtail 
\mathbf{C}\overset{v}\twoheadrightarrow \mathbf{E}$,
where $\mathbf{D},\mathbf{E}$ have zero differential and 
$\mathbf{E}$ is degreewise flat. Then 
\[
\partial(T(\mathbf{C}',u),T(\mathbf{C}',v))=
\Sigma(\mathbf{C}',\mathbf{D})T(\mathbf{C}',\partial(u,v))
\]
\begin{proof}
By \ref{sec:3.1.6}(1) we have the exact sequence 
\begin{displaymath}
\xymatrix{
{T(\mathbf{C}',\mathbf{D})}\ar@{>->}[rr]^{T(\mathbf{C}',u)}&&
{T(\mathbf{C}',\mathbf{C})}\ar@{->>}[rr]^{T(\mathbf{C}',v)}&&{T(\mathbf{C}',\mathbf{E})}
}
\end{displaymath}
where the first and last complexes have zero differential. 
Using the naturality of the morphism $\Sigma$ (\ref{sec:6.1.1}) 
and \ref{sec:7.2.1} we have 
\[
T(\mathbf{C}',u)(-1)\Sigma(\mathbf{C}',\mathbf{D})T(\mathbf{C}',\partial(u,v))T(\mathbf{C}',v)=
\]
\[
\Sigma(\mathbf{C}',\mathbf{C})T(\mathbf{C}',u(-1))T(\mathbf{C}',\partial(u,v))T(\mathbf{C}',v)=
\]
\[
\Sigma(\mathbf{C}',\mathbf{C})T(\mathbf{C}',u(-1))T(\mathbf{C}',m(-1)e)T(\mathbf{C}',v)=
\]
\[
\Sigma(\mathbf{C}',\mathbf{C})T(\mathbf{C}',u(-1)m(-1)ev)=\Sigma(\mathbf{C}',\mathbf{C})T(\mathbf{C}',d)
\]
and 
\[
T(\mathbf{C}',u)(-1)\partial(T(\mathbf{C}',u),T(\mathbf{C}',v))T(\mathbf{C}',v)=
\]
\[
j(\mathbf{C}',\mathbf{C})(-1)i(\mathbf{C}',\mathbf{C})(-1)\delta(\mathbf{C}',\mathbf{C})
=D(\mathbf{C}',\mathbf{C})
\]
By the diagram in \ref{sec:6.1.1} and the fact that $T(\mathbf{C}',u)$ is a
monomorphism and $T(\mathbf{C}',v)$ is an epimorphism we obtain the desired equality.
\end{proof}

\subsubsection{} \label{sec:7.3.9}
Suppose that $\mathcal{A}''$ satisfies AB4 and that there is an exact 
sequence $\mathbf{D}\overset{u}\rightarrowtail 
\mathbf{C}\overset{v}\twoheadrightarrow \mathbf{E}$,
where $\mathbf{D},\mathbf{E}$ have zero differential 
and $\mathbf{E}$ is degreewise flat. Then 
\[
\partial(T(\mathbf{C}',u),T(\mathbf{C}',v))=
H(\Sigma(\mathbf{C}',\mathbf{D}))H(T(\mathbf{C}',\partial(u,v)))
\]
\begin{proof}
By \ref{sec:3.1.6}(1) we have the commutative diagram with exact rows
\begin{displaymath}
\xymatrix{
{T(Z(\mathbf{C}'),\mathbf{D})}\ar@{>->}[rr]^{T(Z(\mathbf{C}'),u)}\ar[d]_{T(j',\mathbf{D})}&&
{T(Z(\mathbf{C}'),\mathbf{C})}\ar[d]_{T(j',\mathbf{C})}\ar@{->>}[rr]^{T(Z(\mathbf{C}'),v)}
&&{T(Z(\mathbf{C}'),\mathbf{E})}\ar[d]_{T(j',\mathbf{E})}\\
{T(\mathbf{C}',\mathbf{D})}\ar@{>->}[rr]^{T(\mathbf{C}',u)}&&
{T(\mathbf{C}',\mathbf{C})}\ar@{->>}[rr]^{T(\mathbf{C}',v)}&&{T(\mathbf{C}',\mathbf{E})}
}
\end{displaymath}
By the naturality of the connecting morphism we have then the commutative diagram
\begin{displaymath}
\xymatrix{
{H(T(Z(\mathbf{C}'),\mathbf{E}))}\ar[rrrr]^{\partial(T(Z(\mathbf{C}'),u),T(Z(\mathbf{C}'),v))} 
\ar[d]_{H(T(j',\mathbf{E}))}&&&&
{H(T(Z(\mathbf{C}'),\mathbf{D}))(-1)} \ar[d]^{H(T(j',\mathbf{D}))(-1)}\\
{H(T(\mathbf{C}',\mathbf{E}))} \ar[rrrr]^{\partial(T(\mathbf{C}',u),T(\mathbf{C}',v))}&&&&
{H(T(\mathbf{C}',\mathbf{D}))(-1)}
}
\end{displaymath}
By the naturality of $\gamma$ we have the commutative diagram
\begin{displaymath}
\xymatrix{
{T(H(Z(\mathbf{C}'))),H(\mathbf{E}))}\ar[rrr]^{\gamma(Z(\mathbf{C}'),\mathbf{E})} 
\ar[d]_{T(H(j'),H(\mathbf{E}))}&&&
{H(T(Z(\mathbf{C}'),\mathbf{E}))} \ar[d]^{H(T(j',\mathbf{E}))}\\
{T(H(\mathbf{C}'),H(\mathbf{E}))} \ar[rrr]^{\gamma(\mathbf{C}',\mathbf{E})}&&&
{H(T(\mathbf{C}',\mathbf{E}))}
}
\end{displaymath}
in which the horizontal arrows are isomorphisms by \ref{sec:7.3.7}. Since $H(j')=p'$ it 
follows that $H(T(j',\mathbf{E}))$ is an epimorphism. Consider now the diagram
\[
\xymatrix{
&{H(T(Z(\mathbf{C}'),\mathbf{D}(-1))}\ar[dr]|-{H(\Sigma(Z(\mathbf{C}'),\mathbf{D}))}\ar[ddr]&\\
{H(T(Z(\mathbf{C}'),\mathbf{E}))}\ar[rr]|-{\partial(T(Z(\mathbf{C}'),u),T(Z(\mathbf{C}'),v))}
\ar@{->>}[ddr]^{H(T(j',\mathbf{E}))} \ar[ur]|-{H(T(Z(\mathbf{C}'),\partial(u,v)))}
&& {H(T(Z(\mathbf{C}'),\mathbf{D})(-1))}\ar[ddr]^{H(T(j',\mathbf{D})(-1))}\\
&&{H(T(\mathbf{C}',\mathbf{D}(-1))}\ar[dr]|-{H(\Sigma(\mathbf{C}',\mathbf{D}))}&\\
& {H(T(\mathbf{C}',\mathbf{E}))} \ar[rr]_{\partial(T(\mathbf{C}',u),T(\mathbf{C}',v))}
\ar[ur]|-{H(T(\mathbf{C}',\partial(u,v)))}&& {H(T(\mathbf{C}',\mathbf{D})(-1))}\\
}
\]
where the unlabeled arrow is $H(T(j',\mathbf{D}(-1)))$. The back face commutes 
by \ref{sec:7.3.8}, the top right face commutes by the naturality of $\Sigma$,
the top left face commutes since $T$ is a bifunctor. It follows that the front face 
commutes since $H(T(j',\mathbf{E}))$ is an epimorphism.
\end{proof}

\subsubsection{} \label{sec:7.3.10}
Suppose that $\mathcal{A}''$ satisfies AB4 and $B(\mathbf{C}), Z(\mathbf{C})$ 
are degreewise flat. Then we have an exact sequence
\begin{displaymath}
\xymatrix{
{T(H(\mathbf{C}'),H(\mathbf{C}))}\ar@{>->}[rr]^{\gamma(\mathbf{C}',\mathbf{C})}&&
{H(T(\mathbf{C}',\mathbf{C}))}\ar@{->>}[r]&{\K T(H(\mathbf{C}'),i(-1))}
}
\end{displaymath}

\begin{proof}
Using the naturality of $\gamma$ and the morphisms $i,j$ (\ref{sec:4.1}) 
we obtain the commutative diagram
\begin{displaymath}
\xymatrix{
{T(H(\mathbf{C}'),H(B(\mathbf{C})))}\ar[rr]^{T(H(\mathbf{C}'),H(i))}
\ar[d]_{\gamma(\mathbf{C}',B(\mathbf{C}))}&&
{T(H(\mathbf{C}'),H(Z(\mathbf{C})))}\ar[d]_{\gamma(\mathbf{C}',Z(\mathbf{C}))}
\ar@{->>}[rr]^{T(H(\mathbf{C}'),H(j))}
&&{T(H(\mathbf{C}'),H(\mathbf{C}))}\ar[d]_{\gamma(\mathbf{C}',\mathbf{C})}\\
{H(T(\mathbf{C}',B(\mathbf{C}))}\ar[rr]^{H(T(\mathbf{C}',i))}&&
{H(T(\mathbf{C}',Z(\mathbf{C}))}\ar[rr]^{H(T(\mathbf{C}',j))}&&{H(T(\mathbf{C}',\mathbf{C}))}
}
\end{displaymath}
The first two vertical arrows are isomorphisms by \ref{sec:7.3.7} and 
the top sequence is exact since $H(i)=i$ and $H(j)=p$. We apply 
\ref{sec:7.3.9} to the exact sequence $(j,\delta)$; 
then \ref{sec:7.2.2} gives
\begin{equation}\label{eq:7.3.10}
\partial(T(\mathbf{C}',j),T(\mathbf{C}',\delta))=
H(\Sigma(\mathbf{C}',Z(\mathbf{C})))H(T(\mathbf{C}',i(-1)))
\end{equation}
This implies that bottom sequence is exact, therefore $\gamma(\mathbf{C}',\mathbf{C})$ 
factors as an isomorphism followed by the monomorphism 
\[
\im H(T(\mathbf{C}',j))\to H(T(\mathbf{C}',\mathbf{C}))
\]
Consider now the commutative diagram
\begin{displaymath}
\xymatrix{
{H(T(\mathbf{C}',\mathbf{C}))}\ar[rr]^{H(T(\mathbf{C}',\delta))}&&
{H(T(\mathbf{C}'),B(\mathbf{C})(-1))}\ar[rr]^{H(T(\mathbf{C}',i(-1)))}
&&{H(T(\mathbf{C}',Z(\mathbf{C})(-1))}\\
{\K T(H(\mathbf{C}'),H(i(-1)))}\ar@{>->}[rr]&&
{T(H(\mathbf{C}'),H(B(\mathbf{C})(-1)))}\ar[u]^{\gamma(\mathbf{C}',B(\mathbf{C})(-1))}
\ar[rr]^{T(H(\mathbf{C}'),H(i(-1)))}&&{T(H(\mathbf{C}'),H(Z(\mathbf{C})(-1)))}
\ar[u]^{\gamma(\mathbf{C}',Z(\mathbf{C})(-1))}
}
\end{displaymath}
The two vertical arrows are isomorphisms by \ref{sec:7.3.7} and the top 
sequence is exact by (\ref{eq:7.3.10}), hence
\[
\K T(H(\mathbf{C}'),H(i(-1)))\cong \K H(T(\mathbf{C}',i(-1)))
\]
Since the sequence $(T(\mathbf{C}',j),T(\mathbf{C}',\delta))$ is exact we obtain 
the desired exact sequence.
\end{proof}

\subsubsection{} \label{sec:7.3.11} \cite[\S4 n\textsuperscript{\scriptsize{o}} 7 thm. 3]{Bou}
Suppose that $\mathcal{A}''$ satisfies AB4 and that  $(\mathcal{F}',\mathcal{F}'^{\perp})$ 
and $(\mathcal{F},\mathcal{F}^{\perp})$ are left complete. Suppose 
moreover that $B(\mathbf{C}), Z(\mathbf{C})$ are degreewise flat.
Then we have an exact sequence
\begin{displaymath}
\xymatrix{
{T(H(\mathbf{C}'),H(\mathbf{C}))}\ar@{>->}[rr]^{\gamma(\mathbf{C}',\mathbf{C})}&&
{H(T(\mathbf{C}',\mathbf{C}))}\ar@{->>}[r]&{\mathsf{Tor}_{1}(H(\mathbf{C}'),H(\mathbf{C}))}
}
\end{displaymath}
where $\mathsf{Tor}_{1}(H(\mathbf{C}'),H(\mathbf{C}))_{n}=\underset{p+q=n}\bigoplus 
\mathsf{Tor}_{1}(H_{p}(\mathbf{C}'),H_{q-1}(\mathbf{C}))$.
\begin{proof}
By \ref{sec:7.3.10} it remains to calculate $\K T(H(\mathbf{C}'),i(-1))$.
We have
\[
(\K T(H(\mathbf{C}'),i(-1)))_{n}=\underset{p+q=n}\bigoplus \K T(H_{p}(\mathbf{C}'),i_{q-1})
\]
Let $p,q$ be integers. Let $\mathbf{P}$ be the complex with
$\mathbf{P}_{1}=B_{q-1}(\mathbf{C}), \mathbf{P}_{0}=Z_{q-1}(\mathbf{C})$
and $\mathbf{P}_{n}=0$ otherwise, so that $\mathbf{P}$ is a flat resolution
of $H_{q-1}(\mathbf{C})$. We have
\[
\K T(H_{p}(\mathbf{C}'),i_{q-1})=H_{1}(T(S^{0}(H_{p}(\mathbf{C}')),\mathbf{P}))\cong
\mathsf{Tor}_{1}(H_{p}(\mathbf{C}'),H_{q-1}(\mathbf{C}))
\]
where the isomorphism exists by the last diagram of \ref{sec:7.1.1}. 
\end{proof}

\subsubsection{} \label{sec:7.3.12}
\cite[\S4 n\textsuperscript{\scriptsize{o}} 7 cor. 2]{Bou} Suppose that 
$B(\mathbf{C}'),B(\mathbf{C})$ are degreewise projective and 
$Z(\mathbf{C})$ is degreewise flat. Then for all integers $n$ the morphism 
\[
\gamma_{n}(\mathbf{C}',\mathbf{C}):T(H(\mathbf{C}'),H(\mathbf{C}))_{n}\to
H_{n}(T(\mathbf{C}',\mathbf{C}))
\] 
has a retraction.

\begin{proof}
The proof is the same as in \emph{loc. cit.}, with the observation that remarque b) 
mentioned in the proof of lemme 4 is on p. 35.
\end{proof}

\subsubsection{} \label{sec:7.3.13}
\cite[\S4 n\textsuperscript{\scriptsize{o}} 7 cor. 3]{Bou} 
Suppose that in $\mathcal{A}'$ and $\mathcal{A}$ 
the class of projective objects is closed under subobjects. If 
$\mathbf{C}'$ and $\mathbf{C}$ are degreewise projective then 
for all integers $n$ the morphism $\gamma_{n}(\mathbf{C}',\mathbf{C})$
has a retraction.

\begin{proof}
It follows from the assumptions that $B(\mathbf{C}'),B(\mathbf{C})$ and
$Z(\mathbf{C})$ are degreewise projective, hence we can apply \ref{sec:7.3.12}.
\end{proof}

\subsubsection{} \label{sec:7.3.14}
\cite[\S4 n\textsuperscript{\scriptsize{o}} 7 cor. 4]{Bou} 
Suppose that $\mathbf{C},H(\mathbf{C})$ are degreewise flat and 
$\mathbf{C}_{n}=0$ for $n<0$; then $\gamma(\mathbf{C}',\mathbf{C})$ 
is an isomorphism.
\begin{proof}
Consider the exact sequences 
\[
B_{n}(\mathbf{C})\rightarrowtail Z_{n}(\mathbf{C})\twoheadrightarrow H_{n}(\mathbf{C})
\]
and
\[
Z_{n}(\mathbf{C})\rightarrowtail \mathbf{C}_{n}\twoheadrightarrow B_{n-1}(\mathbf{C})
\]
Since $Z_{0}(\mathbf{C})=\mathbf{C}_{0}$ is flat we have, using \ref{sec:3.1.5},
the implications ($B_{n-1}(\mathbf{C})$ flat)$\Rightarrow$ ($Z_{n}(\mathbf{C})$ flat)
$\Rightarrow$ ($B_{n}(\mathbf{C})$ flat); we obtain that 
$B(\mathbf{C})$ and $Z(\mathbf{C})$ are degreewise flat. 
We apply $T(H(\mathbf{C}'),-)$ to the exact sequence 
$B(\mathbf{C})(-1)\overset{i(-1)}\rightarrowtail Z(\mathbf{C})(-1)
\twoheadrightarrow H(\mathbf{C})(-1)$; we obtain from \ref{sec:3.1.6}(1)
the exact sequence 
\begin{displaymath}
\xymatrix{
{T(H(\mathbf{C}'),B(\mathbf{C})(-1))}\ar@{>->}[rr]^{T(H(\mathbf{C}'),i(-1))}&& 
{T(H(\mathbf{C}'),Z(\mathbf{C})(-1))} \ar@{->>}[r]&{T(H(\mathbf{C}'),H(\mathbf{C})(-1))}
}
\end{displaymath}
From \ref{sec:7.3.10} we conclude that $\gamma(\mathbf{C}',\mathbf{C})$ 
is an isomorphism.
\end{proof}

\subsubsection{} \label{sec:7.3.15}
\cite[\S4 n\textsuperscript{\scriptsize{o}} 7 cor. 5]{Bou} The functor 
$T(\mathbf{C}',-)$ preserves quasi-isomorphisms between complexes 
that are degreewise flat and zero in negative degrees. 

\begin{proof}
Let $u:\mathbf{C}\to\mathbf{D}$ be a quasi-isomorphism between complexes 
that are degreewise flat and zero in negative degrees. Then $Con(u)$ is exact,
degreewise flat and zero in negative degrees, hence by \ref{sec:7.3.14} 
$\gamma(\mathbf{C}',Con(u))$ is an isomorphism. This implies that 
$T(\mathbf{C}',Con(u))$ is exact, therefore so is $Con(T(\mathbf{C}',u))$ by 
\ref{sec:6.1.2}.
\end{proof}
.%\subsection{} \label{sec:7.3}
%$\mathbf{C}'$ degreewise injective, $Z(\mathbf{C}), B(\mathbf{C})$ degreewise in 
%$\mathcal{I}-\mathcal{F}$ 2.2.7 need something like 2.1.6(1)
%\subsubsection{} \label{sec:7.3.1}

\section{Orthogonality relations for complexes}

\subsection{} \label{sec:8.1}
Let $\mathcal{A}$ be an abelian category, $\mathcal{K}$ a class
of objects of $\mathcal{A}$ and $\mathfrak{K}$ a class of objects of 
$Ch(\mathcal{A})$. We denote by $\mathcal{P}$ (resp. $\mathcal{I}$) 
the class of projective (resp. injective) objects of $\mathcal{A}$. We denote 
by $\mathfrak{K}^{\perp_{e}}$ (resp. $^{\perp_{e}}\mathfrak{K}$)
the class of complexes $\mathbf{C}$ for which 
\[
\mathrm{Homgr}_{\mathcal{A}}(-,\mathbf{C}):
Ch(\mathcal{A})^{op}\to Ch(\mathrm{Ab})
\]
(resp.
\[
\mathrm{Homgr}_{\mathcal{A}}(\mathbf{C},-):
Ch(\mathcal{A})\to Ch(\mathrm{Ab}))
\]
preserves epimorphisms with kernel in $\mathfrak{K}$. We have
$Ob(\mathcal{A})^{\perp_{e}}=Ch(\mathcal{I})$, 
$^{\perp_{e}}Ob(\mathcal{A})=Ch(\mathcal{P})$ and 
$^{\perp_{e}}(Ob(\mathcal{A})^{\perp_{e}})=Ob(\mathcal{A})=
(^{\perp_{e}}Ob(\mathcal{A}))^{\perp_{e}}$. From \ref{sec:4.2}
the classes $\mathfrak{K}^{\perp_{e}}$ and $^{\perp_{e}}\mathfrak{K}$
are closed under translates.

\subsubsection{} \label{sec:8.1.1}
Let $\mathbf{C}$ be a complex in $\mathcal{A}$.

(1) The following are equivalent:

(a) $\mathbf{C}\in Ch(\mathcal{K}^{\perp})$;

(b) $\mathbf{C}\in Ch(\mathcal{K})^{\perp_{e}}$; 

(c) for every exact sequence of complexes
$\mathbf{A}'\overset{u}\rightarrow \mathbf{A}\overset{v}\rightarrow \mathbf{A}''$
with $\mathbf{A}''/\im v\in Ch(\mathcal{K})$ the sequence 
\[
\mathrm{Homgr}_{\mathcal{A}}(\mathbf{A}'',\mathbf{C})\overset{\bar{v}}\to
\mathrm{Homgr}_{\mathcal{A}}(\mathbf{A},\mathbf{C}) \overset{\bar{u}}\to
\mathrm{Homgr}_{\mathcal{A}}(\mathbf{A}',\mathbf{C})
\] 
is exact, where $\bar{u}=\mathrm{Homgr}_{\mathcal{A}}(u,\mathbf{C})$ and 
$\bar{v}=\mathrm{Homgr}_{\mathcal{A}}(v,\mathbf{C})$.

(1bis) The following are equivalent:

(a) $\mathbf{C}\in Ch(^{\perp}\mathcal{K})$;

(b) $\mathbf{C}\in^{\perp_{e}}Ch(\mathcal{K})$; 

(c) for every exact sequence of complexes
$\mathbf{A}'\overset{u}\rightarrow \mathbf{A}\overset{v}\rightarrow \mathbf{A}''$
with $\K u\in Ch(\mathcal{K})$ the sequence 
\[
\mathrm{Homgr}_{\mathcal{A}}(\mathbf{C},\mathbf{A}')\overset{\bar{u}}\to
\mathrm{Homgr}_{\mathcal{A}}(\mathbf{C},\mathbf{A}) \overset{\bar{v}}\to
\mathrm{Homgr}_{\mathcal{A}}(\mathbf{C},\mathbf{A}'')
\] 
is exact, where $\bar{u}=\mathrm{Homgr}_{\mathcal{A}}(\mathbf{C},u)$ and 
$\bar{v}=\mathrm{Homgr}_{\mathcal{A}}(\mathbf{C},v)$.

\begin{proof}
(1) (a)$\Rightarrow$(b) Let $\mathbf{A}\rightarrowtail\mathbf{B}
\twoheadrightarrow\mathbf{X}$ be an exact sequence with $\mathbf{X}$ 
degreewise in $\mathcal{K}$. For all integers $n,p$  we have then the exact sequence
\begin{displaymath}
\mathcal{A}(\mathbf{X}_{p},\mathbf{C}_{p+n})\rightarrowtail 
\mathcal{A}(\mathbf{B}_{p},\mathbf{C}_{p+n})\twoheadrightarrow 
\mathcal{A}(\mathbf{A}_{p},\mathbf{C}_{p+n})
\end{displaymath}
The claim follows from the construction of $\mathrm{Homgr}_{\mathcal{A}}$ (\ref{sec:4.2}).

(b)$\Rightarrow$(c) We factor $v$ as $\mathbf{A}\overset{q}\twoheadrightarrow \im v
\overset{j}\rightarrowtail \mathbf{A}''$. By assumption 
\[
\bar{j}:\mathrm{Homgr}_{\mathcal{A}}(\mathbf{A}'',\mathbf{C})\to
\mathrm{Homgr}_{\mathcal{A}}(\im v,\mathbf{C})
\] 
is an epimorphism. It then suffices to show that the sequence 
$(\bar{q},\bar{u})$ is exact. Since exactness can be checked degreewise, 
the argument is the same as for \ref{sec:1.2.1}(1).

(c)$\Rightarrow$(a) Let $n$ be an integer and $A'\overset{u}\rightarrow 
A\overset{v}\rightarrow A''$ an exact sequence in $\mathcal{A}$ with 
$A''/\im v\in \mathcal{K}$. We have then the exact sequence 
of complexes 
\[
S^{0}(A')\overset{u}\rightarrow S^{0}(A)\overset{v}\rightarrow S^{0}(A'')
\] 
with $S^{0}(A'')/\im v\in Ch(\mathcal{K})$. The claim follows from 
\ref{sec:4.2} and \ref{sec:1.2.1}(1).
\end{proof}

\subsubsection{} \label{sec:8.1.2}
The pair $(\mathcal{K},\mathcal{K}^{\perp})$
is a cotorsion theory if and only if 
$Ch(\mathcal{K})=^{\perp_{e}}(Ch(\mathcal{K})^{\perp_{e}})$.

\subsection{} \label{sec:8.2}

Let $\mathcal{A}$ be an abelian category and $\mathfrak{K}$ 
a class of objects of $Ch(\mathcal{A})$. We define
\[
\mathfrak{K}^{\perp_{h}}=\{\mathbf{D}, 
\ H(\mathrm{Homgr}_{\mathcal{A}}(\mathbf{C},\mathbf{D}))=0 
{\rm \ for\ all} \ \mathbf{C}\in\mathfrak{K}\}
\]
and
\[
^{\perp_{h}}\mathfrak{K}=\{\mathbf{D}, 
\ H(\mathrm{Homgr}_{\mathcal{A}}(\mathbf{D},\mathbf{C}))=0 
{\rm \ for\ all} \ \mathbf{C}\in\mathfrak{K}\}
\]
We denote by $Ob(Ch(\mathcal{A}))_{0}$ the class of complexes with zero differential. 
Then $Ob(Ch(\mathcal{A}))^{\perp_{h}}=Ob(Ch(\mathcal{A}))_{0}^{\perp_{h}}$ 
consists of complexes that are homotopic to zero.
%We denote by $ex(\mathcal{A})$ the class of exact complexes in $\mathcal{A}$. The objects of 
%$ex(\mathcal{A})^{\perp_{h}}$ are called homotopically injective complexes and the objects of
%$^{\perp_{h}}ex(\mathcal{A})$ are called homotopically projective complexes.

\subsubsection{} \label{sec:8.2.1}

(1) If $\mathfrak{K},\mathfrak{L}$ are two classes of objects of $Ch(\mathcal{A})$, 
$\mathfrak{K}\subset \mathfrak{L}$ 
implies $\mathfrak{L}^{\perp_{h}}\subset\mathfrak{K}^{\perp_{h}}$.

(2) $\mathfrak{K}\subset {}^{\perp_{h}}(\mathfrak{K}^{\perp_{h}})$, 
with equality if and only if there is a class $\mathcal{S}$ of complexes such that 
$\mathcal{S}\subset \mathfrak{K}^{\perp_{h}}$ and 
${}^{\perp_{h}}\mathcal{S}\subset \mathfrak{K}$.

(3) $^{\perp_{h}}\mathfrak{K}=^{\perp_{h}}((^{\perp_{h}}\mathfrak{K})^{\perp_{h}})$.

(4) $\mathfrak{K}\cap \mathfrak{K}^{\perp_{h}}$ is the class of complexes 
that are homotopic to zero.

(5) The class $\mathfrak{K}^{\perp_{h}}$ is closed under homotopy 
equivalences, translates and retracts. 

(6) Let $u:\mathbf{D}'\to\mathbf{D}$ be a morphism. 

(a) If any two of $\mathbf{D},\mathbf{D}'$ and $Con(u)$ belong to 
$\mathfrak{K}^{\perp_{h}}$ then so does the third.

(b) If $\mathbf{D},\mathbf{D}'\in\mathfrak{K}^{\perp_{h}}$ and
$Con(u)\in\mathfrak{K}$ then $u$ is a homotopy equivalence.

\begin{proof}
Part 3 follows from 2 by taking $\mathcal{S}=\mathfrak{K}$. 

(6) (a) The exact sequence (\ref{sec:4.1.3}) 
$\mathbf{D}\overset{\pi}\rightarrowtail Con(u)\overset{\delta}
\twoheadrightarrow \mathbf{D}'(-1)$ is degreewise split, therefore 
(\ref{sec:4.2.1}) for each $\mathbf{C}\in\mathfrak{K}$ we have 
the exact sequence
\[
\mathrm{Homgr}_{\mathcal{A}}(\mathbf{C},\mathbf{D})\rightarrowtail 
\mathrm{Homgr}_{\mathcal{A}}(\mathbf{C},Con(u))\twoheadrightarrow
\mathrm{Homgr}_{\mathcal{A}}(\mathbf{C},\mathbf{D}'(-1))
\]
We concluse using \ref{sec:4.2} and part 5. For (b), using (a) and 
part 4 we obtain that $Con(u)$ is homotopic to zero, hence
(\ref{sec:4.1.3}) $u$ is a homotopy equivalence.
\end{proof}

\subsubsection{} \label{sec:8.2.2}
\cite[Proposition 1.4]{Sp} Suppose that $\mathfrak{K}$ is closed 
under translates. Recall (\ref{sec:4.1.3}) that $\sim$ denotes the 
homotopy relation for complexes. For a complex 
$\mathbf{D}$ in $\mathcal{A}$ consider the following statements:

(1) $\mathbf{D}\in\mathfrak{K}^{\perp_{h}}$;

(2) for every morphism $u:\mathbf{D}\to\mathbf{A}$ with $Con(u)\in\mathfrak{K}$
there is a morphism $l:\mathbf{A}\to\mathbf{D}$ such that $lu\sim 1_{\mathbf{D}}$;

(3) for every exact sequence 
$\mathbf{A}\overset{u}\rightarrowtail\mathbf{B}\twoheadrightarrow\mathbf{C}$
with $\mathbf{C}\in\mathfrak{K}$ and every morphism $f:\mathbf{A}\to \mathbf{D}$
there is $l:\mathbf{B}\to\mathbf{D}$ such that $lu\sim f$;

(4) for every morphism $u$ with cone in $\mathfrak{K}$, 
$\mathrm{Homgr}_{\mathcal{A}}(u,\mathbf{D})$ is a quasi-isomorphism.

We have $(1)\Leftrightarrow (4),(1)\Rightarrow (2)$ and $(3)\Rightarrow (2)$. 
In the case $(1)\Rightarrow (2)$ the morphism $l$ is unique up to homotopy.
If the cone of a monomorphism with cokernel in $\mathfrak{K}$ is in $\mathfrak{K}$ 
then the reverse implications hold.

\begin{proof}
$(1)\Rightarrow (4)$
This follows from the natural isomorphism (\ref{sec:4.2})
\[
\Sigma(u,\mathbf{D}):\mathrm{Homgr}_{\mathcal{A}}(Con(u)(1),\mathbf{D})
\cong Con(\mathrm{Homgr}_{\mathcal{A}}(u,\mathbf{D}))
\]
$(4)\Rightarrow (1)$
Let $\mathbf{C}\in\mathfrak{K}$ and consider the zero morphism 
$0:0\to\mathbf{C}$. Then $Con(0)=\mathbf{C}$, hence 
$0:\mathrm{Homgr}_{\mathcal{A}}
(\mathbf{C},\mathbf{D})\to 0$ is a quasi-isomorphism.

$(1)\Rightarrow (2)$ Let $u:\mathbf{D}\to\mathbf{A}$ be a morphism with 
$Con(u)\in\mathfrak{K}$. Consider the exact sequence $\mathbf{A}
\overset{\pi}\rightarrowtail Con(u)\overset{\delta}\twoheadrightarrow
\mathbf{D}(-1)$. Then $\delta$ is homotopic to zero;
let $s$ be a homotopy between $0$ and $\delta$, so that 
$s_{n}:Con(u)_{n}\to\mathbf{D}_{n}$. One can check that 
$s\pi:\mathbf{A}\to\mathbf{D}$ is a morphism. We define 
$t_{n}:\mathbf{D}_{n}\to\mathbf{D}_{n+1}$ as the composite
$\mathbf{D}_{n}\to\mathbf{D}_{n}\oplus\mathbf{A}_{n+1}\overset{-s_{n+1}}
\longrightarrow\mathbf{D}_{n+1}$, where the first arrow is the 
coprojection morphism. One can check that $t$ is a homotopy 
between $(-s\pi)  u$ and $1_{\mathbf{D}}$.

$(1)\Rightarrow (2)$ (uniqueness of $l$) We apply \ref{sec:4.2.2} to 
$u:\mathbf{D}\to\mathbf{A}$ and $\mathbf{C}''=\mathbf{D}$; we obtain 
that $-H(\mathrm{Homgr}_{\mathcal{A}}(u,\mathbf{D})(-1))$ is an isomorphism.
Therefore (\ref{sec:4.2}) for each integer $n$ the morphism
\[
Ch(\mathcal{A})(u,\mathbf{D}(n-1))/\sim:
Ch(\mathcal{A})(\mathbf{A},\mathbf{D}(n-1))/\sim \to
Ch(\mathcal{A})(\mathbf{D},\mathbf{D}(n-1))/\sim
\]
induced by $l\mapsto lu$, is an isomorphism. Taking $n=1$ we obtain the 
uniqueness of $l$.

$(3)\Rightarrow (2)$ Let $u:\mathbf{D}\to\mathbf{A}$ be a morphism with 
$Con(u)\in\mathfrak{K}$. From \ref{sec:4.1.3} we have the factorization
$\mathbf{D}\overset{\tilde{u}}\to Cyl(u)\overset{\beta}\to \mathbf{A}$ 
of $u$ and the exact sequence $\mathbf{D}\overset{\tilde{u}}\rightarrowtail 
Cyl(u)\overset{\tilde{\pi}}\twoheadrightarrow Con(u)$. Then
there is $l:Cyl(u)\to\mathbf{D}$ such that $l\tilde{u}\sim 1_{\mathbf{D}}$. 
If $\alpha$ is a homotopy inverse of $\beta$ (\emph{loc. cit.})
then $l\alpha u=l\alpha\beta\tilde{u}\sim l\tilde{u}\sim 1_{\mathbf{D}}$.

$(2)\Rightarrow (1)$ Let $\mathbf{C}\in\mathfrak{K}$ and 
$u:\mathbf{C}\to\mathbf{D}$ a morphism. Consider the exact 
sequence $\mathbf{D}\overset{\pi}\rightarrowtail Con(u)
\overset{\delta}\twoheadrightarrow\mathbf{C}(-1)$. By 
assumptions $Con(\pi)\in\mathfrak{K}$, therefore there is 
$l:Con(u)\to \mathbf{D}$ such that $l\pi\sim 1_{\mathbf{D}}$.
Since (\ref{sec:4.1.3}) $u=\beta\tilde{u}$ and $\pi=\tilde{\pi}\alpha$ 
we have $u\sim l\pi u=l\tilde{\pi}\alpha\beta\tilde{u}
\sim l\tilde{\pi}\tilde{u}=0$.

$(2)\Rightarrow (3)$ Let $\mathbf{A}\overset{u}\rightarrowtail
\mathbf{B}\twoheadrightarrow\mathbf{C}$ be an exact sequence 
with $\mathbf{C}\in\mathfrak{K}$ and $f:\mathbf{A}\to \mathbf{D}$ 
a morphism. Form the commutative diagram
\[
\xymatrix{
{\mathbf{A}} \ar@{>->}[r]^{u} \ar[d]_{f}& 
{\mathbf{B}}\ar@{->>}[r]\ar[d]^{f'}& {\mathbf{C}}\ar@{=}[d]\\
{\mathbf{D}} \ar@{>->}[r]^{u'} & {PO}\ar@{->>}[r] & {\mathbf{C}}
}
\]
where $PO$ means pushout. By assumption there is 
$l:PO\to\mathbf{D}$ such that $lu'\sim 1_{\mathbf{D}}$; then $lf'u\sim f$.
\end{proof}

\subsubsection{}  \label{sec:8.2.3}
Suppose that $\mathfrak{K}={}^{\perp_{h}}(\mathfrak{K}^{\perp_{h}})$.

(1) Let $u:\mathbf{C}'\to\mathbf{C}$ be a morphism such that
$\mathrm{Homgr}_{\mathcal{A}}(u,\mathbf{D})$ is a quasi-isomorphism
for all $\mathbf{D}\in\mathfrak{K}^{\perp_{h}}$; then $Con(u)\in\mathfrak{K}$.

(2) The class of morphisms with cone in $\mathfrak{K}$ has the two 
out of three property.

\begin{proof}
(1) Let $\mathbf{D}\in\mathfrak{K}^{\perp_{h}}$. From the natural 
isomorphism (\ref{sec:4.2})
\[
\Sigma(u,\mathbf{D}):\mathrm{Homgr}_{\mathcal{A}}(Con(u)(1),\mathbf{D})
\cong Con(\mathrm{Homgr}_{\mathcal{A}}(u,\mathbf{D}))
\]
and the assumption it follows that $Con(u)(1)\in{}^{\perp_{h}}
(\mathfrak{K}^{\perp_{h}})$, hence by \ref{sec:8.2.1}(5) we have 
$Con(u)\in\mathfrak{K}$. Part 2 follows from 1, the equivalence 
$(1)\Leftrightarrow (4)$ of \ref{sec:8.2.2} and the two out 
of three property of quasi-isomorphisms.
\end{proof}

We denote by $ex(\mathcal{A})$ the class of exact complexes in 
an abelian category $\mathcal{A}$.

\subsubsection{} \label{sec:8.2.4}
Suppose that $\mathcal{A}$ has either enough injective objects or a set 
$(E_{i})_{i\in I}$ of injective cogenerators. Then $ex(\mathcal{A})
=^{\perp_{h}}(ex(\mathcal{A})^{\perp_{h}})$.

\begin{proof}
We shall use \ref{sec:8.2.1}(2). Suppose that $\mathcal{A}$ has enough 
injective objects. We put $\mathcal{S}=\{S^{r}(E)\}$, where $r$ is an 
integer and $E$ is injective. We show that $\mathcal{S}\subset 
ex(\mathcal{A})^{\perp_{h}}$. Let $(\mathbf{C}',d')$ be an exact complex 
and $E$ and injective object; then by \ref{sec:4.2} and \ref{sec:1.2.1}(1) 
the complex $\mathcal{A}(\mathbf{C}'_{r-},E)$ is exact.
We show that $^{\perp_{h}}\mathcal{S}\subset ex(\mathcal{A})$.
Let $(\mathbf{C}',d')\in ^{\perp_{h}}\mathcal{S}$; then by \ref{sec:4.2}
the complex $\mathcal{A}(\mathbf{C}'_{r-},E)$ is exact for all injective objects 
$E$ of $\mathcal{A}$ and all integers $r$, hence (\ref{sec:1.1.2}) 
$(\mathbf{C}',d')$ is exact. When $\mathcal{A}$ has a set $(E_{i})_{i\in I}$ 
of injective cogenerators we put $\mathcal{S}=\{S^{0}(E_{i})\}$.
%If we put $d^{-}_{n}=(-1)^{-n}d_{n}$ then $(\mathbf{C},d^{-})$
%is exact, hence so is $\mathcal{A}(\mathbf{C},E_{i})$ since $E_{i}$ is injective.
%We show that $^{\perp_{h}}\mathcal{S}\subset ex(\mathcal{A})$.
%Let $(\mathbf{C},d)\in ^{\perp_{h}}\mathcal{S}$. It suffices to show that 
%$(\mathbf{C},d^{-})$ is exact. Since $\mathcal{A}(\mathbf{C},E_{i})$
%defined above is exact for all $i\in  I$, this follows from \ref{sec:1.2.2}.
\end{proof}

\subsubsection{} \label{sec:8.2.5}
Suppose that $\mathfrak{K}$ is closed under translates.

(1) $\mathfrak{K}^{\perp_{e}}\cap \mathfrak{K}^{\perp_{h}}
\subset\mathfrak{K}^{\perp}\subset \mathfrak{K}^{\perp_{h}}$; when 
$\mathfrak{K}=Ch(\mathcal{K})$, where $\mathcal{K}$ is a class of objects
of $\mathcal{A}$, we have $Ch(\mathcal{K})^{\perp_{e}}\cap 
Ch(\mathcal{K})^{\perp_{h}}= Ch(\mathcal{K})^{\perp}$.

(2) If $^{\perp_{h}}(\mathfrak{K}^{\perp})\subset \mathfrak{K}$ then 
$\mathfrak{K}={}^{\perp_{h}}(\mathfrak{K}^{\perp_{h}})$ and 
$(\mathfrak{K},\mathfrak{K}^{\perp})$ is a cotorsion theory.

(3) Let $\mathfrak{K}'$ consist of those $\mathbf{D}\in\mathfrak{K}^{\perp_{h}}$ 
for which the functor 
\[
\mathrm{Homgr}_{\mathcal{A}}(-,\mathbf{D}):Ch(\mathcal{A})^{op}\to Ch(\mathrm{Ab})
\]
preserves epimorphisms with kernel in $\mathfrak{K}$; then 
$\mathfrak{K}'\subset\mathfrak{K}^{\perp}$.

\begin{proof}
(1) Let $\mathbf{D}\in \mathfrak{K}^{\perp_{e}}\cap \mathfrak{K}^{\perp_{h}}$
and $\mathbf{A}\rightarrowtail \mathbf{B}\twoheadrightarrow\mathbf{C}$ be an 
exact sequence with $\mathbf{C}\in\mathfrak{K}$. We have then the exact sequence 
\[
\mathrm{Homgr}_{\mathcal{A}}(\mathbf{C},\mathbf{D})\rightarrowtail 
\mathrm{Homgr}_{\mathcal{A}}(\mathbf{B},\mathbf{D})\twoheadrightarrow
\mathrm{Homgr}_{\mathcal{A}}(\mathbf{A},\mathbf{D})
\]
Since $\mathrm{Homgr}_{\mathcal{A}}(\mathbf{C},\mathbf{D})$ is exact
and $Z_{0}$ preserves epimorphisms with exact kernel we have the exact sequence
\[
Z_{0}(\mathrm{Homgr}_{\mathcal{A}}(\mathbf{C},\mathbf{D}))\rightarrowtail 
Z_{0}(\mathrm{Homgr}_{\mathcal{A}}(\mathbf{B},\mathbf{D}))\twoheadrightarrow
Z_{0}(\mathrm{Homgr}_{\mathcal{A}}(\mathbf{A},\mathbf{D}))
\]
From \ref{sec:4.2} it follows that $\mathbf{D}\in \mathfrak{K}^{\perp}$.
Let now $\mathbf{D}\in\mathfrak{K}^{\perp},\mathbf{C}\in\mathfrak{K}$ 
and $u:\mathbf{C}\to\mathbf{D}$ a morphism. Then the exact sequence
$\mathbf{D}\overset{\pi}\rightarrowtail Con(u)\overset{\delta}
\twoheadrightarrow\mathbf{C}(-1)$ splits, therefore $\delta$ has a 
section $v$. We define $s_{n-1}:\mathbf{C}_{n-1}\to\mathbf{D}_{n}$
as the composite $\mathbf{C}_{n-1}\overset{v_{n}}\to
\mathbf{C}_{n-1}\oplus\mathbf{D}_{n}\to\mathbf{D}_{n}$, 
where the second arrow is the projection morphism. The fact that 
$v$ is a retract implies that $s$ is a homotopy between $0$ and $u$,
hence $\mathbf{D}\in\mathfrak{K}^{\perp_{h}}$.

Suppose that $\mathfrak{K}=Ch(\mathcal{K})$. Since 
$Ch(\mathcal{K})^{\perp}\subset Ch(\mathcal{K}^{\perp})$
we obtain from \ref{sec:8.1.1}(1) and the first part the desired equality.
Part 2 follows from 1 and \ref{sec:8.2.1}(2) by taking 
$\mathcal{S}=\mathfrak{K}^{\perp}$. Part 3 is similar to 1.
\end{proof}

\subsubsection{} \label{sec:8.2.6}
\cite[Proposition 3.7]{Gi2} Let $(\mathfrak{L},\mathfrak{L}^{\perp})$ be a 
cotorsion theory in $Ch(\mathcal{A})$ with $\mathfrak{L}$ closed under translates.

(1) If $\mathfrak{K}=\mathfrak{L}\cap ^{\perp_{h}}(\mathfrak{K}^{\perp})$ and 
$\mathfrak{K}^{\perp}\subset\mathfrak{L}^{\perp_{e}}$
then $(\mathfrak{K},\mathfrak{K}^{\perp})$ is a cotorsion theory and 
$\mathfrak{K}^{\perp}=\mathfrak{L}^{\perp_{e}}\cap\mathfrak{K}^{\perp_{h}}$.
If $\mathfrak{L}$ is left exact then so is $\mathfrak{K}$.

(2) Suppose that $\mathfrak{L}=^{\perp_{e}}(\mathfrak{L}^{\perp_{e}})$,
$\mathfrak{K}$ is closed under translates, $\mathfrak{K}\subset
\mathfrak{L}$ and $\mathfrak{K}^{\perp}\subset\mathfrak{L}^{\perp_{e}}$.
If $(\mathfrak{K},\mathfrak{K}^{\perp})$ is a cotorsion theory then
\[
\mathfrak{K}=\mathfrak{L}\cap ^{\perp_{h}}(\mathfrak{K}^{\perp})
{\rm \ and}\ \mathfrak{K}^{\perp}=\mathfrak{L}^{\perp_{e}}\cap\mathfrak{K}^{\perp_{h}}
\]
%Let $\mathcal{K}$ be a class of objects of $\mathcal{A}$
%and $\mathfrak{X}$ a class of objects of $Ch(\mathcal{A})$. Suppose that 
%$\mathfrak{X}\subset Ch(\mathcal{K}),\mathfrak{X}^{\perp}\subset Ch(\mathcal{K}^{\perp})$
%and $\mathfrak{X}$ is closed under translates. Then $(\mathfrak{X},\mathfrak{X}^{\perp})$
%is a cotorsion theory if and only if 
%\[
%\mathfrak{X}=Ch(\mathcal{K})\cap ^{\perp_{h}}(\mathfrak{X}^{\perp})
%{\rm \ and}\ \mathfrak{X}^{\perp}=Ch(\mathcal{K}^{\perp})\cap \mathfrak{X}^{\perp_{h}}
%\]
\begin{proof}
(1) One can check that $\mathfrak{K}^{\perp}$ is closed under translates; it 
then follows from \ref{sec:8.2.5}(1) that  $^{\perp}(\mathfrak{K}^{\perp})
\subset ^{\perp_{h}}(\mathfrak{K}^{\perp})$. Also, since 
$\mathfrak{K}\subset \mathfrak{L}$ we have $^{\perp}(\mathfrak{K}^{\perp})
\subset ^{\perp}(\mathfrak{L}^{\perp})=\mathfrak{L}$, hence
$(\mathfrak{K},\mathfrak{K}^{\perp})$ is a cotorsion theory. For the second
part it suffices, by \ref{sec:8.2.5}(1), to show that 
$\mathfrak{L}^{\perp_{e}}\cap\mathfrak{K}^{\perp_{h}}\subset\mathfrak{K}^{\perp}$.
Let $\mathbf{D}\in\mathfrak{L}^{\perp_{e}}\cap\mathfrak{K}^{\perp_{h}}$ 
and $\mathbf{A}\rightarrowtail \mathbf{B}\twoheadrightarrow \mathbf{C}$ 
be an exact sequence with $\mathbf{C}\in\mathfrak{K}$. 
Since $\mathfrak{L}^{\perp_{e}}\subset \mathfrak{K}^{\perp_{e}}$ 
have the exact sequence 
\[
\mathrm{Homgr}_{\mathcal{A}}(\mathbf{C},\mathbf{D})\rightarrowtail 
\mathrm{Homgr}_{\mathcal{A}}(\mathbf{B},\mathbf{D})\twoheadrightarrow
\mathrm{Homgr}_{\mathcal{A}}(\mathbf{A},\mathbf{D})
\]
Since $\mathrm{Homgr}_{\mathcal{A}}(\mathbf{C},\mathbf{D})$ is exact
and $Z_{0}$ preserves epimorphisms with exact kernel we have the exact sequence
\[
Z_{0}(\mathrm{Homgr}_{\mathcal{A}}(\mathbf{C},\mathbf{D}))\rightarrowtail 
Z_{0}(\mathrm{Homgr}_{\mathcal{A}}(\mathbf{B},\mathbf{D}))\twoheadrightarrow
Z_{0}(\mathrm{Homgr}_{\mathcal{A}}(\mathbf{A},\mathbf{D}))
\]
and we conclude by \ref{sec:4.2}. Suppose now that $\mathfrak{L}$ is left exact. 
Let $\mathbf{D}'\rightarrowtail\mathbf{D}\twoheadrightarrow\mathbf{D}''$ 
be an exact sequence with $\mathbf{D},\mathbf{D}''\in\mathfrak{K}$
and let $\mathbf{C}\in\mathfrak{K}^{\perp}$. Since $\mathbf{D}'\in\mathfrak{L}$ 
and $\mathbf{C}\in \mathfrak{L}^{\perp_{e}}$ we have the exact sequence 
\[
\mathrm{Homgr}_{\mathcal{A}}(\mathbf{D}'',\mathbf{C})\rightarrowtail
\mathrm{Homgr}_{\mathcal{A}}(\mathbf{D},\mathbf{C})\twoheadrightarrow
\mathrm{Homgr}_{\mathcal{A}}(\mathbf{D}',\mathbf{C})
\]
Since $\mathfrak{K}\subset^{\perp_{h}}(\mathfrak{K}^{\perp})$
the two left complexes of the previous exact sequence are exact, 
therefore so is the third.

(2) By part 1 it suffices to show the first equality. The inclusion
$\mathfrak{K}\subset\mathfrak{L}\cap ^{\perp_{h}}(\mathfrak{K}^{\perp})$
follows as in the proof of part 1. Let $\mathbf{C}\in
\mathfrak{L}\cap ^{\perp_{h}}(\mathfrak{K}^{\perp})$ and 
$\mathbf{D}\rightarrowtail \mathbf{A}\twoheadrightarrow \mathbf{B}$ 
be an exact sequence with $\mathbf{D}\in\mathfrak{K}^{\perp}$.
Since $\mathfrak{L}=^{\perp_{e}}(\mathfrak{L}^{\perp_{e}})$ 
have the exact sequence 
\[
\mathrm{Homgr}_{\mathcal{A}}(\mathbf{C},\mathbf{D})\rightarrowtail 
\mathrm{Homgr}_{\mathcal{A}}(\mathbf{C},\mathbf{A})\twoheadrightarrow
\mathrm{Homgr}_{\mathcal{A}}(\mathbf{C},\mathbf{B})
\]
Since $\mathrm{Homgr}_{\mathcal{A}}(\mathbf{C},\mathbf{D})$ is exact
and $Z_{0}$ preserves epimorphisms with exact kernel we have the exact sequence
\[
Z_{0}(\mathrm{Homgr}_{\mathcal{A}}(\mathbf{C},\mathbf{D}))\rightarrowtail 
Z_{0}(\mathrm{Homgr}_{\mathcal{A}}(\mathbf{C},\mathbf{A}))\twoheadrightarrow
Z_{0}(\mathrm{Homgr}_{\mathcal{A}}(\mathbf{C},\mathbf{B}))
\]
and we conclude by \ref{sec:4.2}.
\end{proof}

For an abelian category $\mathcal{A}$ and $\mathcal{K}$ a class
of objects of $\mathcal{A}$ we denote by $ex(\mathcal{K})$ the class of 
exact complexes in $\mathcal{A}$ that are degreewise in $\mathcal{K}$
and by $ex[\mathcal{K}]$ the class of exact complexes in $\mathcal{A}$ 
that have cycles in $\mathcal{K}$.

\subsubsection{} \label{sec:8.2.7}
Let $\mathcal{K}$ be a nonempty class of objects of $\mathcal{A}$.

(1) $ex[\mathcal{K}^{\perp}]=ex(\mathcal{K}^{\perp})\cap 
\{S^{0}(X),\ X\in\mathcal{K}\}^{\perp_{h}}$.

(2) \cite[Theorem 2.5]{YD},\cite[Proposition 3.7]{EJX}
Suppose that $\mathcal{A}$ has enough $\mathcal{K}$ objects, meaning 
that for each $A\in Ob(\mathcal{A})$ there is an epimorphism $X\to A$ 
with $X\in \mathcal{K}$, and that $\mathcal{K}$ is closed under extensions. 
Then $\{S^{0}(X),\ X\in\mathcal{K}\}^{\perp_{h}}\subset 
ex(\mathcal{A})\cap ex[\mathcal{K}]^{\perp_{h}}$. If, moreover,
$\mathcal{K}$ is left exact then we have equality and 
\[
ex[\mathcal{K}]\cap Ch(\mathcal{K}^{\perp})\cap ex[\mathcal{K}]^{\perp_{h}}=
ex[\mathcal{K}\cap \mathcal{K}^{\perp}]
\]
\begin{proof}
(1) Let $\mathbf{C}\in ex[\mathcal{K}^{\perp}]$ and $X\in \mathcal{K}$.
Then $\mathrm{Homgr}_{\mathcal{A}}(S^{0}(X),\mathbf{C})$ is exact 
by \ref{sec:1.2.1}((3) and (1bis)). Conversely, let $\mathbf{C}$ be an 
exact complex that is degreewise in $\mathcal{K}^{\perp}$ and such 
that $\mathrm{Homgr}_{\mathcal{A}}(S^{0}(X),\mathbf{C})$ is 
exact for all $X\in\mathcal{K}$. Since $ex(\mathcal{K}^{\perp})\cap 
\{S^{0}(X),\ X\in\mathcal{K}\}^{\perp_{h}}$ is closed under translates
it suffices to show that $Z_{0}(\mathbf{C})\in \mathcal{K}^{\perp}$.
Let $A\rightarrowtail B\twoheadrightarrow X$ be an 
exact sequence with $X\in \mathcal{K}$. Applying the functor 
$\mathrm{Homgr}_{\mathcal{A}}(-,\mathbf{C})$ 
to the exact sequence $S^{0}(A)\rightarrowtail S^{0}(B)
\twoheadrightarrow S^{0}(X)$ and using \ref{sec:8.1.1}(1) 
we obtain the exact sequence  
\[
\mathrm{Homgr}_{\mathcal{A}}(S^{0}(X),\mathbf{C})\rightarrowtail 
\mathrm{Homgr}_{\mathcal{A}}(S^{0}(B),\mathbf{C})\twoheadrightarrow
\mathrm{Homgr}_{\mathcal{A}}(S^{0}(A),\mathbf{C})
\]
Since $\mathrm{Homgr}_{\mathcal{A}}(S^{0}(X),\mathbf{C})$ is exact
and $Z_{0}$ preserves epimorphisms with exact kernel we have the exact sequence
\[
Z_{0}(\mathrm{Homgr}_{\mathcal{A}}(S^{0}(X),\mathbf{C}))\rightarrowtail 
Z_{0}(\mathrm{Homgr}_{\mathcal{A}}(S^{0}(B),\mathbf{C}))\twoheadrightarrow
Z_{0}(\mathrm{Homgr}_{\mathcal{A}}(S^{0}(A),\mathbf{C}))
\]
and we conclude by \ref{sec:4.2} and adjunction (\ref{sec:4.1}).

(2) Let $\mathbf{C}\in\{S^{0}(X),\ X\in\mathcal{K}\}^{\perp_{h}}$.
By \ref{sec:4.2} and \ref{sec:1.1.2}(1) we obtain that $\mathbf{C}$ is 
exact. Let $\mathbf{C}'\in ex[\mathcal{K}]$. Let $i$ be an integer
and consider the complex $\tau_{i}(\mathbf{C}')$ (\ref{sec:4.2.4}).
The complexes $\mathcal{A}(\tau_{i}(\mathbf{C}')_{r},\mathbf{C}_{r+})$ 
and $\mathcal{A}(\tau_{i}(\mathbf{C}')_{i},\mathbf{C})$ (\ref{sec:4.2.3})
are exact, hence $\mathrm{Homgr}_{\mathcal{A}}(\tau_{i}(\mathbf{C}'),\mathbf{C})$
is exact. Using \ref{sec:4.2.4} we conclude that $\mathrm{Homgr}_{\mathcal{A}}
(\mathbf{C}',\mathbf{C})$ is exact. Assume now that $\mathcal{K}$
is left exact. Let $\mathbf{C}\in ex(\mathcal{A})\cap ex[\mathcal{K}]^{\perp_{h}}$
and $f:\mathbf{X}_{0}\to B_{0}(\mathbf{C})$ be a morphism with 
$\mathbf{X}_{0}\in\mathcal{K}$. Choose an epimorphism $\mathbf{X}_{1}\to \mathbf{C}_{1}
\times_{B_{0}(\mathbf{C})}\mathbf{X}_{0}$ with $\mathbf{X}_{1}\in\mathcal{K}$
and put $B_{1}(\mathbf{X})=\K (\mathbf{X}_{1}\twoheadrightarrow\mathbf{X}_{0})$.
Then repeat this with the natural morphism $B_{1}(\mathbf{X})\to B_{1}(\mathbf{C})$
in place of $f$; continuing in this way we obtain $\mathbf{X}\in ex[\mathcal{K}]$ with 
$\mathbf{X}_{n}=0$ for $n<0$ and a morphism $u:\mathbf{X}\to\mathbf{C}$.
Then $u$ is homotopic to zero, therefore there is $s_{0}:\mathbf{X}_{0}\to\mathbf{C}_{1}$
such that $\delta_{1}s_{0}=f$. This implies that the adjoint transpose $f^{\#}$ of
$f$ is homotopic to zero, which is what we wanted to prove. The last part of 2
follows from the preceding considerations and part 1.
\end{proof}

\subsection{} \label{sec:8.3}
Let 
\[
T:\mathcal{A}'\times \mathcal{A}\to \mathcal{A}'', \
H:\mathcal{A}^{op}\times \mathcal{A}''\to \mathcal{A}', \
C:\mathcal{A}'^{op}\times \mathcal{A}''\to \mathcal{A}
\]
be an abelian THC--situation. Consider the 
standard THC--situation (\ref{sec:6.1})
\[
T:Ch(\mathcal{A}')\times Ch(\mathcal{A})\to Ch(\mathcal{A}'')
\]
\[
H:Ch(\mathcal{A})^{op}\times Ch(\mathcal{A}'')\to Ch(\mathcal{A}')
\]
\[
C:Ch(\mathcal{A}')^{op}\times Ch(\mathcal{A}'')\to Ch(\mathcal{A})
\]
associated to it. Let $\mathfrak{K}'$ (resp. $\mathfrak{K}''$) be a 
class of complexes in $\mathcal{A}'$ (resp. $\mathcal{A}''$) that is 
closed under isomorphisms and translates. We define
\[
\mathfrak{K}=T(\mathfrak{K}',-)^{-1}(\mathfrak{K}'')=
\{\mathbf{C}\in Ob(Ch(\mathcal{A})),\ T(\mathbf{C}',\mathbf{C})\in 
\mathfrak{K}'' {\rm \ for\ all}\  \mathbf{C}'\in\mathfrak{K}'\}
\]

\subsubsection{} \label{sec:8.3.1}
(1) If $\mathfrak{K}''$ is invariant under homotopy equivalences then so is $\mathfrak{K}$.

(2) $\mathfrak{K}$ is closed under translates. If, moreover, $\mathfrak{K}''$ is closed under 
extensions then $\mathfrak{K}$ is closed under cones of morphisms between its objects.

(3) For a complex $\mathbf{C}$ in $\mathcal{A}$ the following are equivalent:

(a) $\mathbf{C}\in\mathfrak{K}$;

(b) $T(-,\mathbf{C})$ sends a monomorphism with cone in $\mathfrak{K}'$ to a morphism 
with cone in $\mathfrak{K}''$.

(4) If $\mathfrak{K}''=^{\perp_{h}}(\mathfrak{K}''^{\perp_{h}})$ then 
$\mathfrak{K}=^{\perp_{h}}(\mathfrak{K}^{\perp_{h}})$.

\begin{proof}
(2) The closure of $\mathfrak{K}$ under cones of morphisms between its objects
follows from \ref{sec:6.1.2}.

(3) (a)$\Rightarrow$(b) This follows from \ref{sec:6.1.2}. (b)$\Rightarrow$(a)
Let $\mathbf{C}'\in\mathfrak{K}'$ and consider the zero morphism $0:0\to \mathbf{C}'$;
then $Con(0)=\mathbf{C}'$ and the assertion follows since $T(-,\mathbf{C})$ is additive.

(4) We shall use \ref{sec:8.2.1}(2). We put $\mathcal{S}=\{C(\mathbf{C}',\mathbf{D}''), 
\mathbf{C}'\in\mathfrak{K}', \mathbf{D}''\in\mathfrak{K}''^{\perp_{h}}\}$.
Let $\mathbf{C}'\in\mathfrak{K}'$ and $\mathbf{D}''\in\mathfrak{K}''^{\perp_{h}}$.
We show that $\mathcal{S}\subset \mathfrak{K}^{\perp_{h}}$. Let
$\mathbf{C}\in\mathfrak{K}$. From \ref{sec:6.1.9} we have an isomorphism
\[
\mathrm{Homgr}_{\mathcal{A}''}(T(\mathbf{C}',\mathbf{C}),\mathbf{D}'')\cong 
\mathrm{Homgr}_{\mathcal{A}}(\mathbf{C},C(\mathbf{C}',\mathbf{D}''))
\]
and the inclusion follows from this. We show that 
$^{\perp_{h}}\mathcal{S}\subset \mathfrak{K}$. If 
$\mathbf{C}\in^{\perp_{h}}\mathcal{S}$ then the previous 
isomorphism implies that $T(\mathbf{C}',\mathbf{C})\in^{\perp_{h}}
(\mathfrak{K}''^{\perp_{h}})$, therefore $\mathbf{C}\in\mathfrak{K}$.
\end{proof}

\subsubsection{} \label{sec:8.3.2}
A complex $\mathbf{C}$ (resp. $\mathbf{C}'$) in $\mathcal{A}$ 
(resp. $\mathcal{A}'$) is called $K$-\emph{flat} \cite{Sp}, 
or h-flat, if $T(-,\mathbf{C})$ (resp. $T(\mathbf{C}',-)$)
preserves exact complexes. By \ref{sec:8.3.1}(3) this is 
equivalent to $T(-,\mathbf{C})$ sending monomorphisms 
that are quasi-isomorphisms to quasi-isomorphisms.
Let $u'$ be a quasi-isomorphism in $Ch(\mathcal{A}')$.
Factor it (\ref{sec:4.1.3}) into $u'=\beta'\tilde{u'}$,
where $\beta'$ is a homotopy equivalence and $\tilde{u'}$
is a monomorphism. Since $T(-,\mathbf{C})$ is additive it 
preserves homotopy equivalences, therefore by the 
two out of three property of quasi-isomorphisms we 
have that $\mathbf{C}$ is h-flat if and only if 
$T(-,\mathbf{C})$ preserves quasi-isomorphisms.
We denote by $\mathfrak{F}_{h}$ the class of 
h-flat complexes in $\mathcal{A}$.
%One says that $Ch(\mathcal{A})$ (resp. $Ch(\mathcal{A}')$)
%has \emph{enough h-flat objects} if for every complex 
%$\mathbf{C}$ (resp. $\mathbf{C}'$) in $\mathcal{A}$ 
%(resp. $\mathcal{A}'$) there is a quasi-isomorphism 
%$\mathbf{P}\to\mathbf{C}$ (resp. $\mathbf{P}'\to\mathbf{C}'$)
%with $\mathbf{P}$ (resp. $\mathbf{P}'$) h-flat.

\subsubsection{} \label{sec:8.3.3}

(1) Suppose $\mathcal{A}''$ has either 
enough injectives or a set $(E''_{i})_{i\in I}$ of injective cogenerators. 

(a) $\mathfrak{F}_{h}=^{\perp_{h}}(\mathfrak{F}_{h}^{\perp_{h}})$.

(b) Let $\mathbf{P}$ be a complex in $\mathcal{A}$. When
$\mathcal{A}''$ has enough injectives then $\mathbf{P}$ is h-flat if and only if 
$H(\mathbf{P},S^{0}(E''))\in ex(\mathcal{A}')^{\perp_{h}}$ for all injective 
objects $E''$ of $\mathcal{A}''$ and when $\mathcal{A}''$ has  
a set $(E''_{i})_{i\in I}$ of injective cogenerators then $\mathbf{P}$ is 
h-flat if and only if $H(\mathbf{P},S^{0}(E''_{i}))\in 
ex(\mathcal{A}')^{\perp_{h}}$ for all $i\in I$.

(2) If $\mathcal{A}''$ satisfies AB4 and $\mathcal{A}$ has enough flat objects 
then $\mathfrak{F}_{h}^{\perp_{h}}\subset ex(\mathcal{A})$.

(3)  If $\mathcal{A}''$ hs enough injectives and $\mathcal{A}',\mathcal{A}$
have enough flat objects then 
\[
^{\perp_{h}}ex[\mathcal{F}^{\perp}]
\subset \mathfrak{F}_{h}=^{\perp_{h}}C(ex(\mathcal{A}'),S^{0}(\mathcal{I}''))
\]
where $C(ex(\mathcal{A}'),S^{0}(\mathcal{I}''))$ consists of objects of the form
$C(\mathbf{C}',S^{0}(E''))$, with $\mathbf{C}'$ exact in $\mathcal{A}'$ and 
$E''$ injective object of $\mathcal{A}''$.

\begin{proof}
(1) Part (a) follows from \ref{sec:8.2.4} and \ref{sec:8.3.1}(4). We prove (b).
Let $\mathbf{C}'$ be an exact complex in $\mathcal{A}'$ and $A''\in Ob(\mathcal{A}'')$.
By \ref{sec:6.1.9} and \ref{sec:4.2} we have 
\[
\mathrm{Homgr}_{\mathcal{A}'}(\mathbf{C}',H(\mathbf{P},S^{0}(A'')))\cong
\mathrm{Homgr}_{\mathcal{A}''}(T(\mathbf{C}',\mathbf{P}),S^{0}(A''))=
\mathcal{A}''(T(\mathbf{C}',\mathbf{P})_{0-},A'')
\]
The required equivalences follow from this by putting $A''=E''$ or $A''=E''_{i}$
and using \ref{sec:1.1.2}(1) or \ref{sec:1.1.1}(1).

(2) Let $\mathbf{D}\in\mathfrak{F}_{h}^{\perp_{h}}$. From \ref{sec:6.1.5}
we have that $S^{0}(P)$ is h-flat for every flat object $P$ of $\mathcal{A}$, 
therefore (\ref{sec:4.2}) $\mathcal{A}(P,\mathbf{D})$ is exact.
Since $\mathcal{A}$ has enough flat objects, it follows (\ref{sec:1.1.2}) 
that $\mathbf{D}$ is exact.

(3) We show that $^{\perp_{h}}ex[\mathcal{F}^{\perp}]
\subset \mathfrak{F}_{h}$. Let $\mathbf{P}\in^{\perp_{h}}
ex[\mathcal{F}^{\perp}]$ and $(\mathbf{C}',d')\in ex(\mathcal{A}')$. 
By \ref{sec:1.1.2}(1) it suffices to show that 
$\mathrm{Homgr}_{\mathcal{A}''}(T(\mathbf{C}',\mathbf{P}),S^{0}(E''))$
is exact for all injective objects $E''$ of $\mathcal{A}''$. By \ref{sec:6.1.9}
it suffices to show that $C(\mathbf{C}',S^{0}(E''))=C(\mathbf{C}'_{0-},E'')
\in ex[\mathcal{F}^{\perp}]$, where $C(\mathbf{C}'_{0-},E'')_{n}
=C(\mathbf{C}'_{-n},E'')$ has differential $C((-1)^{n+1}d'_{-n+1},E'')$.
By \ref{sec:3.1.3} the complex $C(\mathbf{C}'_{0-},E'')$ is exact, so 
$B_{n}(C(\mathbf{C}'_{0-},E''))=C(B_{-n-1}(\mathbf{C}'),E'')$ for all
integers $n$. Then from the proof of \ref{sec:3.1.4} we have that 
$C(B_{-n-1}(\mathbf{C}'),E'')\in\mathcal{F}^{\perp}$. The 
displayed equality is straightforward.

\end{proof}

\section{The universal coefficients theorem for complexes in abelian categories}

\subsection{} \label{sec:9.1}
Let $\mathcal{A}$ be an abelian category and $\mathfrak{K}$
a class of complexes in $\mathcal{A}$ such that $\mathfrak{K}=
^{\perp_{e}}(\mathfrak{K}^{\perp_{e}})$ (\ref{sec:8.1}).
Let $(\mathbf{C},d),(\mathbf{C}',d')$ be complexes in $\mathcal{A}$.
Recall (\ref{sec:4.2}) that we denote by $D(\mathbf{C},\mathbf{C}')$
the differential of $\mathrm{Homgr}_{\mathcal{A}}(\mathbf{C},\mathbf{C}')$.
We denote by 
\[
\begin{aligned}
i(\mathbf{C},\mathbf{C}'):B(\mathrm{Homgr}_{\mathcal{A}}(\mathbf{C},\mathbf{C}'))
\to Z(\mathrm{Homgr}_{\mathcal{A}}(\mathbf{C},\mathbf{C}'))\\
j(\mathbf{C},\mathbf{C}'): Z(\mathrm{Homgr}_{\mathcal{A}}(\mathbf{C},\mathbf{C}'))
\to \mathrm{Homgr}_{\mathcal{A}}(\mathbf{C},\mathbf{C}')\\
p(\mathbf{C},\mathbf{C}'):Z(\mathrm{Homgr}_{\mathcal{A}}(\mathbf{C},\mathbf{C}'))
\to H(\mathrm{Homgr}_{\mathcal{A}}(\mathbf{C},\mathbf{C}'))
\end{aligned}
\]
the natural morphisms. We have a natural in both arguments morphism
$h(\mathbf{C},\mathbf{C}'):Z(\mathrm{Homgr}_{\mathcal{A}}(\mathbf{C},\mathbf{C}'))
\to \mathrm{Homgr}_{\mathcal{A}}(H(\mathbf{C}),H(\mathbf{C}'))$
defined as $h(\mathbf{C},\mathbf{C}')(f)=(H_{p}(f))_{p}$. One checks 
\cite[\S5 n\textsuperscript{\scriptsize{o}} 1 p. 82]{Bou} that 
$h(\mathbf{C},\mathbf{C}')i(\mathbf{C},\mathbf{C}')=0$, hence there 
is a unique morphism 
\[
\lambda(\mathbf{C},\mathbf{C}'):H(\mathrm{Homgr}_{\mathcal{A}}(\mathbf{C},\mathbf{C}'))
\to \mathrm{Homgr}_{\mathcal{A}}(H(\mathbf{C}),H(\mathbf{C}'))
\]
such that $\lambda(\mathbf{C},\mathbf{C}')p(\mathbf{C},\mathbf{C}')
=h(\mathbf{C},\mathbf{C}')$. Since $h$ is natural and $p(\mathbf{C},\mathbf{C}')$
is an epimorphism it follows that $\lambda$ is natural in both arguments.

\subsubsection{} \label{sec:9.1.1}
Suppose $d=0$. Then $D(\mathbf{C},\mathbf{C}')=
\mathrm{Homgr}_{\mathcal{A}}(\mathbf{C},d')$ and we identify 
\[
(Z(\mathrm{Homgr}_{\mathcal{A}}(\mathbf{C},\mathbf{C}')),j(\mathbf{C},\mathbf{C}'))
=(\mathrm{Homgr}_{\mathcal{A}}(\mathbf{C},Z(\mathbf{C}')),\bar{j'})
\]
so that we have $h(\mathbf{C},\mathbf{C}')=
\mathrm{Homgr}_{\mathcal{A}}(\mathbf{C},p')(=\bar{p'})$.
Also, the adjoint transpose $D(\mathbf{C},\mathbf{C}')^{\#}$ of 
$D(\mathbf{C},\mathbf{C}')$ is $\mathrm{Homgr}_{\mathcal{A}}(\mathbf{C},d'^{\#})$.
From the commutative diagram
\begin{displaymath}
\xymatrix{
{B(\mathrm{Homgr}_{\mathcal{A}}(\mathbf{C},\mathbf{C}'))}
\ar[rrr]^{i(\mathbf{C},\mathbf{C}')} &&&
{Z(\mathrm{Homgr}_{\mathcal{A}}(\mathbf{C},\mathbf{C}'))} 
\ar[d]_{j(\mathbf{C},\mathbf{C}')}\\
{\mathrm{Homgr}_{\mathcal{A}}(\mathbf{C}),\mathbf{C}'(1))} 
\ar[u]\ar[d]_{\overline{\delta'(1)}}\ar[rrr]^{D(\mathbf{C},\mathbf{C}')^{\#}=
\mathrm{Homgr}_{\mathcal{A}}(\mathbf{C},d'^{\#})}&&&
{\mathrm{Homgr}_{\mathcal{A}}(\mathbf{C},\mathbf{C}')}\\
{\mathrm{Homgr}_{\mathcal{A}}(\mathbf{C},B(\mathbf{C}'))}
\ar[rrr]^{\bar{i'}} &&& {\mathrm{Homgr}_{\mathcal{A}}(\mathbf{C},Z(\mathbf{C}'))} 
\ar[u]_{\bar{j'}}\\
}
\end{displaymath}
we obtain that there is a morphism 
\[
b_{2}(\mathbf{C},\mathbf{C}'):B(\mathrm{Homgr}_{\mathcal{A}}(\mathbf{C},\mathbf{C}'))
\to \mathrm{Homgr}_{\mathcal{A}}(\mathbf{C},B(\mathbf{C}'))
\]
such that $\overline{j'i'}b_{2}(\mathbf{C},\mathbf{C}')=
j(\mathbf{C},\mathbf{C}')i(\mathbf{C},\mathbf{C}')$; it follows
that $\bar{i'}b_{2}(\mathbf{C},\mathbf{C}')=i(\mathbf{C},\mathbf{C}')$.

(1) Suppose $\mathbf{C}\in\mathfrak{K}$ and $Z(\mathbf{C}')\in\mathfrak{K}^{\perp_{e}}$;
then $b_{2}(\mathbf{C},\mathbf{C}')$ is an isomorphism.

(2) Suppose $\mathbf{C}\in\mathfrak{K}$ and $B(\mathbf{C}'),Z(\mathbf{C}')
\in\mathfrak{K}^{\perp_{e}}$; then $\lambda(\mathbf{C},\mathbf{C}')$ 
is an isomorphism.
\begin{proof}
(1) Applying $\mathrm{Homgr}_{\mathcal{A}}(\mathbf{C},-)$ to the exact sequence
$(j'(1),\delta'(1))$ (\ref{sec:4.1}) we obtain that $\overline{\delta'(1)}$
is an epimorphism, hence so is $b_{2}(\mathbf{C},\mathbf{C}')$.

(2) Applying $\mathrm{Homgr}_{\mathcal{A}}(\mathbf{C},-)$ to the exact 
sequence $B(\mathbf{C}')\overset{i'}\rightarrowtail Z(\mathbf{C}')\overset{p'}
\twoheadrightarrow H(\mathbf{C}')$ we obtain the exact sequence
\[
\mathrm{Homgr}_{\mathcal{A}}(\mathbf{C},B(\mathbf{C}'))\overset{\bar{i'}}
\rightarrowtail \mathrm{Homgr}_{\mathcal{A}}(\mathbf{C},Z(\mathbf{C}'))\overset{\bar{p'}}
\twoheadrightarrow \mathrm{Homgr}_{\mathcal{A}}(\mathbf{C},H(\mathbf{C}'))
\]
Using part 1, this implies that $\lambda(\mathbf{C},\mathbf{C}')$ 
is an isomorphism.
\end{proof}

\subsubsection{} \label{sec:9.1.2}
Suppose $B(\mathbf{C})\in\mathfrak{K}$, $d'=0$ and 
$\mathbf{C}'\in\mathfrak{K}^{\perp_{e}}$. Applying 
$\mathrm{Homgr}_{\mathcal{A}}(-,\mathbf{C}')$ to the 
exact sequence $(j,\delta)$ we obtain the exact sequence
\[
\mathrm{Homgr}_{\mathcal{A}}(B(\mathbf{C})(-1),\mathbf{C}')\overset{\bar{\delta}}
\rightarrowtail \mathrm{Homgr}_{\mathcal{A}}(\mathbf{C},\mathbf{C}')\overset{\bar{j}}
\twoheadrightarrow \mathrm{Homgr}_{\mathcal{A}}(Z(\mathbf{C}),\mathbf{C}')
\]
where the outer complexes have zero differential. We want to calculate the connecting 
morphism 
\[
\partial(\bar{\delta},\bar{j}):H(\mathrm{Homgr}_{\mathcal{A}}(Z(\mathbf{C}),\mathbf{C}'))
\to H(\mathrm{Homgr}_{\mathcal{A}}(B(\mathbf{C})(-1),\mathbf{C}'))(-1)
\]
associated to the previous exact sequence. From the standard construction
\cite[\S2 n\textsuperscript{\scriptsize{o}} 3 p. 29]{Bou} of the connecting
morphism it follows that we have the commutative diagram
\begin{displaymath}
\xymatrix{
{H(\mathrm{Homgr}_{\mathcal{A}}(Z(\mathbf{C}),\mathbf{C}'))}
\ar[rrr]^{\partial(\bar{\delta},\bar{j})} 
\ar[ddrrr]_{H(\mathrm{Homgr}_{\mathcal{A}}(i,\mathbf{C}'))}&&&
{H(\mathrm{Homgr}_{\mathcal{A}}(B(\mathbf{C})(-1),\mathbf{C}'))(-1)} 
\ar[d]_{H(\Sigma(B(\mathbf{C}),\mathbf{C}'))(-1)}\\
&&&{H(\mathrm{Homgr}_{\mathcal{A}}(B(\mathbf{C}),\mathbf{C}')(1))(-1)}\ar@{=}[d]\\
&&&{H(\mathrm{Homgr}_{\mathcal{A}}(B(\mathbf{C}),\mathbf{C}'))}
}
\end{displaymath}
where $\Sigma(B(\mathbf{C}),\mathbf{C}')$ has been defined in \ref{sec:4.2}.

\subsubsection{} \label{sec:9.1.3}
Suppose $B(\mathbf{C})\in\mathfrak{K}$ and $Z(\mathbf{C}'),\mathbf{C}'
\in\mathfrak{K}^{\perp_{e}}$. As in \ref{sec:9.1.2} we want to calculate
the connecting morphism 
\[
\partial(\bar{\delta},\bar{j}):H(\mathrm{Homgr}_{\mathcal{A}}(Z(\mathbf{C}),\mathbf{C}'))
\to H(\mathrm{Homgr}_{\mathcal{A}}(B(\mathbf{C})(-1),\mathbf{C}'))(-1)
\]
associated to the exact sequence
\[
\mathrm{Homgr}_{\mathcal{A}}(B(\mathbf{C})(-1),\mathbf{C}')\overset{\bar{\delta}}
\rightarrowtail \mathrm{Homgr}_{\mathcal{A}}(\mathbf{C},\mathbf{C}')\overset{\bar{j}}
\twoheadrightarrow \mathrm{Homgr}_{\mathcal{A}}(Z(\mathbf{C}),\mathbf{C}')
\]
Using the morphism $j':Z(\mathbf{C}')\to \mathbf{C}'$ and the naturality
of the connecting morphism we have the commutative diagram
\begin{displaymath}
\xymatrix{
{H(\mathrm{Homgr}_{\mathcal{A}}(Z(\mathbf{C}),Z(\mathbf{C}')))}
\ar[rrr]^{\partial(\bar{\delta},\bar{j})} 
\ar[d]_{H(\mathrm{Homgr}_{\mathcal{A}}(Z(\mathbf{C}),j'))}&&&
{H(\mathrm{Homgr}_{\mathcal{A}}(B(\mathbf{C})(-1),Z(\mathbf{C}')))(-1)} 
\ar[d]_{H(\mathrm{Homgr}_{\mathcal{A}}(B(\mathbf{C})(-1),j'))(-1)}\\
{H(\mathrm{Homgr}_{\mathcal{A}}(Z(\mathbf{C}),\mathbf{C}'))} 
\ar[rrr]^{\partial(\bar{\delta},\bar{j})}&&&
{H(\mathrm{Homgr}_{\mathcal{A}}(B(\mathbf{C})(-1),\mathbf{C}'))(-1)}
}
\end{displaymath}
By the naturality of the morphism $\Sigma(\mathbf{C}',\mathbf{C})$ 
we have the commutative diagram
\begin{displaymath}
\xymatrix{
{H(\mathrm{Homgr}_{\mathcal{A}}(B(\mathbf{C})(-1),Z(\mathbf{C}')))(-1)}
\ar[rrr]^{H(\Sigma(B(\mathbf{C}),Z(\mathbf{C}')))(-1)} 
\ar[d]^{H(\mathrm{Homgr}_{\mathcal{A}}(B(\mathbf{C})(-1),j'))(-1)}&&&
{H(\mathrm{Homgr}_{\mathcal{A}}(B(\mathbf{C}),Z(\mathbf{C}'))(1))(-1)} 
\ar[d]_{H(\mathrm{Homgr}_{\mathcal{A}}(B(\mathbf{C}),j')(1))(-1)}\\
{H(\mathrm{Homgr}_{\mathcal{A}}(B(\mathbf{C})(-1),\mathbf{C}'))(-1)} 
\ar[rrr]^{H(\Sigma(B(\mathbf{C}),\mathbf{C}'))(-1)}&&&
{H(\mathrm{Homgr}_{\mathcal{A}}(B(\mathbf{C}),\mathbf{C}')(1))(-1)}
}
\end{displaymath}
We claim that $H(\Sigma(B(\mathbf{C}),\mathbf{C}'))(-1)\partial(\bar{\delta},\bar{j})
=H(\mathrm{Homgr}_{\mathcal{A}}(i,\mathbf{C}'))$. Putting the two previous
diagrams side by side we obtain a diagram whose top horizontal composite morphism 
is $H(\mathrm{Homgr}_{\mathcal{A}}(i,Z(\mathbf{C}')))$ by \ref{sec:9.1.2}.
By \ref{sec:9.1.1} we have $j(Z(\mathbf{C}),\mathbf{C}')=
\mathrm{Homgr}_{\mathcal{A}}(Z(\mathbf{C}),j')$, hence
$H(\mathrm{Homgr}_{\mathcal{A}}(Z(\mathbf{C}),j'))=
p(Z(\mathbf{C}),\mathbf{C}')$ is an epimorphism. Using \ref{sec:9.1.2}
we have 
\[
H(\Sigma(B(\mathbf{C}),\mathbf{C}'))(-1)\partial(\bar{\delta},\bar{j})
H(\mathrm{Homgr}_{\mathcal{A}}(Z(\mathbf{C}),j'))=
\]
\[
H(\mathrm{Homgr}_{\mathcal{A}}(B(\mathbf{C}),j'))
H(\Sigma(B(\mathbf{C}),Z(\mathbf{C}')))(-1)\partial(\bar{\delta},\bar{j})=
\]
\[
H(\mathrm{Homgr}_{\mathcal{A}}(B(\mathbf{C}),j'))
H(\mathrm{Homgr}_{\mathcal{A}}(i,Z(\mathbf{C}')))=
\]
\[
H(\mathrm{Homgr}_{\mathcal{A}}(i,\mathbf{C}'))
H(\mathrm{Homgr}_{\mathcal{A}}(Z(\mathbf{C}),j'))
\]
hence the claim.

\subsubsection{} \label{sec:9.1.4}
Suppose that $B(\mathbf{C}),Z(\mathbf{C})\in \mathfrak{K}$
and $B(\mathbf{C}'),Z(\mathbf{C}'),\mathbf{C}'\in \mathfrak{K}^{\perp_{e}}$.
Then we have an exact sequence
\begin{displaymath}
\xymatrix{
{\Co\mathrm{Homgr}_{\mathcal{A}}(i(-1),H(\mathbf{C}'))}\ar@{>->}[r]&
{H(\mathrm{Homgr}_{\mathcal{A}}(\mathbf{C},\mathbf{C}'))}
\ar@{->>}[r]^{\lambda(\mathbf{C},\mathbf{C}')}&
{\mathrm{Homgr}_{\mathcal{A}}(H(\mathbf{C}),H(\mathbf{C}'))}
}
\end{displaymath}
\begin{proof}
Applying $\mathrm{Homgr}_{\mathcal{A}}(-,\mathbf{C}')$ to 
the exact sequence $(j,\delta)$ we obtain the exact sequence
$(\bar{\delta},\bar{j})$ (\ref{sec:9.1.3}). Consider then the top 
diagram on the next page, where the unlabeled top horizontal arrow 
is $H(\Sigma(B(\mathbf{C}),\mathbf{C}'))(-1)$.
\begin{sidewaysfigure}
\xymatrix{
\\
\\
\\
\\
\\
\\
\\
\\
\\
\ar[r]^{H(\bar{\delta})}&{H(\mathrm{Homgr}_{\mathcal{A}}(\mathbf{C},\mathbf{C}'))}
\ar[rr]^{H(\bar{j})} \ar[dd]_{\lambda(\mathbf{C},\mathbf{C}')}&&
{H(\mathrm{Homgr}_{\mathcal{A}}(Z(\mathbf{C}),\mathbf{C}'))} 
\ar[dd]_{\lambda(Z(\mathbf{C}),\mathbf{C}')}
\ar[r]^{\partial(\bar{\delta},\bar{j})}&
{H(\mathrm{Homgr}_{\mathcal{A}}(B(\mathbf{C})(-1),\mathbf{C}'))(-1)}
\ar[r]&{H(\mathrm{Homgr}_{\mathcal{A}}(B(\mathbf{C}),\mathbf{C}'))}
\ar[dd]_{\lambda(B(\mathbf{C}),\mathbf{C}')}\\
\\
&{\mathrm{Homgr}_{\mathcal{A}}(H(\mathbf{C}),H(\mathbf{C}'))} 
\ar[rr]^{\mathrm{Homgr}_{\mathcal{A}}(p,H(\mathbf{C}'))}&&
{\mathrm{Homgr}_{\mathcal{A}}(Z(\mathbf{C}),H(\mathbf{C}'))}
\ar[rr]^{\mathrm{Homgr}_{\mathcal{A}}(i,H(\mathbf{C}'))}&&
{\mathrm{Homgr}_{\mathcal{A}}(B(\mathbf{C}),H(\mathbf{C}'))}\\
\\
&{H(\mathrm{Homgr}_{\mathcal{A}}(Z(\mathbf{C})(-1),\mathbf{C}'))}
\ar[rr] \ar[dd]_{\lambda(Z(\mathbf{C})(-1),\mathbf{C}')}&&
{H(\mathrm{Homgr}_{\mathcal{A}}(Z(\mathbf{C}),\mathbf{C}')(1))} 
\ar[r]^{\partial(\bar{\delta},\bar{j})^{\#}}&
{H(\mathrm{Homgr}_{\mathcal{A}}(B(\mathbf{C})(-1),\mathbf{C}'))}
\ar[r]^{H(\bar{\delta})}\ar[dd]_{\lambda(B(\mathbf{C})(-1),\mathbf{C}')}
&{H(\mathrm{Homgr}_{\mathcal{A}}(\mathbf{C},\mathbf{C}'))}\\
\\
&{\mathrm{Homgr}_{\mathcal{A}}(Z(\mathbf{C})(-1),H(\mathbf{C}'))} 
\ar[rrr]^{\mathrm{Homgr}_{\mathcal{A}}(i(-1),H(\mathbf{C}'))}&&&
{\mathrm{Homgr}_{\mathcal{A}}(B(\mathbf{C})(-1),H(\mathbf{C}'))}
}
\end{sidewaysfigure}
The top sequence is exact from the homology exact sequence associated
to $(\bar{\delta},\bar{j})$ and the bottom sequence is exact by the 
construction of $\mathrm{Homgr}_{\mathcal{A}}$ (\ref{sec:4.2}).
The left square diagram commutes by the naturality of $\lambda$
and the right rectangular diagram commutes by \ref{sec:9.1.3}.
The morphisms $\lambda(Z(\mathbf{C}),\mathbf{C}')$ and 
$\lambda(B(\mathbf{C}),\mathbf{C}')$ are isomorphisms by
\ref{sec:9.1.1}. It follows that $\lambda(\mathbf{C},\mathbf{C}')$
is an epimorphism with $\K\lambda(\mathbf{C},\mathbf{C}')=
\im H(\bar{\delta})$. 

Consider now the bottom diagram on the next page, 
where the unlabeled top horizontal arrow 
is $H(\Sigma(Z(\mathbf{C}),\mathbf{C}'))$.
The top sequence is exact from the homology exact sequence 
associated to $(\bar{\delta},\bar{j})$. By \ref{sec:9.1.3} we 
have $\partial(\bar{\delta},\bar{j})^{\#}H(\Sigma(Z(\mathbf{C}),\mathbf{C}'))
=H(\mathrm{Homgr}_{\mathcal{A}}(i(-1),\mathbf{C}'))$, hence
the rectangular diagram commutes. By \ref{sec:9.1.1} the vertical 
morphisms are isomorphisms. It follows that 
\[ 
\im H(\bar{\delta})=
\Co \mathrm{Homgr}_{\mathcal{A}}(i(-1),H(\mathbf{C}'))
\]
\end{proof}

\subsubsection{} \label{sec:9.1.5}
Let $(\mathcal{K},\mathcal{K}^{\perp})$ be a cotorsion theory in 
$\mathcal{A}$.

(1) \cite[\S5 n\textsuperscript{\scriptsize{o}} 6 cor. 4]{Bou}
Suppose that $B(\mathbf{C}),H(\mathbf{C})\in Ch(\mathcal{K})$
and $B(\mathbf{C}'),H(\mathbf{C}')\in Ch(\mathcal{K}^{\perp})$;
then $\lambda(\mathbf{C},\mathbf{C}')$ is an isomorphism.

(2) The complex 
$\mathrm{Homgr}_{\mathcal{A}}(\mathbf{C},\mathbf{C}')$ is exact
if either (a) $\mathbf{C}\in ex[\mathcal{K}]$ (\ref{sec:4.1.2}) and 
$B(\mathbf{C}'),Z(\mathbf{C}')\in Ch(\mathcal{K}^{\perp})$
or (b) $B(\mathbf{C}),Z(\mathbf{C})\in Ch(\mathcal{K})$ and 
$\mathbf{C}'\in ex[\mathcal{K}^{\perp}]$.
\begin{proof}
(1) Put $\mathfrak{K}=Ch(\mathcal{K})$ in \ref{sec:9.1.4}; then
$\mathfrak{K}^{\perp_{e}}=Ch(\mathcal{K}^{\perp})$ (\ref{sec:8.1.1})
and $\mathfrak{K}=^{\perp_{e}}(\mathfrak{K}^{\perp_{e}})$ (\ref{sec:8.1.2}),
hence we have the exact sequence \ref{sec:9.1.4}. Applying 
$\mathrm{Homgr}_{\mathcal{A}}(-,H(\mathbf{C}'))$ to the 
exact sequence $B(\mathbf{C})(-1)\overset{i(-1)}\rightarrowtail 
Z(\mathbf{C})(-1)\overset{p(-1)}\twoheadrightarrow H(\mathbf{C})(-1)$
we obtain that $\overline{i(-1)}$ is an epimorphism, hence its cokernel
is $0$.

(2) We prove (a). As in part 1, it follows that under the given assumptions 
we have the exact sequence \ref{sec:9.1.4}. Moreover, the target 
of $\lambda(\mathbf{C},\mathbf{C}')$ is $0$ and $\overline{i(-1)}$
is an isomorphism, hence the assertion.
\end{proof}

\section{Some cotorsion theories in categories of complexes}

\subsection{} \label{sec:10.1}

Let $\mathcal{A}$ be an abelian category and $\mathcal{K}$ a class
of objects of $\mathcal{A}$. Recall that $ex(\mathcal{K})$ is the class of 
exact complexes in $\mathcal{A}$ that are degreewise in $\mathcal{K}$
and $ex[\mathcal{K}]$ is the class of exact complexes in $\mathcal{A}$ 
that have cycles in $\mathcal{K}$. We define $dg(\mathcal{K})=
^{\perp}ex[\mathcal{K}^{\perp}]$ and $dg(\mathcal{K}^{\perp})=
ex[\mathcal{K}]^{\perp}$.

\subsubsection{} \label{sec:10.1.1}
\cite[Proposition 3.6]{Gi1} Let $(\mathcal{K},\mathcal{K}^{\perp})$ 
be a cotorsion theory in $\mathcal{A}$.

(1) If $\mathcal{A}$ has enough $\mathcal{K}^{\perp}$ objects 
then $(ex[\mathcal{K}],dg(\mathcal{K}^{\perp}))$ is a 
cotorsion theory and $dg(\mathcal{K}^{\perp})=Ch(\mathcal{K}^{\perp})
\cap ex[\mathcal{K}]^{\perp_{h}}$.

(2) If $\mathcal{A}$ has enough $\mathcal{K}$ objects then 
$(dg(\mathcal{K}),ex[\mathcal{K}^{\perp}])$ is a cotorsion theory 
and $dg(\mathcal{K})=Ch(\mathcal{K})\cap ^{\perp_{h}}ex[\mathcal{K}^{\perp}]$.

\begin{proof}
(1) We only sketch the argument. One uses \ref{sec:1.2.3} with
$\mathcal{S}=\{S^{n}(Y)\}$, where $n$ is an integer and 
$Y\in\mathcal{K}^{\perp}$, together with the exact sequences
$S^{n}(Y)\rightarrowtail D^{n+1}(Y)\twoheadrightarrow S^{n+1}(Y)$ 
and $S^{n}(Y)\rightarrowtail S^{n}(A)\twoheadrightarrow S^{n}(B)$.
The description of $dg(\mathcal{K}^{\perp})$ follows from \ref{sec:8.2.6}(2).
The proof of part 2 is similar.
\end{proof}

\subsubsection{} \label{sec:10.1.2}
In the situation \ref{sec:10.1.1}, if $\mathcal{K}$ is left exact then
(\ref{sec:8.2.7}) $ex[\mathcal{K}]\cap dg(\mathcal{K}^{\perp})=
ex[\mathcal{K}\cap\mathcal{K}^{\perp}]=
dg(\mathcal{K})\cap ex[\mathcal{K}^{\perp}]$.

\subsubsection{} \label{sec:10.1.3}
\cite[Proposition 3.3]{Gi2} Let $\mathcal{A}$ be an abelian 
category and $(\mathcal{K},\mathcal{K}^{\perp})$ a cotorsion theory 
in $\mathcal{A}$. Suppose that $(\mathcal{K},\mathcal{K}^{\perp})$ 
is generated by $\mathcal{S}$ and that $\mathcal{A}$ has enough 
$\underline{\mathcal{K}^{\perp}}$ objects (\ref{sec:1.2.14}), 
meaning that for each $A\in Ob(\mathcal{A})$ there is a monomorphism 
$A\to Y$ with $Y\in \underline{\mathcal{K}^{\perp}}$. Then 
$(ex(\mathcal{K}),ex(\mathcal{K})^{\perp})$ 
is a cotorsion theory and $ex(\mathcal{K})^{\perp}=
Ch(\mathcal{\mathcal{K}^{\perp}})\cap ex(\mathcal{K})^{\perp_{h}}$.

\begin{proof}
We shall use \ref{sec:1.2.3}. We put $\mathcal{S}'=\{D^{m}(Y),S^{n}(\underline{Y})\}$, 
where $Y\in\mathcal{S},\underline{Y}\in \underline{\mathcal{K}^{\perp}}$ 
and $m,n$ are integers.
\paragraph{Step 1} We show that $\mathcal{S}'\subset ex(\mathcal{K})^{\perp}$.
We fix an exact sequence of complexes 
\begin{equation}\label{eq:10.1.3}
\mathbf{C}'\rightarrowtail \mathbf{C}\twoheadrightarrow \mathbf{C}''
\end{equation}
with $\mathbf{C}''\in ex(\mathcal{K})$. Applying the functor 
$Ch(\mathcal{A})(-,D^{m}(Y))$ to (\ref{eq:10.1.3}) we obtain by 
adjunction an exact sequence 
\begin{displaymath}
\mathcal{A}(\mathbf{C}''_{m-1},Y)\rightarrowtail 
\mathcal{A}(\mathbf{C}_{m-1},Y)\twoheadrightarrow 
\mathcal{A}(\mathbf{C}'_{m-1},Y)
\end{displaymath}
with the last morphism an epimorphism since 
$\mathbf{C}''_{m-1}\in \mathcal{K}={}^{\perp}\mathcal{S}$. 
Applying the functor $Ch(\mathcal{A})(-,S^{n}(\underline{Y}))$ 
to the exact sequence (\ref{eq:10.1.3}) we obtain by 
adjunction the exact sequence 
\begin{equation}\label{eq:10.1.3bis}
\mathcal{A}(\mathbf{C}''_{n}/B_{n}(\mathbf{C}''),\underline{Y})\rightarrowtail 
\mathcal{A}(\mathbf{C}_{n}/B_{n}(\mathbf{C}),\underline{Y})\rightarrow 
\mathcal{A}(\mathbf{C}'_{n}/B_{n}(\mathbf{C}'),\underline{Y})
\end{equation}
Let $f:\mathbf{C}'_{n}/B_{n}(\mathbf{C}')\to \underline{Y}$. 
Form the commutative solid arrow diagram
\begin{displaymath}
\xymatrix{
{\mathbf{C}'_{n}/B_{n}(\mathbf{C}')} \ar@{>->}[r] \ar[d]_{f}& 
{\mathbf{C}_{n}/B_{n}(\mathbf{C})}\ar@{->>}[r] \ar[d] \ar@{..>}[dl] & 
{\mathbf{C}''_{n}/B_{n}(\mathbf{C}'')} \ar@{>->}[r] & \mathbf{C}''_{n-1}\ar@{=}[d]\\
{\underline{Y}} \ar@{>->}[r] & {PO}\ar[rr] & &{\mathbf{C}''_{n-1}}
}
\end{displaymath}
where $PO$ means pushout and the two monomorphisms in the top row 
are so since $\mathbf{C}''$ is exact. By assumption the map 
$\underline{Y}\to PO$ has a retraction, therefore there is a dotted arrow 
that makes commutative the upper triangular diagram containing $f$. 
This implies that the right arrow of the exact sequence (\ref{eq:10.1.3bis}) 
is an epimorphism.
\paragraph{Step 2} We show that ${}^{\perp}\mathcal{S}'\subset ex(\mathcal{K})$. 
Let $\mathbf{C}\in{}^{\perp}\mathcal{S}'$ and $Y\rightarrowtail A\twoheadrightarrow B$ 
be an exact sequence with $Y\in \mathcal{S}$. Applying the functor 
$Ch(\mathcal{A})(\mathbf{C},-)$ to the exact sequence 
$D^{m}(Y)\rightarrowtail D^{m}(A)\twoheadrightarrow D^{m}(B)$ 
and using adjunction we obtain the exact sequence  
\begin{displaymath}
\mathcal{A}(\mathbf{C}_{m-1},Y)\rightarrowtail \mathcal{A}(\mathbf{C}_{m-1},A)
\twoheadrightarrow \mathcal{A}(\mathbf{C}_{m-1},B)
\end{displaymath}
This implies that $\mathbf{C}_{m-1}\in {}^{\perp}\mathcal{S}=\mathcal{K}$, 
hence $\mathbf{C}$ is degreewise in $\mathcal{K}$. Let now $n$ be an integer. 
We can find a monomorphism $\mathbf{C}_{n}/B_{n}(\mathbf{C})\to\underline{Y}$ 
with $\underline{Y}\in\underline{\mathcal{K}^{\perp}}$.
Applying the functor $Ch(\mathcal{A})(\mathbf{C},-)$ to the exact sequence 
$S^{n-1}(\underline{Y})\rightarrowtail D^{n}(\underline{Y})
\twoheadrightarrow S^{n}(\underline{Y})$ 
and using adjunction we obtain the exact sequence  
\begin{displaymath}
\mathcal{A}(\mathbf{C}_{n-1}/B_{n-1}(\mathbf{C}),\underline{Y})
\rightarrowtail \mathcal{A}(\mathbf{C}_{n-1},\underline{Y})\twoheadrightarrow 
\mathcal{A}(\mathbf{C}_{n}/B_{n}(\mathbf{C}),\underline{Y})
\end{displaymath}
This implies that there is a dotted arrow in the diagram
\begin{displaymath}
\xymatrix{
{\mathbf{C}_{n}/B_{n}(\mathbf{C})} \ar[r] \ar@{>->}[d]& 
{\mathbf{C}_{n-1}}\ar[r] \ar@{..>}[dl] & 
{\mathbf{C}_{n-1}/B_{n-1}(\mathbf{C})} \\
{\underline{Y}}
}
\end{displaymath}
that makes the triangular diagram commute. Therefore $\mathbf{C}$ is exact.
The description of $ex(\mathcal{K})^{\perp}$ follows from \ref{sec:8.2.6}(2).
\end{proof}

The previous result was proved in \cite{Gi2} under the assumptions 
(1) $\mathsf{Ext}^{n}(X,Y)=0$ for all $n\geqslant 1$ and all 
$X\in \mathcal{K},Y\in\mathcal{K}^{\perp}$ and (2) $\mathcal{K}^{\perp}$ 
contains a cogenerator of finite injective dimension. 

\subsection{} \label{sec:10.2}
Let $\mathcal{A}$ be a complete and cocomplete abelian category.
We consider the natural abelian THC-situation (\ref{sec:2.2.1}) 
associated to $\mathcal{A}$ and the standard THC-situation 
\[
T:Ch(\mathrm{Ab})\times Ch(\mathcal{A})\to Ch(\mathcal{A})
\]
\[
\mathrm{Homgr}_{\mathcal{A}}:Ch(\mathcal{A})^{op}
\times Ch(\mathcal{A})\to Ch(\mathrm{Ab})
\]
\[
C:Ch(\mathrm{Ab})^{op}\times Ch(\mathcal{A})\to Ch(\mathcal{A})
\]
associated to the natural one (\ref{sec:6.1}). Let 
$(\mathcal{K},\mathcal{K}^{\perp})$ be a cotorsion theory 
in $\mathcal{A}$. We consider the cotorsion
theories $(Ch(\mathcal{K}),Ch(\mathcal{K})^{\perp})$ and 
$(^{\perp}Ch(\mathcal{K}^{\perp}),Ch(\mathcal{K}^{\perp}))$ 
in $Ch(\mathcal{A})$ (\ref{sec:4.1.2}). Let $\mathfrak{K}_{2}$ 
be a class that consists of right bounded complexes that are 
degreewise in $\mathcal{K}$.

\subsubsection{} \label{sec:10.2.1} 
We define
$(\mathfrak{K}_{2};\mathcal{K}^{\perp})=Ch(\mathcal{K}^{\perp})
\cap\mathfrak{K}_{2}^{\perp_{h}}$; in the notations \ref{sec:2.2.4} 
we have
\[
(\mathfrak{K}_{2};\mathcal{K}^{\perp})=\mathrm{Homgr}_{\mathcal{A}}
(\mathfrak{K}_{2},-)^{-1}(ex(\mathrm{Ab}))
\cap Ch(\mathcal{K}^{\perp}) 
\]
If $\mathcal{A}$ has enough $\mathcal{K}^{\perp}$ objects then
$(^{\perp}(\mathfrak{K}_{2};\mathcal{K}^{\perp}),(\mathfrak{K}_{2};
\mathcal{K}^{\perp}))$ is a cotorsion theory and 
$^{\perp}(\mathfrak{K}_{2};\mathcal{K}^{\perp})=
Ch(\mathcal{K})\cap ^{\perp_{h}}(\mathfrak{K}_{2};
\mathcal{K}^{\perp})$.
\begin{proof}
We shall use \ref{sec:2.2.4}. We take $\mathcal{K}_{2}=
(\mathfrak{K}_{2};\mathcal{K}^{\perp}), \mathcal{Y}_{3}=Ch(\mathcal{K}^{\perp}),
\mathcal{Y}_{1}=ex(\mathrm{Ab}), \mathcal{S}_{3}=\{D^{m}(X)\}$, where
$m$ is an integer and $X\in\mathcal{K}$, $\mathcal{S}_{1}=S^{n}(\mathbb{Z})$,
where $n$ is an integer. Let $\mathbf{C}'\in(\mathfrak{K}_{2}$.
By \ref{sec:8.1.1}(1bis) the functor $\mathrm{Homgr}_{\mathcal{A}}
(\mathbf{C}',-)$ preserves the epimorphisms with kernel in $Ch(\mathcal{K}^{\perp})$.
We claim that $Ch(\mathcal{A})$ has enough $(\mathfrak{K}_{2};\mathcal{K}^{\perp})$
objects. Let $\mathbf{C}$ be a complex in $\mathcal{A}$. By \ref{sec:4.1.2}
there is a monomorphism $\mathbf{C}\to\mathbf{Y}$ with 
$\mathbf{Y}\in ex[\mathcal{K}^{\perp}]$ and by \ref{sec:4.2.3} 
the complex $\mathrm{Homgr}_{\mathcal{A}}(\mathbf{C}',\mathbf{Y})$
is exact, therefore $\mathbf{Y}\in (\mathfrak{K}_{2};\mathcal{K}^{\perp})$.
The description of $^{\perp}(\mathfrak{K}_{2};\mathcal{K}^{\perp})$ 
follows from \ref{sec:8.2.6}(2).
\end{proof}

\subsubsection{} \label{sec:10.2.2} 
Suppose that $\mathcal{A}$ has enough $\mathcal{K}^{\perp}$ 
objects and enough $\underline{\mathcal{K}}$ objects, the latter
meaning that for each $A\in Ob(\mathcal{A})$ there is an 
epimorphism $X\to A$ with $X\in \underline{\mathcal{K}}$. Then
$(^{\perp}((\mathfrak{K}_{2};\mathcal{K}^{\perp})\cap ex(\mathcal{K})),
(\mathfrak{K}_{2};\mathcal{K}^{\perp})\cap ex(\mathcal{K}))$
is a cotorsion theory and 
$^{\perp}((\mathfrak{K}_{2};\mathcal{K}^{\perp})\cap ex(\mathcal{K}))=
Ch(\mathcal{K})\cap ^{\perp_{h}}((\mathfrak{K}_{2};\mathcal{K}^{\perp})
\cap ex(\mathcal{K}))$.
\begin{proof}
This follows from \ref{sec:1.2.3} applied to the cotorsion
theories $(^{\perp}(\mathfrak{K}_{2};\mathcal{K}^{\perp}),
(\mathfrak{K}_{2};\mathcal{K}^{\perp}))$
and $(^{\perp}ex(\mathcal{K}^{\perp}),ex(\mathcal{K}^{\perp}))$ 
(\ref{sec:10.1.3}).
\end{proof}

\subsection{} \label{sec:10.3}
We recall \cite{AF} that complex in $\mathcal{A}$ is said to be dg-injective  
if it is degreewise injective and it belongs to $ex(\mathcal{A})^{\perp_{h}}$.
Let
\begin{displaymath}
T:\mathcal{A}'\times \mathcal{A}\to \mathcal{A}'', \
H:\mathcal{A}^{op}\times \mathcal{A}''\to \mathcal{A}', \
C:\mathcal{A}'^{op}\times \mathcal{A}''\to \mathcal{A}
\end{displaymath}
be an abelian THC-situation with $\mathcal{A}',\mathcal{A}$ 
complete and $\mathcal{A}''$ satisfying AB4. We consider the 
standard THC-situation associated to it (\ref{sec:6.1}). 
We recall that $\mathcal{F}$ is the class of flat objects of 
$\mathcal{A}$ (\ref{sec:3.1}) and $\mathfrak{F}_{h}$ the 
class of h-flat complexes in $\mathcal{A}$ 
(\ref{sec:8.3.2}). We define $\mathfrak{F}_{dg}=
\mathfrak{F}_{h}\cap Ch(\mathcal{F})$; the objects of 
$\mathfrak{F}_{dg}$ are called dg-flat complexes \cite{AF}.

\subsubsection{} \label{sec:10.3.1}

(1) \cite{AF} For a complex $\mathbf{P}$ in $\mathcal{A}$ 
the following are equivalent:

(a) $\mathbf{P}$ is dg-flat;

(b) $H(\mathbf{P},-)$ preserves dg-injective complexes;

(c) $T(-,\mathbf{P})$ sends monomorphisms with exact 
cokernel to monomorphisms with exact cokernel.

(2) The pair $(\mathfrak{F}_{dg},\mathfrak{F}_{dg}^{\perp})$ 
is a left exact cotorsion theory and 
$\mathfrak{F}_{dg}^{\perp}=Ch(\mathcal{F}^{\perp})
\cap\mathfrak{F}_{dg}^{\perp_{h}}$.

(3) For every exact complex $\mathbf{C}'$ in $\mathcal{A}'$ 
the functor $T(\mathbf{C}',-)$ sends monomorphisms with 
dg-flat cokernel to monomorphisms with exact cokernel.

\begin{proof}
(1) Let $\mathbf{C}'\in ex(\mathcal{A}')$. (a)$\Rightarrow$(b)
Let $\mathbf{C}''$ be a dg-injective complex in $\mathcal{A}''$. 
The object $H(\mathbf{C}'',\mathbf{P})_{n}=\underset{p}\prod 
H(\mathbf{P}_{p},\mathbf{C}''_{p+n})$ is injective by 
(the proof of) \ref{sec:3.1.2} and \ref{sec:1.2.1}(6).
From \ref{sec:6.1.9} we have 
\[
\mathrm{Homgr}_{\mathcal{A}'}(\mathbf{C}',H(\mathbf{P},\mathbf{C}''))\cong
\mathrm{Homgr}_{\mathcal{A}''}(T(\mathbf{C}',\mathbf{P}),\mathbf{C}'')
\]
hence $H(\mathbf{P},\mathbf{C}'')\in ex(\mathcal{A}')^{\perp_{h}}$.
(b)$\Rightarrow$(a) Let $E''$ be an injective object of $\mathcal{A}''$. 
By \ref{sec:6.1.9} and \ref{sec:4.2} we have 
\[
\mathrm{Homgr}_{\mathcal{A}'}(\mathbf{C}',H(\mathbf{P},S^{0}(E'')))\cong
\mathrm{Homgr}_{\mathcal{A}''}(T(\mathbf{C}',\mathbf{P}),S^{0}(E''))=
\mathcal{A}''(T(\mathbf{C}',\mathbf{P})_{0-},E'')
\]
Since $S^{0}(E'')$ is dg-injective we obtain that 
$\mathcal{A}''(T(\mathbf{C}',\mathbf{P})_{0-},E'')$ is exact.
Since $\mathcal{A}''$ has enough injective objects it follows 
(\ref{sec:1.1.2}(1)) that $\mathbf{P}$ is h-flat. Let now
$A'\rightarrowtail B'\twoheadrightarrow C'$ be an 
exact sequence in $\mathcal{A}'$. Applying 
$\mathrm{Homgr}_{\mathcal{A}'}(-,H(\mathbf{P},S^{0}(E'')))$
to the exact sequence $S^{0}(A')\rightarrowtail S^{0}(B')
\twoheadrightarrow S^{0}(C')$ we obtain from \ref{sec:8.1.1}(1) 
and \ref{sec:4.2} the exact sequence
\[
\mathcal{A}'(C',H(\mathbf{P}_{-},E''))\rightarrowtail
\mathcal{A}'(B',H(\mathbf{P}_{-},E''))\twoheadrightarrow 
\mathcal{A}'(A',H(\mathbf{P}_{-},E''))
\]
where $H(\mathbf{P}_{-},E'')_{n}=H(\mathbf{P}_{-n},E'')$ has 
differential $H((-1)^{n+1}d_{-n+1},E'')$, $d$ being the differential of 
$\mathbf{P}$. Using adjunction and \ref{sec:1.1.2}(1) we then have 
that $\mathbf{P}_{n}\in\mathcal{F}$ for all integers $n$.
(a)$\Rightarrow$(c) Let $\mathbf{A}'\rightarrowtail\mathbf{B}'
\twoheadrightarrow\mathbf{C}'$ be an exact sequence in 
$Ch(\mathcal{A}')$ with $\mathbf{C}'$ exact. Then
$T(\mathbf{A}',\mathbf{P})\rightarrowtail T(\mathbf{B}',\mathbf{P})
\twoheadrightarrow T(\mathbf{C}',\mathbf{P})$ is an exact sequence
as required. (c)$\Rightarrow$(a) Applying $T(-,\mathbf{P})$ to the exact 
sequence $0\to\mathbf{C}'=\mathbf{C}'$ we obtain that 
$\mathbf{P}$ is h-flat. Let now $A'\rightarrowtail B'
\twoheadrightarrow C'$ be an exact sequence in 
$\mathcal{A}'$; then we have the exact sequence
$T(D^{1}(A'),\mathbf{P})\rightarrowtail T(D^{1}(B'),\mathbf{P})
\twoheadrightarrow T(D^{1}(C'),\mathbf{P})$. Since 
$T(D^{1}(A'),\mathbf{P})_{n}=T(A',\mathbf{P}_{n-1})\oplus
T(A',\mathbf{P}_{n})$ for all integers $n$ we obtain that 
$T(A',\mathbf{P}_{n})\to T(B',\mathbf{P}_{n})$ is a 
monomorphism.

(2) We shall use \ref{sec:2.2.2}. We take $\mathcal{K}_{1}=ex(\mathcal{A}'),
\mathcal{X}_{2}=Ch(\mathcal{F}), \mathcal{X}_{3}=ex(\mathcal{A}''),
\mathcal{S}_{2}=\{D^{m}(Y)\}$, where $m$ is an integer and 
$Y\in\mathcal{F}^{\perp}$, $\mathcal{S}_{3}=\{S^{n}(E'')\}$, where $n$ 
is an integer and $E''$ is an injective object of $\mathcal{A}''$.
Let $\mathbf{C}'\in ex(\mathcal{A}')$ and 
$\mathbf{A}\rightarrowtail\mathbf{B}\twoheadrightarrow\mathbf{P}$ 
be an exact sequence in $Ch(\mathcal{A})$ with 
$\mathbf{P}\in Ch(\mathcal{F})$. By \ref{sec:3.1.6}(1) we have the exact 
sequence $T(\mathbf{C}',\mathbf{A})\rightarrowtail T(\mathbf{C}',\mathbf{B})
\twoheadrightarrow T(\mathbf{C}',\mathbf{P})$. We now claim that 
$Ch(\mathcal{A})$ has enough $\mathfrak{F}_{dg}$ objects.
Let $\mathbf{C}$ be a complex in $\mathcal{A}$. By \ref{sec:4.1.2}
there is an epimorphism $\mathbf{P}\to\mathbf{C}$ with 
$\mathbf{P}\in ex[\mathcal{F}]$. By \ref{sec:4.2.6} and \ref{sec:8.3.3}(3)
we have $ex[\mathcal{F}]\subset ^{\perp_{h}}ex[\mathcal{F}^{\perp}]
\subset \mathfrak{F}_{h}$, hence the claim. Left exactness of 
$\mathfrak{F}_{dg}$ follows from \ref{sec:2.2.3} and \ref{sec:3.1.5}.
The description of $\mathfrak{F}_{dg}^{\perp}$ follows from \ref{sec:8.2.6}(2).

(3) Let $\mathbf{A}\rightarrowtail \mathbf{B}\twoheadrightarrow\mathbf{P}$
be an exact sequence in $Ch(\mathcal{A})$ with $\mathbf{P}$ dg-flat.
By \ref{sec:3.1.6}(1) we have the exact sequence 
$T(\mathbf{C}',\mathbf{A})\rightarrowtail T(\mathbf{C}',\mathbf{B})
\twoheadrightarrow T(\mathbf{C}',\mathbf{P})$ with 
$T(\mathbf{C}',\mathbf{P})$ exact.
\end{proof}

\subsubsection{} \label{sec:10.3.2}
We define (\ref{sec:2.2.2})
\[
tex(\mathcal{F})=T(S^{0}(\mathcal{I}'),-)^{-1}(ex(\mathcal{A}''))
\cap ex(\mathcal{F}) 
\]
where $\mathcal{I}'$ is the class of injective objects of 
$\mathcal{A}'$; the objects of $tex(\mathcal{F})$ 
are often called $F$-totally acyclic complexes. If $\mathcal{A}$ 
has enough $\underline{\mathcal{F}^{\perp}}$ objects then 
$(tex(\mathcal{F}),tex(\mathcal{F})^{\perp})$ is a left exact 
cotorsion theory and $tex(\mathcal{F})^{\perp}=
Ch(\mathcal{F}^{\perp})\cap tex(\mathcal{F})^{\perp_{h}}$.
\begin{proof}
We shall use \ref{sec:2.2.2}. We take $\mathcal{K}_{1}=\{S^{0}(E')\}$,
where $E'$ is an injective object of $\mathcal{A}'$, $\mathcal{X}_{2}
=ex(\mathcal{F})$, $\mathcal{X}_{3}=ex(\mathcal{A}''),
\mathcal{S}_{2}=\{D^{m}(C(A',E'')),S^{n}(\underline{Y'})\}$,
where $m,n$ are integers, $A'\in Ob(\mathcal{A}')$, $E''$ is an 
injective object of $\mathcal{A}''$, $\underline{Y}\in
\underline{\mathcal{F}^{\perp}}$, $\mathcal{S}_{3}=\{S^{p}(I'')\}$,
where $p$ is an integer and $I_{3}$ is an injective object of $\mathcal{A}''$.
Let $E'$ be an injective object of $\mathcal{A}'$ and 
$\mathbf{A}\rightarrowtail\mathbf{B}\twoheadrightarrow\mathbf{P}$ 
an exact sequence in $Ch(\mathcal{A})$ with $\mathbf{P}\in ex(\mathcal{F})$. 
By \ref{sec:3.1.6}(1) we have the exact sequence 
$T(S^{0}(E'),\mathbf{A})\rightarrowtail T(S^{0}(E'),\mathbf{B})
\twoheadrightarrow T(S^{0}(E'),\mathbf{P})$. We now claim that 
$Ch(\mathcal{A})$ has enough $tex(\mathcal{F})$ objects.
Let $\mathbf{C}$ be a complex in $\mathcal{A}$. By \ref{sec:4.1.2}
there is an epimorphism $\mathbf{P}\to\mathbf{C}$ with 
$\mathbf{P}\in ex[\mathcal{F}]$, therefore it suffices to show 
that $\mathbf{P}\in tex(\mathcal{F})$. We have $T(S^{0}(E'),\mathbf{P})
=T(E',\mathbf{P}_{-0})$, where $T(E',\mathbf{P}_{-0})_{n}
=T(E',\mathbf{P}_{n})$ has differential $T(E',d)$, $d$ being the 
differential of $\mathbf{P}$. But then $T(E',\mathbf{P}_{-0})$ is exact
by \ref{sec:3.1.6}(1). Left exactness of $tex(\mathcal{F})$ 
follows from \ref{sec:2.2.3} and \ref{sec:3.1.5}. The description of 
$tex(\mathcal{F})^{\perp}$ follows from \ref{sec:8.2.6}(2).
\end{proof}
%Recall \cite{EGR} that a complex $\mathbf{C}$ in $\mathcal{A}$ is called dg-cotorsion 
%if it is exact and $Z(\mathbf{C})\in Ch(\mathcal{F})^{\perp}$.

\section{Cotorsion theories and corner morphisms for categories
of complexes}

\subsection{} \label{sec:11.1}
Let $\mathcal{A}$ be an abelian category and let 
$\mathbf{A}\overset{u}\rightarrowtail\mathbf{B}
\overset{\pi}\twoheadrightarrow\mathbf{C}$,
$\mathbf{K}\overset{i}\rightarrowtail\mathbf{E}
\overset{v}\twoheadrightarrow\mathbf{F}$ be
exact sequences in $Ch(\mathcal{A})$. We denote 
by $uHv$ the natural morphism
\begin{equation}\label{eq:11.1}
\mathrm{Homgr}_{\mathcal{A}}(\mathbf{B},\mathbf{E})\to
\mathrm{Homgr}_{\mathcal{A}}(\mathbf{A},\mathbf{E})
\times_{\mathrm{Homgr}_{\mathcal{A}}(\mathbf{A},\mathbf{F})}
\mathrm{Homgr}_{\mathcal{A}}(\mathbf{B},\mathbf{F})
\end{equation}
Let $\mathcal{K}$ be a class of objects of $\mathcal{A}$ and assume
$\mathbf{C}\in dg(\mathcal{K}), \mathbf{K}\in dg(\mathcal{K}^{\perp})$ 
(\ref{sec:10.1}). Then we have an exact sequence 
\[
\mathrm{Homgr}_{\mathcal{A}}(\mathbf{C},\mathbf{K})\overset{\bar{\pi}\bar{i}}
\rightarrowtail\mathrm{Homgr}_{\mathcal{A}}(\mathbf{B},\mathbf{E})\overset{uHv}
\twoheadrightarrow\mathrm{Homgr}_{\mathcal{A}}(\mathbf{A},\mathbf{E})
\times_{\mathrm{Homgr}_{\mathcal{A}}(\mathbf{A},\mathbf{F})}
\mathrm{Homgr}_{\mathcal{A}}(\mathbf{B},\mathbf{F})
\]
If, moreover, $\mathbf{C}\in ex[\mathcal{K}]$ or $\mathbf{K}\in ex[\mathcal{K}^{\perp}]$
then $uHv$ is a quasi-isomorphism.
\begin{proof}
Consider the commutative diagram
\begin{displaymath}
\xymatrix{
{\mathrm{Homgr}_{\mathcal{A}}(\mathbf{C},\mathbf{K})}\ar@{>->}[rr]^{\bar{i}}
\ar@{>->}[d]_{\bar{\pi}}&&{\mathrm{Homgr}_{\mathcal{A}}(\mathbf{C},\mathbf{E})}
\ar@{>->}[d]^{\bar{\pi}}\\
{\mathrm{Homgr}_{\mathcal{A}}(\mathbf{B},\mathbf{K})}\ar@{>->}[rr]^{\bar{i}}&&
{\mathrm{Homgr}_{\mathcal{A}}(\mathbf{B},\mathbf{E})}\ar@{->>}[rr] \ar@{->>}[d]
&&{\mathrm{Homgr}_{\mathcal{A}}(\mathbf{B},\mathbf{F})} \ar@{->>}[d]\\
& &{\mathrm{Homgr}_{\mathcal{A}}(\mathbf{A},\mathbf{E})} \ar@{->>}[rr]
&&{\mathrm{Homgr}_{\mathcal{A}}(\mathbf{A},\mathbf{F})}
}
\end{displaymath}
There is a natural morphism $\alpha:\mathrm{Homgr}_{\mathcal{A}}
(\mathbf{C},\mathbf{K})\to\K uHv$ such that the triangular diagram
\begin{displaymath}
\xymatrix{
{\K uHv}\ar@{>->}[r]^{k}&{\mathrm{Homgr}_{\mathcal{A}}(\mathbf{B},\mathbf{E})}
\ar[r]&{\mathrm{Homgr}_{\mathcal{A}}(\mathbf{A},\mathbf{E})
\times_{\mathrm{Homgr}_{\mathcal{A}}(\mathbf{A},\mathbf{F})}
\mathrm{Homgr}_{\mathcal{A}}(\mathbf{C},\mathbf{K})}\\
{\mathrm{Homgr}_{\mathcal{A}}(\mathbf{C},\mathbf{K})}
\ar[u]_{\alpha}\ar@{>->}[ur]_{\bar{\pi}\bar{i}}
}
\end{displaymath}
commutes.
\paragraph{Step 1} We show that $\alpha$ is an epimorphism. Let $n,p$
be integers and $f_{p}:\mathbf{B}_{p}\to\mathbf{E}_{p+n}$ a morphism
in $\mathcal{A}$ such that $f_{p}u_{p}=0$ and $v_{p+n}f_{p}=0$.
Then there are morphisms $g_{p}:\mathbf{C}_{p}\to\mathbf{E}_{p+n}$
and $h_{p}:\mathbf{B}_{p}\to\mathbf{K}_{p+n}$ such that $g_{p}\pi_{p}=f_{p}$
and $i_{p+n}h_{p}=f_{p}$. We obtain the commutative diagram
\begin{displaymath}
\xymatrix{
{\mathbf{B}_{p}}\ar[r]^{h_{p}}\ar@{->>}[d]_{\pi_{p}}&{\mathbf{K}_{p+n}} 
\ar@{>->}[d]^{i_{p+n}}\\
{\mathbf{C}_{p}}\ar[r]^{g_{p}}&{\mathbf{E}_{p+n}}
}
\end{displaymath} 
This diagram has a diagonal filler $w_{p}:\mathbf{C}_{p}\to\mathbf{K}_{p+n}$,
therefore $\alpha_{n}((w_{p}))=(f_{p})$.
\paragraph{Step 2} We show that $uHv$, defined as 
$uHv((f_{p}))=(f_{p}u_{p},v_{p+n}f_{p})$, is an epimorphism.  
Let $n,p$ be integers and 
$f_{p}:\mathbf{A}_{p}\to\mathbf{E}_{p+n},g_{p}:\mathbf{B}_{p}
\to\mathbf{F}_{p+n}$ morphisms in $\mathcal{A}$ such that 
$g_{p}u_{p}=v_{p+n}f_{p}$. By \ref{sec:1.2.1}(4) the 
commutative diagram
\begin{displaymath}
\xymatrix{
&{\mathbf{K}_{p+n}}\ar@{>->}[d]\\
{\mathbf{A}_{p}}\ar[r]^{f_{p}}\ar@{>->}[d]_{u_{p}}&
{\mathbf{E}_{p+n}} \ar@{->>}[d]^{v_{p+n}}\\
{\mathbf{B}_{p}}\ar@{->>}[d] \ar[r]^{g_{p}}&{\mathbf{F}_{p+n}}\\
{\mathbf{C}_{p}}
}
\end{displaymath}
has a diagonal filler $h_{p}:\mathbf{B}_{p}\to\mathbf{E}_{p+n}$.
But this means that $(uHv)_{n}((h_{p}))=((f_{p}),(g_{p}))$.
\paragraph{Step 3} If  $\mathbf{C}\in ex[\mathcal{K}]$ or 
$\mathbf{K}\in ex[\mathcal{K}^{\perp}]$ then 
$\mathrm{Homgr}_{\mathcal{A}}(\mathbf{C},\mathbf{K})$ is 
exact, hence $uHv$ is a quasi-isomorphism.
\end{proof}

\subsection{} \label{sec:11.2}
Let 
\[
T:\mathcal{A}'\times \mathcal{A}\to \mathcal{A}'', \
H:\mathcal{A}^{op}\times \mathcal{A}''\to \mathcal{A}', \
C:\mathcal{A}'^{op}\times \mathcal{A}''\to \mathcal{A}
\]
be an abelian THC--situation. Assume that $\mathcal{A}''$ 
satisfies AB4 and has enough injective objects. Consider 
the standard THC--situation (\ref{sec:6.1})
\[
T:Ch(\mathcal{A}')\times Ch(\mathcal{A})\to Ch(\mathcal{A}'')
\]
\[
H:Ch(\mathcal{A})^{op}\times Ch(\mathcal{A}'')\to Ch(\mathcal{A}')
\]
\[
C:Ch(\mathcal{A}')^{op}\times Ch(\mathcal{A}'')\to Ch(\mathcal{A})
\]
associated to the given one. Let $\mathbf{A}'\overset{i'}\rightarrowtail
\mathbf{B}'\twoheadrightarrow\mathbf{C}'$ and $\mathbf{A}\overset{j}
\rightarrowtail\mathbf{B}\twoheadrightarrow\mathbf{C}$ be exact sequences
$Ch(\mathcal{A}')$ and $Ch(\mathcal{A})$, respectively. We denote by $i'Tj$ 
the natural morphism
%\begin{equation}\label{eq:11.2}
\[
T(\mathbf{A}',\mathbf{B})\bigcup_{T(\mathbf{A}',\mathbf{A})} 
T(\mathbf{B}',\mathbf{A})\to T(\mathbf{B}',\mathbf{B})
\]
%\end{equation}

\subsubsection{} \label{sec:11.2.1}
If $\mathbf{C}'\in dg(\mathcal{F}')$ and $\mathbf{C}\in dg(\mathcal{F})$
then $i'Tj$ is a monomorphism that is a quasi-isomorphism 
if either $\mathbf{C}'\in ex[\mathcal{F}']$ or $\mathbf{C}\in ex[\mathcal{F}]$.
\begin{proof}
By \ref{sec:6.1.4} and \ref{sec:3.1.7} the morphism $i'Tj$ is a monomorphism.
By \ref{sec:7.3.5} we have the exact sequence
\[
T(\mathbf{A}',\mathbf{B})\bigcup_{T(\mathbf{A}',\mathbf{A})} 
T(\mathbf{B}',\mathbf{A})\overset{i'Tj}\rightarrowtail 
T(\mathbf{B}',\mathbf{B})\twoheadrightarrow T(\mathbf{C}',\mathbf{C})
\]
By \ref{sec:8.3.3}(3) we have $dg(\mathcal{F})\subset\mathfrak{F}_{dg}$.
Suppose $\mathbf{C}'\in ex[\mathcal{F}']$. By \ref{sec:10.3.1}(3) applied
to $\mathbf{C}'$ and the exact sequence $\mathbf{A}
\rightarrowtail\mathbf{B}\twoheadrightarrow\mathbf{C}$ we obtain that
$T(\mathbf{C}',\mathbf{C})$ is exact, which is what we wanted.
The case $\mathbf{C}\in dg(\mathcal{F})$ is dealt with similarly.
\end{proof}

\subsubsection{} \label{sec:11.2.2}
Let $\mathbf{K}''\rightarrowtail \mathbf{C}''\overset{p''}
\twoheadrightarrow\mathbf{D}'$ be an exact sequence in $Ch(\mathcal{A}'')$.
We denote by $jHp''$ the natural morphism
%\begin{equation}\label{eq:11.2.2a}
\[
H(\mathbf{B},\mathbf{C}'')\to H(\mathbf{B},\mathbf{D}'')\times_{
H(\mathbf{A},\mathbf{D}'')} H(\mathbf{A},\mathbf{C}'')
\]
%\end{equation}
and by $i'Cp''$ the natural morphism 
%\begin{equation}\label{eq:11.2.2b}
\[
C(\mathbf{B}',\mathbf{C}'')\to C(\mathbf{B}',\mathbf{D}'')\times_{
C(\mathbf{A}',\mathbf{D}'')} C(\mathbf{A}',\mathbf{C}'')
\]
%\end{equation}
(1) If $\mathbf{C}\in dg(\mathcal{F})$ and $\mathbf{K}''$ is 
dg-injective then $jHp''$ is an epimorphism with kernel in 
$dg(\mathcal{F}'^{\perp})$. If, moreover, $\mathbf{C}\in 
ex[\mathcal{F}]$ or $\mathbf{K}''$ is injective 
then $\K jHp'' \in ex[\mathcal{F}']$.

(2) If $\mathbf{C}'\in dg(\mathcal{F}')$ and $\mathbf{K}''$ is dg-injective 
then $i'Cp''$ is an epimorphism with kernel in $dg(\mathcal{F}^{\perp})$. 
If, moreover, $\mathbf{C}'\in ex[\mathcal{F}']$ or $\mathbf{K}''$ is 
injective then $\K i'Cp'' \in ex[\mathcal{F}^{\perp}]$.

\begin{proof}
Both 1 and 2 are equivalent to \ref{sec:11.2.1}.
\end{proof}

\subsubsection{} \label{sec:11.2.3} \cite[Theorem 5.1]{Gi4}
Suppose now that $\mathcal{A}'=\mathcal{A}=\mathcal{A}''$ and for all 
$X,Y,Z\in Ob(\mathcal{A})$ there is a natural isomorphism $T(X,T(Y,Z))
\cong T(T(X,Y),Z)$. Then (\ref{sec:3.3.1}) $(\mathcal{F},\mathcal{F}^{\perp})$
is $T$-closed. The cotorsion theory $(dg(\mathcal{F}),ex[\mathcal{F}^{\perp}])$
is $T$-closed. Moreover, if $\mathbf{C}'\in ex[\mathcal{F}]$ and 
$\mathbf{C}\in dg(\mathcal{F})$ then $T(\mathbf{C}',\mathbf{C})\in 
ex[\mathcal{F}]$.
\begin{proof}
For the first part, by \ref{sec:3.3.2} it suffices to show 
that if $(\mathbf{C}',d')\in dg(\mathcal{F})$ and 
$(\mathbf{C}'',d'')\in ex[\mathcal{F}^{\perp}]$ then 
$C(\mathbf{C}',\mathbf{C}'')\in ex[\mathcal{F}^{\perp}]$.
Let $P$ be a flat object of $\mathcal{A}$.
\paragraph{Step 1} We show that $H(S^{r}(P),\mathbf{C}'')\in 
ex[\mathcal{F}^{\perp}]$ for all integers $r$. We have 
$H(S^{r}(P),\mathbf{C}'')=H(P,\mathbf{C}''_{r+})$, where
$H(P,\mathbf{C}''_{r+})_{n}=H(P,\mathbf{C}''_{r+n})$ has
differential $H(P,d''_{r+n})$. By \ref{sec:2.1.1}(1) the functor
$H(P,-)$ preserves epimorphisms with kernel in 
$\mathcal{F}^{\perp}$, hence $H(S^{r}(P),\mathbf{C}'')$ is exact.
By \ref{sec:2.1.1}(2bis) we have $B_{n}(H(S^{r}(P),\mathbf{C}''))
=H(P,B_{r+n}(\mathbf{C}''))\in\mathcal{F}^{\perp}$.
\paragraph{Step 2} We show that 
$C(\mathbf{C}',\mathbf{C}'')\in ex(\mathcal{F}^{\perp})$.
By \ref{sec:2.1.1}(2bis) and \ref{sec:1.2.1}(6) we have 
$C(\mathbf{C}',\mathbf{C}'')\in Ch(\mathcal{F}^{\perp})$.
Since (\ref{sec:6.1.9})
\[
\mathrm{Homgr}_{\mathcal{A}}(S^{r}(P),C(\mathbf{C}',\mathbf{C}''))\cong 
\mathrm{Homgr}_{\mathcal{A}}(\mathbf{C}',H(S^{r}(P),\mathbf{C}''))
\]
we obtain from Step 1 that the complex 
$\mathcal{A}(P,C(\mathbf{C}',\mathbf{C}'')_{r+})$ (\ref{sec:4.2}) is 
exact, therefore (\ref{sec:1.1.2}(1)) $C(\mathbf{C}',\mathbf{C}'')$ is 
exact.
\paragraph{Step 3} We show that $Z_{n}(C(\mathbf{C}',\mathbf{C}''))
\in\mathcal{F}^{\perp}$ for all integers $n$. Let  
$A\rightarrowtail B\twoheadrightarrow P$ be an exact sequence in
$\mathcal{A}$. Applying the functor 
$\mathrm{Homgr}_{\mathcal{A}}(-,C(\mathbf{C}',\mathbf{C}''))$
to the exact sequence $S^{n}(A)\rightarrowtail S^{n}(B)\twoheadrightarrow 
S^{n}(P)$ we obtain from \ref{sec:8.1.1}(1) the exact sequence 
\[
\mathrm{Homgr}_{\mathcal{A}}(S^{n}(P),C(\mathbf{C}',\mathbf{C}''))
\rightarrowtail\mathrm{Homgr}_{\mathcal{A}}(S^{n}(B),C(\mathbf{C}',\mathbf{C}''))
\twoheadrightarrow\mathrm{Homgr}_{\mathcal{A}}(S^{n}(A),C(\mathbf{C}',\mathbf{C}''))
\]
By Step 2 the complex $\mathrm{Homgr}_{\mathcal{A}}(S^{n}(P),
C(\mathbf{C}',\mathbf{C}''))$
is exact, hence, upon applying the functor $Z_{0}$ to the previous 
exact sequence, we obtain the exact sequence 
\[
Ch(\mathcal{A})(S^{n}(P),C(\mathbf{C}',\mathbf{C}''))\rightarrowtail
Ch(\mathcal{A})(S^{n}(B),C(\mathbf{C}',\mathbf{C}''))\twoheadrightarrow
Ch(\mathcal{A})(S^{n}(A),C(\mathbf{C}',\mathbf{C}''))
\]
By adjunction it follows that $Z_{n}(C(\mathbf{C}',\mathbf{C}''))
\in\mathcal{F}^{\perp}$. Therefore the cotorsion theory
$(dg(\mathcal{F}),ex[\mathcal{F}^{\perp}])$
is $T$-closed.
Suppose now that $(\mathbf{C}',d')\in ex[\mathcal{F}]$ 
and $\mathbf{C}\in dg(\mathcal{F})$. Since 
$dg(\mathcal{F})\subset\mathfrak{F}_{dg}$ we have that
$T(\mathbf{C}',\mathbf{C})$ is exact. Since  
$(\mathcal{F},\mathcal{F}^{\perp})$ is $T$-closed we have 
$T(\mathbf{C}',\mathbf{C})\in  Ch(\mathcal{F})$. Let $n$
be an integer and $Y\rightarrowtail A\twoheadrightarrow B$ 
an exact sequence in $\mathcal{A}$ with $Y\in\mathcal{F}^{\perp}$.
By \ref{sec:8.1.1}(1bis) we have the exact sequence 
\[
\mathrm{Homgr}_{\mathcal{A}}(T(\mathbf{C}',\mathbf{C}),S^{n}(Y))
\rightarrowtail\mathrm{Homgr}_{\mathcal{A}}(T(\mathbf{C}',\mathbf{C}),S^{n}(A))
\twoheadrightarrow\mathrm{Homgr}_{\mathcal{A}}(T(\mathbf{C}',\mathbf{C}),S^{n}(B))
\]
From \ref{sec:6.1.9} we have 
\[
\mathrm{Homgr}_{\mathcal{A}}(T(\mathbf{C}',\mathbf{C}),S^{n}(Y))\cong 
\mathrm{Homgr}_{\mathcal{A}}(\mathbf{C},C(\mathbf{C}',S^{n}(Y)))
\]
We claim that $C(\mathbf{C}',S^{n}(Y))\in ex[\mathcal{F}^{\perp}]$.
One has $C(\mathbf{C}',S^{n}(Y))=C(\mathbf{C}'_{n-},Y)$, where
$C(\mathbf{C}'_{n-},Y)_{m}=C(\mathbf{C}'_{n-m},Y)$ has 
differential $C((-1)^{m+1}d'_{n-m+1},Y)$. By \ref{sec:2.1.1}((2bis) and (1bis))
the functor $C(-,Y)$ preserves monomorphisms with cokernel in
$\mathcal{F}$, therefore $C(\mathbf{C}',S^{n}(Y))$ is exact.
By an argument similar to the proof of \ref{sec:6.2.3} we have 
$B_{m}(C(\mathbf{C}',S^{n}(Y)))=C(B_{n-m}(\mathbf{C}'),Y)$.
By \ref{sec:2.1.1}(2bis) we have $C(B_{n-m}(\mathbf{C}'),Y)
\in\mathcal{F}^{\perp}$, hence the claim is proved. Consequently, the complex
$\mathrm{Homgr}_{\mathcal{A}}(\mathbf{C},C(\mathbf{C}',S^{n}(Y)))$
is exact. Applying the functor $Z_{0}$ to the last displayed exact sequence
we obtain the exact sequence
\[
Ch(\mathcal{A})(T(\mathbf{C}',\mathbf{C}),S^{n}(Y))\rightarrowtail
Ch(\mathcal{A})(T(\mathbf{C}',\mathbf{C}),S^{n}A))\twoheadrightarrow
Ch(\mathcal{A})(T(\mathbf{C}',\mathbf{C}),S^{n}(B))
\]
hence by adjunction the exact sequence
\[
\mathcal{A}(T_{n}/B_{n},Y)\rightarrowtail\mathcal{A}(T_{n}/B_{n},A)
\twoheadrightarrow\mathcal{A}(T_{n}/B_{n},B)
\]
where $T_{n}/B_{n}=T(\mathbf{C}',\mathbf{C})_{n}/B_{n}(T(\mathbf{C}',\mathbf{C}))$.
It follows that $Z_{n-1}(T(\mathbf{C}',\mathbf{C}))\in\mathcal{F}$.
\end{proof}

\section{Some resolutions of complexes}

\subsection{} \label{sec:12.1}
Let $\mathcal{A}$ be an abelian category. Recall (\ref{sec:5.1.6}) 
that a cotorsion theory $(\mathcal{K},\mathcal{K}^{\perp})$ in 
$\mathcal{A}$ is complete if it is left and right complete. This is 
equivalent to $(\mathcal{K},\mathcal{K}^{\perp})$ being left 
complete and $\mathcal{A}$ having enough $\mathcal{K}^{\perp}$ 
objects and to $(\mathcal{K},\mathcal{K}^{\perp})$ 
being right complete and $\mathcal{A}$ having enough 
$\mathcal{K}$ objects.

\subsubsection{} \label{sec:12.1.1}
\cite[Proposition 5.4]{Ho} Let $(\mathcal{K},\mathcal{K}^{\perp})$ be 
a complete cotorsion theory in $\mathcal{A}$. Every morphism of 
$\mathcal{A}$ can be factored into a monomorphism with cokernel in 
$\mathcal{K}$ followed by an epimorphism with kernel
in $\mathcal{K}^{\perp}$.
\begin{proof}
Let $f:A\to B$ be a morphism. Choose a monomorphism $A\to Y$
with $Y\in\mathcal{K}^{\perp}$ and form the pushout diagram
\begin{displaymath}
\xymatrix{
{A}\ar[r]^{f} \ar[d]&{B} \ar[d]\\
{Y} \ar[r]&{PO}
}
\end{displaymath} 
This diagram is exact (\ref{sec:1.2.8}) hence we have an exact
sequence $A\rightarrowtail B\oplus Y\twoheadrightarrow PO$
and $f$ is the composite $A\rightarrowtail B\oplus Y
\twoheadrightarrow B$, where the second arrow is a projection
morphism. We can find an exact sequence 
$Y'\rightarrowtail X\twoheadrightarrow PO$ with $Y'\in\mathcal{K}^{\perp}$
and $X\in\mathcal{K}$. Form then the commutative diagram
\begin{displaymath}
\xymatrix{
&{Y'}\ar@{=}[r]\ar@{>->}[d]&{Y'}\ar@{>->}[d]\\
{A}\ar@{>->}[r] \ar@{=}[d]&
{PB} \ar@{->>}[d]\ar@{->>}[r]\ar[d]&{X}\ar@{->>}[d]\\
{A}\ar@{>->}[r]\ar[dr]_{f}&{B\oplus Y}\ar@{->>}[r]\ar@{->>}[d]&{PO}\\
&{B}
}
\end{displaymath}
where $PB$ means pullback. Since the epimorphisms with kernel in 
$\mathcal{K}^{\perp}$ are closed under composition, we obtain
the desired factorization.
\end{proof}

%If $\mathcal{K}$ is a class of objects of $\mathcal{A}$ such that
%for every $A\in Ob(\mathcal{A})$ there is an epimorphism $X\to A$
%with $X\in\mathcal{K}$, then every morphism of $\mathcal{A}$ can
%be factored into a monomorphism with cokernel in $\mathcal{K}$
%followed by an epimorphism.
\subsubsection{} \label{sec:12.1.2}
\cite[Lemma 2.3]{YD} (1) Let $(\mathcal{K},\mathcal{K}^{\perp})$ 
be a left exact and left complete cotorsion theory 
in $\mathcal{A}$ such that $\mathcal{A}$ has enough 
$\mathcal{K}^{\perp}$ objects. Then $ex(\mathcal{A})=
\Co(ex[\mathcal{K}^{\perp}]\rightarrowtail ex[\mathcal{K}])$.

(1bis) Let $(\mathcal{K},\mathcal{K}^{\perp})$ be a right exact and right
complete cotorsion theory in $\mathcal{A}$ such that $\mathcal{A}$ has enough 
$\mathcal{K}$ objects. Then 
$ex(\mathcal{A})=\K(ex[\mathcal{K}^{\perp}]\twoheadrightarrow ex[\mathcal{K}])$.

\begin{proof}
(1) Let $\mathbf{C}$ be an exact complex and $n$ an integer. We can find an 
exact sequence $B_{n}(\mathbf{Y})\rightarrowtail B_{n}(\mathbf{X})
\twoheadrightarrow B_{n}(\mathbf{C})$ with $B_{n}(\mathbf{Y})\in\mathcal{K}^{\perp},
B_{n}(\mathbf{X})\in\mathcal{K}$. We shall complete the solid arrows diagram
\[
\xymatrix{
{B_{n+1}(\mathbf{Y})}\ar@{>..>}[r] \ar@{>..>}[d]& {\mathbf{Y}_{n+1}}
\ar@{..>>}[r]\ar@{>..>}[d] & {B_{n}(\mathbf{Y})}\ar@{>->}[d]
\ar@{>..>}[r] \ar@{>..>}[d]& {\mathbf{Y}_{n}}
\ar@{..>>}[r]\ar@{>..>}[d] & {B_{n-1}(\mathbf{Y})}\ar@{>..>}[d]\\
{B_{n+1}(\mathbf{X})} \ar@{>..>}[r] \ar@{..>>}[d]& {\mathbf{X}_{n+1}}\ar@{..>>}[r]
\ar@{..>>}[d] & {B_{n}(\mathbf{X})}  \ar@{->>}[d]
\ar@{>..>}[r] \ar@{..>>}[d]& {\mathbf{X}_{n}}\ar@{..>>}[r]
\ar@{..>>}[d] & {B_{n-1}(\mathbf{X})}  \ar@{..>>}[d]\\
{B_{n+1}(\mathbf{C})}\ar@{>->}[r] & {\mathbf{C}_{n+1}}\ar@{->>}[r] & {B_{n}(\mathbf{C})} 
\ar@{>->}[r] & {\mathbf{C}_{n}}\ar@{->>}[r] & {B_{n-1}(\mathbf{C})}
}
\]
with the dotted arrows, first to the left of the vertical sequence and then to the right.
Continuing we obtain the desired exact sequence $\mathbf{Y}\rightarrowtail \mathbf{X}
\twoheadrightarrow\mathbf{C}$.
\paragraph{Step 1} 
Form the commutative diagram
\begin{displaymath}
\xymatrix{
{B_{n+1}(\mathbf{Y})}\ar@{>->}[d]\ar@{=}[r]&{B_{n+1}(\mathbf{Y})}\ar@{>->}[d]\\
{PB'_{n+1}}\ar@{>->}[r] \ar@{->>}[d]& {\mathbf{X}_{n+1}}\ar@{->>}[r]\ar@{->>}[d] 
& {B_{n}(\mathbf{X})}\ar@{=}[d]\\
{B_{n+1}(\mathbf{C})} \ar@{>->}[r] \ar@{=}[d]& {PB_{n+1}}\ar@{->>}[r]\ar@{->>}[d] & 
{B_{n}(\mathbf{X})}\ar@{->>}[d]\\
{B_{n+1}(\mathbf{C})} \ar@{>->}[r] & {\mathbf{C}_{n+1}}\ar@{->>}[r] & {B_{n}(\mathbf{C})}
}
\end{displaymath}
where $PB_{n+1}$ and $PB'_{n+1}$ mean pullback and $B_{n+1}(\mathbf{Y})\rightarrowtail 
\mathbf{X}_{n+1}\twoheadrightarrow PB_{n+1}$ is an exact sequence with
$\mathbf{X}_{n+1}\in\mathcal{K}$ and $B_{n+1}(\mathbf{Y})\in\mathcal{K}^{\perp}$.
We put $B_{n+1}(\mathbf{X})=PB'_{n+1}$ and we have $B_{n+1}(\mathbf{X})\in \mathcal{K}$.
Consider now the commutative solid arrows diagram
\begin{displaymath}
\xymatrix{
{B_{n+1}(\mathbf{Y})}\ar@{>..>}[r] \ar@{>->}[d]& {\mathbf{Y}_{n+1}}\ar@{..>>}[r]\ar@{>->}[d] 
& {B_{n}(\mathbf{Y})}\ar@{>->}[d]\\
{B_{n+1}(\mathbf{X})} \ar@{>->}[r] \ar@{->>}[d]& {\mathbf{X}_{n+1}}\ar@{->>}[r]
\ar@{->>}[d] & {B_{n}(\mathbf{X})}\ar@{->>}[d]\\
{B_{n+1}(\mathbf{C})} \ar@{>->}[r] & {\mathbf{C}_{n+1}}\ar@{->>}[r] & {B_{n}(\mathbf{C})}
}
\end{displaymath}
where $\mathbf{Y}_{n+1}$ is the kernel of the composite morphism 
$\mathbf{X}_{n+1}\to PB_{n+1}\to \mathbf{C}_{n+1}$.
By the universal property of kernel we have the induced dotted arrows 
that make the resulting square diagrams commute; moreover, by the snake 
diagram the dotted sequence of the previous diagram is a short exact sequence.
\paragraph{Step 2} Factor (\ref{sec:12.1.1}) the composite morphism
$B_{n}(\mathbf{X})\to B_{n}(\mathbf{C})\to\mathbf{C}_{n}$ into a
monomorphism $B_{n}(\mathbf{X})\to\mathbf{X}_{n}$ with cokernel
$B_{n-1}(\mathbf{X})\in\mathcal{K}$ followed by an epimorphism
$\mathbf{X}_{n}\to\mathbf{C}_{n}$ with kernel $\mathbf{Y}_{n}\in\mathcal{K}^{\perp}$.
We obtain the commutative solid arrows diagram
\begin{displaymath}
\xymatrix{
{B_{n}(\mathbf{Y})}\ar@{>..>}[r] \ar@{>->}[d]& {\mathbf{Y}_{n}}\ar@{..>>}[r]\ar@{>->}[d] 
& {B_{n-1}(\mathbf{Y})}\ar@{>..>}[d]\\
{B_{n}(\mathbf{X})} \ar@{>->}[r] \ar@{->>}[d]& {\mathbf{X}_{n}}\ar@{->>}[r]
\ar@{->>}[d] & {B_{n-1}(\mathbf{X})}\ar@{..>>}[d]\\
{B_{n}(\mathbf{C})} \ar@{>->}[r] & {\mathbf{C}_{n}}\ar@{->>}[r] & {B_{n-1}(\mathbf{C})}
}
\end{displaymath}
The dotted arrows $B_{n}(\mathbf{Y})\to\mathbf{Y}_{n}$ and 
$B_{n-1}(\mathbf{X})\to B_{n-1}(\mathbf{C})$ are induced by the universal property 
of kernel and cokernel. Letting $B_{n-1}(\mathbf{Y})$ be the kernel of
$B_{n-1}(\mathbf{X})\to B_{n-1}(\mathbf{C})$, we obtain as before that the 
dotted horizontal sequence of the previous diagram is a short exact sequence.
By \ref{sec:1.2.10} we then have that $B_{n-1}(\mathbf{Y})\in\mathcal{K}^{\perp}$.
\end{proof}

\subsubsection{} \label{sec:12.1.3}
\cite[Theorem 3.4]{Sp} %cite[Application 2.4]{BN}
Suppose that $\mathcal{A}$ satisfies AB4. Let $\mathfrak{K},\mathfrak{L}$ 
be classes of complexes in $\mathcal{A}$ such that 

(a) $\mathfrak{L}$ is closed under extensions;

(b) for every right bounded complex $\mathbf{C}$ in $\mathcal{A}$ there are 
a complex $\mathbf{X}\in\mathfrak{L}$ and a quasi-isomorphism
$\mathbf{X}\twoheadrightarrow \mathbf{C}$;

(c) $\mathfrak{L}\subset ^{\perp_{h}}\mathfrak{K}$;

(d) $\mathfrak{K}\subset \mathfrak{L}^{\perp_{e}}$.

Then for every complex $\mathbf{C}$ in $\mathcal{A}$ 
there are $\mathbf{X}\in ^{\perp_{h}}\mathfrak{K}$ and 
a quasi-isomorphism $\mathbf{X}\twoheadrightarrow \mathbf{C}$.
\begin{proof}
Let $k$ be an integer and consider the complex $\tau_{k}(\mathbf{C})$
(\ref{sec:4.2.4}). We can find a complex $\mathbf{X}^{k}\in\mathfrak{L}$ 
and a quasi-isomorphism $v^{k}:\mathbf{X}^{k}\twoheadrightarrow 
\tau_{k}(\mathbf{X})$. Let $u$ be the composite 
$\mathbf{X}^{k}\overset{v^{k}}\to\tau_{k}(\mathbf{X})\to\tau_{k-1}(\mathbf{X})$;
then $u$ is a morphism of right bounded complexes, therefore 
$Con(u)$ and $Cyl(u)$ are right bounded. Consider the diagram
\begin{displaymath}
\xymatrix{
{\mathbf{X}^{k}}\ar@{>->}[r]^{u^{k}} \ar@{=}[d]&
{\mathbf{X}^{k-1}}\ar@{->>}[d]^{v''}\ar@{->>}[r]&
{\mathbf{X}'}\ar@{->>}[d]^{v'}\\
{\mathbf{X}^{k}}\ar@{>->}[r]^{\tilde{u}}\ar[dr]_{u}&{Cyl(u)}
\ar@{->>}[r]^{\tilde{\pi}}\ar@{->>}[d]^{\beta}&{Con(u)}\\
&{ \tau_{k-1}(\mathbf{C})}
}
\end{displaymath}
where $\tilde{u},\tilde{\pi},\beta$ are as in \ref{sec:4.1.3},
$\mathbf{X}'\in\mathfrak{L}$, $v'$ is a quasi-isomorphism and
$\mathbf{X}^{k-1}$ is the pullback of $v'$ along $\tilde{\pi}$.
Then $v''$ is a quasi-isomorphism and $\mathbf{X}^{k-1}\in\mathfrak{L}$.
Since $\tilde{u}$ is degreewise split, so is $u^{k}$. We put 
$v^{k-1}=\beta v''$ and we obtain the commutative diagram
\begin{displaymath}
\xymatrix{
{\mathbf{X}^{k}}\ar@{>->}[r]^{u^{k}} \ar@{->>}[d]_{v^{k}}
&{\mathbf{X}^{k-1}} \ar@{->>}[d]^{v^{k-1}}\\
{\tau_{k}(\mathbf{C})} \ar[r]&{\tau_{k-1}(\mathbf{C})}
}
\end{displaymath}
The morphism $H(\tau_{k}(\mathbf{C}))\to H(\tau_{k-1}(\mathbf{C}))$
is $0\to H_{k}(\mathbf{C})$ in degree $k-1$ and the identity in degrees 
$\geqslant k$ and $\leqslant k-2$, therefore in the commutative diagram
\begin{displaymath}
\xymatrix{
{H(\mathbf{X}^{k})}\ar[r]^{H(u^{k})} \ar[d]_{H(v^{k})}
&{H(\mathbf{X}^{k-1})} \ar[d]^{H(v^{k-1})}\\
{H(\tau_{k}(\mathbf{C}))}\ar[r]&{H(\tau_{k-1}(\mathbf{C}))}
}
\end{displaymath}
the morphism $H(u^{k})$ is degreewise split.
Using the same procedure for $v^{k-1}$ in place of $v^{k}$
we obtain a commutative diagram
\begin{displaymath}
\xymatrix{
{\mathbf{X}^{k-1}}\ar@{>->}[r]^{u^{k-1}} \ar@{->>}[d]_{v^{k-1}}
&{\mathbf{X}^{k-2}} \ar@{->>}[d]^{v^{k-2}}\\
{\tau_{k-1}(\mathbf{C})} \ar[r]&{\tau_{k-2}(\mathbf{C})}
}
\end{displaymath}
and so on. The resulting direct system of right bounded complexes
\[
\mathbf{X}^{k}\overset{u^{k}}\rightarrowtail \mathbf{X}^{k-1}\overset{u^{k-1}}
\rightarrowtail \mathbf{X}^{k-2}\rightarrowtail  ...
\]
has the properties that $u^{i}$ is degreewise split,
$\Co u^{i}\in\mathfrak{L}$ and $H(u^{i})$ is degreewise 
split for all $i\leqslant k$. Let $\mathbf{X}=\underset{i\in(\mathbb{Z}_{\leqslant k})^{op}}
\cl \mathbf{X}^{i}$; we then obtain an epimorphism $\mathbf{X}\to \mathbf{C}$
and a commutative diagram
\begin{displaymath}
\xymatrix{
{\underset{i\in(\mathbb{Z}_{\leqslant k})^{op}}
\cl H(\mathbf{X}^{i})}\ar[rr]\ar[d] &&{H(\mathbf{X})} \ar[d]\\
{\underset{i\in(\mathbb{Z}_{\leqslant k})^{op}}
\cl H(\tau_{i}(\mathbf{C}))} \ar[rr]&&{H(\mathbf{C})}
}
\end{displaymath}
It is standard that the bottom horizontal arrow is an isomorphism. 
The left vertical arrow is an isomorphism since so is $H(v^{i})$. 
By \ref{sec:4.2.7} the top horizontal arrow is an isomorphism, 
therefore $\mathbf{X}\to \mathbf{C}$ is a quasi-isomorphism.
Let $\mathbf{D}\in\mathfrak{K}$; then we have the exact sequence
\[
\mathrm{Homgr}_{\mathcal{A}}(\Co u^{i},\mathbf{D})\rightarrowtail 
\mathrm{Homgr}_{\mathcal{A}}(\mathbf{X}^{i-1},\mathbf{D})
\twoheadrightarrow\mathrm{Homgr}_{\mathcal{A}}(\mathbf{X}^{i},\mathbf{D})
\]
hence the inverse system 
$(\mathrm{Homgr}_{\mathcal{A}}(\mathbf{X}^{i},\mathbf{D}))_{i\geqslant k}$
satisfies the Mittag-Leffler condition. It follows that 
$\mathbf{X}\in ^{\perp_{h}}\mathfrak{K}$.
\end{proof}
We note that \ref{sec:12.1.3} admits a dual formulation
whose proof uses the `mapping path factorization' of a 
morphism of complexes instead of the factorization given
in \ref{sec:4.1.3}.

\subsubsection{} \label{sec:12.1.4}
Suppose that $\mathcal{A}$ satisfies AB4. Let $\mathcal{K}$
be a class of objects of $\mathcal{A}$ such that for every 
$A\in Ob(\mathcal{A})$ there is an epimorphism $X\to A$
with $X\in^{\perp}\mathcal{K}$. For every complex $\mathbf{C}$
in $\mathcal{A}$ ther are $\mathbf{X}\in Ch(^{\perp}\mathcal{K})
\cap ^{\perp_{h}}ex[(^{\perp}\mathcal{K})^{\perp}]$ and a 
quasi-isomorphism $\mathbf{X}\twoheadrightarrow \mathbf{C}$.
In particular, let $(\mathcal{K},\mathcal{K}^{\perp})$ be a 
cotorsion theory in $\mathcal{A}$ such that $\mathcal{A}$ has 
enough $\mathcal{K}$ objects. For every complex $\mathbf{C}$
in $\mathcal{A}$ there are $\mathbf{X}\in dg(\mathcal{K})$ and a 
quasi-isomorphism $\mathbf{X}\twoheadrightarrow \mathbf{C}$.

\begin{proof}
In \ref{sec:12.1.3} we take $\mathfrak{K}=
ex[(^{\perp}\mathcal{K})^{\perp}]$ and $\mathfrak{L}$
to be the class of right bounded complexes that are degreewise
in $^{\perp}\mathcal{K}$. Condition (a) follows from \ref{sec:1.2.1}(5),
condition (b) is standard and follows from assumption, condition (c)
follows from \ref{sec:4.2.3} and condition (d) from \ref{sec:8.1.1}(1)
since $\mathfrak{K}\subset Ch((^{\perp}\mathcal{K})^{\perp})$.
The construction of $\mathbf{X}$ and \ref{sec:1.2.1}(7) show that
$\mathbf{X}\in Ch(^{\perp}\mathcal{K})$.
 \end{proof}

\subsubsection{} \label{sec:12.1.5}
\cite[Section 5 A]{Sp} 
Let 
\[
T:\mathcal{A}'\times \mathcal{A}\to \mathcal{A}'', \
H:\mathcal{A}^{op}\times \mathcal{A}''\to \mathcal{A}', \
C:\mathcal{A}'^{op}\times \mathcal{A}''\to \mathcal{A}
\]
be an abelian THC--situation. Suppose that 
$\mathcal{A}',\mathcal{A}$ satisfy AB4 and have 
enough flat objects and that $\mathcal{A}''$ has enough 
injective objects. Consider the 
standard THC--situation (\ref{sec:6.1})
\[
T:Ch(\mathcal{A}')\times Ch(\mathcal{A})\to Ch(\mathcal{A}'')
\]
\[
H:Ch(\mathcal{A})^{op}\times Ch(\mathcal{A}'')\to Ch(\mathcal{A}')
\]
\[
C:Ch(\mathcal{A}')^{op}\times Ch(\mathcal{A}'')\to Ch(\mathcal{A})
\]
associated to the given one. 

(1) For every complex $\mathbf{C}$ in $\mathcal{A}$ there 
are a dg-flat complex $\mathbf{P}$ and quasi-isomorphism 
$\mathbf{P}\twoheadrightarrow \mathbf{C}$.

(2) For every complex $\mathbf{C}'$ in $\mathcal{A}'$

(a) $T(\mathbf{C}',\mathbf{P})$ is exact whenever $\mathbf{P}\in 
\mathfrak{F}_{h}\cap ex(\mathcal{A})$;

(b) the functor $T(\mathbf{C}',-)$ preserves quasi-isomorphisms 
between h-flat objects.

\begin{proof}
(1) Since (\ref{sec:8.3.3}(3)) $dg(\mathcal{F})\subset\mathfrak{F}_{dg}$,
the assertion follows from \ref{sec:12.1.4} applied to 
$\mathcal{K}=\mathcal{F}$.

(2) (a) By part 1 we can find a dg-flat complex $\mathbf{P}'$
in $\mathcal{A}'$ and a quasi-isomorphism $u':\mathbf{P}'\to\mathbf{C}'$.
The assertion follows then from \ref{sec:8.3.2}.

(b) Let $u:\mathbf{P}_{1}\to\mathbf{P}_{2}$ be a quasi-isomorphism
between h-flat objects in $Ch(\mathcal{A})$ and let 
$u':\mathbf{P}'\to\mathbf{C}'$ be a quasi-isomorphism with 
$\mathbf{P}'$ dg-flat; we have then the commutative diagram
\begin{displaymath}
\xymatrix{
{T(\mathbf{P}',\mathbf{P}_{1})}\ar[r]^{T(\mathbf{P}',u)} 
\ar[d]_{T(u',\mathbf{P}_{1})}&{T(\mathbf{P}',\mathbf{P}_{2})} 
\ar[d]^{T(u',\mathbf{P}_{2})}\\
{T(\mathbf{C}',\mathbf{P}_{1})} \ar[r]^{T(\mathbf{C}',u)}
&{T(\mathbf{C}',\mathbf{P}_{2})}
}
\end{displaymath}
By \ref{sec:8.3.2} the top horizontal and the vertical arrows are 
quasi-isomorphisms, therefore so is the bottom horizontal arrow 
by the two out of three property of quasi-isomorphisms.
\end{proof}

\subsubsection{} \label{sec:12.1.6}
Let $A$ be a left noetherian ring that has finite injective dimension 
as a left $A$-module and $\mathcal{I}$ be the class of injective 
left $A$-modules. For every complex $\mathbf{C}$ of left $A$-modules 
there are $\mathbf{X}\in  ^{\perp_{h}}(Ch(\mathcal{I})\cap\{S^{0}(E),
E\in\mathcal{I}\}^{\perp_{h}})$ and a quasi-isomorphism 
$\mathbf{X}\twoheadrightarrow \mathbf{C}$.
\begin{proof}
Let $_{A}Mod$ be the category of left $A$-modules.
We denote $\mathrm{Homgr}_{_{A}Mod}$ by
$\mathrm{Homgr}_{A}$.
In \ref{sec:12.1.3} we take $\mathfrak{K}=
Ch(\mathcal{I})\cap\{S^{0}(E), E\in\mathcal{I}\}^{\perp_{h}}$ 
and $\mathfrak{L}$ to be the class of right bounded complexes of
left $A$-modules that are degreewise of finite injective dimension.
Condition (a) follows from horseshoe lemma, condition (b) follows 
from assumptions on $A$ and condition (d) from \ref{sec:8.1.1}(1).
We show (c). Let $\mathbf{C}'\in\mathfrak{L}, \mathbf{C}\in\mathfrak{K}$
and let $i$ be an integer such that $\mathbf{C}'_{n}=0$
for $n<i$. By \ref{sec:4.2.3} it suffices to show that 
$_{A}Mod(\mathbf{C}'_{r},\mathbf{C}_{r+})$ and 
$_{A}Mod(\mathbf{C}'_{i},\mathbf{C})$ are exact.
Let 
\[
\xymatrix{
{\mathbf{C}'_{i}}\ar@{>->}[r] & {E_{0}}\ar[r]^{d_{0}} & {E_{1}}\ar[r]^{d_{1}} 
&{...} \ar@{->>}[r]^{d_{n-1}} & {E_{n}}
}
\]
be an exact sequence with $E_{0},E_{1},...,E_{n}$ injective.
By \ref{sec:8.1.1}(1) we have exact sequences
\[
\mathrm{Homgr}_{A}(S^{0}(E_{n}),\mathbf{C})\rightarrowtail 
\mathrm{Homgr}_{A}(S^{0}(E_{n-1}),\mathbf{C})\twoheadrightarrow
\mathrm{Homgr}_{A}(S^{0}(\K d_{n-1}),\mathbf{C}),
\]
\[
\mathrm{Homgr}_{A}(S^{0}(\K d_{n-1}),\mathbf{C})\rightarrowtail 
\mathrm{Homgr}_{A}(S^{0}(E_{n-2}),\mathbf{C})\twoheadrightarrow
\mathrm{Homgr}_{A}(S^{0}(\K d_{n-2}),\mathbf{C}),
\]
\[
..., \mathrm{Homgr}_{A}(S^{0}(\K d_{1}),\mathbf{C})\rightarrowtail 
\mathrm{Homgr}_{A}(S^{0}(E_{0}),\mathbf{C})\twoheadrightarrow
\mathrm{Homgr}_{A}(S^{0}(\mathbf{C}'_{i}),\mathbf{C})
\]
hence the implication ($\mathrm{Homgr}_{A}
(S^{0}(\K d_{i}),\mathbf{C})$ exact) $\Rightarrow$ ($\mathrm{Homgr}_{A}
(S^{0}(\K d_{i-1}),\mathbf{C})$ exact) for $2\leqslant i\leqslant n-1$.
Since $\mathrm{Homgr}_{A}(S^{0}(\K d_{n-1}),\mathbf{C})$
is exact we conclude that $\mathrm{Homgr}_{A}
(S^{0}(\mathbf{C}'_{i}),\mathbf{C})$ is exact. The proof that 
$_{A}Mod(\mathbf{C}'_{r},\mathbf{C}_{r+})$ is exact is similar.
\end{proof}

\subsubsection{} \label{sec:12.1.7}
Suppose that $\mathcal{A}$ satisfies AB4. Let $(\mathcal{K},
\mathcal{K}^{\perp})$ be a right exact and right complete cotorsion 
theory in $\mathcal{A}$ such that $\mathcal{A}$ has enough 
$\mathcal{K}$ objects. Then the cotorsion theory 
$(dg(\mathcal{K}),ex[\mathcal{K}^{\perp}])$ is left complete.

\begin{proof}
This follows from \ref{sec:12.1.4}, \ref{sec:12.1.2}(1bis), a standard 
pushout argument and \ref{sec:9.1.5}(2).
\end{proof}

\subsubsection{} \label{sec:12.1.8}
Suppose that $\mathcal{A}$ satisfies AB4$^{\ast}$ or is a 
Grothendieck category. Let $(\mathcal{K},\mathcal{K}^{\perp})$ 
be a left exact and left complete cotorsion theory 
in $\mathcal{A}$ such that $\mathcal{A}$ has enough 
$\mathcal{K}^{\perp}$ objects. Then the cotorsion theory
$(ex[\mathcal{K}],dg(\mathcal{K}^{\perp}))$ is right complete.
\begin{proof}
If $\mathcal{A}$ satisfies AB4$^{\ast}$ the proof is dual to the proof
of \ref{sec:12.1.7}. Suppose $\mathcal{A}$ is a Grothendieck category
and let $\mathbf{C}$ be a complex in $\mathcal{A}$. 
By \cite[Th\'{e}or\`{e}me 2]{Jo} there are a dg-injective (\ref{sec:10.3})
complex $\mathbf{Y}$ and a quasi-isomorphism $\mathbf{C}\rightarrowtail 
\mathbf{Y}$. In particular $\mathbf{Y}\in dg(\mathcal{K}]^{\perp})$.
The rest is a standard pushout argument using \ref{sec:12.1.2}(1) and 
\ref{sec:9.1.5}(2).
\end{proof}

\subsubsection{} \label{sec:12.1.9}
\cite[Lemma I.4.6 2)]{Ha} 
Let $\mathcal{K}$ be a class of objects of $\mathcal{A}$ such that 
(i) $\mathcal{A}$ has enough $\mathcal{K}$ objects, meaning that 
for each $A\in Ob(\mathcal{A})$ there is an epimorphism $X\to A$ 
with $X\in \mathcal{K}$, and (ii) $\mathcal{K}$ is closed under
extensions and left exact. For every complex $\mathbf{C}$ in 
$\mathcal{A}$ there are $\mathbf{X}\in Ch(\mathcal{K})$ 
and a quasi-isomorphism $\mathbf{X}\twoheadrightarrow \mathbf{C}$. 
\begin{proof}
Choose an epimorphism $\mathbf{C}^{0}_{0}\to \mathbf{C}_{0}$
with $\mathbf{C}^{0}_{0}\in \mathcal{K}$ and then choose, 
successively, epimorphisms $\mathbf{C}^{0}_{1}\to \mathbf{C}_{1}
\times_{\mathbf{C}_{0}}\mathbf{C}^{0}_{0}$,
...,$\mathbf{C}^{0}_{n}\to \mathbf{C}_{n}\times_{\mathbf{C}_{n-1}}
\mathbf{C}^{0}_{n-1}$,...
with $\mathbf{C}^{0}_{1},...,\mathbf{C}^{0}_{n}\in \mathcal{K}$.
We obtain a complex $\mathbf{C}^{0}$ and a quasi-isomorphism (\ref{sec:4.1.1})
$u^{0}:\mathbf{C}^{0}\twoheadrightarrow \mathbf{C}$ such that 
$\mathbf{C}^{0}_{n}=\mathbf{C}_{n}$ and $u^{0}_{n}=1_{\mathbf{C}_{n}}$
for $n<0$. Next choose an epimorphism $\mathbf{C}^{1}_{-1}
\to \mathbf{C}^{0}_{-1}=\mathbf{C}_{-1}$ with 
$\mathbf{C}^{1}_{-1}\in \mathcal{K}$, then choose an epimorphism
$\mathbf{C}^{1}_{0}\to \mathbf{C}^{0}_{0}\times_{\mathbf{C}_{-1}}
\mathbf{C}^{1}_{-1}$ with $\mathbf{C}^{1}_{0}\in \mathcal{K}$ and put 
$\mathbf{C}^{1}_{1}=\mathbf{C}^{0}_{1}\times_{\mathbf{C}^{0}_{0}}
\mathbf{C}^{1}_{0}$. The square diagram 
\begin{displaymath}
\xymatrix{
{\mathbf{C}^{1}_{1}}\ar[r] \ar@{->>}[d]&{\mathbf{C}^{1}_{0}} \ar@{->>}[d]\\
{\mathbf{C}^{0}_{1}} \ar[r]&{\mathbf{C}^{0}_{0}}
}
\end{displaymath} 
is exact hence (\ref{sec:1.2.8}) $\mathbf{C}^{1}_{1}\in \mathcal{K}$.
We obtain a complex $\mathbf{C}^{1}$ and a quasi-isomorphism 
$u^{1}:\mathbf{C}^{1}\twoheadrightarrow \mathbf{C}^{0}$
such that $\mathbf{C}^{1}_{n}=\mathbf{C}^{0}_{n}$ for $n\geqslant 2$,
$\mathbf{C}^{1}_{n}=\mathbf{C}_{n}$ for $n\leqslant -2$ and 
$u^{1}_{n}$ is an identity morphism except for $n\in\{-1,0,1\}$.
Continuing like in the preceding step we obtain, for every integer $k>0$, 
a complex $\mathbf{C}^{k}$ and a quasi-isomorphism 
$u^{k}:\mathbf{C}^{k}\twoheadrightarrow \mathbf{C}^{k-1}$
such that $\mathbf{C}^{k}_{n}=\mathbf{C}^{k-1}_{n}$ for $n\geqslant -k+3$,
$\mathbf{C}^{k}_{n}=\mathbf{C}_{n}$ for $n\leqslant -k-1$ and 
$u^{k}_{n}$ is an identity morphism except for $n\in\{-k,-k+1,-k+2\}$.
Below is the diagram $\mathbf{C}^{3}\overset{u^{3}}
\to\mathbf{C}^{2}\overset{u^{2}}\to
\mathbf{C}^{1}\overset{u^{1}}\to\mathbf{C}^{0}
\overset{u^{0}}\to\mathbf{C}$.
\[
\xymatrix{
{\mathbf{C}^{0}_{3}}\ar@{=}[d]\ar[r]&{\mathbf{C}^{0}_{2}}\ar[r]\ar@{=}[d]& 
{\mathbf{C}^{1}_{1}}\ar[r]\ar@{=}[d] & {\mathbf{C}^{2}_{0}}
\ar[r]\ar@{=}[d]& {\mathbf{C}^{3}_{-1}}\ar[r]\ar@{->>}[d] & {\mathbf{C}^{3}_{-2}}
\ar@{->>}[d]\ar[r]&{\mathbf{C}^{3}_{-3}}\ar@{->>}[d] \\
{\mathbf{C}^{0}_{3}}\ar@{=}[d]\ar[r]&{\mathbf{C}^{0}_{2}}\ar[r]\ar@{=}[d]& 
{\mathbf{C}^{1}_{1}}\ar[r]\ar@{=}[d] & {\mathbf{C}^{2}_{0}}
\ar[r] \ar@{->>}[d]& {\mathbf{C}^{2}_{-1}}\ar[r]\ar@{->>}[d] & 
{\mathbf{C}^{2}_{-2}}\ar@{->>}[d]\ar[r]&{\mathbf{C}_{-3}}\ar@{=}[d]\\
{\mathbf{C}^{0}_{3}}\ar[r]\ar@{=}[d]&{\mathbf{C}^{0}_{2}}\ar@{=}[d]\ar[r]& 
{\mathbf{C}^{1}_{1}}\ar[r]\ar@{->>}[d] & {\mathbf{C}^{1}_{0}}\ar@{->>}[d]
\ar[r] & {\mathbf{C}^{1}_{-1}}\ar[r] \ar@{->>}[d]& 
{\mathbf{C}_{-2}}\ar[r]\ar@{=}[d]&{\mathbf{C}_{-3}}\ar@{=}[d]\\
{\mathbf{C}^{0}_{3}}\ar[r]\ar@{->>}[d]&{\mathbf{C}^{0}_{2}}\ar[r]\ar@{->>}[d]& 
{\mathbf{C}^{0}_{1}}\ar[r]\ar@{->>}[d]& {\mathbf{C}^{0}_{0}}
\ar[r]\ar@{->>}[d]& {\mathbf{C}_{-1}}\ar[r]\ar@{=}[d]& {\mathbf{C}_{-2}}
\ar[r]\ar@{=}[d]& {\mathbf{C}_{-3}}\ar@{=}[d]\\
{\mathbf{C}_{3}}\ar[r]&{\mathbf{C}_{2}}\ar[r]& 
{\mathbf{C}_{1}}\ar[r]& {\mathbf{C}_{0}}
\ar[r]& {\mathbf{C}_{-1}}\ar[r]& {\mathbf{C}_{-2}}
\ar[r]& {\mathbf{C}_{-3}}
}
\]
The resulting inverse system $...\twoheadrightarrow \mathbf{C}^{2}
\twoheadrightarrow \mathbf{C}^{1}\twoheadrightarrow \mathbf{C}^{0}$ 
has the property that for every integer $n$ the inverse system 
$...\twoheadrightarrow \mathbf{C}^{2}_{n}\twoheadrightarrow 
\mathbf{C}^{1}\twoheadrightarrow \mathbf{C}^{0}_{n}$
is eventually constant and eventually in $\mathcal{K}$. Thus,
if we put $\mathbf{X}=\underset{k\in\mathbb{N}^{op}}\li \mathbf{C}^{k}$,
then $\mathbf{X}\in Ch(\mathcal{K})$ and there is a quasi-isomorphism
$\mathbf{X}\twoheadrightarrow \mathbf{C}$.
\end{proof}

\subsection{} \label{sec:12.2}
Let $\mathcal{A}$ be an abelian category. We denote by $Ch_{\geqslant 0}(\mathcal{A})$
the full subcategory of $Ch(\mathcal{A})$ with objects $\mathbf{C}$ such that 
$\mathbf{C}_{n}=0$ for $n<0$. The objects of $Ch_{\geqslant 0}(\mathcal{A})$
will be called nonnegative complexes. 

The inclusion functor $e_{0}:Ch_{\geqslant 0}(\mathcal{A})\to Ch(\mathcal{A})$ 
has a right adjoint $\tau_{0}$, where $\tau_{0}(\mathbf{C})$ is the complex
defined in \ref{sec:4.2.4}. The functor $e_{0}$ has a left adjoint $t_{0}$ 
that sends $\mathbf{C}$ to the complex $...\to\mathbf{C}_{1}\to
\mathbf{C}_{0}\to 0\to...$. The functor $t_{0}$ has a left adjoint 
$q_{0}$ that sends a nonnegative complex $\mathbf{C}$ to
$...\to\mathbf{C}_{1}\to\mathbf{C}_{0}\to\mathbf{C}_{0}/B_{0}(\mathbf{C})\to 0\to...$.

For each integer $n\geqslant 0$ we denote by $(-)_{n}:Ch_{\geqslant 0}(\mathcal{A})
\to \mathcal{A}$ the evaluation at $n$ functor; for $n>0$ the functor
$(-)_{n-1}$ has the $n$th disk functor $D^{n}$ as right adjoint and $D^{n}$ has
$(-)_{n}$ as right adjoint. 

We put $ex_{\geqslant 0}(\mathcal{A})=ex(\mathcal{A})\cap Ch_{\geqslant 0}(\mathcal{A})$
and for a class $\mathcal{K}$ of objects of $\mathcal{A}$ we put
$Ch_{\geqslant 0}(\mathcal{K})=Ch(\mathcal{K})\cap Ch_{\geqslant 0}(\mathcal{A})$
and $ex_{\geqslant 0}[\mathcal{K}]=ex[\mathcal{K}]\cap Ch_{\geqslant 0}(\mathcal{A})$.

\subsubsection{} \label{sec:12.2.1}
\cite[Lemma I.4.6 1)]{Ha} Let $(\mathcal{K},\mathcal{K}^{\perp})$ be a cotorsion 
theory in $\mathcal{A}$ such that $\mathcal{A}$ has enough $\mathcal{K}$ objects. 
Then $(Ch_{\geqslant 0}(\mathcal{K}),ex_{\geqslant 0}[\mathcal{K}^{\perp}])$ 
is a cotorsion theory in $Ch_{\geqslant 0}(\mathcal{A})$ that is left exact (left complete) 
if $(\mathcal{K},\mathcal{K}^{\perp})$ is left exact (left complete). If $\mathcal{A}$
has enough $\mathcal{K}^{\perp}$ objects then $Ch_{\geqslant 0}(\mathcal{A})$
has enough $ex_{\geqslant 0}[\mathcal{K}^{\perp}]$ objects. If 
$\mathcal{K}$ is left exact then $ex_{\geqslant 0}[\mathcal{K}]=
Ch_{\geqslant 0}(\mathcal{K})\cap ex_{\geqslant 0}(\mathcal{A})$.

\begin{proof}
We apply \ref{sec:2.1.2}(1bis) to the adjoint pair $(t_{0},e_{0})$ 
and $(dg(\mathcal{K}),ex[\mathcal{K}^{\perp}])$.
We have $e_{0}^{-1}(ex[\mathcal{K}^{\perp}])=
ex_{\geqslant 0}[\mathcal{K}^{\perp}]$. Let $\mathbf{C}\in 
^{\perp}(ex_{\geqslant 0}[\mathcal{K}^{\perp}])$. Applying the functor 
$Ch_{\geqslant 0}(\mathcal{A})(\mathbf{C},-)$ to the exact sequence 
$D^{n}(Y)\rightarrowtail D^{n}(A)\twoheadrightarrow 
D^{n}(B)$ with $n>0$ and $Y\in\mathcal{K}^{\perp}$
we obtain by adjunction and the fact that 
$D^{n}(Y)\in ex_{\geqslant 0}[\mathcal{K}^{\perp}]$ 
the exact sequence 
$\mathcal{A}(\mathbf{C}_{n-1},Y)\rightarrowtail \mathcal{A}(\mathbf{C}_{n-1},A)
\twoheadrightarrow\mathcal{A}(\mathbf{C}_{n-1},B)$.
This implies that $\mathbf{C}_{n-1}\in\mathcal{K}$ for $n>0$. Conversely,
we have $Ch_{\geqslant 0}(\mathcal{K})\subset 
^{\perp_{h}}ex_{\geqslant 0}[\mathcal{K}^{\perp}]$
by \ref{sec:4.2.2}, so it suffices (\ref{sec:8.2.5}(3)) to show 
that if $\mathbf{D}\in Ch_{\geqslant 0}(\mathcal{K})$ then
$\mathrm{Homgr}_{\mathcal{A}}(\mathbf{D},-)$ preserves 
epimorphisms with kernel in $ex_{\geqslant 0}[\mathcal{K}^{\perp}]$. 
But this is true by \ref{sec:8.1.1}(1bis). Therefore 
$(Ch_{\geqslant 0}(\mathcal{K}),ex_{\geqslant 0}[\mathcal{K}^{\perp}])$ 
is a cotorsion theory. Left exactness is straightforward.
We show that $(Ch_{\geqslant 0}(\mathcal{K}),ex[\mathcal{K}^{\perp}])$
is left complete. Let $\mathbf{C}$ be a nonnegative complex.
We use induction over $n\geqslant 0$ to construct a commutative diagram
\begin{displaymath}
\xymatrix{
{B_{n}(\mathbf{Y})}\ar@{>->}[r] \ar@{>->}[d]& {\mathbf{Y}_{n}}\ar@{>->}[d]\\
{Z_{n}(\mathbf{X})} \ar@{>->}[r] \ar@{->>}[d]& {\mathbf{X}_{n}}\ar@{->>}[d]\\
{Z_{n}(\mathbf{C})} \ar@{>->}[r] & {\mathbf{C}_{n}}
}
\end{displaymath}
with exact columns and with $\mathbf{X}_{n}\in
\mathcal{K},B_{n}(\mathbf{Y})\in\mathcal{K}^{\perp}$.
When $n=0$ we can find an exact sequence 
$\mathbf{Y}_{0}\rightarrowtail \mathbf{X}_{0}
\twoheadrightarrow \mathbf{C}_{0}$ with $\mathbf{X}_{0}\in\mathcal{K},
\mathbf{Y}_{0}\in\mathcal{K}^{\perp}$ and we put 
$Z_{0}(\mathbf{X})=\mathbf{X}_{0},B_{0}(\mathbf{Y})=\mathbf{Y}_{0}$.
When $n>0$ and given the previous diagram, we put 
$B_{n}(\mathbf{X})=B_{n}(\mathbf{X})\times_{Z_{n}
(\mathbf{C})}Z_{n}(\mathbf{X})$
and then we form the commutative diagram
\begin{displaymath}
\xymatrix{
{B_{n+1}(\mathbf{Y})}\ar@{>->}[d]\ar@{=}[r]&{B_{n+1}(\mathbf{Y})}\ar@{>->}[d]\\
{PB'_{n+1}}\ar@{>->}[r] \ar@{->>}[d]& {\mathbf{X}_{n+1}}\ar@{->>}[r]\ar@{->>}[d] 
& {B_{n}(\mathbf{X})}\ar@{=}[d]\\
{Z_{n+1}(\mathbf{C})} \ar@{>->}[r] \ar@{=}[d]& {PB_{n+1}}\ar@{->>}[r]\ar@{->>}[d] & 
{B_{n}(\mathbf{X})}\ar@{->>}[d]\\
{Z_{n+1}(\mathbf{C})} \ar@{>->}[r] & {\mathbf{C}_{n+1}}\ar@{->>}[r] & {B_{n}(\mathbf{C})}
}
\end{displaymath}
where $PB_{n+1}$ and $PB'_{n+1}$ mean pullback and $B_{n+1}(\mathbf{Y})
\rightarrowtail \mathbf{X}_{n+1}\twoheadrightarrow PB_{n+1}$ is an exact 
sequence with $\mathbf{X}_{n+1}\in\mathcal{K}$ and 
$B_{n+1}(\mathbf{Y})\in\mathcal{K}^{\perp}$.
We put $Z_{n+1}(\mathbf{X})=PB'_{n+1}$ and $\mathbf{Y}_{n+1}$ 
to be the kernel of the composite morphism 
$\mathbf{X}_{n+1}\to PB_{n+1}\to \mathbf{C}_{n+1}$.
We obtain the commutative solid arrows diagram
\begin{displaymath}
\xymatrix{
{B_{n+1}(\mathbf{Y})}\ar@{>..>}[r] \ar@{>->}[d]& 
{\mathbf{Y}_{n+1}}\ar@{..>>}[r]\ar@{>->}[d] 
& {B_{n}(\mathbf{Y})}\ar@{>->}[d]\\
{Z_{n+1}(\mathbf{X})} \ar@{>->}[r] \ar@{->>}[d]& 
{\mathbf{X}_{n+1}}\ar@{->>}[r]
\ar@{->>}[d] & {B_{n}(\mathbf{X})}\ar@{->>}[d]\\
{Z_{n+1}(\mathbf{C})} \ar@{>->}[r] & 
{\mathbf{C}_{n+1}}\ar@{->>}[r] & {B_{n}(\mathbf{C})}
}
\end{displaymath}
By the universal property of kernel we have the induced dotted arrows that 
make the resulting square diagrams commute; moreover, by the snake diagram 
the dotted sequence of the previous diagram is a short exact sequence.
This finishes the inductive step. We have constructed an exact sequence 
$\mathbf{Y}\rightarrowtail \mathbf{X}\twoheadrightarrow\mathbf{C}$
with $\mathbf{X}\in Ch_{\geqslant 0}(\mathcal{K})$ and 
$\mathbf{Y}\in ex_{\geqslant 0}[\mathcal{K}^{\perp}]$.
The fact that $Ch_{\geqslant 0}(\mathcal{A})$ has enough 
$ex_{\geqslant 0}[\mathcal{K}^{\perp}]$ objects follows in 
the same way as for $Ch(\mathcal{A})$ (\ref{sec:4.1.2}).
Finally, we show that $ex_{\geqslant 0}[\mathcal{K}]=
Ch_{\geqslant 0}(\mathcal{K})\cap ex_{\geqslant 0}(\mathcal{A})$. 
For this, let $\mathbf{C}\in Ch_{\geqslant 0}(\mathcal{K})\cap 
ex_{\geqslant 0}(\mathcal{A})$. For every $n\geqslant 1$ we have 
the exact sequence $Z_{n}(\mathbf{C})\rightarrowtail \mathbf{C}_{n}
\twoheadrightarrow Z_{n-1}(\mathbf{C})$. Since $Z_{0}(\mathbf{C})=
\mathbf{C}_{0}\in\mathcal{K}$, we obtain inductively that 
$Z_{n}(\mathbf{C})\in\mathcal{K}$ for all $n\geqslant 0$.
\end{proof}
%(2) NOT CLEAR Let $(\mathcal{K},\mathcal{K}^{\perp})$ be a cotorsion 
%theory in $\mathcal{A}$ such that $\mathcal{A}$ has enough $\mathcal{K}^{\perp}$ objects. 
%Then $(ex_{\geqslant 0}[\mathcal{K}],Ch_{\geqslant 0}(\mathcal{K}^{\perp}))$ 
%is a cotorsion theory in $Ch_{\geqslant 0}(\mathcal{A})$ that is right exact (right complete) 
%if $(\mathcal{K},\mathcal{K}^{\perp})$ is right exact (right complete). If $\mathcal{A}$
%has enough $\mathcal{K}$ objects then $Ch_{\geqslant 0}(\mathcal{A})$
%has enough $ex_{\geqslant 0}[\mathcal{K}]$ objects. 

\end{document}